\documentclass[12pt,a4paper,english]{book}
\usepackage[latin1]{inputenc}
\usepackage[T1]{fontenc}
\usepackage{amsmath,amssymb,amscd}
\usepackage{babel,epsfig}
\usepackage{a4dutch}
\usepackage{theorem,alltt}
\usepackage{mathrsfs}

\newcommand{\id}{\mathit{id}}
\newcommand{\pr}{\mathit{pr}}

\newcommand{\op}{\operatorname{op}}

\newcommand{\C}{\mathbb{C}}
\newcommand{\R}{\mathbb{R}}
\newcommand{\Z}{\mathbb{Z}}
\newcommand{\U}{\mathbb{U}}
\newcommand{\bbP}{\mathbb{P}}

\newcommand{\Ind}{\operatorname{Ind}}
\newcommand{\Orb}{\operatorname{Orb}}

\newcommand{\Ens}{\mathcal{E}\!\mathit{ns}}
\newcommand{\Top}{\mathcal{T}\!\!\mathit{op}}

\newcommand{\cy}{\operatorname{cy}}
\newcommand{\ev}{\operatorname{ev}}
\newcommand{\sd}{\operatorname{sd}}
\newcommand{\sk}{\operatorname{sk}}

\newcommand{\All}{\mathscr{A}\ell\ell}
\newcommand{\Tr}{\operatorname{Tr}}
\newcommand{\incl}{\operatorname{incl}}
\newcommand{\hocolim}{\operatornamewithlimits{hocolim}}
\newcommand{\holim}{\operatornamewithlimits{holim}}
\newcommand{\colim}{\operatornamewithlimits{colim}}
\newcommand{\catDelta}{\mathbf{\Delta}}

\newcommand{\catDeltaT}{\mathbf{\Delta T}}
\newcommand{\catDeltaC}{\mathbf{\Delta C}}
\newcommand{\catDeltaG}{\mathbf{\Delta G}}
\newcommand{\catDeltaD}{\mathbf{\Delta D}}

\newcommand{\Sing}{\operatorname{Sing}_{\bullet}}

\newcommand{\Fix}{\operatorname{Fix}}

\newcommand{\mb}[1]{\mathbf{#1}}

\theoremstyle{break}
\newtheorem{Def}{Definition}[section]

\newtheorem{Thm}[Def]{Theorem}
\newtheorem{Prop}[Def]{Proposition}
\newtheorem{Lem}[Def]{Lemma}

\newtheorem{Cor}[Def]{Corollary}
{\theorembodyfont{\rmfamily} \newtheorem{Exa}[Def]{Example}}
{\theorembodyfont{\rmfamily} \newtheorem{Rem}[Def]{Remark}}
{\theorembodyfont{\rmfamily} \newtheorem{Constr}[Def]{Construction}}
{\theorembodyfont{\rmfamily} \newtheorem{Conj}[Def]{Conjecture}}

\newenvironment{proof}{\mbox{}\newline\textbf{Proof:} }{\nopagebreak\hfill$\square$\\}

\title{Involutions on $S[\Omega M]$}
\author{Tore August Kro\\ \texttt{toreak@math.uio.no}}

\begin{document}

\pagestyle{empty}

\newpage

\maketitle

\newpage
\mbox{}
\newpage

\vspace*{\stretch{1}}
\LARGE
\hfill{}To Tordis and Andreas
\vspace*{\stretch{1}}
\vspace*{\stretch{1}}
\normalsize
\newpage

\chapter*{Introduction}
\thispagestyle{empty}

An ultimate goal is to calculate the homotopy type of automorphism spaces
of compact manifolds $M$. A promising approach uses the parametrized
surgery theory introduced by Hsiang and Sharpe in~\cite{HsiangSharpe:76}.
Historically, one studied the fiber of map from parametrized surgery to surgery
using the involution on concordance theory, see~\cite{HsiangJahren:82} and~\cite{Burghelea:78}.
Two step computations along these lines are, for example, found 
in~\cite{HsiangSharpe:76},~\cite{HsiangJahren:83},~\cite{HsiangJahren:82b} and~\cite{FarrellHsiang:78}.
However, the experts knew that one should try to combine surgery theory and concordance theory 
to get direct computations of parametrized surgery. Recent advances in homotopy theory has made
it possible for Weiss and Williams to define this
$LA$-theory, see~\cite{WeissWilliams:01}. 
Still, $LA$-theory is not easily computed, but it is
related to algebraic $K$-theory via a map called $\Xi$,
and algebraic $K$-theory can be studied via  
trace maps into $TC$ or $THH$, see~\cite{Madsen:94}.

Surgery theory classifies manifolds within a given simple homotopy type. A basic
ingredient in this theory is the $L$-groups, whose definition depends on the
ring $\Z[\pi_1 M]$ together with an involution, see~\cite{Wall:99}. 
In concordance theory, turning a concordance upside down gives an involution
on the homotopy groups of the stabilized concordance space, $\mathcal{C}(M)$.
Hatcher's spectral sequence, see~\cite{Hatcher:78} proposition~2.2, 
has in a stable range an $E^2$-page given by
$E^2_{pq}=H_p(\Z/2;\pi_q\mathcal{C}(M))$, and converges to give information
relating surgery theory to automorphism spaces at the level of homotopy groups.
Weiss and Williams strengthen Hatcher's result to the 
level of spaces in~\cite{WeissWilliams:88}.
In~\cite{Vogell:85}, Vogell shows that the involution on concordance spaces
corresponds to his canonical involution on algebraic $K$-theory of spaces. Also Steiner,~\cite{Steiner:81},
and Hsiang and Jahren,\cite{HsiangJahren:82},~\cite{HsiangJahren:preprint}, considers involutions on $A$-theory.

The input for $LA$, $K$, $TC$ and $THH$ should be a ``ring up to homotopy''.
We choose orthogonal ring spectra, as defined in~\cite{MandellMaySchwedeShipley:01},
to be our model for such rings. Hence, the theories above need to be redefined for this setting.
Waldhausen's work on algebraic $K$-theory tells us that $S[\Omega M]$ is the correct
ring to consider when relating to concordances.
Geometry, as discussed above, demands an involution on our ring, and Vogell
shows that interesting involutions come from bundles $\xi$ over $M$.
However, it is not immediately clear that such $\xi$ give involutions
on the orthogonal ring spectrum $S[\Omega M]$.

The contribution of this thesis is to construct for each vector bundle $\xi$ 
an orthogonal ring spectrum $R$, weakly equivalent to $S[\Omega M]$,
together with an involution on $R$. On the homotopy groups $\pi_*S[\Omega M]$
our involution corresponds to parallel transportation in $\xi$, and reversing loops in $M$.
The main result is theorem~\ref{Thm:main}. The orthogonal ring spectrum $R$ with involution
is intended as input for the theories mentioned above. 
We take a first step in this direction by considering 
the definition and a few basic properties of $TC(L)$ and $THH(L)$ 
for arbitrary orthogonal ring spectra $L$ (with involution).
However, in the future the author hopes to show 
that $LA(R)$ will yield information about the automorphism space of $M$.

Chapter~\ref{Chapt:simp} recalls from the literature various simplicial techniques.
Much of this should be well known to the reader. The reason the author  
included this material is mainly to point out certain viewpoints and to introduce
notation. 

In order to get strong results in stable homotopy theory, one can prove
theorems in a category of spectra with symmetric smash product. We
find orthogonal spectra particularly convenient for our purposes.
In chapter~\ref{Chapt:ort} we give an exposition of the relevant theory 
and develop the techniques needed to work within this category.
 
We study the equivariant homotopy theory of orthogonal spectra for two reasons;
the definition of an operad involves actions of symmetric groups, and $THH$ comes
with an $S^1$-action, or $O(2)$-action in the involutive setting. 
Chapter~\ref{Chapt:equiv} is an introduction
to equivariant orthogonal spectra and
provides the results necessary for our applications.

Chapter~\ref{Cha:invoper} proves the main theorem. 
An important ingredient is the concept of operads, and we introduce the operad
$\mathcal{H}$, which encodes multiplication together with an anti-commutative involution.
Moreover, we explain the notion of an operad in orthogonal spectra.
The main idea of the proof is to start out with a vector bundle $\xi$ over our compact manifold
$M$ and then attempt to construct a related involution on $S[\Omega M]$.
Doing this involves many choices, in fact there are orthogonal spectra $\mathcal{D}_n(j)$
parameterizing this. Could the collection $\mathcal{D}_n$ be an operad 
where $S[\Omega M]$ is its algebra by the parametrization? The answer is yes, 
and the formulas for the composition operations of the operad are 
forced by the algebra structure.
Furthermore, $\mathcal{D}_n$ is ``up to homotopy'' sufficiently equal to $\mathcal{H}$.
Using May's two-sided bar construction, we therefore can replace $S[\Omega M]$ by a weakly
equivalent $\mathcal{H}$-algebra. This gives our orthogonal ring spectrum $R$ with involution.
Unfortunately, the logic demands that we reverse this argument when writing out the proof.
 
Chapter~\ref{Chapt:trace} ends this thesis. It contains some theory regarding
$THH$ and $TC$ of orthogonal ring spectra with involution.

I am very grateful to my advisor Bj\o{}rn Jahren for many enlightening discussions, and for
his constructive feedback during my writing of the manuscript. You have always been available, and 
you have a keen eye for the beauty in mathematics. 
Furthermore, I would like to thank 
Sverre Lun\o{}e-Nielsen, Christian Schlichtkrull, Halvard Fausk, and John Rognes
for helpful conversations.
The support from my wife, Tordis Fuskeland, has been very valuable to me.  
You have carefully proofread the final manuscript, but the errors that remain are mine.
Together with our son, Andreas, you are the most important part of my life.

\newpage
\tableofcontents
\addtocontents{toc}{\protect\thispagestyle{empty}}
\newpage
\mbox{}
\newpage

\newpage
\pagestyle{headings}
\setcounter{page}{1}
\chapter{Simplicial techniques}\label{Chapt:simp}

The theory of simplicial sets and simplicial spaces is classical.
Simplicial sets were first defined in~\cite{EilenbersZilber:50}. The geometric realization
was defined in~\cite{Milnor:57}. Other references to simplicial techniques
are~\cite{May:67},~\cite{GoerssJardine:99} and~\cite{DwyerHenn:01}. 
For the theory of simplicial spaces see~\cite{May:72},~\cite{Segal:74}
and~\cite{Madsen:94}.
This chapter will recall from the literature the simplicial techniques which
are relevant to this thesis. One reason for including this material is for completeness,
but also important is the viewpoint and the notation.

\section{The category $\catDelta$ and its relatives}

\begin{Def}
Let $[n]$ be the ordered set $\{0<1<\ldots<n\}$. The category $\catDelta$ has one object $[n]$ for each
non-negative integer $n$, and the morphisms are ordering preserving functions $\phi:[m]\rightarrow[n]$.
\end{Def}

It is customary to let $\delta_i:[n-1]\rightarrow[n]$ be the order preserving function that misses $i$,
and $\sigma_i:[n+1]\rightarrow[n]$ be the order preserving function that hits $i$ twice. The $\delta$'s and $\sigma$'s
generate all morphisms in $\catDelta$.  

\begin{Def}
A \textit{simplicial set} is a functor $X_{\bullet}:\catDelta^{\op}\rightarrow\Ens$.
A \textit{simplicial space} is a functor $X_{\bullet}:\catDelta^{\op}\rightarrow\Top$.
More generally, one can define simplicial objects in any category.
\end{Def}

Observe that simplicial sets can be considered as simplicial spaces by giving each $X_n$ the 
discrete topology. Hence, in most cases we can do our constructions
for simplicial spaces, the corresponding results for simplicial sets follow implicitly.

Given a simplicial space $X_{\bullet}$, the following notation and terminology is standard:
The space $X_n$ is called the \textit{$n$-simplices} of $X_{\bullet}$. $\delta_i:[n-1]\rightarrow[n]$
induces a map $d_i:X_n\rightarrow X_{n-1}$ called the \textit{$i$'th face map}, and
$\sigma_i:[n+1]\rightarrow[n]$ induces the \textit{$i$'th degeneracy map}, $s_i:X_n\rightarrow X_{n+1}$.
A simplex $x$ in $X_n$ is said to be \textit{degenerate} if $x=s_ix'$ for some $i$ and $x'\in X_{n-1}$.
We denote by $sX_{n-1}$ the subspace of $X_n$ consisting of the degenerate simplices.

\subsection{Geometric realization of simplicial spaces}

Simplicial spaces are combinatorial models for topological spaces, and
geometric realization is the functor which turns a simplicial space
into the topological space for which it is a model.
The geometric realization, due to Milnor~\cite{Milnor:57}, has several good properties,
it commutes with products and it commutes with all colimits.
See~\cite{May:67} or~\cite{DwyerHenn:01}.
Furthermore, every point in the geometric realization is uniquely determined
as the interior point of a non-degenerate simplex. We give a modern 
formulation of this result; giving a filtration for the geometric realization.

Geometric realization of simplicial spaces is defined
using a functor $\Delta:\catDelta\rightarrow\Top$. We send $[n]$ to the space 
$\Delta^n=\{(t_0,\ldots,t_n)\in\R^{n+1}\;|\; \sigma_i t_i=1, t_i\geq 0\}$.
And we call $\Delta^n$ the \textit{topological $n$-simplex}.
On morphisms the functor is defined by sending $\delta_i:[n-1]\rightarrow [n]$ to the map
\[
\delta_i(t_0,\ldots,t_{n-1})=(t_0,\ldots,t_{i-1},0,t_i,\ldots,t_{n-1})\quad,
\]
and $\sigma_i:[n+1]\rightarrow [n]$ to
\[
\sigma_i(t_0,\ldots,t_{n+1})=(t_0,\ldots,t_{i-1},t_i+t_{i+1},t_{i+2},\ldots,t_{n+1})\quad.
\]
Now we define the \textit{geometric realization} of a simplicial space $X_{\bullet}$ as the coend
\[
|X_{\bullet}|=\int^{[n]\in\catDelta}X_n\times \Delta^n\quad.
\]
Coends are defined in section~IX.5 in~\cite{MacLane:98}. The space
$|X_{\bullet}|$ is isomorphic to the quotient of $\coprod X_n\times \Delta^n$ 
where we identify $(x,\phi(\mathbf{t}))$ with $(\phi^*x,\mathbf{t})$
for all morphisms $\phi$ in $\catDelta$.

\begin{Rem}\label{Rem:prodformulaforgeomreal}
In order to have convenient technical properties, one should form
the geometric realization in the category of compactly generated spaces (=weak Hausdorff $k$-spaces),
see~\cite{McCord:69}. This ensures, for example, that the product theorem holds.
\end{Rem}

There is also a presimplicial realization, defined using only the injective morphisms in $\catDelta$.
The injective morphisms are those generated by the $\delta$'s. Let $i\catDelta$ denote this subcategory.
We define the \textit{presimplicial realization} as the coend
\[
\|X_{\bullet}\|=\int^{[n]\in i\catDelta}X_n\times \Delta^n\quad.
\]
This space is the quotient of $\coprod X_n\times \Delta^n$
where we identify $(x,\delta_i(\mathbf{t}))$ with $(d_i(x),\mathbf{t})$
for all $\delta_i$'s.

Whereas the geometrical realization, $|X_{\bullet}|$, has better formal properties, it is often
easier to prove results about the homotopy of the presimplicial realization, $\|X_{\bullet}\|$.
And we can compare the two realizations via a natural map
\[
\|X_{\bullet}\|\rightarrow |X_{\bullet}|\quad.
\]
It is natural to ask when this is a weak homotopy equivalence. This question is answered by
Segal in~\cite{Segal:74} and by May in~\cite{May:72}. We follow Segal and define:

\begin{Def}
A simplicial space $X_{\bullet}$ is \textit{good} if for all $n$ and $i$, the map
$s_i:X_n\rightarrow X_{n+1}$ is a closed cofibration.
\end{Def}

We refer to our remark~\ref{Rem:cofibdef} or one of the articles~\cite{Steenrod:67} or~\cite{Strom:66}
for the definition of a closed cofibration. Observe that any simplicial set automatically is good, since
an injective map between discrete spaces is a closed cofibration. Now, Segal shows in his proposition~A.1(iv) that:

\begin{Prop}\label{Prop:geomandpresimplrealization}
If $X_{\bullet}$ is a good simplicial space, then the natural map
$\|X_{\bullet}\|\rightarrow |X_{\bullet}|$ is a weak equivalence.
\end{Prop}

Let us now describe the realizations more carefully. We have already mentioned Milnor's result,
that a point in $|X_{\bullet}|$ is uniquely given as the interior point of a non-degenerate simplex.
There is a similar description for the presimplicial realization. We now give a
modern formulation of these statements: 

\begin{Constr}\label{Constr:presimplfilt}
First we consider the case of the presimplicial realization.
Recall that $\|X_{\bullet}\|$ is the quotient space formed from $\coprod X_n\times \Delta^n$ 
by identifying $(x,\delta^i(\mathbf{t}))$ with $(d_i x,\mathbf{t})$ for all morphisms $d_i$.
We define a filtration by letting the $q$'th space, $F_q\|X_{\bullet}\|$, $q\geq0$,
be the image of $\coprod_{n\leq q} X_n\times \Delta^n$ in $\|X_{\bullet}\|$.
Notice that $F_q\|X_{\bullet}\|$ is the pushout of
\[
X_q\times\Delta^q\leftarrow X_q\times\partial\Delta^q\rightarrow F_{q-1}\|X_{\bullet}\|\quad.
\]
Now observe that
\begin{itemize}
\item[] $\colim F_q\| X_{\bullet}\|$ is equal to $\| X_{\bullet}\|$, and
\item[] each $X_q\times\partial\Delta^q\rightarrow X_q\times\Delta^q$ is a closed cofibration. It follows that also
$F_{q-1}\| X_{\bullet}\|\rightarrow F_q\| X_{\bullet}\|$ is a closed cofibration.
\end{itemize}
The last observation explains why the presimplicial realization behaves so well homotopically;
it is easy to give inductive arguments using the pushout diagram relating $F_{q-1}\| X_{\bullet}\|$ 
to $F_q\| X_{\bullet}\|$.
\end{Constr}

\begin{Constr}\label{Constr:filtergeomrealization}
Next consider the geometric realization.
Recall that $|X_{\bullet}|$ is the quotient of $\coprod X_n\times \Delta^n$ 
where we identify $(x,\phi(\mathbf{t}))$ with $(\phi^* x,\mathbf{t})$ for all morphisms $\phi$ in $\catDelta$.
We define $F_q|X_{\bullet}|$, $q\geq 0$, to be the image of $\coprod_{n\leq q} X_n\times \Delta^n$.
Recall that the degenerate simplices, $sX_{q-1}$, are the
points in $X_q$ which are in the image of some map $s_i:X_{q-1}\rightarrow X_q$.
Notice that the following diagram is pushout:
\[
\begin{CD}
X_q\times\partial\Delta^q\cup sX_{q-1}\times\Delta^q @>>> F_{q-1}|X_{\bullet}|\\
@V{j}VV @VVV\\
X_q\times \Delta^q@>>> F_q|X_{\bullet}|
\end{CD}\quad.
\]
Here the map $j$ comes from the square
\[
\begin{CD}
sX_{q-1}\times\partial\Delta^q @>>> sX_{q-1}\times\Delta^q\\
@VVV @VVV\\
X_{q}\times\partial\Delta^q @>>> X_{q}\times\Delta^q
\end{CD}\quad.
\]
By Lillig's union theorem, see~\cite{Lillig:73}, we have that
$j$ is a closed cofibration whenever $sX_{q-1}\subset X_q$ is a closed cofibration.
Now observe that
\begin{itemize}
\item[] $\colim F_q| X_{\bullet}|$ is equal to $| X_{\bullet}|$, and
\item[] if each $sX_{q-1}\subset X_q$ is a closed cofibration, then
$X_q\times\partial\Delta^q\cup sX_{q-1}\times\Delta^q\rightarrow X_q\times\Delta^q$ and 
$F_{q-1}| X_{\bullet}|\rightarrow F_q| X_{\bullet}|$ are closed cofibrations.
\end{itemize}
The last observation explains why good simplicial spaces behave well with respect to homotopy.
\end{Constr}

These filtrations are extremely useful when proving results about realizations.
We illustrate this by proving a few well known facts:

\begin{Prop}
Let $f_{\bullet}:X_{\bullet}\rightarrow Y_{\bullet}$ be a map of simplicial spaces
such that each $f_q$ is a weak homotopy equivalence. Then the induced map
\[
\|f_{\bullet}\|:\|X_{\bullet}\|\rightarrow \|Y_{\bullet}\|
\]
also is a weak homotopy equivalence.
\end{Prop}

\begin{proof}
We use the filtration from construction~\ref{Constr:presimplfilt}, and prove inductively
that $F_q\|X_{\bullet}\|\rightarrow F_q\|Y_{\bullet}\|$ is a weak homotopy equivalence.
It will follow that $\|f\|$ is a weak homotopy equivalence.

$F_0\|X_{\bullet}\|=X_0\rightarrow Y_0=F_0\|Y_{\bullet}\|$ is a weak homotopy equivalence by assumption.
Now consider the inductive step. We have the diagram
\[
\begin{CD}
X_q\times\Delta^q @<{j}<< X_q\times\partial\Delta^q @>>> F_{q-1}\|X_{\bullet}\|\\
@V{f_q\times\id}VV @V{f_q\times\id}VV @VVV\\
Y_q\times\Delta^q @<{j'}<< Y_q\times\partial\Delta^q @>>> F_{q-1}\|Y_{\bullet}\|
\end{CD}\quad.
\]
Here $j$ and $j'$ are closed cofibrations, while all vertical maps are weak homotopy equivalences.
By proposition~\ref{Prop:gluewe}, the map of the row-wise pushouts,
\[
F_q\|X_{\bullet}\|\rightarrow F_q\|Y_{\bullet}\|
\]
also is a weak homotopy equivalence. This finishes the proof.
\end{proof}

A straight forward corollary of this proposition together with 
proposition~\ref{Prop:geomandpresimplrealization} is:

\begin{Cor}
If $f_{\bullet}:X_{\bullet}\rightarrow Y_{\bullet}$ is a map between good
simplicial spaces and each $f_q$ is a weak homotopy equivalence, then
\[
|f_{\bullet}|:|X_{\bullet}|\rightarrow |Y_{\bullet}|
\]
also is a weak homotopy equivalence.
\end{Cor}

Let us also prove the product theorem. 
Given two simplicial spaces, $X_{\bullet}$ and $Y_{\bullet}$, we define their product
$X_{\bullet}\times Y_{\bullet}$ to be the simplicial space with $n$-simplices
$X_n\times Y_n$. We have a natural map
$\eta:|X_{\bullet}\times Y_{\bullet}|\rightarrow|X_{\bullet}|\times |Y_{\bullet}|$ 
defined using the natural projections from $X_{\bullet}\times Y_{\bullet}$
into $X_{\bullet}$ and $Y_{\bullet}$. The product theorem states that $\eta$ is a 
homeomorphism. It is hard to find a proof in the literature,
which is of the generality suggested in remark~\ref{Rem:prodformulaforgeomreal}.
This is the reason for including a proof here:

\begin{Prop}
Since the geometric realization is formed in the category of compactly generated spaces,
the natural map
\[
\eta:|X_{\bullet}\times Y_{\bullet}|\rightarrow|X_{\bullet}|\times |Y_{\bullet}|
\]
is a homeomorphism.
\end{Prop}

\begin{proof}
It is well known that $\eta$ is a continuous bijection, see theorem~2 in~\cite{Milnor:57} or
theorem~11.5 in~\cite{May:72}. The hard part is to check that $\eta^{-1}$ is continuous.
May proves this when ``spaces'' is the category of compactly generated Hausdorff spaces,
although his proof for continuity of $\eta^{-1}$ is not particularly clear. The author 
hopes that the argument below will be more understandable, but in essence the proofs are the same.

We use the filtration of the geometrical realization given in construction~\ref{Constr:filtergeomrealization}.
The product $|X_{\bullet}|\times|Y_{\bullet}|$ inherits a filtration given by
\[
F_q(|X_{\bullet}|\times|Y_{\bullet}|)=\bigcup_{m+n=q}F_m|X_{\bullet}|\times F_n|Y_{\bullet}|\quad.
\]
And by the constructions one can see that $\eta$ restricts to a continuous bijection
$F_q|X_{\bullet}\times Y_{\bullet}|\cong F_q(|X_{\bullet}|\times|Y_{\bullet}|)$.

A continuous bijection $\eta$ between compactly generated spaces is a homeomorphism,
if $\eta^{-1}(K)$ is compact whenever $K$ is compact. This follows from the definition
of compactly generated by lemma~2.1 in~\cite{McCord:69}. 

We will now try to apply lemma~2.8 in~\cite{McCord:69}.
Suppose that $K\subset |X_{\bullet}|\times|Y_{\bullet}|$ is a compact subset. Then $K$
is contained in $F_q(|X_{\bullet}|\times|Y_{\bullet}|)$ for some $q$. 
$q$ is now fixed.
Below we will specify $Z_{\alpha}$'s such that for all $\alpha$ we have commutative diagrams
\[
\begin{CD}
\coprod_{p\leq q} X_p\times Y_p\times \Delta^p @<<< 
Z_{\alpha}\subset \coprod_{n,m}X_n\times \Delta^n\times Y_m\times \Delta^m\\
@VVV @VV{\pi}V\\
F_q|X_{\bullet}\times Y_{\bullet}| @>{\eta_q}>> |X_{\bullet}|\times |Y_{\bullet}|
\end{CD}\quad.
\]
Here the vertical maps are surjective
and $\eta_q$ is injective. Furthermore, the target of the map $\pi$
has the quotient topology. 
The $Z_{\alpha}$'s depend on the standard triangulation of $\Delta^n\times\Delta^m$.
For given $m$ and $n$ let $i_{\alpha}:\Delta^{n+m}\rightarrow \Delta^n\times\Delta^m$
be the inclusion of a maximal topological simplex in this triangulation. We set
$Z_{\alpha}=X_n\times Y_m\times \Delta^{n+m}$, included in 
$\coprod_{n,m}X_n\times \Delta^n\times Y_m\times \Delta^m$
via $i_{\alpha}$. Moreover, via the appropriate degeneracy maps, there are
maps $Z_{\alpha}\rightarrow X_p\times Y_p\times \Delta^p$, $p=n+m$, such that the diagram
above commutes. These maps are explicitly constructed in May's proof.

Now observe that the collection of $Z_{\alpha}$'s with $n+m\leq q$ covers the image
of $F_q|X_{\bullet}\times Y_{\bullet}|$. This collection is finite. Thus all conditions
of McCord's lemma~2.8 are satisfied. This implies that $\eta_q$ is an embedding, and
consequently we have that $\eta^{-1}(K)=\eta_q^{-1}(K)$ is compact. And we are done.
\end{proof}

\subsection{Crossed simplicial categories}

Techniques involving simplicial sets or simplicial spaces are extremely useful when working with
topological spaces. However, if we want to consider involutions, $S^1$- or $O(2)$-actions on our spaces,
it is handy to replace $\catDelta$ by other categories; $\catDeltaT$, $\catDeltaC$ and $\catDeltaD$.
We will recall the notion of a crossed simplicial group from~\cite{FiedorowiczLoday:91}.
The categories mentioned above are examples of such, they will be defined below 
and we will introduce notation for their morphisms.

\begin{Def}
A sequence of groups $\{G_n\}$, $n\geq0$, is a \textit{crossed simplicial group}
if it is equipped with the following structure. There is a small category $\catDeltaG$,
which is part of the structure, such that
\begin{itemize}
\item[] the objects of $\catDeltaG$ are $[n]$, $n\geq0$,
\item[] $\catDeltaG$ contains $\catDelta$ as a subcategory,
\item[] the automorphisms of $[n]$ in $\catDeltaG$ is the opposite group of $G_n$, and
\item[] any morphism in $\catDeltaG([m],[n])$ can be uniquely written as a composite $\phi\circ g$,
where $\phi\in\catDelta([m],[n])$ and $g\in G^{\op}_m$.
\end{itemize}
\end{Def}

\begin{Rem}\label{Rem:Gissimplset}
The last axiom implies that for any $g\in G_n$ and $\phi\in\catDelta([m],[n])$ there exist unique
$\phi^*(g)\in G_m$ and $g^*(\phi)\in \catDelta([m],[n])$ such that
\[
g\circ \phi= g^*(\phi)\circ \phi^*(g)\quad.
\]
The functor that sends $[n]$ to $G_n$ and $\phi$ to $\phi^*:G_n\rightarrow G_m$ gives $G_{\bullet}$ the
structure of a simplicial set. 
\end{Rem}

Unlike~\cite{FiedorowiczLoday:91}, our focus will not be these simplicial sets, but
rather the categories $\catDeltaG$ and their analogue of simplicial sets and spaces, i.e. 
functors from $\catDeltaG^{\op}$ into sets and spaces. We will therefore refer to $\catDeltaG$
as a \textit{crossed simplicial category}.

Here are the crossed simplicial categories relevant for this thesis, 
they are taken from the examples~2, 4, 5 and~7 in~\cite{FiedorowiczLoday:91}:

\begin{Def}
Define $\catDeltaT$ to be the crossed simplicial category with
the automorphism group of $[n]$ cyclic of order $2$.
Let $\rho_n$ be the generator of the automorphism group and put
$\rho_{n}\delta_i=\delta_{n-i}\rho_{n-1}$ and  
$\rho_{n}\sigma_i=\sigma_{n-i}\rho_{n+1}$.
\end{Def}

\begin{Def}\label{Def:catDeltaC}
Let $\catDeltaC$ be the crossed simplicial category 
where the automorphism group of $[n]$ is cyclic of order $(n+1)$.
We name the preferred generator $\tau_n$, and introduce the relations:
\begin{align*}
\tau_{n}\delta_0=\delta_{n}\quad&\text{and}\quad \tau_n\delta_i=\delta_{i-1}\tau_{n-1}\text{, for $1\leq i\leq n$, and}\\ 
\tau_{n}\sigma_0=\sigma_{n}\tau_{n+1}^2\quad&\text{and}\quad \tau_{n}\sigma_i=\sigma_{i-1}\tau_{n+1}\text{, for $1\leq i\leq n$.}
\end{align*}
\end{Def}

\begin{Def}\label{Def:catDeltaD}
Let $\catDeltaD$ be crossed simplicial category 
where the automorphism group of $[n]$ is the dihedral group of order $2(n+1)$. 
We name the preferred generators $\rho_n$ and $\tau_n$, where $\rho_n^2=\tau_n^{n+1}=\id$ and $\rho_n\tau_n=\tau_n^{-1}\rho_n$,
and introduce the relations:
\begin{align*}
\rho_{n}\delta_i&=\delta_{n-i}\rho_{n-1},\\
\rho_{n}\sigma_i&=\sigma_{n-i}\rho_{n+1},\\
\tau_{n}\delta_0=\delta_{n}\quad&\text{and}\quad \tau_n\delta_i=\delta_{i-1}\tau_{n-1}\text{, for $1\leq i\leq n$, and}\\ 
\tau_{n}\sigma_0=\sigma_{n}\tau_{n+1}^2\quad&\text{and}\quad \tau_{n}\sigma_i=\sigma_{i-1}\tau_{n+1}\text{, for $1\leq i\leq n$.}
\end{align*}
\end{Def}

\begin{Def}
Let $\catDeltaC_r$, $r\geq 1$, be the crossed simplicial category 
where the automorphism group of $[n]$ is cyclic of order $r(n+1)$.
We name the preferred generator $\tau_n$, where $\tau_n^{r(n+1)}=\id$, and introduce the same relations
as in definition~\ref{Def:catDeltaC}.
\end{Def}

\begin{Def}
Let $\catDeltaD_r$, $r\geq 1$, be crossed simplicial category 
where the automorphism group of $[n]$ is the dihedral group of order $2r(n+1)$. 
We name the preferred generators $\rho_n$ and $\tau_n$, where $\rho_n^2=\tau_n^{r(n+1)}=\id$ and $\rho_n\tau_n=\tau_n^{-1}\rho_n$,
and introduce the same relations as in definition~\ref{Def:catDeltaD}.
\end{Def}

We now give names to these crossed simplicial categories, and call $\catDeltaT$, $\catDeltaC$, $\catDeltaD$, $\catDeltaC_r$ and $\catDeltaD_r$
the \textit{involutive simplicial category}, the \textit{cyclic category}, the \textit{dihedral category}, 
the \textit{$r$-cyclic category} and the \textit{$r$-dihedral category} respectively.
Notice that $\catDeltaC_1=\catDeltaC$ and $\catDeltaD_1=\catDeltaD$. We see that $\catDeltaC_r$ is a subcategory of $\catDeltaD_r$,
and that $\catDeltaT$ is a subcategory of $\catDeltaD_r$, for any $r\geq1$.

Our reason for introducing crossed simplicial categories is to 
study $G_{\bullet}$-objects in some category $\mathscr{C}$:

\begin{Def}
Let $\catDeltaG$ be a crossed simplicial category
and $\mathscr{C}$ any category. 
A \textit{$G_{\bullet}$-object} in $\mathscr{C}$
is a functor $\catDeltaG^{\op}\rightarrow\mathscr{C}$. A \textit{$G_{\bullet}$-map}
between $G_{\bullet}$-objects is a natural transformation of functors.
\end{Def}

If $\catDeltaG$ is one of the crossed simplicial categories above and $\mathscr{C}$ is $\Top$, the category of 
(compactly generated) spaces,
then we call $G_{\bullet}$-objects for \textit{involutive simplicial spaces}, 
\textit{cyclic spaces}, \textit{dihedral spaces}, \textit{$r$-cyclic spaces} and \textit{$r$-dihedral spaces} accordingly,
and similarly we replace the word ``spaces'' by ``sets'' when $\mathscr{C}=\Ens$, the category of sets.

Given an $r$-dihedral space $X_{\bullet}$ we have the following notation: The map induced by
$\delta_i$ is denoted by $d_i:X_{n}\rightarrow X_{n-1}$ and called the \textit{$i$'th face map}.
The map induced by $\sigma_i$ is denoted $s_i:X_{n}\rightarrow X_{n+1}$ and called the 
\textit{$i$'th degeneracy map}. The map induced by $\rho_n$ is denoted by $r_n:X_n\rightarrow X_n$
and called the \textit{involutive operator}. And the map induced by $\tau_n$ is denoted by
$t_n:X_n\rightarrow X_n$ and called the \textit{cyclic operator}.

For an $r$-cyclic space, we use the same notation and terminology, but in this case there are
no involutive operators. Analogously, there are no cyclic operators for involutive simplicial spaces.

\subsection{Geometric realization of $G_{\bullet}$-spaces}

We now turn toward the geometric realization of
$G_{\bullet}$-spaces. Via the inclusion $j:\catDelta\rightarrow\catDeltaG$ we associate to
any $G_{\bullet}$-space $X_{\bullet}$ its underlying simplicial space $j^*X_{\bullet}$, which is given as 
the composition
$\catDelta^{\op}\xrightarrow{j}\catDeltaG^{\op}\xrightarrow{X_{\bullet}}\Top$.
And we define:

\begin{Def}
The geometric realization of a $G_{\bullet}$-space $X_{\bullet}$ is the geometric realization
of its underlying simplicial space $j^*X_{\bullet}$.
\end{Def}

From the article~\cite{FiedorowiczLoday:91} we now summarize 
results about the geometric realization of a $G_{\bullet}$-space.

\begin{Thm}\label{Thm:MainThmOfCrossedSimplCat}
Let $\catDeltaG$ be a crossed simplicial category, and $X_{\bullet}$ a simplicial space.
We have:
\begin{itemize}
\item[] The functor $j^*$ from $G_{\bullet}$-spaces to simplicial spaces
has a left adjoint, denoted by $F_G$, and there are projection maps
$p_1:|F_G(X_{\bullet})|\rightarrow |G_{\bullet}|$ and $p_2:|F_G(X_{\bullet})|\rightarrow |X_{\bullet}|$.
\item[] The map $(p_1,p_2):|F_G(X_{\bullet})|\rightarrow |G_{\bullet}|\times |X_{\bullet}|$
is a homeomorphism.
\item[] For any simplicial map $f_{\bullet}:X_{\bullet}\rightarrow Y_{\bullet}$ the following diagrams commute:
\[
\begin{CD}
|F_G X_{\bullet}| @>{|F_G f_{\bullet}|}>> |F_G Y_{\bullet}|\\
@V{p_2}VV @VV{p_2}V\\
|X_{\bullet}| @>{|f_{\bullet}|}>> |Y_{\bullet}|
\end{CD}\quad\text{and}\quad
\begin{CD}
|F_G X_{\bullet}| @>{|F_G f_{\bullet}|}>> |F_G Y_{\bullet}|\\
@V{p_1}VV @VV{p_1}V\\
|G_{\bullet}| @= |G_{\bullet}|
\end{CD}\quad.
\]
\item[] Since $F_G$ is a left adjoint, there are canonical natural transformations
$\mu_{\bullet}:F_G(F_G(X_{\bullet}))\rightarrow F_G(X_{\bullet})$ and $\iota_{\bullet}:X_{\bullet}\rightarrow F_G(X_{\bullet})$.
And the following diagrams commute:
\[
\begin{CD}
|F_G(F_G( X_{\bullet}))| @>{|\mu_{\bullet}|}>> |F_G (X_{\bullet})|\\
@V{p_2}VV @VV{p_2}V\\
|F_G(X_{\bullet})| @>{p_2}>> |X_{\bullet}|
\end{CD}\quad\text{and}\quad
\begin{CD}
|X_{\bullet}| @>{|\iota_{\bullet}|}>> |F_G X_{\bullet}|\\
@V{=}VV @VV{p_2}V\\
|X_{\bullet}| @= |X_{\bullet}|
\end{CD}\quad.
\]
\item[] There is a canonical isomorphism $G_{\bullet}\cong F_G(*)$ and the composition 
$|G_{\bullet}|\cong |F_G(*)|\xrightarrow{p_1}|G_{\bullet}|$ is the identity.
\item[] Let $1$ denote the point in $|G_{\bullet}|$ determined by the unit in $G_0$. The following diagram commutes:
\[
\begin{CD}
|X_{\bullet}| @>{|\iota_{\bullet}|}>> |F_G(X_{\bullet})|\\
@VVV @VV{p_1}V\\
\{1\} @>>> |G_{\bullet}|\\
\end{CD}\quad.
\]
\item[] $|G_{\bullet}|$ is a topological group.
\item[] If $X_{\bullet}$ is a $G_{\bullet}$-space, then there is an 
induced action $|G_{\bullet}|\times|X_{\bullet}|\rightarrow|X_{\bullet}|$.
\item[] $(p_1,p_2):|F_G(X_{\bullet})|\rightarrow |G_{\bullet}|\times |X_{\bullet}|$ is an equivariant homeomorphism.
\item[] For every $n$ there is an inclusion of $G_n$ in $|G_{\bullet}|$ as a discrete subgroup.
\end{itemize}
\end{Thm}

For a proof see propositions~4.4,~5.1,~5.3 and~5.13 in~\cite{FiedorowiczLoday:91}.

\begin{Rem}\label{Rem:sescsg}
Chasing Fiedorowicz and Loday's proof of theorem~\ref{Thm:MainThmOfCrossedSimplCat} above, 
it is not hard to see that
all results are natural with respect to a morphism $\catDeltaG\rightarrow\catDeltaG'$
of crossed simplicial categories. In particular we get an induced homomorphism of topological groups
$|G_{\bullet}|\rightarrow|G_{\bullet}'|$. Furthermore, it is possible to consider short exact sequences
of crossed simplicial categories. It is more convenient to write such a sequence in terms of
the corresponding crossed simplicial groups. The sequence
\[
0\rightarrow G''_{\bullet}\rightarrow G_{\bullet}\rightarrow G'_{\bullet}\rightarrow 0
\]
is short exact if the evaluation at each $[n]$ is. Taking the geometric realization
one gets a sequence
\[
|G''_{\bullet}|\rightarrow |G_{\bullet}|\xrightarrow{f} |G'_{\bullet}|\quad,
\]
which an extension of topological groups. 
\end{Rem}

Let us now determine what the group $|G_{\bullet}|$ is for our crossed simplicial categories.

\begin{Exa}
Consider the involutive simplicial category, $\catDeltaT$. The automorphism group, $G_n^{\op}$,
of $[n]$ in $\catDeltaT$ is isomorphic to $\Z/2$. 
Recall that $G_{\bullet}$ is a 
simplicial set, the face and degeneracy maps are given by the formula in remark~\ref{Rem:Gissimplset}. 
The degeneracy map $s_0$ is always injective. By counting the order of $G_n$, we immediately see
that the only non-degenerate simplices lie in degree $0$. Hence, we have that
$|G_{\bullet}|$, in this case, is the group $\Z/2$. This means that the geometric realization of an involutive
simplicial space is a topological space with involution.
\end{Exa}

\begin{Exa}
Next consider the cyclic category, $\catDeltaC$. 
Using the formula from remark~\ref{Rem:Gissimplset}, we find that the non-degenerate
simplices are $\tau_0\in G_0$ and $\tau_1\in G_1$. Hence, $|G_{\bullet}|\cong S^1$.
We now determine the group structure.
A theorem by von Neumann says that any compact, locally Euclidean topological group
is a Lie group, see theorem~57 in~\cite{Pontrjagin:39}. 
The theory of Lie groups now tells us that the only topological group structure on $S^1$
is the ordinary group structure.
\end{Exa}

\begin{Exa}
Now look at the $r$-cyclic category, $\catDeltaC_r$. Let $G_{\bullet}$
be the associated crossed simplicial group. 
To determine $|G_{\bullet}|$ as a topological space, we find the non-degenerate simplices.
The $0$-simplices, $G_0=C_r$, are non-degenerate. Recall from remark~\ref{Rem:Gissimplset}
the formula defining the simplicial structure on $G_{\bullet}$. The relation
\[
\tau_0\sigma_0=\sigma_0\tau_1^2
\]
implies that $s_0(\tau_0^i)=\tau_1^{2i}$. Hence, $\tau_1$, $\tau_1^3$,$\ldots$, $\tau_1^{2r-1}$
are the non-degenerate simplices in $G_1$. Playing with the relations in $\catDeltaC_r$,
we see that there are no more non-degenerate simplices. Furthermore, we have that
$d_0(\tau_1^{2i-1})=\tau_0^{i-1}$ and $d_1(\tau_1^{2i-1})=\tau_0^{i}$. Hence,
$|G_{\bullet}|\cong S^1$. And $S^1$ has a unique structure as a topological group.
\end{Exa}

\begin{Exa}
Let us now study the $r$-dihedral category, $\catDeltaD_r$.
We can use the definition of the category and the formula from remark~\ref{Rem:Gissimplset}
to determine the simplicial structure on the associated simplicial group $G_{\bullet}$.
Finding non-degenerate simplices and calculating the face maps, we see that
\[
|G_{\bullet}|\cong S^1\times \Z/2
\]
as topological spaces. 
Hence there are two possibilities for the group structure on $|G_{\bullet}|$:
it is isomorphic either to $S^1\times\Z/2$ or $O(2)$. 
By the last statement of theorem~\ref{Thm:MainThmOfCrossedSimplCat},
$|G_{\bullet}|$ contains dihedral subgroups.
This excludes $S^1\times\Z/2$, so $|G_{\bullet}|=O(2)$.
\end{Exa}

The theorem~\ref{Thm:MainThmOfCrossedSimplCat} above tells us that the geometric 
realization of a $G_{\bullet}$-space has a $|G_{\bullet}|$ action.
However, it is usually the case that the action takes one out of the topological simplex one starts in. 
In particular, the $q$'th space of the filtration $F_q|X_{\bullet}|$ is seldom
$|G_{\bullet}|$-equivariant.
In many situations it would be easier if the action stayed inside the topological simplices and
the filtration had $|G_{\bullet}|$-action.
We can achieve this by defining the topological $|G_{\bullet}|$-simplices according to
the crossed simplicial category under consideration.

Let $\catDeltaG$ be a crossed simplicial category. Consider the representable functors
\[
\catDeltaG(-,[n]):\catDeltaG^{\op}\rightarrow \Ens\quad.
\]

\begin{Def}
Let $\Delta G:\catDeltaG\rightarrow \Top$ be the functor with $\Delta G^n=|\catDeltaG(-,[n])|$.
The \textit{topological $|G_{\bullet}|$-simplices} are the spaces $\Delta G^n$, $n\geq0$.
\end{Def}

Observe that the representable functor $\catDeltaG(-,[n])$ is $F_G(\catDelta(-,[n]))$, hence we have
homeomorphisms $\Delta G^n=|\catDeltaG(-,[n])|\cong|G_{\bullet}|\times \Delta^n$. 
So the $|G_{\bullet}|$-action does not take points outside $\Delta G^n$.

Using the functor $\Delta G^{\bullet}$ we can now define a geometric realization of 
$G_{\bullet}$-spaces $X$ given by:
\[
|X|_{\Delta G}=\int^{[n]\in\catDeltaG}X_n\times \Delta G^n\quad.
\]
This space is isomorphic to the quotient of $\coprod X_n\times \Delta G^n$ 
where we identify $(x,\phi(\mathbf{t}))$ with $(\phi^*x,\mathbf{t})$
for all morphisms $\phi$ in $\catDeltaG$.

\begin{Lem}\label{Lem:equalityofrealizations}
There is a natural homeomorphism $|X|_{\Delta G}\cong|X|$ for $\catDeltaG^{\op}$-spaces $X$.
\end{Lem}

\begin{proof}
Consider the functor 
$F:(\catDelta\times\catDeltaG)^{\op}\times (\catDelta\times\catDeltaG)\rightarrow \Top$
given by
\[
F([n_o],[m_o],[n],[m])=X_{m_o}\times\catDeltaG(j(n_o),m)\times\Delta^n\quad.
\]
We have that
\[
\int^{[n]\in\catDelta}F([n],[m_o],[n],[m])\cong X_{m_o}\times\Delta G^{m}
\]
and
\[
\int^{[m]\in\catDeltaG}F([n_o],[m],[n],[m])\cong X(j(n_o))\times\Delta^{n}\quad.
\]
The result now follows from the Fubini theorem for coends, see \S{}IX.8 in~\cite{MacLane:98}:
coends can be interchanged.
\end{proof}

To achieve full control of the $|G_{\bullet}|$-action on $|X_{\bullet}|$, it suffices to have
an explicit description of the functor $\Delta G^{\bullet}$. This description should
specify the map $\Delta G^n\rightarrow\Delta G^m$ induced by a morphism $\phi:[n]\rightarrow[m]$ in $\catDeltaG$.
In the case $\catDeltaC$, this description is given implicitly in proposition~2.7 in~\cite{DwyerHopkinsKan:85},
and more explicitly in theorem~3.4 in~\cite{Jones:87}. For the $r$-cyclic case a formula is given by
lemma~1.6 in~\cite{BokstedtHsiangMadsen:93}, and by formula~(2.1.3) in~\cite{Madsen:94}. In general it is just
a question about writing out the equivariant homeomorphism 
$(p_1,p_2):|F_G(\Delta^n_{\bullet})|\rightarrow |G_{\bullet}|\times |\Delta^n_{\bullet}|$
from theorem~\ref{Thm:MainThmOfCrossedSimplCat}. Here $\Delta^n_{\bullet}$ is the simplicial $n$-simplex
$\catDelta(-,[n])$.

Explicitly we have in our cases:

\begin{Exa}
For the involutive simplicial category $\catDeltaT$ we define the functor $\Delta T^{\bullet}$ by
sending $[n]$ to $\Z/2\times\Delta^n$. We write $\Z/2$ multiplicatively. The generators of $\catDeltaT$
induce the following maps:
\begin{align*}
\delta_i(\epsilon;t_0,\ldots,t_n)&= (\epsilon; t_0,\ldots, t_{i-1},0,t_i,\ldots,t_n)\quad,\\
\sigma_i(\epsilon;t_0,\ldots,t_n)&= (\epsilon; t_0,\ldots, t_{i-1},t_i+t_{i+1},t_{i+2},\ldots,t_n)\quad\text{, and}\\
\rho_n(\epsilon;t_0,\ldots,t_n)&= (-\epsilon; t_n,t_{n-1},\ldots, t_1,t_0)\quad.
\end{align*}
\end{Exa}

\begin{Exa}
For the $r$-cyclic category $\catDeltaC_r$ we define the functor $\Delta C_r^{\bullet}$ by
sending $[n]$ to $S^1\times\Delta^n$. We identify $S^1$ with the quotient $\R/\Z$. 
The generators of $\catDeltaC_r$
induce the following maps:
\begin{align*}
\delta_i(\theta;t_0,\ldots,t_n)&= (\theta; t_0,\ldots, t_{i-1},0,t_i,\ldots,t_n)\quad,\\
\sigma_i(\theta;t_0,\ldots,t_n)&= (\theta; t_0,\ldots, t_{i-1},t_i+t_{i+1},t_{i+2},\ldots,t_n)\quad\text{, and}\\
\tau_n(\theta;t_0,\ldots,t_n)&= (\theta-\frac{1}{r}t_0; t_1,t_{2},\ldots, t_n,t_0)\quad.
\end{align*}
\end{Exa}

\begin{Exa}\label{Exa:dihedraltopsimpl}
For the $r$-dihedral category $\catDeltaD_r$ we define the functor $\Delta D_r^{\bullet}$ by
sending $[n]$ to $O(2)\times\Delta^n$. $O(2)$ is the space of orthogonal $2\times 2$-matrices.
For $t\in\R/\Z$ let $R(t)$ denote the rotation matrix 
$\begin{pmatrix}\cos(2\pi t)&\sin (2\pi t)\\ -\sin(2\pi t)& \cos(2\pi t)\end{pmatrix}$, and let
$T$ be the matrix $\begin{pmatrix} 0 & 1\\ 1& 0\end{pmatrix}$.
The generators of $\catDeltaD_r$
induce the following maps:
\begin{align*}
\delta_i(M;t_0,\ldots,t_n)&= (M; t_0,\ldots, t_{i-1},0,t_i,\ldots,t_n)\quad,\\
\sigma_i(M;t_0,\ldots,t_n)&= (M; t_0,\ldots, t_{i-1},t_i+t_{i+1},t_{i+2},\ldots,t_n)\quad,\\
\tau_n(M;t_0,\ldots,t_n)&= (M R(-\frac{1}{r}t_0); t_1,t_{2},\ldots, t_n,t_0)\quad\text{, and}\\
\rho_n(M;t_0,\ldots,t_n)&= (M T;t_n,t_{n-1},\ldots, t_1,t_0)\quad.
\end{align*}
\end{Exa}

\subsection{Filtering the geometric realization}

Similar to the constructions~\ref{Constr:presimplfilt} and~\ref{Constr:filtergeomrealization},
we now design a filtration of $|X_{\bullet}|$, when $X_{\bullet}$ is a $G_{\bullet}$-space.
This filtration is $|G_{\bullet}|$-equivariant.

\begin{Constr}\label{Constr:filtercrossedsimplreal}
Let $\catDeltaG$ be a crossed simplicial category and $X_{\bullet}$ a $G_{\bullet}$-space.
The drawback of using the filtration above to study $|X_{\bullet}|$ is that
$F_q|X_{\bullet}|$ has no $|G_{\bullet}|$ action. Therefore we define another
filtration $F_q^G|X_{\bullet}|$. 
Recall that $|X_{\bullet}|$ can be described as the quotient of
$\coprod X_n\times\Delta G^n$ where we identify $(x,\phi(\mathbf{t}))$ 
with $(\phi^* x,\mathbf{t})$ for all morphisms $\phi$ in $\catDeltaG$.
Define $F^G_q|X_{\bullet}|$ to be the image of  $\coprod_{n\leq q} X_n\times \Delta G^n$.
We define the \textit{$G_{\bullet}$-degenerate simplices} of 
$X_{q}$ to be the subspace $s^GX_{q-1}$ consisting of all points which lie in the image of
some map $\phi^*:X_{q-1}\rightarrow X_q$, $\phi\in\catDeltaG([q],[q-1])$.
Recall that the opposite group of $G_q$ is the automorphisms of $[q]$ in $\catDeltaG$. 
Hence $X_q$ and $sX_{q-1}$ have $G_q$ actions, while 
$\Delta G^q$  and $\partial \Delta G^q$ have $G_q^{\op}$ actions.
Let $X_q\times_{G_q} \Delta G^q$ denote the quotient of the product
where we have identified $(gx,\mb{t})$ with $(x,g^*\mb{t})$ for every $g$ in $G_q$.
We now have a pushout diagram
\[
\begin{CD}
X_q\times_{G_q}\partial\Delta G^q\cup sX_{q-1}\times_{G_q}\Delta G^q @>>> F_{q-1}^G|X_{\bullet}|\\
@V{i}VV @VVV\\
X_q\times_{G_q} \Delta G^q @>>> F_q^G|X_{\bullet}|
\end{CD}\quad.
\]
\end{Constr}

\begin{Rem}
Here is a warning: In general it is not true that natural map
$X_0\times_{G_0}|G_{\bullet}|\rightarrow F_0^G|X_{\bullet}|$ is an homeomorphism,
but it is always an equivariant quotient map. 
\end{Rem}

\subsection{Edgewise subdivision}

Above we have seen that both cyclic and $r$-cyclic spaces yield $S^1$-spaces
after geometric realization. Similarly both dihedral and $r$-dihedral spaces
realize to $O(2)$-spaces. So why do we bother with the $r$-cyclic and $r$-dihedral
categories? Observe that neither the $S^1$- nor the $O(2)$-action is simplicial. 
Let $C$ be a finite cyclic group. Notice that $C$ embeds as
a normal subgroup of both $S^1$ and $O(2)$. The answer to the question is that
$C$-fixed points can be studied simplicially whenever the order of $C$ divides $r$.

After making precise the observations above, we shall define the $c$'th
edgewise subdivision, $c\geq 1$.
This is a functor $\sd_c$ from $r$-cyclic spaces to $rc$-cyclic spaces,
and similarly from $r$-dihedral spaces to $rc$-dihedral spaces.
The edgewise subdivisions come with natural equivariant homeomorphisms 
$D_c:|\sd_c X_{\bullet}|\rightarrow |X_{\bullet}|$.
In particular we can replace a cyclic space
with an $r$-cyclic space for the purpose of studying its restricted $C_r$-action.

Let $C$ be a finite cyclic group $C$ of order $c$. 
Recall from example~\ref{Exa:dihedraltopsimpl} that $R(t)\in O(2)$ denotes a rotation by $2\pi t$, 
while $T\in O(2)$ is a reflection. We identify $C$ as the normal subgroup of $O(2)$
generated by $R(\frac{1}{c})$.
Now we construct homomorphisms
\[
\rho_C: O(2)\rightarrow O(2)/C
\]
by letting $\rho_C(R(t))=R(\frac{t}{c})$ and $\rho_C(T)=T$. Observe that $\rho_C$
is an isomorphism. The restriction of $\rho_C$ to $S^1$ is the ``$c$-th root map''
$S^1\xrightarrow{\cong}S^1/C$. 

Two basic facts are: The $C$-fixed point space of an $O(2)$-space $Y$ is an $O(2)/C$-space $Y^C$,
and an $O(2)/C$-space $Z$ can be viewed as an $O(2)$-space $\rho_C^*Z$ via the isomorphism
$\rho_C$.

After these preliminaries we show:

\begin{Prop}
Assume that $X_{\bullet}$ is an $r$-dihedral space and $C$ a finite cyclic group
of order $c$. Assume that $c$ divides $r$ and let $cs=r$. 
Each $X_n$ has a $C$-action and $X_{\bullet}^C$
is an $s$-dihedral space. Furthermore, there is a natural $O(2)$-equivariant homeomorphism
\[
\rho_C^*|X_{\bullet}|^C\cong |X_{\bullet}^C|\quad.
\]
A similar result holds for $r$-cyclic spaces.
\end{Prop}

\begin{proof}
The $C$ action on $X_n$ is given by the map $t_n^{s(n+1)}$. 
Observe that all the operators $d_i$, $s_i$, $t_n$ and $r_n$
preserve $C$-fixed points. Hence, $X_{\bullet}^C$ is an
$r$-dihedral space. But since $t_n^{s(n+1)}$
is the identity when restricted to $X_n^C$, we see that 
$X_{\bullet}^C$ satisfies the identities for an $s$-dihedral space.

To define the natural $O(2)$-homeomorphism we use the filtration
from construction~\ref{Constr:filtercrossedsimplreal}. Assume inductively
that we have an $O(2)$-homeomorphism
\[
\rho_C^*F_{n-1}^{\catDeltaD_r}|X_{\bullet}|^C\cong F_{n-1}^{\catDeltaD_s}|X_{\bullet}^C|\quad.
\]

Recall that the automorphism group
of $[n]$ in $\catDeltaD_r$ is the dihedral group $D_{2r(n+1)}$ of order $2r(n+1)$.
If $Y$ is a $D_{2r(n+1)}$-space, then we may form the induced $O(2)$-space
$Y\times_{D_{2r(n+1)}} O(2)$. It is a basic fact about induced $O(2)$-spaces
and $C$-fixed points, compare lemma~\ref{Lem:fixedpointsofinducedspaces},
that 
\[
\left(Y\times_{D_{2r(n+1)}} O(2)\right)^C\cong Y^C\times_{D_{2r(n+1)}/C} O(2)/C\quad.
\]

For the induction step we inspect the $n$-simplices, and calculate:
\begin{align*}
\rho_C^*\left(X_n\times_{D_{2r(n+1)}}\Delta D_r^n\right)^C
&\cong \rho_C^*\left( (X_n\times\Delta^n)\times_{D_{2r(n+1)}} O(2)\right)^C\\
&\cong \rho_C^*\left( (X_n\times\Delta^n)^C\times_{D_{2r(n+1)}/C} O(2)/C\right)\\
&\cong (X_n\times\Delta^n)^C\times_{D_{2s(n+1)}}\rho_C^{-1}(O(2)/C)\\
&\cong (X_n^C\times\Delta^n)\times_{D_{2s(n+1)}}O(2)\\
&\cong X_n\times_{D_{2s(n+1)}}\Delta D_s^n\quad.
\end{align*}
Similarly, we have an $O(2)$-equivariant homeomorphism for the degenerate points. And these
$O(2)$-homeomorphisms fit into a diagram
\scriptsize
\[
\begin{CD}
\rho_C^*\left(X_n\times_{D_{2r(n+1)}}\Delta D_r^n\right)^C
&\leftarrow& \rho_C^*\left(X_n\times_{D_{2r(n+1)}}\partial\Delta D_r^n\cup sX_{n-1}\times_{D_{2r(n+1)}}\Delta D_r^n\right)^C
&\rightarrow& \rho_C^*F_{n-1}^{\catDeltaD_r}|X_{\bullet}|^C\\
@V{\cong}VV @V{\cong}VV @VV{\cong}V\\
X_n\times_{D_{2s(n+1)}}\Delta D_s^n
&\leftarrow& X_n^C\times_{D_{2s(n+1)}}\partial\Delta D_s^n\cup sX_{n-1}^C\times_{D_{2s(n+1)}}\Delta D_s^n
&\rightarrow& F_{n-1}^{\catDeltaD_s}|X_{\bullet}^C|
\end{CD}\quad.
\]
\normalsize
By construction~\ref{Constr:filtercrossedsimplreal} we see that the map
of the row-wise pushouts is $F_{n}^{\catDeltaD_s}|X_{\bullet}^C|\cong \rho_C^*F_{n}^{\catDeltaD_r}|X_{\bullet}|^C$.

The statement for $r$-cyclic spaces is proved similarly.
\end{proof}

We now define the \textit{edgewise subdivision functor} $\sd_c:\catDeltaD_{rc}\rightarrow \catDeltaD_r$. 
The idea behind $\sd_c$ is to send the ordered set $[q]$ to the disjoint union of $c$ copies of $[q]$:
\[
\sd_c[q]=[q]\amalg\cdots\amalg[q]=[c(q+1)-1]\quad.
\]
This yields the following formulas for $\sd_c$ of the generators in the dihedral case:
\begin{align*}
\sd_c(\delta_i)&=\delta_{i+(c-1)(q+1)}\cdots\delta_{i+(q+1)}\delta_i\quad,\\
\sd_c(\sigma_i)&=\sigma_i\sigma_{i+(q+2)}\cdots\sigma_{i+(c-1)(q+2)}\quad,\\
\sd_c(\tau_q)&=\tau_{c(q+1)-1}\quad\text{, and}\\
\sd_c(\rho_q)&=\rho_{c(q+1)-1}\quad.
\end{align*}
Observe that $\sd_c$ restricts to functors $\catDeltaC_{rc}\rightarrow \catDeltaC_r$,
$\catDeltaT\rightarrow \catDeltaT$ and $\catDelta\rightarrow \catDelta$.

\begin{Def}
Let $X_{\bullet}$ be an $r$-dihedral space. Its \textit{$c$'th edgewise subdivision}, $\sd_cX_{\bullet}$ is the
composition $\catDeltaD_{rc}^{\op}\xrightarrow{\sd_c}\catDeltaD_r^{\op}\xrightarrow{X_{\bullet}}\Top$.
Similarly, we also define the $c$'th edgewise subdivision of $r$-cyclic, involutive simplicial and
simplicial spaces.
\end{Def}

To compare the geometric realization of $\sd_cX_{\bullet}$ and $X_{\bullet}$, we first
define a diagonal map from the topological $rc$-dihedral $q$-simplex $\Delta D_{rc}^q$
to the topological $r$-dihedral $(c(q+1)-1)$-simplex $\Delta D_r^q$. This map is given by
\[
(M;t_0,\ldots,t_q)\mapsto 
(M;\frac{1}{c}t_0,\ldots,\frac{1}{c}t_q,\frac{1}{c}t_0,\ldots,\frac{1}{c}t_q,\ldots,\frac{1}{c}t_0,\ldots,\frac{1}{c}t_q)
\quad.
\]
This map is $O(2)$-equivariant.
Varying $q$, we get a natural transformation $\Delta D_{rc}^{\bullet}\rightarrow\Delta D_{r}^{\bullet}\circ\sd_c$.
Using a trick with coends we define a natural $O(2)$-map $D_c:|\sd_cX_{\bullet}|\rightarrow|X_{\bullet}|$.
Consider
\[
\int^{[p]\in\catDeltaD_{rc}}\int^{[q]\in\catDeltaD_r} X_q\times\catDeltaD_{r}(\sd_c[p],[q])\times\Delta D_{rc}^p\quad.
\]
Observe that the evaluation
\[
\int^{[q]\in\catDeltaD_r} X_q\times\catDeltaD_{r}(\sd_c[p],[q])\rightarrow (\sd_c X)_p
\]
is a homeomorphism. (The identity in $\catDeltaD_{r}(\sd_c[p],\sd_c[p])$ gives an inverse map.)
It follows that the double coend above equals
\[
\int^{[p]\in\catDeltaD_{rc}} (\sd_cX)_p\times\Delta D_{rc}^p=|\sd_cX_{\bullet}|\quad.
\]
On the other hand, by the Fubini theorem for coends, we can consider the coend over $[p]\in\catDeltaD_{rc}$ first.
Via the diagonal map given above, we get a natural $O(2)$-map
\begin{align*}
\int^{[q]\in\catDeltaD_r}&\int^{[p]\in\catDeltaD_{rc}} X_q\times\catDeltaD_{r}(\sd_c[p],[q])\times\Delta D_{rc}^p\\
&\cong \int^{[q]\in\catDeltaD_r}X_q\times
\left(\int^{[p]\in\catDeltaD_{rc}}\catDeltaD_{r}(\sd_c[p],[q])\times\Delta D_{rc}^p\right)\\
&\xrightarrow{\text{diagonal}} \int^{[q]\in\catDeltaD_r}X_q\times
\left(\int^{[p]\in\catDeltaD_{rc}}\catDeltaD_{r}(\sd_c[p],[q])\times\Delta D_{r}^{\sd_c[p]}\right)\\
&\xrightarrow{\text{evaluate}} \int^{[q]\in\catDeltaD_r}X_q\times\Delta D_r^q\\
&=|X_{\bullet}|\quad.
\end{align*}
Putting this together we see that the diagonal map on topological simplices gives a natural $O(2)$-map
\[
D_c:|\sd_cX_{\bullet}|\rightarrow|X_{\bullet}|\quad.
\]
Similarly, for the cyclic, the involutive simplicial and the simplicial categories we have a
natural $S^1$-map, $\Z/2$-map and map respectively. 

\begin{Prop}
Let $X_{\bullet}$ be an $r$-dihedral space, an $r$-cyclic space, an involutive simplicial space
or a simplicial space.
In all cases, the (equivariant) map $D_c:|\sd_cX_{\bullet}|\rightarrow|X_{\bullet}|$
is a homeomorphism.
\end{Prop}

\begin{proof}
Recall that we can compute the geometric realization either over the crossed simplicial category or
over $\catDelta$. Because both methods yield the same space, lemma~\ref{Lem:equalityofrealizations},
it is enough to inspect the map in the simplicial case. 

The proof for simplicial sets, lemma~1.1 in~\cite{BokstedtHsiangMadsen:93} applies also to
the case of simplicial spaces: One first checks by explicit computation that 
$D_c$ is a homeomorphism when $X_{\bullet}$ is the simplicial $1$-simplex $\catDelta(-,[1])$.
It follows that $D_c$ also is a homeomorphism for products $\catDelta(-,[1])^{\times q}$.
Then it holds for the simplicial $q$-simplex because of the retraction 
$\catDelta(-,[q])\xrightarrow{i}\catDelta(-,[1])^{\times q}\xrightarrow{r}\catDelta(-,[q])$.
Let $\eta_q$ denote the inverse of $D_c:|\sd_c\catDelta(-,[q])|\rightarrow|\catDelta(-,[q])|$.
For general simplicial spaces $X_{\bullet}$ we now define the inverse as follows:
\begin{align*}
|X_{\bullet}| &= \int^{[q]\in\catDelta}X_q\times\Delta^q\\
&= \int^{[q]\in\catDelta}X_q\times|\catDelta(-,[q])|\\
&\xrightarrow{\id\times\eta_q}\int^{[q]\in\catDelta}X_q\times|\sd_c\catDelta(-,[q])|\\
&=\int^{[q]\in\catDelta}X_q\times\left(\int^{[p]\in\catDelta}\catDelta(\sd_c[p],[q])\times\Delta^p\right)\\
&=\int^{[q]\in\catDelta}\int^{[p]\in\catDelta} X_q\times\catDelta(\sd_c[p],[q])\times\Delta^p\\
&=|\sd_cX_{\bullet}|\quad.
\end{align*}
\end{proof}

\section{Homotopy colimits over topological categories}

In this short section we will define the homotopy colimit of a continuous functor
over a topological category. Also, we give a condition
on $F$ such that $\hocolim_{\mathscr{C}}F\rightarrow B\mathscr{C}$ is a $\lambda$-quasi cofibration.

Assume that $\mathscr{C}$ is a small topological category; we have a discrete set of objects,
while for each pair of objects, $a,b\in\mathscr{C}$, we have a topological space $\mathscr{C}(a,b)$
of morphisms from $a$ to $b$.
For continuous functors $F:\mathscr{C}\rightarrow\Top$ we would like to define a homotopy colimit.

\begin{Def}
We define $\hocolim_{\mathscr{C}}F$ as the realization of a simplicial space. Its $q$-simplices are
\[
X_q=\coprod_{a_0,\ldots,a_q\in\mathscr{C}} \mathscr{C}(a_{q-1},a_q)
\times\cdots\times \mathscr{C}(a_{0},a_1)\times F(a_0)\quad.
\]
Face and degeneracy maps are given by
\begin{align*}
d_i(f_{q-1},\ldots,f_0;x)&=\begin{cases}
(f_{q-1},\ldots,f_1;f_0(x))&\text{for $i=0$,}\\
(f_{q-1},\ldots,f_{i+1},f_{i}\circ f_{i-1},f_{i-2},\ldots,f_0;x)&\text{for $0<i<q$,}\\
(f_{q-2},\ldots,f_0;x)&\text{for $i=q$, and}
\end{cases}\\
s_i(f_{q-1},\ldots,f_0;x)&=(f_{q-1},\ldots,f_i,\id_{a_i},f_{i-1},\ldots,f_0;x)\quad.
\end{align*}
\end{Def}

$\hocolim$ is functorial. If $\tau:F\rightarrow F'$ is a natural transformation, then there is an induced map
\[
\hocolim_{\mathscr{C}}F\rightarrow \hocolim_{\mathscr{C}}F'\quad.
\]
Furthermore, if $j:\mathscr{D}\rightarrow\mathscr{C}$ is a functor, then there is an induced map
\[
\hocolim_{\mathscr{D}}j^*F\rightarrow \hocolim_{\mathscr{C}}F\quad,
\]
where $j^*f$ is the composite $f\circ F:\mathscr{D}\rightarrow\Top$.

\begin{Prop}\label{prop:hocolimhomotopy}
If $\tau:j_0\rightarrow j_1$ is a natural transformation between continuous
functors $\mathscr{D}\rightarrow\mathscr{C}$,
then there is a simplicial homotopy between 
\[
\hocolim_{\mathscr{D}}j_0^*F\xrightarrow{(j_0)_*} \hocolim_{\mathscr{C}}F
\] 
and 
\[
\hocolim_{\mathscr{D}}j_0^*F\xrightarrow{\tau_*}\hocolim_{\mathscr{D}}j_1^*F\xrightarrow{(j_1)_*}\hocolim_{\mathscr{C}}F
\]
for any continuous functor $F:\mathscr{C}\rightarrow\Top$.
\end{Prop}

\begin{proof}
We define a simplicial homotopy. It is given by maps
\begin{align*}
h_i:\coprod &\mathscr{D}(b_{q-1},b_q)\times\cdots\times \mathscr{D}(b_{0},b_1)\times F(j_0(b_0))\\
&\rightarrow \coprod \mathscr{C}(a_{q},a_{q+1})\times\mathscr{C}(a_{q-1},a_q)\times\cdots\times \mathscr{C}(a_{0},a_1)\times F(a_0)
\end{align*}
for $0\leq i\leq q$.
To define the $h_i$'s we consider the diagram in $\mathscr{C}$:
\[
\begin{CD}
j_0(b_0)@>{j_0(f_0)}>> j_0(b_1) @>{j_0(f_1)}>> \cdots @>{j_0(f_{q-1})}>> j_0(b_q)\\
@V{\tau_{b_0}}VV @V{\tau_{b_1}}VV &\cdots& @VV{\tau_{b_q}}V\\
j_1(b_0)@>{j_1(f_0)}>> j_1(b_1) @>{j_1(f_1)}>> \cdots @>{j_1(f_{q-1})}>> j_1(b_q)
\end{CD}\quad.
\]
$h_i$ is now given by the formula:
\[
h_i(f_{q-1},\ldots,f_0;x)=(j_1(f_{q-1}),\ldots,j_1(f_i),\tau_{b_i},j_0(f_{i-1}),\ldots,j_0(f_0);x)\quad.
\]
It is easily checked that this is the required simplicial homotopy.
\end{proof}

We now define $\lambda$-quasi fibrations:

\begin{Def}
A map $p:E\rightarrow B$ is a \textit{$\lambda$-quasi fibration} if for any $b\in B$ the induced map
$\pi_i(E,p^{-1}(b))\rightarrow \pi_i (B,b)$ is an isomorphism for $0\leq i< \lambda$ and a surjection
for $i=\lambda$.
\end{Def}

\begin{Prop}\label{prop:gluelqf}
Consider the diagram:
\[
\begin{CD}
E_2@<{F}<< E_0 @>>> E_1\\
@V{p_2}VV @VV{p_0}V @VV{p_1}V\\
B_2@<{f}<< B_0 @>{i}>> B_1
\end{CD}\quad.
\]
Assume that the $p_i$'s are $\lambda$-quasi fibrations with $p_i^{-1}(b)$ path-connected for all $i$ and $b\in B_i$.
If $i$ is a cofibration, the right square pullback, and $p_0^{-1}(b)\rightarrow p_2^{-1}(f(b))$ $\lambda$-connected
for all $b\in B_0$, then the induced map of pushouts $p:E\rightarrow B$ is a $\lambda$-quasi fibration.
\end{Prop}

\begin{proof}
We can assume that $f$ is a cofibration, if not one can replace $B_2$ by the mapping cylinder $M_f$, and
$E_2$ by the pullback $r^*E_2$ over the retraction $r:M_f\rightarrow B_2$.
Moreover, we can assume that $F$ is a cofibration, if not we can replace $E_2$ by $M_F$. Using that $f$ is injective
it follows that $M_F\rightarrow B_2$ is a $\lambda$-quasi fibration.

Now compare the long exact sequences of homotopy groups
for the triples $(E_1,E_0,p_0^{-1}(b))$ and $(B_1,B_0,b)$, where $b\in B_0$.
Since $p_0^{-1}(b)=p_1^{-1}(b)$, remember that the right square is pullback, and using that $p_0$ and $p_1$ are 
$\lambda$-quasi fibrations, we get that $\pi_i(E_1,E_0)\rightarrow \pi_i(B_1,B_0)$ is an isomorphism 
for $0\leq i<\lambda$ and surjective for $i=\lambda$.

Regarding the connectedness of $\pi_i(E_2,p_0^{-1}(b))\rightarrow \pi_i(B_2,b)$, we reason as follows:
Since $p_0^{-1}(b)\rightarrow p_2^{-1}(f(b))$ is $\lambda$-connected, we get that $\pi_i(p_2^{-1}(f(b)),p_0^{-1}(b))$
is the trivial group when $i\leq\lambda$. Now consider the long 
exact sequence of homotopy groups for $(E_2,p_2^{-1}(f(b)),p_0^{-1}(b))$. The homomorphism
$\pi_i(E_2,p_0^{-1}(b))\rightarrow \pi_i(E_2,p_2^{-1}(f(b)))$ is an isomorphism for $i<\lambda$ and 
surjective for $i=\lambda$.
Using that $p_2$ is a $\lambda$-quasi fibration, the composed map
\[
\pi_i(E_2,p_0^{-1}(b))\rightarrow \pi_i(E_2,p_2^{-1}(f(b)))\rightarrow \pi_i(B_2,b)
\]
is also an isomorphism for $i<\lambda$ and surjective for $i=\lambda$.

Comparing the long exact sequences of homotopy groups 
for $(E_2,E_0,p_0^{-1}(b))$ and $(B_2,B_0,b)$, we see that 
$\pi_i(E_2,E_0)\rightarrow \pi_i(B_2,B_0)$ is an isomorphism for $0\leq i<\lambda$ and 
surjective for $i=\lambda$.

Since the maps under consideration are cofibrations, the Mayer-Vietoris property for homotopy groups holds as stated
in~\cite{Hatcher:02} proposition~4K.1. Therefore, we have that $\pi_i(E,E_1)\rightarrow \pi_i(B,B_1)$ 
is an isomorphism for $i\leq \lambda$ and surjective for $i=\lambda$. The same is also true for
$\pi_i(E,E_2)\rightarrow \pi_i(B,B_2)$.

At last we can check whether $p:E\rightarrow B$ is a $\lambda$-quasi fibration. If $b\in B_2$ we
compare the long exact sequences of homotopy groups for $(E,E_2,p^{-1}(b))$ and $(B,B_2,b)$. By the five lemma
we see that $\pi_i(E,p^{-1}(b))\rightarrow \pi_i(B,b)$ is an isomorphism 
for $i\leq \lambda$ and surjective for $i=\lambda$.
When $b\in B_1\smallsetminus B_0$, we compare long exact sequences of homotopy 
groups for $(E,E_1,p^{-1}(b))$ and $(B,B_1,b)$.
The same conclusion holds.
\end{proof}

Observe that for any functor $F:\mathscr{C}\rightarrow\Top$ there is a natural map
\[
\hocolim_{\mathscr{C}}F\rightarrow B\mathscr{C}\quad.
\]
Here $B\mathscr{C}$ is the bar construction (=geometric realization of the nerve).
In some cases, this map is a $\lambda$-quasi fibration:

\begin{Prop}\label{prop:hocolimqf}
If the induced map $F(a)\rightarrow F(b)$ is $\lambda$-connected for all morphisms of $\mathscr{C}$,
$\mathscr{C}$ is well-pointed
and all $F(a)$'s are path-connected, then
\[
\hocolim_{\mathscr{C}}F\rightarrow B\mathscr{C}
\]
is a $\lambda$-quasi fibration.
\end{Prop}

\begin{proof}
As above let $X_{\bullet}$ denote the simplicial space whose realization is $\hocolim_{\mathscr{C}}F$.
Now compare the presimplicial realization with the geometric realization:
\[
\begin{CD}
F(a_0)@>>> \| X_{\bullet} \| @>>> \|B_{\bullet}\mathscr{C} \|\\
@V{=}VV @VV{\simeq}V @VV{\simeq}V\\
F(a_0)@>>> | X_{\bullet} | @>>> |B_{\bullet}\mathscr{C} |\\
\end{CD}\quad.
\]
Here $F(a_0)$ is the fiber over some point $b$ in $\|B_{\bullet}\mathscr{C}\|$.
The fiber over $b$'s image in $|B_{\bullet}\mathscr{C}|$ is identical. This can be seen by inspecting the
definition of the degeneracy maps.

Since $\mathscr{C}$ is well-pointed, it follows that $X_{\bullet}$ and $B_{\bullet}\mathscr{C}$ are
good simplicial spaces. Hence, the vertical maps are weak equivalences. Therefore it is enough to
show that $\| X_{\bullet} \| \rightarrow \|B_{\bullet}\mathscr{C} \|$ is a $\lambda$-quasi fibration.

Following Quillen, we now consider the skeletal filtration of the presimplicial realization.
\[
\begin{CD}
F_{q-1}\|X_{\bullet}\|@<<< X_q\times \partial\Delta^q @>>> X_q\times \Delta^q\\
@VVV @VVV @VVV\\
F_{q-1}\|B_{\bullet}\mathscr{C}\|@<<< B_{q}\mathscr{C}\times \partial\Delta^q @>>> B_{q}\mathscr{C}\times \Delta^q
\end{CD}\quad.
\]
This diagram satisfies the conditions of proposition~\ref{prop:gluelqf}, so the map of pushouts
$F_{q}\|X_{\bullet}\|\rightarrow F_{q}\|B_{\bullet}\mathscr{C}\|$ is a $\lambda$-quasi fibration.

Now the result follows since the direct limit of $\lambda$-quasi fibrations is a $\lambda$-quasi fibration.
\end{proof}

\chapter{Orthogonal spectra}\label{Chapt:ort}

This chapter will introduce the relevant results about orthogonal spectra. The main reference for these
results is the article~\cite{MandellMaySchwedeShipley:01}. The aim of that article is to compare
different constructions of a category of spectra with an associative and commutative smash product.
Their theorem~0.1 says that the categories of $\mathcal{N}$-spectra, symmetric spectra, orthogonal spectra,
and $\mathscr{W}$-spaces are Quillen equivalent. However, the aim of this thesis is
to study involutions on certain ring spectra related to geometry of manifolds, see chapter~\ref{Cha:invoper}.
Therefore we are free to choose the category of spectra most convenient for our purposes.
This is the category of orthogonal spectra, and we will focus on how to work within this
category.
 
Below we will give an exposition of the theory of orthogonal spectra. 
All relevant definitions are included here. For completeness we also reprove some of the results
of~\cite{MandellMaySchwedeShipley:01}. However, there are also new results here: 
We introduce l-cofibrations, definition~\ref{Def:lcof}, in order to study simplicial orthogonal spectra,
propositions~\ref{Prop:sos:comparerealizarions} and~\ref{Prop:piisoofrealiz}.
We consider induced functors, corollary~\ref{Cor:indfunct}. And we construct cofibrant and fibrant 
replacement functors with additional properties, theorem~\ref{Thm:cofrepl} and theorem~\ref{Thm:fibrrepl}
respectively.

We use the convention that topological spaces mean compactly generated spaces (=weak Hausdorff $k$-spaces).
This category satisfies Steenrod's convenient technical properties as defined in~\cite{Steenrod:67}. In addition
the category is closed under the operation of passing to the quotient $X/A$ 
of any closed pair $(X,A)$, and under the 
operation of taking the union of an expanding sequence of closed subspaces.
We refer to \S{}2 of~\cite{McCord:69} for the definition and further properties of
compactly generated spaces.
We let $\Top$ denote this category, and $\Top_*$ based compactly generated spaces.

\section{Basic definitions}

In this section we will define the category of orthogonal spectra, $\mathscr{IS}$.
It is a topological category. 
To define $\mathscr{IS}$ we
introduce the topological category $\mathscr{I}$ of finite dimensional real inner product spaces
and linear isometric isomorphisms. The morphism spaces $\mathscr{I}(V,W)$ are empty when
$V$ and $W$ have different dimensions, and homeomorphic to the orthogonal group $O(n)$ when
$n=\dim V=\dim W$. Direct sum gives $\mathscr{I}$ the structure of a symmetric monoidal category,
and one-point-compactification gives a functor $S$ from $\mathscr{I}$ to compactly generated spaces.

\begin{Def}\label{Def:ortsp}
The category $\mathscr{IS}$ of orthogonal spectra has as its objects
continuous functors $L$ from $\mathscr{I}$
to based compactly generated spaces together with maps
$\sigma:L(V)\wedge S^W\rightarrow L(V\oplus W)$, natural in $V$ and $W$, such that
the composite
\[
L(V)\cong L(V)\wedge S^0\xrightarrow{\sigma}L(V\oplus 0)\cong L(V)
\]
is the identity and $\sigma$ is associative in the sense that the following diagram commutes
\[
\begin{CD}
L(U)\wedge S^V\wedge S^W @>{\sigma\wedge\id}>> L(U\oplus V)\wedge S^W\\
@V{\cong}VV @VV{\sigma}V\\
L(U)\wedge S^{V\oplus W} @>{\sigma}>> L(U\oplus V\oplus W)
\end{CD}.
\]
A map of orthogonal spectra is a natural transformation $f:K\rightarrow L$
of functors such that the following diagram commutes
\[
\begin{CD}
K(V)\wedge S^W @>{\sigma}>> K(V\oplus W)\\
@V{f_V\wedge \id}VV @VV{f_{V\oplus W}}V\\
L(V)\wedge S^W @>{\sigma}>> L(V\oplus W)
\end{CD}\quad.
\]
\end{Def}

We call $\sigma$ the right assembly map. There is also a unique left assembly $\bar{\sigma}$
corresponding to $\sigma$ via the symmetry of $\wedge$ and $\oplus$.

There are several interesting examples of orthogonal spectra. First observe that 
the functor $S$ is an example by letting $\sigma:S^V\wedge S^W\rightarrow S^{V\oplus W}$
be the natural homeomorphism. We call $S$ the sphere spectrum. For based topological spaces
$X$ the suspension spectrum is defined by $V\mapsto X\wedge S^V$. We can also
define Thom spectra by letting $TO(V)$, for an $n$-dimensional $V$, be the 
Thom space of the tautological $n$-plane bundle over the Grassmannian of $n$-planes
in $V\oplus V$.

\subsection{Shift desuspension functors}

There is a shift desuspension functor
$F_V$ from based compactly generated spaces to $\mathscr{IS}$ for any $V$.
It is defined by the formula
\[
(F_V A)(W)=\mathscr{I}(V\oplus\R^d,W)_+\wedge_{O(d)}(A\wedge S^{d})\quad,
\]
where $A$ is a based space and $d=\dim W-\dim V$. We let $(F_V A)(W)=*$ for $\dim W<\dim V$.
The right assembly $\sigma:(F_V A)(W)\wedge S^U\rightarrow (F_V A)(W\oplus U)$ is defined
by choosing an isomorphism $g:\R^{d'}\cong U$, and is well defined since we divide out by
$O(d+d')$ in the definition of $(F_V A)(W\oplus U)$.

$F_V$ is left adjoint to the evaluation at level $V$:
\[
\mathscr{IS}(F_V A,L)\cong \Top_*(A,L(V))\quad,
\]
for all $V$, $A$ and $L$.

\subsection{Notions of equivalence} 

For orthogonal spectra there are two different notions of equivalence:

\begin{Def}
A map $f:K\rightarrow L$ of orthogonal spectra is a \textit{level equivalence}
if for every $V$ the map $f_V:K(V)\rightarrow L(V)$ is a weak equivalence.
\end{Def}

To define the other kind of equivalence we make use of a forgetful functor
$\U$ from $\mathscr{IS}$ to prespectra. (The theory of prespectra can be found in
chapter~II of~\cite{Rudyak:98}.) The $n$'th
space of $\U L$ is $L(\R^n)$ and the suspension map $s_n:(\U L)_n\wedge S^1\rightarrow (\U L)_{n+1}$
comes from the right assembly by identifying $\R^n\oplus\R$ with $\R^{n+1}$.
Recall that the homotopy groups of a prespectrum $X$ is defined as
\[
\pi_q(X)=\colim_{n}\pi_{q+n}(X_{n})\quad.
\]
We now define:

\begin{Def}\label{Def:piiso}
A map $f:K\rightarrow L$ of orthogonal spectra is a \textit{$\pi_*$-isomorphism}
if the underlying map of prespectra $\U f:\U K\rightarrow \U L$ induces an 
isomorphism on all homotopy groups.
\end{Def}

A level equivalence $K\rightarrow L$ induces an isomorphism $\pi_{q+n} (\U K)_n\rightarrow \pi_{q+n} (\U L)_n$
for all $q$ and $n$, thus we have:

\begin{Lem}
Any level equivalence is a $\pi_*$-isomorphism.
\end{Lem}

We also have a notion of $\Omega$-spectra:

\begin{Def}
An orthogonal spectrum $E$ is an $\Omega$-spectrum if the adjoint of $\sigma$,
\[
E(V)\rightarrow\Omega^W E(V\oplus W)\quad,
\]
is a weak equivalence for all $V$ and $W$.
\end{Def}

Notice that $E$ is an $\Omega$-spectrum if and only if $\U E$ is an $\Omega$-prespectrum.

\begin{Rem}
For general diagram spectra, and symmetric spectra in particular, there is a third
notion of equivalence, namely stable equivalence. Let $[L,E]$ denote the set of maps in the level homotopy category.
If $L$ is cofibrant (see section~\ref{sect:cellullar} below for a definition of cofibrant), then
$[L,E]$ is isomorphic to the set of path components of the topological space $\mathscr{IS}(L,E)$. 
We say that $f:K\rightarrow L$ is a stable equivalence
if $f^*:[L,E]\rightarrow [K,E]$ is a bijection for all $\Omega$-spectra $E$.
However, for orthogonal spectra there is no difference between
$\pi_*$-isomorphisms and stable equivalences. See proposition~8.7 in~\cite{MandellMaySchwedeShipley:01}.
\end{Rem}

\subsection{l-cofibrations}\label{subsect:lcof}

There are many conditions on maps $i:A\rightarrow L$ that could be taken as the definition of some kind of cofibration of orthogonal spectra.
In section~\ref{sect:cellullar} below we are going to study q-cofibrations. They depend on cellular techniques;
this is a strong condition. However, in this section we will consider a very weak way of
defining a notion of cofibrancy:

\begin{Def}\label{Def:lcof}
A map $i:A\rightarrow L$ of orthogonal spectra is an \textit{l-cofibration} if for every $V$ the map
\[
A(V)\rightarrow L(V)
\]
is an unbased closed cofibration of topological spaces. We call $L$ \textit{well-pointed} if $*\rightarrow L$
is an l-cofibration. 
\end{Def}

\begin{Rem}\label{Rem:cofibdef}
Recall that any unbased cofibration of topological spaces is a homeomorphism 
onto its image (Theorem~1 in~\cite{Strom:66}).
Therefore, we can always assume without loss of generality that 
any unbased closed cofibration is an inclusion of a closed subspace.
Furthermore, if $i:A\subseteq X$ is the inclusion of a subspace, 
then the following are equivalent ways to define that $i$ is a cofibration: 
\begin{itemize}
\item[] For any map $f:X\rightarrow Y$ and any homotopy 
$\bar{F}:A\times I\rightarrow Y$ with $\bar{F}(a,0)=fi(a)$ for all $a\in A$,
there exists a homotopy $F:X\times I\rightarrow Y$ 
such that $F$ restricts to $\bar{F}$ on $A\times I$ and $F(x,0)=f(x)$ for all $x\in X$.
\item[] The subspace $X\times 0\cup A\times I$ is a 
retract of $X\times I$. (Theorem~2 in~\cite{Strom:68}.)
\item[] There exists a continuous function $\phi:X\rightarrow I$ 
and a homotopy $H:X\times I\rightarrow X$ such that $A\subseteq \phi^{-1}(0)$,
$H(x,0)=x$ for all $x\in X$, $H(a,t)=a$ for all $a\in A$ 
and $t\in I$, and $H(x,t)\in A$ whenever $t>\phi(x)$. (Lemma~4 in~\cite{Strom:68}.)
\end{itemize}
Notice that the subspace topology on $X\times 0\cup A\times I$ does not always 
coincides with the mapping cylinder topology, but in two important cases
these topologies are identical: 1) If $A$ is a closed 
subspace of $X$, or 2) if $A\subseteq X$ is a cofibration.
\end{Rem}

Let us now look at some properties of l-cofibrations of orthogonal spectra:

\begin{Prop}\label{Prop:lcofseqpiiso}
If we are given a map $i$ between sequences of l-cofibrations 
$A_0\rightarrow A_1\rightarrow \cdots$ and $L_0\rightarrow L_1\rightarrow \cdots$ and
each $i_m:A_m\rightarrow L_m$ is a $\pi_*$-iso, then the induced map of colimits $i:A\rightarrow L$
is also a $\pi_*$-iso.  
\end{Prop}

\begin{proof}
Since spheres $S^q$ and disks $D^{q+1}$ are compact, we have
\[
\pi_q A(V)=\colim_m \pi_q A_m(V)\quad,
\]
and similar for $L$. Thus
\[
\pi_q A=\colim_{n,m}\pi_{q+n} A_m(\R^n)\xrightarrow{\cong}\colim_{n,m}\pi_{q+n}L_m(\R^n)=\pi_q L 
\]
is an isomorphism because each $i_m:A_m\rightarrow L_m$ is a $\pi_*$-iso.
\end{proof}

From the category of spaces we immediately inherit union and gluing theorems for l-cofibrations:

\begin{Prop}
If $A\rightarrow L$, $B\rightarrow L$ and $A\cap B\rightarrow L$ are l-cofibrations and $A\cup B\rightarrow L$ an
inclusion, then $A\cup B\rightarrow L$ is also an l-cofibration.
\end{Prop}

\begin{proof}
Notice that intersection and union are level-wise constructions on orthogonal spectra. Now the result follows
directly from the definition of l-cofibration and Lillig's union theorem~\cite{Lillig:73}.
\end{proof}

\begin{Prop}
If we have a diagram
\[
\begin{CD}
B @<<< A @>{i_1}>> L\\
@V{f_2}VV @V{f_0}VV @VV{f_1}V\\
B' @<<< A' @>{i_2}>> L'
\end{CD}
\]
of orthogonal spectra, 
where $i_1$, $i_2$, $f_0$, $f_1$ and $f_2$ are l-cofibrations and the right square is pullback, then the map
of the row-wise pushouts is an l-cofibration.
\end{Prop}

\begin{proof}
Pushouts and pullbacks are level-wise constructions, therefore the result follows from proposition~2.5 in~\cite{Lewis:82}.
\end{proof}

If $X$ is a based space and $L$ an orthogonal spectrum, then we may form the function spectrum
\[
F(X,L)
\]
level-wise. To be precise we let $F(X,L)(V)$ be the space of based maps $X\rightarrow L(V)$. This is again an orthogonal spectrum, 
see example~\ref{exa:lwconstr}. We now apply $F(X,-)$ to an l-cofibration: 

\begin{Prop}\label{Prop:lcofexpthm}
If $C$ is a compact based space, and $i:A\rightarrow L$ is an l-cofibration of orthogonal spectra, then
\[
F(C,A)\rightarrow F(C,L)
\]
is also an l-cofibration.
\end{Prop}

\begin{proof}
Fix $V$. Then we have $H:L(V)\times I\rightarrow L(V)$ and $\phi:L(V)\rightarrow I$ satisfying Str\o{}m's criterion.
Define $\bar{H}:F(C,L(V))\times I\rightarrow F(C,L(V))$ and $\bar{\phi}:F(C,L(V))\rightarrow I$ by
\[
\bar{H}(f,t)(c)=H(f(c),t)\quad\text{and}\quad\bar{\phi}(f)=\sup_{c\in C}\phi f(c)\quad,
\]
for $f:C\rightarrow L(V)$ and $t\in I$. Then $(\bar{H},\bar{\phi})$ shows that $F(C,A)(V)\rightarrow F(C,L)(V)$
is a cofibration.
\end{proof}

\begin{Rem}
One can define h-cofibrations as the maps $i:A\rightarrow L$ having the homotopy extension property, 
see \S{}5 in~\cite{MandellMaySchwedeShipley:01}. These should behave more or less like based cofibrations of spaces.
Therefore, we run into problems if we try to prove union and gluing theorems for h-cofibrations without introducing extra conditions.
\end{Rem}

\subsection{A symmetric monoidal smash product}

The main advantage of orthogonal spectra compared to prespectra 
is the existence of a symmetric monoidal smash product.
To define this we follow~\cite{MandellMaySchwedeShipley:01}. 
Define the category of $\mathscr{I}$-spaces, $\mathscr{I}\Top_*$,
to be functors $\mathscr{I}\rightarrow \Top_*$. It is a 
topological category, the morphisms being the space of natural transformations.

Before defining $\wedge$ on $\mathscr{IS}$, we define the smash product, $\tilde{\wedge}$, of $\mathscr{I}$-spaces.
This is given by
\[
(X\tilde{\wedge} Y)(V)= 
\bigvee_{d_1,d_2}\mathscr{I}(\R^{d_1}\oplus\R^{d_2},V)_+\wedge_{O(d_1)\times O(d_2)}(X(\R^{d_1})\wedge Y(\R^{d_2}))\quad.
\]
If $L$ and $K$ are orthogonal spectra, notice that the 
assembly induces a map of $\mathscr{I}$-spaces $\sigma:L\tilde{\wedge}S\rightarrow L$,
and similarly the left assembly induces $\bar{\sigma}:S\tilde{\wedge}K\rightarrow K$. We now define the smash product, 
$L\wedge K$, of orthogonal spectra by the coequalizer diagram of $\mathscr{I}$-spaces:
\[
L\tilde{\wedge}S\tilde{\wedge}K\overset{\sigma\wedge\id}{\underset{\id\wedge\bar{\sigma}}{\rightrightarrows}}
L\tilde{\wedge}K\rightarrow L\wedge K\quad.
\]
Coequalizers in $\mathscr{I}\Top_*$ are formed level-wise. Recall that $\Top_*$ is
the category of compactly generated spaces. It is cocomplete, and hence
the topology of $(L\wedge K)(V)$ is given as a coequalizer in this category.
As explained in~\cite{MandellMaySchwedeShipley:01}, the smash product is symmetric monoidal.

Having the smash product we define the $\square$ product of maps:

\begin{Def}\label{Def:squareproduct}
Let $f:A\rightarrow L$ and $g:B\rightarrow K$ be maps of orthogonal spectra, then we define $f\square g$
as the map
\[
f\square g: A\wedge K\cup L\wedge B \rightarrow L\wedge K\quad.
\]
\end{Def}

We also have internal function objects. Again we start by defining the internal function object, $\tilde{F}(-,-)$, on
$\mathscr{I}$-spaces. This is given by
\[
\tilde{F}(X,Y)(V)=\mathscr{I}\Top_*(X,Y(V\oplus -))\quad.
\]
And we have an adjunction for $\mathscr{I}$-spaces $X$, $Y$ and $Z$:
\[
\mathscr{I}\Top_*(X\tilde{\wedge}Y,Z)\cong \mathscr{I}\Top_*(X,\tilde{F}(Y,Z))\quad.
\]
If $L$ and $K$ are orthogonal spectra, the assembly induces a map $\sigma^*:\tilde{F}(L,K)\rightarrow\tilde{F}(L\tilde{\wedge}S,K)$.
By the adjunction above there is an evaluation map $\epsilon:\tilde{F}(L,K)\tilde{\wedge}L\rightarrow K$. Now consider the composite
\[
\tilde{F}(L,K)\tilde{\wedge}L\tilde{\wedge}S\xrightarrow{\epsilon\tilde{\wedge}\id} K\tilde{\wedge}S
\xrightarrow{\sigma} K\quad,
\]
let $\omega$ be its adjoint. Define the internal function spectrum, $F(L,K)$, by the equalizer diagram of $\mathscr{I}$-spaces:
\[
F(L,K)\rightarrow \tilde{F}(L,K)\overset{\sigma^*}{\underset{\omega}{\rightrightarrows}}\tilde{F}(L\tilde{\wedge}S,K)\quad.
\]
We immediately get an adjunction for orthogonal spectra $L$, $K$, $X$:
\[
\mathscr{IS}(L\wedge K,X)\cong\mathscr{IS}(L,F(K,X))\quad.
\]

\begin{Lem}
There is an adjunction for the internal hom objects:
\[
F(X\wedge Y,Z)\cong F(X,F(Y,Z))
\]
for $X$, $Y$ and $Z$ orthogonal spectra.
\end{Lem}

\begin{proof}
Assume first that $X$, $Y$ and $Z$ are $\mathscr{I}$-spaces.
Then we check by the definitions that
\[
\tilde{F}(X\tilde{\wedge} Y,Z)(V)\cong \tilde{F}(X,\tilde{F}(Y,Z))(V)
\]
for all $V$.

By the coequalizer defining $\wedge$ and the equalizer defining $F(-,-)$, the
adjunction also holds for internal hom objects in $\mathscr{IS}$.
\end{proof}

\subsection{The external viewpoint}\label{sect:external}

We have presented the symmetric monoidal category $\mathscr{IS}$ of orthogonal spectra. The formal properties
are nice, but when we actually want to do constructions things usually are a bit harder. For example it is not
easy to define a map $L\wedge K\rightarrow X$ directly using the definition of $L\wedge K$. Therefore it is
useful to have alternative viewpoints.

Orthogonal spectra may be described as diagram spaces, see~II.4 in~\cite{MandellMay:02} and
\S{}23 in~\cite{MandellMaySchwedeShipley:01}: There is a topological category $\mathscr{J}$ such that
continuous functors $\mathscr{J}\rightarrow\Top_*$ corresponds to orthogonal spectra. The objects of $\mathscr{J}$
are the same as the objects of $\mathscr{I}$, finite dimensional real inner product spaces $V$. 
Let $\mathscr{E}(V,W)$ be the space of linear isometries $V\hookrightarrow W$. And let $E(V,W)$ consist
of pairs $(f,w)$ where $f:V\rightarrow W$ is a linear isometry and $w\in W$ is orthogonal to $f(V)$.
$E(V,W)$ is a vector bundle over $\mathscr{E}(V,W)$, and we define the space of morphisms $\mathscr{J}(V,W)$
to be the Thom space of $E(V,W)$. (First apply fiber-wise one-point-compactification to $E(V,W)$,
then identify the points at $\infty$.)
Composition
\[
\circ:\mathscr{J}(W,U)\wedge\mathscr{J}(V,W)\rightarrow\mathscr{J}(V,U)
\]
is defined by the formula $(g,u)\circ(f,w)=(g\circ f,g(w)+u)$. The identity of $V$ in $\mathscr{J}$ is represented by $(\id_V,0)$. 
Direct sum gives $\mathscr{J}$ a symmetric monoidal structure: Here
\[
\oplus:\mathscr{J}(V,W)\wedge\mathscr{J}(V',W')\rightarrow \mathscr{J}(V\oplus V,W\oplus W')
\]
is defined by $(f,w)\oplus(f',w')=(f\oplus f',(w,w'))$. Observe that when $V\subseteq W$, we have the identification:
\[
\mathscr{J}(V,W)\cong O(W)_+\wedge_{O(W-V)} S^{W-V}\quad.
\]

\begin{Thm}\label{thm:orthspasdiagr}
The category $\mathscr{IS}$ of orthogonal spectra is isomorphic to the category of $\mathscr{J}$-spaces
as symmetric monoidal categories.
\end{Thm}

\begin{proof}
This is the special case $R=S$ of theorem~2.2 in~\cite{MandellMaySchwedeShipley:01}.
Given an orthogonal spectrum $L$, the corresponding $\mathscr{J}$-space $L'$ is defined by
\[
L'(V)=L(V)
\]
and the map $L'(V)\rightarrow L'(W)$ induced by $(f,w)$ is the composition
\[
L'(V)=L(V)\xrightarrow{x\mapsto(x,w)} L(V)\wedge S^{W-f(V)}\xrightarrow{\sigma} L(V\oplus(W-f(V)))\cong L(W)=L'(W)\quad.
\]
\end{proof}

\begin{Exa}[Level-wise constructions]\label{exa:lwconstr}
Given a continuous endofunctor $F$ on $\Top_*$, we may apply this level-wise to
an orthogonal spectrum $L$. This yields a new orthogonal spectrum $F(L)$. To see this we view $L$ as a functor $\mathscr{J}\rightarrow \Top_*$
and consider the composition
\[
\mathscr{J}\xrightarrow{L}\Top_*\xrightarrow{F}\Top_*\quad.
\]
Examples of such endofunctors are: Based loops $\Omega(-)$, suspension $\Sigma(-)$, 
function spaces $F(X,-)$, the Barratt-Eccles functor $\Gamma^+(-)$.
\end{Exa}

\begin{Exa}[External description of the smash product]\label{exa:extsmashdescr}
Given $\mathscr{J}$-spaces $K$ and $L$
we have an external smash product $\bar{\wedge}$. This produces a functor $\mathscr{J}\times \mathscr{J}\rightarrow\Top_*$
defined by
\[
(K\bar{\wedge}L)(V_1,V_2)=K(V_1)\wedge L(V_2)\quad.
\]
Recall that $\oplus$ is a functor $\mathscr{J}\times \mathscr{J}\rightarrow\mathscr{J}$. Left Kan extension, 
see section~X.4 in~\cite{MacLane:98},
of $K\bar{\wedge}L$ along $\oplus$ gives an internal product. Theorem~\ref{thm:orthspasdiagr} says that this internal product
is equal to $K\wedge L$ defined above. Adjunction for left Kan extensions now says that for orthogonal spectra $K$, $L$ and $X$ 
there is a homeomorphism between
\[
\text{the space of natural transformations}\quad (K\bar{\wedge}L)(V_1,V_2)\rightarrow X(V_1\oplus V_2)
\]
where $V_1,V_2\in\mathscr{J}$ and
\[
\text{the space of maps of orthogonal spectra}\quad K\wedge L\rightarrow X\quad.
\]

This adjunction is useful when defining maps $K\wedge L\rightarrow X$. All we have to do is to provide maps
\[
K(V_1)\wedge L(V_2)\rightarrow X(V_1\oplus V_2)
\]
for all $V_1$ and $V_2$ such that the following diagrams commute
\[
\begin{CD}
S^W\wedge K(V_1)\wedge L(V_2) @>>> S^W\wedge X(V_1\oplus V_2)\\
@V{\bar{\sigma}\wedge\id_L}VV @VV{\bar{\sigma}}V\\
K(W\oplus V_1)\wedge L(V_2) @>>> X(W\oplus V_1\oplus V_2)
\end{CD}
\]
and
\[
\begin{CD}
K(V_1)\wedge L(V_2) \wedge S^W@>>> X(V_1\oplus V_2)\wedge S^W\\
@V{\id_K\wedge\sigma}VV @VV{\sigma}V\\
K(V_1)\wedge L(V_2\oplus W) @>>> X(V_1\oplus V_2\oplus W)
\end{CD}
\]
for all $V_1$, $V_2$ and $W$.
\end{Exa}

\subsection{Orthogonal ring spectra; $S$-algebras}

Having a symmetric smash product $\wedge$ of orthogonal spectra we define orthogonal ring spectra, also called
$S$-algebras, as follows:

\begin{Def}
An orthogonal ring spectrum, or $S$-algebra, is an orthogonal spectrum $L$ together with
maps $\eta:S\rightarrow L$ and $\mu:L\wedge L\rightarrow L$ such that the following diagrams commute:
\[
\begin{CD}
S\wedge L @>{\text{$S$ is the unit for $\wedge$}}>> L\\
@V{\eta\wedge\id}VV @VV{=}V\\
L\wedge L @>{\mu}>> L
\end{CD}\quad,
\]
\[
\begin{CD}
L\wedge S @>{\text{$S$ is the unit for $\wedge$}}>> L\\
@V{\id\wedge\eta}VV @VV{=}V\\
L\wedge L @>{\mu}>> L
\end{CD}
\]
and
\[
\begin{CD}
L\wedge L\wedge L @>{\mu\wedge\id}>> L\wedge L\\
@V{\id\wedge\mu}VV @VV{\mu}V\\
L\wedge L @>{\mu}>> L
\end{CD}\quad.
\]
\end{Def}

\begin{Def}\label{Def:ringLwithinvolution}
An involution on an orthogonal ring spectrum is a map $\iota:L\rightarrow L$ such that
the following diagram commutes:
\[
\begin{CD}
L\wedge L @>{\iota\wedge\iota}>> L\wedge L @>{\text{twist}}>> L\wedge L\\
@V{\mu}VV && @VV{\mu}V\\
L @>{\iota}>> L @= L
\end{CD}\quad.
\]
\end{Def}

\begin{Rem}\label{Rem:extring}
We can externalize the definition of an orthogonal ring spectrum. What we then get is a continuous
functor $L:\mathscr{I}\rightarrow \Top_*$ together with natural transformations
\[
\eta:S^V\rightarrow L(V)\quad\text{and}\quad\mu:L(V)\wedge L(W)\rightarrow L(V\oplus W)
\]
satisfying certain conditions. $L$ has involution if we in addition have $\iota:L(V)\rightarrow L(V)$.
We often call an orthogonal ring spectrum an $\mathscr{I}$-FSP when we view it externally.
\end{Rem}

\section{Cellular techniques; q-cofibrations}\label{sect:cellullar}

According to~\cite{MandellMaySchwedeShipley:01} the q-cofibrations are the retracts of
relative cellular maps. Let us therefore see what a relative cellular map is:
We start by defining our set of cells.

\begin{Def}
Let $FI$ be the set of all maps $F_{\R^m}S^{n-1}_+\rightarrow F_{\R^m}D^n_+$, where $m\geq0$ and $n\geq0$.
\end{Def}

We think about $F_{\R^m}D^n_+$ as a cell with boundary $F_{\R^m}S^{n-1}_+$. In the case $n=0$ the boundary
is $*$.

\begin{Rem}[Symmetries of cells]\label{Rem:cellsymmetries}
If we inspect a cell $F_{\R^m}S^{n-1}_+\rightarrow F_{\R^m}D^n_+$,
we see that it has internal symmetries. We will particularly be interested in two different actions.
First we have an action of the permutation group $\Sigma_n$. It acts on $\R^n$ by permuting the factors.
Any permutation preserves the subspaces $D^n$ and $S^{n-1}$. Therefore, $\sigma\in\Sigma_n$ induces a map of pairs
\[
\sigma_*:(D^n_+,S^{n-1}_+)\rightarrow (D^n_+,S^{n-1}_+)\quad.
\]
Applying the shift desuspension functor $F_{\R^m}(-)$ we get self-maps of the cell $F_{\R^m}S^{n-1}_+\rightarrow F_{\R^m}D^n_+$.
We denote this map by $F(\sigma)$.

There is also another action. $\Sigma_m$ acts on $\R^m$ by permuting the factors. These maps are isometries, so for any
space $A$ a permutation $\sigma\in\Sigma_m$ induces a map
\[
F_{\sigma}: F_{\R^m} A\rightarrow F_{\R^m} A
\]
natural in $A$. This gives another action on $F_{\R^m}S^{n-1}_+\rightarrow F_{\R^m}D^n_+$.
\end{Rem}

\begin{Def}\label{def:FIcell}
A map $i:A\rightarrow L$ of orthogonal spectra is relative $FI$-cellular if:
\begin{itemize}
\item[] $i(A)$ is a subspectrum of $L$.
\item[] There is a set $C$ of subspectra $L_{\alpha}$ such that each $L_{\alpha}$ contains $i(A)$ and $\bigcup_{\alpha\in C}L_{\alpha}=L$.
\item[] $C$ is partially ordered by inclusion. We write $\beta\leq\alpha$ if $L_{\beta}\subseteq L_{\alpha}$. And for all $\alpha$ the
set $P_{\alpha}=\{ \beta\in C\;|\; \beta<\alpha\}$ is finite.
\item[] For every $\alpha\in C$ there is a pushout diagram
\[
\begin{CD}
F_{\R^m}S^{n-1}_+@>>> F_{\R^m}D^n_+\\
@VVV @VVV\\
\bigcup_{\beta<\alpha} L_{\beta} @>>> L_{\alpha}
\end{CD}\quad.
\]
\end{itemize}
\end{Def}

Recall that in the category of orthogonal spectra pushouts are formed level-wise.
If $\alpha$ is minimal, then the union $\bigcup_{\beta<\alpha} L_{\beta}$ is indexed over the empty set.
If this is the case we interpret the union as $i(A)$.
$C$ is the set of cells in the given relative cellular decomposition of $i:A\rightarrow L$, and
$\beta<\alpha$ if the cell $\beta$ is attached prior to $\alpha$. Observe that we allow
some redundancy in this definition, since we do not insist that the map 
$F_{\R^m}S^{n-1}_+\rightarrow \bigcup_{\beta<\alpha} L_{\beta}$
meets each $L_{\beta}$ non-trivially.

\begin{Rem}
There is also a notion of relative CW-orthogonal spectra: What we do is to put on the extra condition that
the cells are attached to cells of lower dimension only. The dimension of a cell $F_{\R^m}D^n_+$ is $n-m$.
In other words: A relative $FI$-cell structure of $i:A\rightarrow L$ is 
CW if the map $\dim:C\rightarrow \Z$ is strictly increasing.
\end{Rem}

\begin{Def}
A map $i:A\rightarrow L$ is a \textit{q-cofibration} if it is a retract of a relative $FI$-cellular map.
We call $L$ \textit{cofibrant} if $*\rightarrow L$ is a q-cofibration.
\end{Def}

This definition says that there exists a relative $FI$-cellular map $B\rightarrow K$ and a diagram
\[
\begin{CD}
A @>>> B @>>> A\\
@V{i}VV @VVV @VV{i}V\\
L @>>> K @>>> L
\end{CD}
\]
such that the horizontal compositions are the identity. Now observe that there is no loss
of generality if we assume that $B=A$. This follows from the elementary fact that relative $FI$-cellular maps are
closed under cobase change.

Observe that all q-cofibrations are both l-cofibrations and h-cofibrations.

\begin{Rem}
Alternatively one could define q-cofibrations as the maps which has the left lifting property with respect to
all level acyclic fibrations. See \S{}6 in~\cite{MandellMaySchwedeShipley:01}.
Recall that a level acyclic fibration $f:E\rightarrow B$ is by definition
a map such that for each $V$ the map $f_V:E(V)\rightarrow B(V)$ is both a weak equivalence and
a Serre fibration.
\end{Rem}

Here is an example of an orthogonal spectrum which is not cofibrant:

\begin{Exa}
Let $S'$ be the orthogonal spectrum given by
\[
S'(V)=\begin{cases}
S^V&\text{if $\dim V>0$, and}\\
*&\text{if $\dim V=0$.}
\end{cases} 
\]
This is a subspectrum of $S$, and the
assembly maps are inherited.
We will now show that $S'$ is not cofibrant.

Assume that $S'$ is cofibrant.
For contradiction we now construct a level acyclic fibration $f:E\rightarrow B$
and a diagram 
\[
\begin{CD}
* @>>> E\\
@VVV @VVV\\
S'@>>> B
\end{CD}
\] 
such that no lifting $S'\rightarrow E$ exists. 
Hence, $S'$ cannot be cofibrant.

Consider the map
\[
p:S^{\infty}_+\wedge S^1\rightarrow S^1
\]
given by collapsing $S^{\infty}$ to a point. This is a weak equivalence since $S^{\infty}$ is contractible.
The map $p$ is involutive. Here $\Z/2$ acts on $S^1$ by reflecting the circle across
a line, and $\Z/2$ acts freely on $S^{\infty}$. The smash product is given the diagonal action.
However $p$ is not a fibration, so we use the standard trick: Let $E_p$ be the $\Z/2$-space
of pairs $(x,\gamma)$, where $x\in S^{\infty}_+\wedge S^1$ and $\gamma$ is a path in $S^1$ such that
$p(x)=\gamma(0)$. The natural map 
\[
(x,\gamma)\mapsto \gamma(1):E_p\rightarrow S^1
\]
is again a weak equivalence.
Taking $\Z/2$-fixed points we see that
\[
*=(E_p)^{\Z/2}\rightarrow (S^1)^{\Z/2}=S^0\quad.
\]
Now define the level acyclic fibration $E\rightarrow B$ of orthogonal spectra by
\[
E(V)=\begin{cases}
\mathscr{I}(\R^1,V)\wedge_{O(1)}E_p&\text{if $\dim V=1$,}\\
*&\text{otherwise,}
\end{cases}
\]
and
\[
B(V)=\begin{cases}
\mathscr{I}(\R^1,V)\wedge_{O(1)} S^1&\text{if $\dim V=1$, and}\\
*&\text{otherwise.}
\end{cases}
\]

There is a map $S'\rightarrow B$ defined by letting the evaluation on level $\R^1$ be the identity.
If $S'$ is cofibrant, there exists a lift to $E$ and at level $\R^1$ we have
\[
S'(\R^1)\rightarrow E(\R^1)\rightarrow B(\R^1)
\]
where the composite is a $\Z/2$-equivariant homeomorphism. Taking $\Z/2$-fixed points we get
\[
S^0\rightarrow *\rightarrow S^0\quad,
\]
and since the composition cannot be a homeomorphism, this yields the contradiction.
\end{Exa}

For topological spaces we know that a compact subset of a $CW$-complex only meets finitely many cells. The same
is true for relative $FI$-cellular orthogonal spectra, as the following lemma shows:

\begin{Lem}\label{lem:compactmeetsfinitelymanycells}
If $K$ is a compact space, $i:A\rightarrow L$ a relative $FI$-cellular map of orthogonal spectra and $f:K\rightarrow L(\R^m)$ a map,
then there exists a finite set of cells $P$ such that $f$ factors through $\bigcup_{\alpha\in P}L_{\alpha}(\R^m)$.
\end{Lem}

\begin{proof}
We say that $K$ meets a cell $\alpha$ non-trivially if there exists a point $x$ in $K$ such that 
$f(x)\in L_{\alpha}(\R^m)$, but $f(x)\not\in L_{\beta}(\R^m)$ for any $\beta<\alpha$. We have to prove that 
$K$ only meets finitely many cells non-trivially. For a contradiction assume that $S=\{x_1,x_2,\ldots\}$
is a countable subset of $K$ such that each $x_i$ meets a distinct cell non-trivially. Then we can show that
$f(S)\cap L_{\alpha}(\R^m)$ is closed for all $\alpha$ by induction on the number of elements in $P_{\alpha}$.
The induction step uses the pushout diagram in definition~\ref{def:FIcell}. Since $L$ has the topology of
$\colim_{\alpha\in C}L_{\alpha}$, it follows that $f(S)$ is closed in $L(\R^m)$. The same argument shows that any
subset of $f(S)$ also is closed. Hence $f(K)$ contains an infinite discrete set. This contradicts compactness of $K$.

Let $P$ be the set of cells which $K$ meets non-trivially.
\end{proof}

We have the following reformulation of cellularity.

\begin{Prop}\label{Prop:seqcelldef}
A map $j:A\rightarrow L$ of orthogonal spectra is relative $FI$-cellular if and only if 
there exists a sequence
$L_0\rightarrow L_1\rightarrow\cdots$ of orthogonal spectra such that:
\begin{itemize}
\item[] $A=L_0$.
\item[] $L=\colim_i L_i$ and $j$ equals the natural map $L_0\rightarrow\colim_i L_i$.
\item[] For each $i$ there is a pushout diagram
\[
\begin{CD}
\bigvee_{\alpha\in C_i} F_{V_{\alpha}}S^{n_{\alpha}-1}_+ @>>> \bigvee_{\alpha\in C_i} F_{V_{\alpha}}D^{n_{\alpha}}_+\\
@VVV @VVV\\
L_{i} @>>> L_{i+1}
\end{CD}\quad.
\]
Here $C_i$ is the set of cells attached to $L_i$.
\end{itemize}
\end{Prop}

\begin{proof}
Assume that $j$ is relative $FI$-cellular. 
Let $C_i$ be a collection of subsets of $C$ such that
\begin{itemize}
\item[] the set of indexes $i$ is the non-negative natural numbers, 
\item[] the $C_i$'s are disjoint, and $\bigcup_{i=0}^{\infty} C_i=C$,
\item[] if $\alpha\in C_i$ and $\beta<\alpha$ then $\beta\in C_k$ for some $k<i$.
\end{itemize}
For example, we could let $C_i$ be the set of all cells $\alpha$ such that $P_{\alpha}$ contains
exactly $i$ elements. 

We set $L_0=A$, and let $L_{i+1}$ be the union of all $L_{\alpha}$ when $\alpha$ runs through $C_0\cup C_1\cup \cdots\cup C_i$.
Then $\bigcup_{\alpha\in C} L_{\alpha}=L$ and $j:A\rightarrow L$ is the natural map $L_0\rightarrow\colim_i L_i$.
To get the pushout diagram of the proposition, consider the functor $D$ from $C_0\cup C_1\cup \cdots\cup C_i$
to pushout diagrams which sends $\alpha\in C_i$ to 
\[
\begin{CD}
F_{\R^{m_{\alpha}}}S^{n_{\alpha}-1}_+@>>> F_{\R^{m_{\alpha}}}D^{n_{\alpha}}_+\\
@VVV @VVV\\
\bigcup_{\beta<\alpha} L_{\beta} @>>> L_{\alpha}
\end{CD}\quad,
\]
and $\alpha \in C_0\cup C_1\cup \cdots\cup C_{i-1}$ to
\[
\begin{CD}
* @>>> *\\
@VVV @VVV\\
L_{\alpha} @>>> L_{\alpha}
\end{CD}\quad.
\]
If $\beta<\alpha$, then we clearly have a map of pushout diagrams $D(\beta)\rightarrow D(\alpha)$. Taking the colimit
yields the desired diagram, and this is a pushout, since forming pushouts commutes with forming colimits.

Now assume that $j$ satisfies the properties of the proposition. We will inductively construct relative
$FI$-cell structures on $L_{i+1}$ with cells $C_0\cup\cdots\cup C_{i}$,
such that $L_i$ is a subcomplex. To start the induction we regard $A\xrightarrow{=}L_0$
as a relative $FI$-cell complex with the set of cells being the empty set.
Assume that $L_i$ already has been given a relative $FI$-cellular structure. The set of cells in $L_{i+1}$
should be $C_0\cup\cdots\cup C_{i}$, but we need to extend the partial ordering from $C_0\cup\cdots\cup C_{i-1}$.
We do this by specifying subspectra $L_{\alpha}$ for all $\alpha$ in $C_i$.
Recall that for each such $\alpha$ we have a diagram
\[
\begin{CD}
F_{V_{\alpha}}S^{n_{\alpha}-1}_+ @>>> F_{V_{\alpha}}D^{n_{\alpha}}_+\\
@VVV @VVV\\
L_{i} @>>> L_{i+1}
\end{CD}\quad.
\]
Look at the attaching map $S^{n_{\alpha}-1}_+\rightarrow L_i(V_{\alpha})$. By lemma~\ref{lem:compactmeetsfinitelymanycells}
this map factors through $\bigcup_{\beta\in P}L_{\beta}(V_{\alpha})$ for some finite subset $P$ of $C_0\cup\cdots\cup C_{i-1}$.
Define $L_{\alpha}$ by the pushout diagram
\[
\begin{CD}
F_{V_{\alpha}}S^{n_{\alpha}-1}_+ @>>> F_{V_{\alpha}}D^{n_{\alpha}}_+\\
@VVV @VVV\\
\bigcup_{\beta\in P}L_{\beta} @>>> L_{i+1}
\end{CD}\quad.
\]
Here the $L_{\beta}$'s are already defined since $\beta$ is a cell in $L_i$.

Letting $i$ go to $\infty$ we get a relative $FI$-cellular structure for $j:A\rightarrow L$.
\end{proof}

\subsection{The smash product of relative cellular maps}

Assume that $A\rightarrow L$ and $B\rightarrow K$ are relative $FI$-cellular maps. We
will now describe the relative $FI$-cellular structure of $L\wedge B\cup A\wedge K\rightarrow L\wedge K$.

We will need a technical lemma:

\begin{Lem}\label{Lem:polem}
Consider a diagram of spaces:
\[
\begin{CD}
B_0 @<<< A_0 @>>> X_0\\
@VVV @VVV @VVV\\
B_1 @<<< A_1 @>>> X_1
\end{CD}\quad.
\]
Let $Y_0$ and $Y_1$ be the pushout of the top and bottom row respectively.
Then the diagram
\[
\begin{CD}
A_1\cup_{A_0} X_0 @>>> X_1\\
@VVV @VVV\\
B_1\cup_{B_0} Y_0 @>>> Y_1
\end{CD}
\]
is pushout.
\end{Lem}

Here we use the notation $B\cup_A X$ for the pushout of $B\leftarrow A\rightarrow B$, even 
when neither of the two maps are injective.

\begin{proof}
To see this, take a look at the diagram
\[
\begin{CD}
A_1 @>>> A_1\cup_{A_0} X_0 @>>> X_1\\
@VVV @VVV @VVV\\
B_1 @>>> B_1\cup_{B_0} Y_0 @>>> Y_1
\end{CD}\quad.
\]
By the observation that $B_1\cup_{B_0} Y_0=B_1\cup_{B_0}(B_0\cup_{A_0} X_0)=B_1\cup_{A_0}X_0=B_1\cup_{A_1}(A_1\cup_{A_0}X_0)$,
we get that the left square is pushout.
Since the outer square is pushout, it now follows by cancellation that the right square also is.
\end{proof}

\begin{Lem}\label{Lem:wedgepushout}
Assume that
\[
\begin{CD}
A_0 @>{i_0}>> L_0\\
@VVV @VVV\\
A_1 @>>> L_1
\end{CD}\quad\text{and}\quad
\begin{CD}
B_0 @>{j_0}>> K_0\\
@VVV @VVV\\
B_1 @>>> K_1
\end{CD}
\]
are pushout squares of orthogonal spectra, where $i_0$ and $j_0$ are l-cofibrations.
Then the diagram
\[
\begin{CD}
A_0\wedge K_0 \cup L_0\wedge B_0 @>>> L_0\wedge K_0\\
@VVV @VVV\\
A_1\wedge K_1 \cup L_1\wedge B_1 @>>> L_1\wedge K_1
\end{CD}
\]
is also pushout.
\end{Lem}

Using the $\square$ product, definition~\ref{Def:squareproduct}, we can say that the bottom diagram is the row-wise
$\square$ of the upper diagrams.

\begin{proof}
It is enough to prove the result in the case where $B_1=B_0$ and $K_1=K_0$.

Since the functor $-\wedge X$ preserve pushout diagrams for any orthogonal spectrum $X$, we have that
the row-wise pushout of
\[
\begin{CD}
A_1\wedge B_0 @<<< A_0\wedge B_0 @>>> L_0\wedge B_0\\
@VVV @VVV @VVV\\
A_1\wedge K_0 @<<< A_0\wedge K_0 @>>> L_0\wedge K_0
\end{CD}
\]
is $L_1\wedge B_0\rightarrow L_1\wedge K_0$. Since the property of being a pushout diagram
is level-wise, we can apply lemma~\ref{Lem:polem}. Thus we get that
\[
\begin{CD}
A_0\wedge K_0\cup L_0\wedge B_0 @>>> L_0\wedge K_0\\
@VVV @VVV\\
A_1\wedge K_0 \cup L_1\wedge B_0 @>>> L_1\wedge K_0
\end{CD}
\]
is pushout. This completes the proof.
\end{proof}

\begin{Prop}\label{Prop:prodcell}
Assume that $i:A\rightarrow L$ and $j:B\rightarrow K$ are relative $FI$-cellular maps. Then
$i\square j:L\wedge B\cup A\wedge K\rightarrow L\wedge K$ is also relative $FI$-cellular.
\end{Prop}

\begin{proof}
We will describe the relative $FI$-cellular structure of $i\square j$.
Let $C$ and $D$ be the sets of cells of $i$ and $j$ respectively. The set of cells for 
$i\square j$ will be $C\times D$, and we define $(L\wedge K)_{(\alpha_1,\alpha_2)}$ to be
\[
L\wedge B\cup L_{\alpha_1}\wedge K_{\alpha_2} \cup A\wedge K\quad.
\]
Observe that $P_{(\alpha_1,\alpha_2)}=P_{\alpha_1}\times P_{\alpha_2}\cup\{\alpha_1\}\times P_{\alpha_2} \cup P_{\alpha_1}\times \{\alpha_2\}$,
thus it is finite.
It is clear that $\bigcup_{C\times D} (L\wedge K)_{(\alpha_1,\alpha_2)}=L\wedge K$. It remains to show that
for each cell $(\alpha_1,\alpha_2)$ there are pushout diagrams:
\[
\begin{CD}
F_{V_{\alpha_1}\oplus V_{\alpha_2} }S^{n_{\alpha_1}+n_{\alpha_2}-1}_+ @>>> F_{V_{\alpha_1}\oplus V_{\alpha_2}}D^{n_{\alpha_1}+n_{\alpha_2}}_+\\
@VVV @VVV\\
\bigcup_{(\beta_1,\beta_2)\in P_{(\alpha_1,\alpha_2)}}(L\wedge K)_{(\beta_1,\beta_2)}@>>> (L\wedge K)_{(\alpha_1,\alpha_2)}
\end{CD}\quad.
\]
To see this, apply lemma~\ref{Lem:wedgepushout} to the diagrams:
\[
\begin{CD}
F_{V_{\alpha_1}}S^{n_{\alpha_1}-1}_+ @>>> F_{V_{\alpha_1}}D^{n_{\alpha_1}}_+\\
@VVV @VVV\\
\bigcup_{\beta_1\in P_{\alpha_1}}L_{\beta_1}@>>> L_{\alpha_1}
\end{CD}\quad\text{and}\quad
\begin{CD}
F_{V_{\alpha_2}}S^{n_{\alpha_2}-1}_+ @>>> F_{V_{\alpha_2}}D^{n_{\alpha_2}}_+\\
@VVV @VVV\\
\bigcup_{\beta_2\in P_{\alpha_2}}K_{\beta_2}@>>> K_{\alpha_2}
\end{CD}\quad.
\]
Then apply $(L\wedge B\cup A\wedge K)\cup_{(L_{\alpha_1}\wedge B\cup A\wedge K_{\alpha_2})}-$ to the lower map.
This concludes the proof.
\end{proof}

\subsection{Cofibrant replacement functor}

Now we shall introduce a cofibrant replacement functor. Among other uses, we want to apply this functor
to orthogonal ring spectra with involution to get cofibrant orthogonal ring spectra with involution. 
Therefore we want the functor to
be lax skew-symmetric with respect to $\wedge$. See appendix~\ref{sect:MonCat} for definitions of 
(symmetric) (co)monoidal categories and lax/strong (symmetric) (co)monoidal functors.
Skew-symmetry will be defined below.
Our structure theorem is:

\begin{Thm}\label{Thm:cofrepl}
There is an endofunctor $\Gamma$ on orthogonal spectra having the following properties:
\begin{itemize}
\item[] $\Gamma L$ is cofibrant for all $L$.
\item[] If $K\rightarrow L$ is the inclusion of a subspectrum, then $\Gamma K\rightarrow \Gamma L$ is a q-cofibration.
\item[] $\Gamma$ comes with a natural level-wise acyclic fibration $\gamma_L:\Gamma L\rightarrow L$.
\item[] $\Gamma$ comes with a natural transformation $\phi_{L,K}:\Gamma L\wedge \Gamma K\rightarrow \Gamma (L\wedge K)$.
\item[] $\Gamma$ comes with an involution $\iota:\Gamma L\rightarrow\Gamma L$, and $\iota^2=\id$.
\item[] There is a canonical level-wise acyclic q-cofibration $\lambda:S\rightarrow \Gamma S$.
\item[] With $\lambda$, $\phi$ and $\iota$ the functor $\Gamma$ is lax skew-symmetric monoidal with respect to $\wedge$.
\item[] $\Gamma$ comes with a natural level equivalence $\rho_{L,K}:\Gamma(L\times K)\rightarrow \Gamma L\times \Gamma K$.
\item[] With $\gamma_*:\Gamma*\rightarrow *$ and $\rho$ the functor $\Gamma$ is lax symmetric comonoidal with respect to $\times$.
\end{itemize}
\end{Thm}

Note that $\phi$ is not always a $\pi_*$-iso. If $L$ and $K$ both are sufficiently bad, the smash product $L\wedge K$
can have homotopy unrelated to $L$ and $K$. Thus $\Gamma (L\wedge K)$ will have the same bad homotopy, 
whereas $\Gamma L\wedge \Gamma K$
will have the correct homotopy. However, if either $L$ or $K$ is cofibrant, then 
$\phi$ will be a $\pi_*$-isomorphism, since smashing
with a cofibrant orthogonal spectrum preserves $\pi_*$-isomorphisms, see proposition~\ref{Prop:os:wedge}.

\begin{proof}
We will break the proof of the structure theorem into several propositions, 
and the proof will span the rest of this subsection.
However, we first show that the comonoidality with respect to $\times$ is a formal consequence of the other properties:

To define a map into a $\times$-product, it is enough to define one map into each factor. $\Gamma$ applied to
the projections gives
\[
\Gamma(L\times K)\rightarrow \Gamma L\quad\text{and}\quad \Gamma(L\times K)\rightarrow \Gamma K\quad,
\]
and $\rho_{L,K}$ is determined by these maps. It is elementary to see that 
the cross product of two level equivalences is a level equivalence.
Hence, we have the diagram
\[
\begin{CD}
\Gamma(L\times K) @>{\rho}>> \Gamma L\times \Gamma K\\
@V{\gamma}V{\simeq}V @V{\simeq}V{\gamma\times\gamma}V\\
L\times K @= L\times K
\end{CD}\quad,
\]
and it follows that $\rho$ is a level equivalence. 

We know that $\times$ is the categorical product on $\mathscr{IS}$, hence lemma~\ref{Lem:catprodimpllaxsymfunct}
implies that $\Gamma$ is lax symmetric comonoidal with respect to $\times$.
\end{proof}

Let us now look at the construction of $\Gamma$. 
The idea is to apply Quillen's small object argument, see~7.12 in~\cite{DwyerSpalinski:95}.
We proceed as follows:

Suppose that $p:A\rightarrow L$ is a map of orthogonal spectra.
We now construct an orthogonal spectrum $G(p)$ and a factorization of $p$:
\[
A\rightarrow G(p)\rightarrow L\quad.
\]
Let $C$ be the set of all diagrams
\[
\begin{CD}
F_{\R^{m}}S^{n-1}_+ @>{F_{\R^{m}}i_n}>> F_{\R^{m}}D^{n}_+\\
@V{f}VV @VV{g}V\\
A @>{p}>> L
\end{CD}\quad,
\]
where $n,m\geq 0$. Now define $G(p)$ by the pushout diagram
\[
\begin{CD}
\bigvee_{\alpha\in C} F_{\R^{m_{\alpha}}}S^{n_{\alpha}-1}_+ @>>> 
\bigvee_{\alpha\in C}F_{\R^{m_{\alpha}}}D^{n_{\alpha}}_+\\
@V{\bigvee_{\alpha\in C}f_{\alpha}}VV @VV{\bigvee_{\alpha\in C}g_{\alpha}}V\\
A @>{p}>> G(p)
\end{CD}\quad.
\]

By construction we see that $A$ is a subspectrum of $G(p)$, and the natural map $G(p)\rightarrow L$ comes 
from the universal property of the pushout. Observe that $G(p)$ is a functor; a morphism between maps
$p_1:A_1\rightarrow L_1$ and $p_2:A_2\rightarrow L_2$ is a commutative diagram
\[
\begin{CD}
A_1 @>{p_1}>> L_1\\
@VVV @VVV\\
A_2 @>{p_2}>> L_2
\end{CD}\quad.
\]

Without a proof we observe that:

\begin{Lem}
There exists a relative $FI$-cellular structure on the map $A\rightarrow G(p)$.
\end{Lem}

To define $\Gamma L$ we iterate the gluing construction. 
Start with $p_0:*\rightarrow L$.
Apply the construction above and set $G^1(L)=G(p_0)$ to get a diagram:
\[
*\xrightarrow{j_0} G^1(L)\xrightarrow{p_1} L\quad.
\]
Iterate to get diagrams
\[
G^i(L)\xrightarrow{j_i} G^{i+1}(L)\xrightarrow{p_{i+1}} L\quad.
\]

\begin{Def}
Define the q-cofibrant replacement of $L$, $\Gamma L$, to be the colimit of the $G^i(L)$'s
constructed above. 
\end{Def}

As a colimit of a sequence relative $FI$-cellular maps starting from $*$, we see that $\Gamma L$ is cofibrant.
$\Gamma$ is a functor and the natural transformation $\gamma_L:\Gamma L\rightarrow L$
is induced by the natural maps $G^i(L)\rightarrow L$.
Let us now begin to prove the various statements of theorem~\ref{Thm:cofrepl}.

\begin{Prop}\label{Prop:Gammaincltoqcof}
If $j:K\rightarrow L$ is the inclusion of a subspectrum, then $\Gamma K\rightarrow \Gamma L$ is a q-cofibration.
\end{Prop}

\begin{proof}
Inspect the gluing construction. By induction it is enough to consider a diagram
\[
\begin{CD}
B @>{j_0}>> A\\
@V{p}VV @VV{q}V\\
K @>{j}>> L
\end{CD}\quad,
\]
where $j_0$ is a relative $FI$-cellular. We must show that $G(p)\cup_B A$
maps homeomorphically onto a subcomplex of $G(q)$. Compare with proposition~\ref{Prop:qcofincolim}.

Consider a cell $\alpha$ of $G(p)$:
\[
\begin{CD}
F_{\R^{m}}S^{n-1}_+ @>{F_{\R^{m}}i_n}>> F_{\R^{m}}D^{n}_+\\
@V{f}VV @VV{g}V\\
B @>{p}>> K
\end{CD}\quad.
\]
Composing $f$ with $j_0$ and $g$ with $j$ we get a new diagram $\beta$ representing a cell in $G(q)$.
Since $j_0$ and $j$ are injective, we see that different cells $\alpha$ and $\alpha'$ in $G(p)$
gives different cells $\beta$ and $\beta'$ in $G(q)$. It follows that 
\[
G(p)\cup_B A\rightarrow G(q)
\]
is relative $FI$-cellular.
\end{proof}

It is easy to construct the map $\lambda: S\rightarrow \Gamma S$:
The diagram 
\[
\begin{CD}
* @>>> F_0 D^0_+\\
@VVV @VV{\cong}V\\
* @>>> S
\end{CD}
\]
determines a cell $\alpha$ of $\Gamma S$. We have $\Gamma S_{\alpha}\cong S$, and
define $\lambda$ to be the composition
\[
S\xrightarrow{\cong}\Gamma S_{\alpha}\rightarrow \Gamma S\quad.
\]

\begin{Prop}\label{Prop:Gammalevelacyclfibr}
$\gamma:\Gamma L\rightarrow L$ is a level acyclic fibration.
\end{Prop}

\begin{proof}
Fix some level $V=\R^m$. We have to show that for every diagram
\[
\begin{CD}
S^{n-1} @>{f}>> \Gamma L(V)\\
@VVV @VV{\gamma}V\\
D^n @>>> L(V)
\end{CD}
\]
there is a lift $D^n\rightarrow \Gamma L(V)$. Since $S^{n-1}$ is compact there is an $i$ such
that $f$ factors through $G^i(L)(V)$. Then we get the diagram
\[
\begin{CD}
S^{n-1} @>{f}>> G^i(L)(V)\\
@VVV @VV{p_i}V\\
D^n @>>> L(V)
\end{CD}\quad,
\]
but this is exactly what determines a new cell $\alpha$ in $G(p_i)$. 
And we see that $D^n$ lifts into $G^{i+1}(L)(V)$.
\end{proof}

\begin{Constr}\label{Constr:GammaMult}
Next we construct the natural 
transformation $\phi:\Gamma L\wedge \Gamma K\rightarrow \Gamma(L\wedge K)$.
Inductively we define maps
\[
\phi_{i,j}:G^i(L)\wedge G^j(K)\rightarrow G^{i+j-1}(L\wedge K)
\]
such that the following diagrams commute for all $i$ and $j$:
\[
\begin{CD}
G^i(L)\wedge G^j(K) @>{\phi_{i,j}}>> G^{i+j-1}(L\wedge K)\\
@V{\subseteq}VV @VV{\subseteq}V\\
G^{i+1}(L)\wedge G^j(K) @>{\phi_{i+1,j}}>> G^{i+j}(L\wedge K)
\end{CD}\quad,
\]
\[
\begin{CD}
G^i(L)\wedge G^j(K) @>{\phi_{i,j}}>> G^{i+j-1}(L\wedge K)\\
@V{\subseteq}VV @VV{\subseteq}V\\
G^{i}(L)\wedge G^{j+1}(K) @>{\phi_{i,j+1}}>> G^{i+j}(L\wedge K)
\end{CD}\quad\text{and}
\]
\[
\begin{CD}
G^i(L)\wedge G^j(K) @>{\phi_{i,j}}>> G^{i+j-1}(L\wedge K)\\
@VVV @VVV\\
L\wedge K @= L\wedge K
\end{CD}\quad.
\]
By taking the colimit as both $i$ and $j$ tend to infinity, we get our natural transformation
$\Gamma L\wedge \Gamma K\rightarrow \Gamma(L\wedge K)$.

Let $G^i(L)$ be $*$ for $i\leq0$.
If $i\leq0$ or $j\leq0$, then $\phi_{i,j}$ is trivially defined. 
We construct the $\phi$'s by induction on $i+j$.
Let $\alpha$ and $\beta$ be the diagrams
\[
\begin{CD}
F_{\R^{m}}S^{n-1}_+ @>>> F_{\R^{m}}D^{n}_+\\
@V{f}VV @VV{g}V\\
G^{i-1}(L) @>>> L
\end{CD}\quad\text{and}\quad
\begin{CD}
F_{\R^{m'}}S^{n'-1}_+ @>>> F_{\R^{m'}}D^{n'}_+\\
@V{f'}VV @VV{g'}V\\
G^{j-1}(K) @>>> K
\end{CD}
\]
respectively.
By the construction of $G^i(L)$ and $G^j(K)$ there are unique lifts of $\alpha$ and $\beta$ to diagrams
$\bar{\alpha}$ and $\bar{\beta}$:
\[
\begin{CD}
F_{\R^{m}}S^{n-1}_+ @>>> F_{\R^{m}}D^{n}_+\\
@V{f}VV @VV{\bar{g}}V\\
G^{i-1}(L) @>>> G^i(L)
\end{CD}\quad\text{and}\quad
\begin{CD}
F_{\R^{m'}}S^{n'-1}_+ @>>> F_{\R^{m'}}D^{n'}_+\\
@V{f'}VV @VV{\bar{g}'}V\\
G^{j-1}(K) @>>> G^j(K)
\end{CD}\quad.
\]
Recall the definition of $\square$ of two maps, see definition~\ref{Def:squareproduct}.
Now consider the diagram
\[
\begin{CD}
F_{\R^{m+m'}}S^{n+n'-1}_+ @>>> F_{\R^{m+m'}}D^{n+n'}_+\\
@V{f\wedge \bar{g}'\cup \bar{g}\wedge f'}VV @VV{\bar{g}\wedge \bar{g}'}V\\
G^{i-1}(L)\wedge G^j(K) \cup G^i(L)\wedge G^{j-1}(K) @>>> G^i(L) \wedge G^j(K)\\
@V{\phi_{i-1,j}\cup\phi_{i,j-1}}VV @VVV\\
G^{i+j-2}(L\wedge K) @>>> L\wedge K
\end{CD}\quad,
\]
where the upper part is row-wise $\square$ of $\bar{\alpha}$ and $\bar{\beta}$.
The map $\phi_{i-1,j}\cup\phi_{i,j-1}$ exists by our induction hypothesis.
The outer square is a diagram in $C$ for the gluing construction applied to 
$G^{i+j-2}(L\wedge K) \rightarrow L\wedge K$. Call this diagram $\delta$.
And by the cell $\delta$ we get a map
\[
G^i(L)_{\alpha}\wedge G^j(K)_{\beta} \rightarrow G^{i+j-1}(L\wedge K)_{\delta}\subseteq G^{i+j-1}(L\wedge K)\quad.
\]
Letting $\alpha$ and $\beta$ run through all cells of $G^i(L)$ and $G^j(K)$ respectively,
we get our map
\[
\phi_{i,j}:G^i(L)\wedge G^j(K)\rightarrow G^{i+j-1}(L\wedge K)\quad.
\]
This finishes the construction of $\phi:\Gamma L\wedge \Gamma K\rightarrow \Gamma(L\wedge K)$.
\end{Constr}

\begin{Lem}\label{Lem:twistcells}
Let $A\rightarrow L$ and $B\rightarrow K$ be maps of orthogonal spectra. 
Consider diagrams
\[
\begin{CD}
F_{\R^{m_1}}S^{n_1-1}_+@>>> F_{\R^{m_1}}D^{n_1}_+\\
@VVV @VVV\\
A@>>> L
\end{CD}\quad\text{and}\quad
\begin{CD}
F_{\R^{m_2}}S^{n_2-1}_+@>>> F_{\R^{m_2}}D^{n_2}_+\\
@VVV @VVV\\
B@>>> K
\end{CD}\quad,
\]
representing cells called $\alpha$ and $\beta$ respectively. 
We can compare the cells $\alpha\square\beta$ and $\beta\square\alpha$
via the twist map $L\wedge K\rightarrow K\wedge L$. And we have
\[
\text{twist}\circ(\alpha\square\beta)=(\beta\square\alpha)\circ F(\sigma) F_{\rho}\quad,
\]
where $F(\sigma)$ is the cell symmetry  
permuting coordinates of the $(n_1+n_2)$-disk as indicated by
the map $D^{n_1}\times D^{n_2}\cong D^{n_2}\times D^{n_1}$,
and $F_{\rho}$ is the cell symmetry permuting the
indexing spaces as indicated by $\R^{m_1}\oplus \R^{m_2}\cong \R^{m_2}\oplus \R^{m_1}$.
\end{Lem}

\begin{proof}
We write out the proof only for the disks. The boundary of the cells can be treated similarly.

Consider the diagram
\[
\begin{CD}
F_{\R^{m_1+m_2}}D^{n_1+n_2}_+ @>{F(\sigma)F_{\rho}}>> F_{\R^{m_2+m_1}}D^{n_2+n_1}_+ \\
@A{\cong}AA @AA{\cong}A\\
F_{\R^{m_1}}D^{n_1}_+\wedge F_{\R^{m_2}}D^{n_2}_+ @>{\text{twist}}>> F_{\R^{m_2}}D^{n_2}_+\wedge F_{\R^{m_1}}D^{n_1}_+ \\
@VVV @VVV\\
L\wedge K @>{\text{twist}}>> K\wedge L
\end{CD}\quad.
\]
The bottom part clearly commutes, so it remains to check that the top part also does.
We evaluate at level $\R^{m_1+m_2}$ and get:
\footnotesize
\[
\begin{CD}
\mathscr{I}(\R^{m_1+m_2},\R^{m_1+m_2})_+\wedge D^{n_1+n_2}_+ @>{\sigma_*\rho_*}>> 
\mathscr{I}(\R^{m_2+m_1},\R^{m_1+m_2})_+\wedge D^{n_2+n_1}_+ \\
@A{\cong}AA @AA{\cong}A\\
\mathscr{I}(\R^{m_1}\oplus \R^{m_2},\R^{m_1+m_2})_+\wedge D^{n_1}_+\wedge D^{n_2}_+ 
@>{\text{twist}}>> 
\mathscr{I}(\R^{m_2}\oplus \R^{m_1},\R^{m_1+m_2})_+\wedge D^{n_2}_+\wedge D^{n_1}_+ \\
\end{CD}\quad.
\]
\normalsize
The map $\text{twist}$ swaps both the indexing spaces $\R^{m_1}$ and $\R^{m_2}$, and
the disks $D^{n_1}$ and $D^{n_2}$. $\sigma_*$ is the map $D^{n_1+n_2}\cong D^{n_2+n_1}$
permuting the factors, while $\rho_*$ is the linear map $\R^{m_1+m_2}\cong\R^{m_2+m_1}$
applied to the first factor of $\mathscr{I}(\R^{m_1+m_2},\R^{m_1+m_2})$.
The left vertical identification is defined via
$\R^{m_1}\oplus \R^{m_2}\cong\R^{m_1+m_2}$ and $D^{n_1}_+\wedge D^{n_2}_+\cong D^{n_1+n_2}_+$,
and the right vertical identification is given by
$\R^{m_2}\oplus \R^{m_1}\cong\R^{m_2+m_1}$ and $D^{n_2}_+\wedge D^{n_1}_+\cong D^{n_2+n_1}_+$.
And we see that the diagram commutes.
\end{proof}

Now we are ready to define skew-symmetry and to show that $\Gamma$ satisfies this.

\begin{Def}
A lax monoidal functor $F:M\rightarrow B$ is skew-symmetric 
with respect to a product $\square$ if there exists a
natural transformation $\iota:F(a)\rightarrow F(a)$ with $\iota^2=\id$, such
that the following diagram commutes:
\[
\begin{CD}
F(a)\square F(b)@>{\iota\square\iota}>> F(a)\square F(b) @>{\gamma}>> F(b)\square F(a)\\
@V{\phi}VV && @VV{\phi}V\\
F(a\square b) @>{\iota}>> F(a\square b) @>{F(\gamma)}>> F(b\square a)
\end{CD}\quad.
\]
\end{Def}

We can say that $\iota$ measures the failure of $F$ being symmetric.
In our case we have:

\begin{Prop}\label{Prop:Gammaskew}
$\Gamma$ is a lax skew-symmetric monoidal functor with respect to $\wedge$. 
\end{Prop}

\begin{proof}
We begin by constructing $\iota$. It is a cellular map and can be defined on the gluing construction.
Let $\alpha$ be the cell
\[
\begin{CD}
F_{\R^{m}}S^{n-1}_+@>>> F_{\R^{m}}D^{n}_+\\
@VVV @VVV\\
A@>{p}>> L
\end{CD}
\]
of $G(p)$. Let $\tau_n$ and $\tau_m$ denote the order reversing permutations in $\Sigma_n$ and $\Sigma_m$
respectively. $\iota:G(p)\rightarrow G(p)$ is the map that sends $\alpha$ to the cell $F(\tau_n)F_{\tau_m}(\alpha)$.
By construction of $\Gamma L$ as the iterated gluing construction, we get the natural transformation
\[
\iota:\Gamma L\rightarrow \Gamma L\quad.
\] 
Clearly $\iota^2=\id$.

Inspecting the construction of $\phi:\Gamma L\wedge \Gamma K\rightarrow \Gamma (L\wedge K)$,
we see that $\Gamma$ is lax monoidal. In order to check skew-symmetry, it
remains to check that the following diagram commutes:
\[
\begin{CD}
\Gamma L\wedge \Gamma K @>{\iota\wedge\iota}>> \Gamma L\wedge \Gamma K 
@>{\text{twist}}>> \Gamma K\wedge \Gamma L\\
@V{\phi}VV && @VV{\phi}V\\
\Gamma (L\wedge K) @>{\iota}>> \Gamma (L\wedge K)
@>{\Gamma(\text{twist})}>> \Gamma (K\wedge L)
\end{CD}\quad.
\]
By induction on $i+j$ we prove that
\[
\begin{CD}
G^i(L)\wedge G^j(K) @>{\iota\wedge\iota}>> G^i(L)\wedge G^j(K) @>{\text{twist}}>> G^j(K)\wedge G^i(L)\\
@V{\phi_{i,j}}VV && @VV{\phi_{j,i}}V\\
G^{i+j-1}(L\wedge K) @>{\iota}>> G^{i+j-1}(L\wedge K) @>{G^{i+j-1}(\text{twist})}>> G^{i+j-1}(K\wedge L)
\end{CD}
\]
commutes. Let $\alpha$ and $\beta$ be the diagrams
\[
\begin{CD}
F_{\R^{m}}S^{n-1}_+ @>>> F_{\R^{m}}D^{n}_+\\
@V{f}VV @VV{g}V\\
G^{i-1}(L) @>>> L
\end{CD}\quad\text{and}\quad
\begin{CD}
F_{\R^{m'}}S^{n'-1}_+ @>>> F_{\R^{m'}}D^{n'}_+\\
@V{f'}VV @VV{g'}V\\
G^{j-1}(K) @>>> K
\end{CD}
\]
respectively, and let $\bar{\alpha}$ and $\bar{\beta}$ be liftings as defined in the construction of $\phi_{i,j}$.
By the previous lemma the diagrams
$\text{twist}\circ(\bar{\alpha}\square\bar{\beta})$
and $\bar{\beta}\square\bar{\alpha}$ differ by the permutation of coordinates $\sigma_*:D^{n+n'}\cong D^{n'+n}$
and the permutation of indexing space $\rho_*:\R^{m+m'}\cong\R^{m'+m}$. 
Computing with permutations in $\Sigma_{n+n'}$ we have that
\[
(\tau_n\amalg\tau_{n'})\sigma=\sigma\tau_{n+n'}\quad.
\]
And similarly $(\tau_m\amalg\tau_{m'})\rho=\rho\tau_{m+m'}$ in $\Sigma_{m+m'}$.
Therefore the diagram
\[
\begin{CD}
G^i(L)_{\alpha}\wedge G^j(K)_{\beta} @>{\iota}>> 
G^i(L)_{F(\tau_n)F_{\tau_m}(\alpha)}\wedge G^j(K)_{F(\tau_{n'})F_{\tau_{m'}}(\beta)}\\
@V{\phi_{i,j}}VV @V{\text{twist}}VV \\
G^{i+j-1}(L\wedge K) &&G^j(K)_{F(\tau_{n'})F_{\tau_{m'}}(\beta)}\wedge G^i(L)_{F(\tau_n)F_{\tau_m}(\alpha)}\\
@V{\iota}VV @VV{\phi_{j,i}}V\\
G^{i+j-1}(L\wedge K) @>{G^{i+j-1}(\text{twist})}>> G^{i+j-1}(K\wedge L)
\end{CD}
\]
commutes. And the result follows.
\end{proof}

\section{Boundedness}

When doing constructions with orthogonal spectra, it can be useful to consider those spectra
which are bounded below, or those which in addition have highly connected assembly maps.
However cellular orthogonal spectra does not in general have these properties. An example
is $\bigvee_{m=0}^{\infty} F_{\R^m}S^0$. But we may approximate any cellular orthogonal spectra by 
spectra satisfying these properties.

Let us begin with some definitions:

\begin{Def}
Let $L$ be an orthogonal spectrum.
\begin{itemize}
\item[]$L$ is \textit{strictly $c$-connected} if there exists an integer $N$ such that $L(\R^n)$ is $(n+c)$-connected for all $n\geq N$.
\item[]We call $L$ \textit{strictly connected} if $L$ is strictly $(-1)$-connected.
\item[]$L$ is \textit{strictly bounded below} if there exists a $c$ such that $L$ is strictly $c$-connected.
\end{itemize}
\end{Def}

We can simplify the definition of strictly bounded below:

\begin{Lem}
$L$ is strictly bounded below if and only if there exists a $c$ such that $L(\R^n)$ is $(n+c)$-connected for all $n$.
\end{Lem}

\begin{proof}
The ``if'' direction is obvious. For the ``only if'' direction assume that $L$ is strictly $c$-connected.
Let $N$ be such that $L(\R^n)$ is $(n+c)$-connected for all $n\geq N$. Now set $c'=\min(c,-N)$. When $n\geq N$, the
space $L(\R^n)$ is $(n+c)$-connected, hence also $(n+c')$-connected.
Since every based space is $(-1)$-connected, we have that $L(\R^n)$ is $(n+c')$-connected also for $n<N$.
\end{proof}

By the definition of the homotopy groups of an orthogonal spectrum we immediately get:

\begin{Prop}
Let $L$ be an orthogonal spectrum.
\begin{itemize}
\item[]If $L$ is strictly $c$-connected, then $\pi_q L=0$ for $q\leq c$.
\item[]If $L$ is strictly bounded below, then there exists a $c$ such that $\pi_q L=0$ for all $q\leq c$.
\end{itemize}
\end{Prop}

The converse is not true as the following examples show:

\begin{Exa}
Let $m$ be an integer and consider $L$ defined by
\[
L(V)=\begin{cases}
\mathscr{I}(\R^{n},V)_+\wedge S^{n-m} &\text{if $\dim V=n\geq m$, and}\\
* &\text{otherwise}
\end{cases}
\]
with trivial assembly maps. Assume $n\geq m$. By the Freudenthal suspension theorem we have that
$L(\R^n)$ is $(n-m-1)$-connected, and by homology calculations we get that $L(\R^n)$ is not
$(n-m)$-connected. Consequently, we have that $L$ is strictly $(-m-1)$-connected, but not strictly $(-m)$-connected. 
However, the homotopy groups $\pi_q L$ are trivial for all $q$. 
\end{Exa}

\begin{Exa}
Let $L$ be given by
\[
L(V)=\bigvee_{m=0}^{\infty}\mathscr{I}(\R^m,V)_+
\]
with trivial assembly maps. Then for all $n\geq0$ we have that $L(\R^n)$ is not $0$-connected.
Hence, there exists no $c$ such that $L(\R^n)$ is $(n+c)$-connected for all $n$. It follows that
$L$ is not strictly bounded below. But since the assembly maps are trivial, it follows that
$\pi_q L=0$ for all $q$.
\end{Exa}

\begin{Lem}
If $A\rightarrow X$ is an l-cofibration, $A$, $B$ and $X$ strictly $c$-connected orthogonal spectra, then the pushout, $Y$, of
$B\leftarrow A\rightarrow X$ is also strictly $c$-connected.
\end{Lem}

\begin{proof}
The spaces $A(\R^n)$, $B(\R^n)$ and $X(\R^n)$ all are $(n+c)$-connected.
We want to show that $Y(\R^n)$ also is $(n+c)$-connected.
Consider the diagram
\[
\begin{CD}
A(\R^n) @>>> X(\R^n)\\
@VVV @VVV\\
B(\R^n) @>>> Y(\R^n)
\end{CD}\quad.
\]
For $(n+c)=-1$ and $(n+c)=0$, it is obvious that $Y(\R^n)$ is $(n+c)$-connected. Assume
$(n+c)>0$. Blakers-Massey applies and shows that
\[
\pi_q(X(\R^n),A(\R^n))\rightarrow \pi_q(Y(\R^n),B(\R^n))
\]
is an isomorphism when $q<2(n+c)$. By the long exact sequences in homotopy for the pairs
$(X(\R^n),A(\R^n))$ and $(Y(\R^n),B(\R^n))$, the result follows.
\end{proof}

\begin{Cor}\label{Cor:posbb}
Assume $A$, $B$ and $X$ strictly bounded below. Then the pushout of $B\leftarrow A\rightarrow X$
is also strictly bounded below if at least one of the maps is an l-cofibration.
\end{Cor}

\begin{Lem}\label{Lem:VFAconn}
Let $L=\bigvee_{\alpha} F_{V_{\alpha}}A_{\alpha}$. If $A_{\alpha}$ is well-pointed
and $\dim V_{\alpha}\leq k$ for all $\alpha$, then $L$ is strictly $(-k-1)$-connected.
\end{Lem}

\begin{proof}
First consider the case with only one wedge summand: $A=A_{\alpha}$ is a well-pointed space and 
$V_{\alpha}=\R^k$.
We want to calculate the connectivity of $F_{\R^k} A(\R^n)$. If $n<k$, then by definition
$F_{\R^k} A(\R^n)=*$, so assume that $k\leq n$. Then
\[
F_{\R^k} A(\R^n)= O(n)_+\wedge_{O(n-k)}(A\wedge S^{n-k})\cong A\wedge \left( O(n)_+\wedge_{O(n-k)}S^{n-k}\right)\quad.
\]
Now consider the diagram
\[
\begin{CD}
&&O(n-k)&&\\
&&@VVV&&\\
&&O(n)\times S^{n-k}&&\\
&&@VVV&&\\
O(n)/O(n-k)@>>> O(n)\times_{O(n-k)} S^{n-k} @>>> O(n)_+\wedge_{O(n-k)}S^{n-k}
\end{CD}\quad.
\]
Here the vertical sequence is a fibration, and the horizontal sequence a cofibration.
Since $O(n-k)\rightarrow O(n)$ is $(n-k-1)$-connected, the long exact sequence of the fibration yields
that $O(n)\times_{O(n-k)} S^{n-k}$ is $(n-k-1)$-connected. Furthermore,
we know that $O(n)/O(n-k)$ is $(n-k-1)$-connected.
Using proposition~4.28 in~\cite{Hatcher:02}, 
we see that $O(n)_+\wedge_{O(n-k)}S^{n-k}$ is also $(n-k-1)$-connected.

By Blakers-Massey or CW-approximation, one can prove that for well-pointed spaces $X$ and $Y$ which are $r$- and $s$-connected
respectively, the smash product $X\wedge Y$ is $(r+s+1)$-connected. Since any space is $(-1)$-connected, applying
this to $A\wedge \left( O(n)_+\wedge_{O(n-k)}S^{n-k}\right)$ yields that $F_{\R^k} A(\R^n)$ is $(n-k-1)$-connected.

Now consider the wedge
\[
L=\bigvee_{\alpha} F_{V_{\alpha}}A_{\alpha}\quad.
\]
By the calculation above $L(\R^n)$ is a wedge of well-pointed $(n-k-1)$-connected spaces. 
Using CW-approximation or Blakers-Massey we can prove that the wedge of well-pointed $l$-connected
spaces again is $l$-connected. And it follows that $L(\R^n)$ is $(n-k-1)$-connected for all $n$.
\end{proof}

For some purposes we need a stronger condition than strictly bounded below:

\begin{Def}
An orthogonal spectrum $L$ is \textit{meta-stable} if it is strictly bounded below 
and there exists an integer $d$ such that
$\sigma:L(\R^n)\wedge S^1\rightarrow L(\R^{n+1})$ is $(2n+d)$-connected for all $n$.
\end{Def}

\begin{Lem}\label{lem:pomstab}
Let $K$ be the pushout of $L\leftarrow A\xrightarrow{i} B$.
Assume that $A$, $B$ and $L$ are meta-stable orthogonal spectra, that $i$
is an l-cofibration and $A$ is well-pointed. Then $K$ is also meta-stable.
\end{Lem}

\begin{proof}
Consider the diagram
\[
\begin{CD}
L(\R^n)\wedge S^1 @<<< A(\R^n)\wedge S^1 @>{i}>> B(\R^n)\wedge S^1\\
@VVV @VVV @VVV\\
L(\R^{n+1}) @<<< A(\R^{n+1}) @>{i}>> B(\R^{n+1})
\end{CD}\quad.
\]
Since $A$, $B$ and $L$ are bounded below, we may increase $n$ until all spaces in the diagram above are simply connected.
Therefore it is enough to consider homology when calculating connectedness. Comparing Mayer-Vietoris sequences for the two
rows, we get that $K(\R^n)\wedge S^1\rightarrow K(\R^{n+1})$ is $(2n+d)$-connected, 
assuming that all three vertical maps in
the diagram above also have this connectedness.
\end{proof}

\begin{Lem}\label{Lem:VFAmstab}
Let $L=\bigvee_{\alpha} F_{V_{\alpha}}A_{\alpha}$. If all $A_{\alpha}$ are well-pointed
and there exists a $k$ such that $\dim V_{\alpha}\leq k$ for all $\alpha$, then
$L$ is meta-stable.
\end{Lem}

\begin{proof}
First we consider the case with a single wedge summand.
Assume that $n>k$. Then we have by lemma~\ref{Lem:VFAconn} that both 
$F_{\R^k} A(\R^n)\wedge S^1$ and $F_{\R^k} A(\R^{n+1})$ are simply connected.
Thus we can calculate the connectivity of $\sigma:F_{\R^k} A(\R^n)\wedge S^1\rightarrow F_{\R^k}(\R^{n+1})$
using homology.
By the definition of $F_{\R^k}$ the map under consideration is
\[
O(n)_+\wedge_{O(n-k)}(A\wedge S^{n-k})\wedge S^1\rightarrow O(n+1)_+\wedge_{O(n-k+1)}(A\wedge S^{n-k+1})\quad.
\]
We see that this map is $l$-connected if and only if the map
\[
O(n)_+\wedge_{O(n-k)}(A\wedge S^{n-k})\wedge S^1\wedge S^k\rightarrow O(n+1)_+\wedge_{O(n-k+1)}(A\wedge S^{n-k+1})\wedge S^k
\]
is $(l+k)$-connected. Now observe that the $O(n-k+1)$-action on $(A\wedge S^{n-k+1})\wedge S^k$ can be extended to a
$O(n+1)$-action, and thus we have
\[
O(n+1)_+\wedge_{O(n-k+1)}(A\wedge S^{n-k+1})\wedge S^k\cong O(n+1)/O(n-k+1)_+\wedge (A\wedge S^{n+1})\quad,
\]
and similarly for the source space. Thus we are considering the connectivity of
\[
O(n)/O(n-k)_+\wedge (A\wedge S^{n+1}) \rightarrow O(n+1)/O(n-k+1)_+\wedge (A\wedge S^{n+1})\quad.
\]
This map is easily seen to be $(2n-k)$-connected. And it follows that $\sigma$ is $(2n-2k)$-connected for $n>k$.
An inspection of the case $n=k$ shows that $F_{\R^k} A(\R^{k+1})$ is $0$-connected, and it follows that
$\sigma$ is $0$-connected. Hence we can take $d$ in the definition of meta-stability 
to be $(-2k)$.

In the general case we observe that the suspension map is the composition
\[
L(\R^n)\wedge S^1\cong \bigvee_{\alpha} (F_{V_{\alpha}}A_{\alpha}(\R^n)\wedge S^1)\xrightarrow{\bigvee\sigma}
\bigvee_{\alpha} F_{V_{\alpha}}A_{\alpha}(\R^{n+1})\quad.
\]
But the wedge of $(2n-2k)$-connected maps are $(2n-2k)$-connected. And it follows that $L$ is meta-stable.
\end{proof}

The following result is a useful property of relative $FI$-cellularity.

\begin{Prop}\label{Prop:mstabcellfiltr}
Assume that $A\rightarrow L$ is a relative $FI$-cellular map of orthogonal spectra. If $A$ is strictly bounded below,
then there exists a sequence $A=L_0\rightarrow L_1\rightarrow\cdots$ of l-cofibrations with colimit $L$, and
such that each $L_i$ is strictly bounded below.
If $A$ is meta-stable, then each $L_i$ can also be assumed meta-stable.
\end{Prop}

\begin{proof}
Let $C$ be the poset of cells. For a cell $\alpha$ let $m_{\alpha}$ denote the desuspension degree, i.e.
the dimension of $V$ in $F_VS^{n-1}_+\rightarrow F_V D^n_+$. Now define $C_i$ inductively:
\begin{itemize}
\item[] Let $C_0$ be those $\alpha\in C$ such that $m_{\alpha}=0$ and $P_{\alpha}=\emptyset$.
\item[] Given $C_0,\ldots,C_{i-1}$, let $C_i$ be the cells $\alpha$ not in $C_k$ for any $k<i$ such that
$m_{\alpha}\leq i$ and for any $\beta<\alpha$ there is an $l<i$ such that $\beta\in C_l$.
\end{itemize}
Since each $P_{\alpha}$ is finite, we see that $\bigcup_i C_i=C$.
Thus, the collection of $C_i$'s satisfies the conditions in the proof of proposition~\ref{Prop:seqcelldef}.
And by the proof we can then construct the sequence $A=L_0\rightarrow L_1\rightarrow\cdots$ with colimit $L$,
inductively. $L_i$ is defined by the pushout diagram
\[
\begin{CD}
\bigvee_{\alpha\in C_i}F_{\R^{m_{\alpha}}}S^{n_{\alpha}-1}_+@>>> \bigvee_{\alpha\in C_i}F_{\R^{m_{\alpha}}}D^{n_{\alpha}}_+\\
@VVV @VVV\\
L_{i-1} @>>> L_{i}
\end{CD}\quad.
\]
Since $m_{\alpha}\leq i$ for all $\alpha$ in $C_i$, it follows
by the lemma~\ref{Lem:VFAconn} and corollary~\ref{Cor:posbb} that each $L_i$ is 
strictly bounded below. If $A$ in addition is meta-stable, 
we can use the lemmas~\ref{lem:pomstab} and~\ref{Lem:VFAmstab} 
to show that each $L_i$ is meta-stable.
\end{proof}

Let us prove the following property of meta-stable orthogonal spectra:

\begin{Lem}\label{Lem:sigmastab}
If $L$ is meta-stable and well-pointed, then there exists a constant $e$ such that
the assembly induces a $(2n+k+e)$-connected map
\[
\Omega^k L(\R^{n+k})\rightarrow \Omega^{k+l} L(\R^{n+k+l})
\]
for all $n$, $k$ and $l$.
\end{Lem}

\begin{proof}
By induction we may reduce to the case where $l=1$.

Since $L(\R^{n+k})\wedge S^1\rightarrow L(\R^{n+k+1})$ is $(2n+2k+d)$-connected for 
some constant $d$ independent of
$n$ and $k$, 
it follows that
\[
\Omega \left( L(\R^{n+k})\wedge S^1\right)\rightarrow \Omega L(\R^{n+k+1})
\] 
is $(2n+2k+d-1)$-connected. There exists a constant $c$, also independent of $n$ and $k$, such that
$L(\R^{n+k})$ is $(n+k+c)$-connected. By Freudenthal's suspension theorem we have that
\[
L(\R^{n+k})\rightarrow \Omega \left( L(\R^{n+k})\wedge S^1\right)
\]
is $(2n+2k+2c+1)$-connected. Hence, there exists an $e$ such that the composite
\[
L(\R^{n+k})\rightarrow\Omega L(\R^{n+k+1})
\]
is $(2n+2k+e)$-connected. Applying $\Omega^k-$ we get that
\[
\Omega^k L(\R^{n+k})\rightarrow \Omega^{k+1} L(\R^{n+k+1})
\]
is $(2n+k+e)$-connected.
\end{proof}

\subsection{Induced functors on orthogonal spectra}\label{subsect:indfunct}

Consider a continuous endofunctor $F$ defined on based spaces. We already know from
example~\ref{exa:lwconstr} that applying $F$ level-wise to an orthogonal spectrum $L$
yields a new orthogonal spectrum $F(L)$.
What can be said about $F(L)$? The results we are looking for will compare
$F(L)$ and $G(L)$ whenever we have a natural transformation $f:F\rightarrow G$
of continuous endofunctors.

\begin{Lem}\label{Lem:Fpresbonuded}
Assume that there exists an integer $d$ such that $f_X:F(X)\rightarrow G(X)$
is $(2n+d)$-connected when $X$ is $n$-connected and well-pointed.
If $L$ is strictly bounded below and well-pointed, then
$L\rightarrow F(L)$ is a $\pi_*$-iso. 
\end{Lem}

\begin{proof}
The natural transformation $f$ induces a map of orthogonal spectra
\[
F(L)\rightarrow G(L)\quad.
\] 
We take the underlying prespectra, and look at the $q$'th homotopy groups:
\[
\pi_q F(L)=\colim_n\pi_{q+n}F(L(\R^n))\rightarrow \colim_n\pi_{q+n}G(L(\R^n))=\pi_q G(L)\quad.
\]
For a fixed $q$ and under the given assumptions on $f$ and $L$, 
the map $\pi_{q+n}F(L(\R^n))\rightarrow \pi_{q+n}G(L(\R^n))$ is eventually an
isomorphism: There exist integers $N$ and $c$ such that $L(\R^n)$ is $(n+c)$-connected
for $n\geq N$. Furthermore, there is a $d$ such that $F(X)\rightarrow G(X)$ is $(2n+d)$-connected
when $X$ is $n$-connected. Therefore, $F(L(\R^n))\rightarrow G(L(\R^n))$ is
$(2n+2c+d)$-connected for $n\geq N$, so $\pi_{q+n}F(L(\R^n))\rightarrow \pi_{q+n}G(L(\R^n))$
is an iso for $n>\max(q-2c-d,N)$.
\end{proof}

However, the condition we usually want to assume is cofibrancy, not strictly bounded below.
Therefore we use proposition~\ref{Prop:mstabcellfiltr} to transform the lemma above.

\begin{Cor}\label{Cor:indfunct}
Let $f:F\rightarrow G$ be a natural transformation of endofunctors on $\Top_*$.
Assume that
\begin{itemize}
\item[] there exists an integer $d$ such that $f_X$ is $(2n+d)$-connected when $X$ is $n$-connected and well-pointed,
\item[] $F$ and $G$ preserves cofibrations of spaces, and
\item[] if $X_0\rightarrow X_1\rightarrow X_2\rightarrow \cdots$ is any sequence of cofibrations of spaces, then
the natural maps $\colim_i F(X_i)\rightarrow F(\colim_i X_i)$ and $\colim_i G(X_i)\rightarrow G(\colim_i X_i)$
are weak equivalences.
\end{itemize}
If $L$ is cofibrant, then the induced map
\[
F(L) \rightarrow G(L)
\]
is a $\pi_*$-isomorphism.
\end{Cor}

\begin{proof}
First suppose that $L$ is $FI$-cellular. By proposition~\ref{Prop:mstabcellfiltr} there is 
a sequence $*=L_0\rightarrow L_1\rightarrow$ of l-cofibrations with colimit $L$ such that each
$L_i$ is strictly bounded below.
Apply $f:F\rightarrow G$ to the sequence and compare:
\[
\begin{CD}
F(L_0)@>>> F(L_1) @>>> F(L_2) @>>> \cdots\\
@VVV @VVV @VVV\\
G(L_0)@>>> G(L_1) @>>> G(L_2) @>>> \cdots
\end{CD}\quad.
\]
All horizontal maps are l-cofibrations, 
while lemma~\ref{Lem:Fpresbonuded} implies that all vertical maps are $\pi_*$-isomorphisms.
It follows that $\colim_i F(L_i)\rightarrow \colim_i G(L_i)$ is a $\pi_*$-iso. 
By the last assumption, we get that $F(L)\rightarrow G(L)$ is a $\pi_*$-iso.

Now suppose that $L$ is any cofibrant orthogonal spectrum. Then $L$ is a retract of an
$FI$-cellular orthogonal spectrum $K$. It follows that $F(L)\rightarrow G(L)$ is a retract of a $\pi_*$-iso
$F(K)\rightarrow G(K)$, hence the first map is also a $\pi_*$-iso.
\end{proof}

The technique of induced functors can be extended to multi-functors $\Top_*^n\rightarrow\Top_*$.
But instead of giving the most general statement, we will illustrate this by
considering the example:
\[
i_{X,Y}:X\vee Y\rightarrow X\times Y\quad,\text{ for spaces $X$ and $Y$.}
\]

\begin{Prop}\label{Prop:veepiisotimes}
If $L$ and $K$ are cofibrant orthogonal spectra, then $L\vee K\rightarrow L\times K$
is a $\pi_*$-isomorphism.
\end{Prop}

\begin{proof}
Observe that $i_{X,Y}$ is $(n+m+1)$-connected if $X$ is $n$-connected, $Y$ is $m$-connected
and both spaces are well-pointed. Now assume that $L$ and $K$ are well-pointed and strictly bounded below.
Then there exists a $c$ such that $L(\R^n)$ and $K(\R^n)$ are $(n+c)$-connected. It follows that
\[
\pi_{q+n}(L(\R^n)\vee K(\R^n))\rightarrow \pi_{q+n}(L(\R^n)\times K(\R^n))
\]
is an isomorphism when $n\geq q-2c$. Hence $L\vee K\rightarrow L\times K$ is a $\pi_*$-iso.

To get the result in the general case, we use proposition~\ref{Prop:mstabcellfiltr} and the following observations:
\begin{itemize}
\item[] $L\mapsto L\vee K$ and $L\mapsto L\times K$ preserves l-cofibrations when $K$ is well-pointed,
\item[] $(\colim_i L_i)\vee K=\colim_i (L_i\vee K)$, and
\item[] $(\colim_i L_i)\times K=\colim_i (L_i\times K)$.
\end{itemize}
First assume that $L$ is cofibrant, and $K$ is strictly bounded below and well-pointed. Filtrating $L$ with each
$L_i$ being strictly bounded below and well-pointed, we see that $L\vee K\rightarrow L\times K$ is a $\pi_*$-iso.
Next, assume that $L$ and $K$ are both cofibrant, filtrate $K$ and use the previous sentence to finish the proof.
\end{proof}

\section{$\mathscr{IS}$ as a model category}

Model categories (=closed model categories)
were introduced by Quillen~\cite{Quillen:67} and~\cite{Quillen:69} as
an axiomatization of homotopy theory. See also the survey article~\cite{DwyerSpalinski:95} or
the book~\cite{Hirschhorn:03}.
We will recall the definition of a model category below.
Mandell, May, Schwede and Shipley~\cite{MandellMaySchwedeShipley:01}
show that the category of orthogonal spectra has several model structures.
We will explain this. The section ends by listing various results
concerning the model category theory of orthogonal spectra. 

\begin{Def}
A \textit{model category} is a category $\mathscr{C}$ with three
distinguished classes of maps: \textit{weak equivalences}, 
\textit{fibrations} and \textit{cofibrations}. Each of these classes
is closed under composition and contains all identity maps.
A map which is both a fibration (resp. cofibration) and
a weak equivalence is called an \textit{acyclic fibration} 
(resp. \textit{acyclic cofibration}). And we have the following
axioms:
\begin{itemize}
\item[\textbf{MC1}] Finite limits and colimits exist in $\mathscr{C}$.
\item[\textbf{MC2}] If $A\xrightarrow{f} B\xrightarrow{g} C$ are composable
maps in $\mathscr{C}$, and if two of the three maps $f$, $g$ and $gf$ are weak
equivalences, then so is the third.
\item[\textbf{MC3}] If $f$ is a retract of $g$ and $g$ is a fibration, cofibration
or weak equivalence, then so is $f$.
\item[\textbf{MC4}] Given a commutative diagram
\[
\begin{CD}
A @>{f}>> X\\
@V{i}VV @VV{p}V\\
B @>{g}>> Y
\end{CD}\quad,
\]
then there exists a lift $h:B\rightarrow X$ such that $hi=f$ and $ph=g$ in the 
following two situations: when $i$ is a cofibration and $p$ is an acyclic fibration,
or when $i$ is an acyclic cofibration and $p$ is a fibration.
\item[\textbf{MC5}] Any map $f$ can be factored in two ways: as $f=pi$ where
$i$ is a cofibration and $p$ is an acyclic fibration, and as $f=pi$ where
$i$ is an acyclic cofibration and $p$ is a fibration.
\end{itemize}
A \textit{model structure} on a category $\mathscr{C}$ is a model category
with $\mathscr{C}$ as its underlying category.
\end{Def}

Axiom \textbf{MC4} gives liftings of certain diagrams. Since liftings are important
in model category theory, we have the following standard terminology:
A map $i:A\rightarrow B$ has the \textit{left lifting property} with respect to
another map $p:X\rightarrow Y$ if for any diagram of the same form as the diagram in
\textbf{MC4}, there exists a lift $h:B\rightarrow X$ such that $hi=f$ and $ph=g$.
Dually, we say that $p:X\rightarrow Y$ has the \textit{right lifting property}
with respect to $i:A\rightarrow B$.

A basic result about model categories, proposition~3.13 in~\cite{DwyerSpalinski:95},
says that $i$ is a cofibration if and only if it has the left lifting property with respect to all acyclic
fibrations. Dually, $p$ is a fibration if and only if it has the right lifting property
with respect to all acyclic cofibrations. Hence, when specifying a model category
it is enough to define the weak equivalences and either cofibrations or fibrations.
The remaining class is determined by lifting properties.

Model categories often come with extra structure. For example we have \textit{simplicial model categories},
see \S{}II.3 in~\cite{GoerssJardine:99}. More relevant
to us are \textit{topological model categories}, see \S{}5 in~\cite{MandellMaySchwedeShipley:01}.

There are several model structures on the category of orthogonal spectra.
The first model structure is:

\begin{Def}
The \textit{level model structure} on orthogonal spectra is given by setting
\begin{itemize}
\item[] $f:K\rightarrow L$ is a weak equivalence if $f$ is a level equivalence,
\item[] $f:K\rightarrow L$ is a cofibration if $f$ is a q-cofibration (=retract of relative $FI$-cellular map), and
\item[] $f:K\rightarrow L$ is a fibration if for each level $V$ the map $f:K(V)\rightarrow L(V)$ is a Serre fibration.
\end{itemize}
\end{Def}

Theorem~6.5 in~\cite{MandellMaySchwedeShipley:01} says that the level model structure on orthogonal spectra
is a model structure. Next we have:

\begin{Def}\label{Def:osstablemodelstr}
The \textit{stable model structure} on orthogonal spectra is given by setting
\begin{itemize}
\item[] $f:K\rightarrow L$ is a weak equivalence if $f$ is a $\pi_*$-isomorphism,
\item[] $f:K\rightarrow L$ is a cofibration if $f$ is a q-cofibration (=retract of relative $FI$-cellular map), and
\item[] $f:K\rightarrow L$ is a fibration if $f$ has the right lifting property with respect to the 
acyclic cofibrations (=maps which are both $\pi_*$-isomorphisms and q-cofibrations).
\end{itemize}
\end{Def}

Theorem~9.2 in~\cite{MandellMaySchwedeShipley:01} says that the stable model structure on orthogonal spectra
is a model structure. 

\begin{Rem}
There is also a \textit{positive stable model structure} on the category of
orthogonal spectra. The weak equivalences of this model structure are the $\pi_*$-isomorphisms. 
There are fewer cofibrations than the previous model structures, because in the positive stable structure
one does not allow cells $F_{\R^m}S^{n-1}_+\rightarrow F_{\R^m}D^n_+$ with $m=0$.
The fibrations are defined via the right lifting property. Theorem~14.2 in~\cite{MandellMaySchwedeShipley:01}
verifies that the positive stable model structure is a model structure.

The purpose of the positive stable models structure is to study commutative orthogonal ring spectra.
However, in this thesis we study orthogonal ring spectra with involution, and they are rarely commutative.
Hence we do not need the positive stable model structure.
\end{Rem}

We now list miscellaneous results:

\begin{Prop}\label{Prop:os:gluepiiso}
Consider the diagram
\[
\begin{CD}
L @<{i}<< A @>>> K\\
@V{\simeq}VV @VV{\simeq}V @VV{\simeq}V\\
L' @<{i'}<< A' @>>> K'
\end{CD}
\]
where $i$ and $i'$ are h-cofibrations and the vertical maps are $\pi_*$-isos. Then
the map of the row-wise pushouts is also a $\pi_*$-iso. 
\end{Prop}

For a proof see theorem~8.12(iv) in~\cite{MandellMaySchwedeShipley:01}.
The result also holds when $i$ and $i'$ are q-cofibrations, since any q-cofibration is an h-cofibration.

\begin{Prop}\label{Prop:os:gluepiisolcof}
The conclusion of proposition~\ref{Prop:os:gluepiiso} also holds if we assume that $i$ and $i'$ are
l-cofibrations instead of h-cofibrations.
\end{Prop}

\begin{proof}
Apply the cofibrant replacement functor to the diagram of proposition~\ref{Prop:os:gluepiiso}. In the resulting diagram,
\[
\begin{CD}
\Gamma L @<{\Gamma i}<< \Gamma A @>>> \Gamma K\\
@V{\simeq}VV @VV{\simeq}V @VV{\simeq}V\\
\Gamma L' @<{\Gamma i'}<< \Gamma A' @>>> \Gamma K'
\end{CD}\quad,
\]
we observe that $\Gamma i$ and $\Gamma i'$ are q-cofibrations, by theorem~\ref{Thm:cofrepl}, and the vertical maps
are $\pi_*$-isomorphisms. From the proposition above it follows that
\[
\Gamma L\cup_{\Gamma A}\Gamma K\rightarrow \Gamma L'\cup_{\Gamma A'}\Gamma K'
\]
is a $\pi_*$-iso. Now inspect the diagram
\[
\begin{CD}
\Gamma L @<{\Gamma i}<< \Gamma A @>>> \Gamma K\\
@V{\simeq_l}VV @VV{\simeq_l}V @VV{\simeq_l}V\\
L @<{i}<< A @>>> K
\end{CD}\quad.
\]
Here the vertical maps are level-equivalences. And since the pushout is formed level-wise, we may evaluate at some level $V$
and use the gluing theorem for weak equivalences between spaces, proposition~\ref{Prop:gluewe},
to conclude that
\[
\Gamma L\cup_{\Gamma A}\Gamma K\rightarrow L\cup_A K
\]
is a level-equivalence. Similarly the map $\Gamma L'\cup_{\Gamma A'}\Gamma K'\rightarrow L'\cup_{A'} K'$
is also a level-equivalence. Now consider the commutative square
\[
\begin{CD}
\Gamma L\cup_{\Gamma A}\Gamma K @>{\simeq_l}>> L\cup_A K\\
@V{\simeq}VV @VVV\\
\Gamma L'\cup_{\Gamma A'}\Gamma K' @>{\simeq_l}>> L'\cup_{A'} K'
\end{CD}\quad.
\]
We see that the last map must be a $\pi_*$-iso, and we are done.
\end{proof}

\begin{Prop}\label{Prop:os:wedge}
If $X$ is cofibrant and $K\rightarrow L$ is a $\pi_*$-iso, then also $K\wedge X\rightarrow L\wedge X$ is a
$\pi_*$-iso.
\end{Prop}

For a proof see proposition~12.3 in~\cite{MandellMaySchwedeShipley:01}.

\begin{Prop}\label{Prop:os:qcofsquare}
If $f:A\rightarrow L$ and $g:B\rightarrow K$ are q-cofibrations, then
\[
f\square g: A\wedge K\cup L\wedge B\rightarrow L\wedge K
\]
is also a q-cofibration. Furthermore, if $f$ or $g$ is in addition a $\pi_*$-iso, then $f\square g$ is also a $\pi_*$-iso.
\end{Prop}

The first part is a corollary of proposition~\ref{Prop:prodcell}, this is also lemma~6.6 in~\cite{MandellMaySchwedeShipley:01}.
The last part is the pushout-product axiom, proposition~12.6 in~\cite{MandellMaySchwedeShipley:01}.

\begin{Prop}\label{Prop:squarepiiso}
If $f:A\rightarrow L$, $g:B\rightarrow K$ and $g':B'\rightarrow K'$ are q-cofibrations, and $A$ is cofibrant 
and there is a commutative diagram
\[
\begin{CD}
B @>{g}>> K\\
@V{\simeq}VV @VV{\simeq}V\\
B' @>{g'}>> K'
\end{CD}
\]
where the vertical maps are $\pi_*$-isomorphisms, then the vertical maps in the diagram
\[
\begin{CD}
A\wedge K\cup L\wedge B @>{f\square g}>> L\wedge K\\
@V{\simeq}VV @VV{\simeq}V\\
A\wedge K'\cup L\wedge B' @>{f\square g'}>> L\wedge K'
\end{CD}
\]
are also $\pi_*$-isomorphisms.
\end{Prop}

\begin{proof}
Since $L$ is cofibrant,
the map $L\wedge K\rightarrow L\wedge K'$ is a $\pi_*$-iso by
proposition~\ref{Prop:os:wedge}.

Let $h$ be the map $A\wedge K\cup L\wedge B\rightarrow A\wedge K'\cup L\wedge B'$.
Notice that $h$ is the row-wise pushout of
\[
\begin{CD}
A\wedge K @<<< A\wedge B @>>> L\wedge B\\
@VVV @VVV @VVV\\
A\wedge K' @<<< A\wedge B' @>>> L\wedge B'
\end{CD}\quad.
\]
Since $A$ and $L$ are cofibrant, the vertical maps are $\pi_*$-isos by proposition~\ref{Prop:os:wedge}.
The maps $A\wedge B\rightarrow A\wedge K$ and $A\wedge B'\rightarrow A\wedge K'$ are q-cofibrant by
proposition~\ref{Prop:os:qcofsquare}. Now we apply proposition~\ref{Prop:os:gluepiiso}
and conclude that $h$ also is a $\pi_*$-iso.
\end{proof}

\begin{Prop}\label{Prop:qcofincolim}
Assume that we have a map between two sequences of orthogonal spectra:
\[
\begin{CD}
K_0 @>>> K_1 @>>> K_2 @>>> \cdots\\
@VVV    @VVV     @VVV\\
L_0 @>>> L_1 @>>> L_2 @>>> \cdots
\end{CD}\quad.
\]
If $K_0\rightarrow L_0$ is a q-cofibration, and
$K_i\cup_{K_{i-1}}L_{i-1}\rightarrow L_i$ is a q-cofibration for every $i\geq0$, then
\[
\colim_i K_i\rightarrow \colim_i L_i
\]
also is a q-cofibration. In particular, the case where $K_i$ is constant equal to $L_0$ yields that
$L_0\rightarrow \colim_i L_i$ is a q-cofibration if each $L_{i-1}\rightarrow L_i$ is. 
\end{Prop}

\begin{proof}
This statement holds in all model categories, we give an abstract proof.
Recall that a map in a model category is a cofibration if and only if it has the left lifting property with respect to acyclic
fibrations, see~\cite{DwyerSpalinski:95}. Consider a diagram
\[
\begin{CD}
\colim_i K_i @>>> X\\
@VVV @VVV\\
\colim_i L_i @>>> B
\end{CD}\quad,
\]
where $X\rightarrow B$ is an acyclic fibration. Define $f_0:L_0\rightarrow X$ to be a lift in
\[
\begin{CD}
K_0 @>>> X\\
@VVV @VVV\\
L_0 @>>> B
\end{CD}\quad.
\]
Inductively, choose a lift $f_i:L_i\rightarrow X$ in the diagram
\[
\begin{CD}
K_i\cup_{K_{i-1}}L_{i-1} @>>> X\\
@VVV @VVV\\
L_i @>>> B
\end{CD}\quad.
\]
Now observe that $\colim_i f_i:\colim_i L_i\rightarrow X$ lifts the original left lifting problem.
\end{proof}

\begin{Lem}\label{Lem:FnSntoSpiiso}
The maps $F_{\R^n}S^n\rightarrow S$, adjoint to the homeomorphisms $S^n\rightarrow S(\R^n)$, 
are $\pi_*$-isomorphisms for all $n\geq0$.
\end{Lem}

This follows from lemma~8.6 in~\cite{MandellMaySchwedeShipley:01}.

\section{Simplicial orthogonal spectra}

We need to discuss simplicial orthogonal spectra and we will use the theory
of simplicial spaces as our guideline. 
In~\cite{Segal:74} Segal defines
what it means for a simplicial space to be good, and shows that good 
simplicial spaces behaves well with respect to geometric realization.
May has a similar definition, proper, in~\cite{May:72}. 
Using Lillig's union theorem~\cite{Lillig:73} one can prove that 
proper and good are equivalent notions.

Let us now define simplicial orthogonal spectra.

\begin{Def}
A \textit{simplicial orthogonal spectrum} is a functor
$L_{\bullet}:\catDelta^{\op}\rightarrow \mathscr{IS}$.
It is \textit{good} if each $s_i:L_q\rightarrow L_{q+1}$ is an l-cofibration.
\end{Def}

As usual there is a geometric realization functor $|-|$ from
simplicial orthogonal spectra to $\mathscr{IS}$. 
A quick definition of $|L_{\bullet}|$ 
is given by a coend: 
\[
|L_{\bullet}|=\int^{[q]\in\catDelta}
L_q\wedge \Delta^q_+\quad.
\]
It is easy to see that $|-|$ is the same as applying the geometric realization
of simplicial spaces level-wise. We therefore have the formula $|L_{\bullet}|(V)\cong |L_{\bullet}(V)|$.
There also is a presimplicial realization, $\|-\|$, given by identifying 
along injective maps of $\catDelta$ only. As above this construction is also
level-wise: We have that $\|L_{\bullet}\|(V)\cong \|L_{\bullet}(V)\|$.

There is a natural map $\|L_{\bullet}\|\rightarrow |L_{\bullet}|$ and we have the
following standard result:

\begin{Prop}\label{Prop:sos:comparerealizarions}
For a good simplicial orthogonal spectrum $L_{\bullet}$, the natural map
$\|L_{\bullet}\|\rightarrow |L_{\bullet}|$ is a level equivalence.
\end{Prop}

\begin{proof}
We evaluate at $V$ and apply the corresponding result
for simplicial spaces, proposition~A.1(iv) in~\cite{Segal:74}.
\end{proof}

Our sufficient criterion for $|K_{\bullet}|\rightarrow |L_{\bullet}|$
to be a $\pi_*$-isomorphism, is a bit harder to prove:

\begin{Prop}\label{Prop:piisoofrealiz}
Let $f:K_{\bullet}\rightarrow L_{\bullet}$ be a map of simplicial orthogonal spectra.
If $K_{\bullet}$ and $L_{\bullet}$ are good and the map 
$f_q:K_q\rightarrow L_q$ is a $\pi_*$-isomorphism for any $q$, then the induced map
$|f|:|K_{\bullet}|\rightarrow |L_{\bullet}|$ is also a $\pi_*$-isomorphism.
\end{Prop}

\begin{proof}
By the previous proposition it is enough to prove
that $\|f\|:\|K_{\bullet}\|\rightarrow \|L_{\bullet}\|$ is a $\pi_*$-isomorphism.

We have a filtration $F_q\|K_{\bullet}\|$ of $\|K_{\bullet}\|$ by skeleta, and 
pushout diagrams
\[
\begin{CD}
K_q\wedge \partial \Delta_+^{q} @>>> F_{q-1}\|K_{\bullet}\|\\
@VVV @VVV\\
K_q\wedge \Delta_+^{q} @>>> F_{q}\|K_{\bullet}\|
\end{CD}
\]
for each $q\geq 1$. It can be checked directly that the left vertical map is an l-cofibration,
and consequently the right vertical map is also an l-cofibration. There is a similar 
filtration for $\|L_{\bullet}\|$. We compare the two filtrations. By 
proposition~\ref{Prop:lcofseqpiiso} it is enough to show that each map
$F_q\|K_{\bullet}\|\rightarrow F_q\|L_{\bullet}\|$ is a $\pi_*$-isomorphism.
This is proved by induction:

$F_0\|K_{\bullet}\|=K_0\xrightarrow{f_0}L_0=F_0\|L_{\bullet}\|$
is a $\pi_*$-isomorphism by assumption.

For the induction step we consider the diagram
\[
\begin{CD}
K_q\wedge \Delta_+^{q} @<<< K_q\wedge \partial \Delta_+^{q} @>>> F_{q-1}\|K_{\bullet}\|\\
@VVV @VVV @VVV\\
L_q\wedge \Delta_+^{q} @<<< L_q\wedge \partial \Delta_+^{q} @>>> F_{q-1}\|L_{\bullet}\|
\end{CD}\quad.
\]
Here the first two vertical maps are $\pi_*$-isomorphisms by proposition~\ref{Prop:os:wedge},
and the last vertical map is a $\pi_*$-iso by induction. By the
gluing lemma, proposition~\ref{Prop:os:gluepiisolcof}, we get that
$F_q\|K_{\bullet}\|\rightarrow F_q\|L_{\bullet}\|$ is a $\pi_*$-iso.
\end{proof}

An important feature of simplicial sets and simplicial spaces is that realization
commutes with products. This also holds for simplicial orthogonal spectra:

\begin{Lem}\label{Lem:sosrealizprod}
The category of simplicial orthogonal spectra, $s\mathscr{IS}$, is symmetric monoidal under the
product sending $K_{\bullet}$ and $L_{\bullet}$ to $[q]\mapsto K_q\wedge L_q$.
And geometric realization is strong symmetric monoidal. In particular there
is a natural isomorphism
\[
|K_{\bullet}|\wedge |L_{\bullet}| \rightarrow |K_{\bullet}\wedge L_{\bullet}|
\quad.
\]
\end{Lem}

\begin{proof}
It is clear that $s\mathscr{IS}$ is symmetric monoidal, with unit $[q]\mapsto S$.
To check that $|-|$ is strong symmetric monoidal, we first check the corresponding
statement for the geometric realization of simplicial $\mathscr{I}$-spaces. Here the product is $\tilde{\wedge}$.
Let $V$ be a finite dimensional real vector space, and evaluate.
We have:
\[
\left( |K_{\bullet}|\tilde{\wedge}|L_{\bullet}|\right)(V)\cong
\mathscr{I}(\R^{d}\oplus\R^{d'},V)_+\wedge_{O(d)\times O(d')}(|K_{\bullet}(\R^d)|\wedge |L_{\bullet}(\R^{d'})|)
\]
and
\[
|K_{\bullet}\tilde{\wedge}L_{\bullet}|(V)\cong
|\mathscr{I}(\R^{d}\oplus\R^{d'},V)_+\wedge_{O(d)\times O(d')}(K_{\bullet}(\R^d)\wedge L_{\bullet}(\R^{d'}))|
\quad.
\]
Since realization is a strong symmetric monoidal functor from simplicial spaces to spaces,
these formulas imply that $|-|$ is strong symmetric monoidal on simplicial $\mathscr{I}$-spaces.
Now the result also follows for simplicial orthogonal spectra by inspecting the coequalizer
definition of $\wedge$. We have
\[
\begin{CD}
|K_{\bullet}|\tilde{\wedge}|S|\tilde{\wedge} |L_{\bullet}| &\rightrightarrows & |K_{\bullet}|\tilde{\wedge} |L_{\bullet}| 
&\rightarrow& |K_{\bullet}|\wedge |L_{\bullet}|\\
@V{\cong}VV @V{\cong}VV @VVV\\
|K_{\bullet}\tilde{\wedge}S\tilde{\wedge} L_{\bullet}| &\rightrightarrows & |K_{\bullet}\tilde{\wedge} L_{\bullet}| 
&\rightarrow& |K_{\bullet}\wedge L_{\bullet}|
\end{CD}\quad,
\]
and it follows that the last map is an isomorphism.
\end{proof}

\section{The fibrant replacement functor $Q$}

We will need the underlying infinite loop space of an orthogonal spectrum 
in order to define $\widehat{GL}$ and consequently also $K$-theory.
Getting the underlying infinite loop space is a two step process. 
First there is a functor, which we will call $Q$,
that tries to turn orthogonal spectra into $\Omega$-spectra. To get the underlying infinite loop space of $L$, one then 
picks out the $0$'th space of $QL$.

The classical idea for constructing $QL$ is to take the homotopy colimit of $\Omega^n L(\R^n\oplus -)$
as $n$ grows to infinity. This construction would give the correct homotopy, at least when $L$ is suitably nice,
but the monoidal properties with respect to $\wedge$ are bad. B{\"o}kstedt solved this problem by instead
considering a homotopy colimit over the category of finite sets and injections. See the proof of lemma~2.3.7 in~\cite{Madsen:94}.
When $n$ lives in this category the functor $n\mapsto \Omega^n L(\R^n\oplus -)$ has monoidal properties, and consequently
also its homotopy colimit. However, for the purpose of constructing a fibrant replacement functor of orthogonal spectra
with monoidal properties, it is more natural to let the indexing category be finite dimensional real inner product spaces
and isometric embeddings.

When reading the proof of proposition~8.8 in~\cite{MandellMaySchwedeShipley:01} or the proof of theorem~3.1.11 
in~\cite{HoveyShipleySmith:00}, one can get the impression that the construction indicated above does not
yield a fibrant replacement functor. But their problem is closely tied to symmetric spectra, 
rather than with the construction. In that category of spectra it is not true that the $FI$-cells are meta-stable. 
For example consider the symmetric spectrum $F_1S^1$, see example~3.1.10 in~\cite{HoveyShipleySmith:00}.

We now state the structure theorem for $Q$:

\begin{Thm}\label{Thm:fibrrepl}
There is an endofunctor $Q$ on orthogonal spectra having the following properties:
\begin{itemize}
\item[] $QL$ is an $\Omega$-spectrum if $L$ is well-pointed.
\item[] $Q$ preserves l-cofibrations of well-pointed orthogonal spectra.
\item[] $Q$ commutes with sequential colimits. 
\item[] If $K\rightarrow L$ is a $\pi_*$-isomorphism and $L$ and $K$ are well-pointed, then $QK\rightarrow QL$ is a level-equivalence.
\item[] There is a natural inclusion $\eta_L:L\rightarrow QL$, this is a $\pi_*$-iso if $L$ is well-pointed.
\item[] There is a natural map $\mu_{L,K}:QL\wedge QK\rightarrow Q(L\wedge K)$ such
that $\mu_{L,K}\circ(\eta_L\wedge\eta_K)=\eta_{L\wedge K}$.
\item[] With $\eta_S$ and $\mu$ the functor $Q$ is lax monoidal with respect to $\wedge$.
\item[] There is a natural transformation $\iota_L:QL\rightarrow QL$ such that $\iota^2=\id$, $\iota\eta=\eta$, $\iota$ is 
level equivalent to $\id$ when $L$ is well-pointed and the following diagram commutes:
\[
\begin{CD}
QL\wedge QK @>{\iota\wedge\iota}>> QL\wedge QK @>{\text{twist}}>> QK\wedge QL\\
@VVV && @VVV\\
Q(L\wedge K) @>{\iota}>> Q(L\wedge K) @>{\text{twist}}>> Q(K\wedge L)
\end{CD}\quad.
\]
\item[] There is a natural map $\alpha_{L,K}:QL\times QK\rightarrow Q(L\times K)$, this is a level equivalence if $L$ and $K$
are well-pointed.
\item[] With $\eta_*$ and $\alpha$ the functor $Q$ is lax monoidal with respect to $\times$.
\end{itemize}
\end{Thm}

\begin{Rem}\label{Rem:ISnotconvenient}
\textbf{Warning:} $Q$ is not symmetric. That would lead to a contradiction: The $0$'th space of $QS$ would then 
be a commutative topological monoid with unit and zero elements, and have the homotopy type of $\Omega^{\infty}\Sigma^{\infty}S^0$.
By a result of Moore, \cite{Moore:58}, this would imply that the path component of the unit of $\Omega^{\infty}\Sigma^{\infty}S^0$ 
has the homotopy type of a cross product of
Eilenberg-MacLane spaces. This is not true.

The fact that $Q$ is not symmetric is the precise point where the category of orthogonal spectra fails to 
be a ``convenient category of spectra'' as defined by Lewis, \cite{Lewis:91}.

Our theorem above is therefore the best possible result regarding the fibrant replacement functor: It is lax
skew-symmetric monoidal, and the involution $\iota$ is homotopic to $\id$.
\end{Rem}

Let us now look into the construction of $QL$. First we let $\mathscr{E}$ be the topological category
of the finite dimensional real inner product spaces $\R^n$, $n\geq0$, and isometric embeddings. 
Given an orthogonal spectrum $L$, we have a continuous functor from $\mathscr{E}$ to $\mathscr{IS}$
given by
\[
W\mapsto \Omega^W L(W\oplus -)\quad.
\]
For morphisms in $\mathscr{E}$ from $W$ to $U$ we define the map
\[
\mathscr{E}(W,U)_+\wedge \Omega^W L(W\oplus V)\rightarrow \Omega^U L(U\oplus V)
\]
as follows: Assume $f:W\rightarrow U$ is an isometric embedding and $\alpha:S^W\rightarrow L(W\oplus V)$
represent a point in $\Omega^W L(W\oplus V)$. Let $d$ be the codimension of $W$ in $U$. Now consider the composition
\[
S^U\cong S^d\wedge S^W\xrightarrow{\id\wedge\alpha} S^d\wedge L(W\oplus V)\xrightarrow{\bar{\sigma}}
L(\R^d\oplus W\oplus V)\cong L(U\oplus V)\quad,
\]
where the first and last map is induced by $f$, and $\bar{\sigma}$ is the left assembly.
This composition represents a point in $\Omega^U L(U\oplus V)$.

We have chosen to work with the Euclidean 
spaces $\R^n$, $n\geq 0$, instead of all finite dimensional real inner product spaces. 
The reason is that we would like to take a ``homotopy colimit'' over $\mathscr{E}$.
Therefore, the objects should be a set.
Define $Q_{\bullet}L$ to be the simplicial orthogonal spectrum given by
\[
Q_qL(V)=\bigvee_{n_0,\ldots,n_q\geq0}\mathscr{E}(\R^{n_{q-1}},\R^{n_q})_+\wedge\cdots\wedge\mathscr{E}(\R^{n_{0}},\R^{n_1})_+
\wedge\Omega^{n_0} L(\R^{n_0}\oplus V)\quad.
\]
The face and degeneracy maps are given by
\begin{align*}
d_i(f_{q-1},\ldots,f_0;\alpha)&=\begin{cases}
(f_{q-1},\ldots,f_1;f_0(\alpha))&\text{for $i=0$,}\\
(f_{q-1},\ldots,f_{i+1},f_{i}\circ f_{i-1},f_{i-2},\ldots,f_0;\alpha)&\text{for $0<i<q$,}\\
(f_{q-2},\ldots,f_0;\alpha)&\text{for $i=q$, and}
\end{cases}\\
s_i(f_{q-1},\ldots,f_0;\alpha)&=(f_{q-1},\ldots,f_i,\id_{\R^{n_i}},f_{i-1},\ldots,f_0;\alpha)\quad.
\end{align*}

We now define:

\begin{Def}
The functor $Q$ is an endofunctor on orthogonal spectra, given on $L$ as 
the geometric realization of $Q_{\bullet}L$.
\end{Def}

Clearly there is a natural inclusion $\eta_L:L\rightarrow QL$. This comes from the inclusion of
$L(V)$ as the wedge summand of $Q_0L(V)=\bigvee_{n_0\geq0}\Omega^{n_0} L(\R^{n_0}\oplus V)$
corresponding to $n_0=0$.

\begin{Lem}
If $L$ is well-pointed, then $Q_{\bullet}L$ is a good simplicial orthogonal spectrum.
\end{Lem}

\begin{proof}
The space $\mathscr{E}(\R^n,\R^n)$ is well-pointed at $\id_{\R^n}$, it is even a smooth manifold.
By assumption $\Omega^{n_0} L(\R^{n_0}\oplus V)$ is well-pointed. By applying the smash product theorem for
well-pointed cofibrations in $\Top_*$, proposition~12 in~\cite{Strom:72}, we get that each $s_i$ is an l-cofibration.
\end{proof}

\begin{Lem}\label{lem:Qpresle}
$Q$ preserves level equivalences between well-pointed orthogonal spectra.
\end{Lem}

\begin{proof}
If $K\rightarrow L$ is a level equivalence, it follows that 
$\Omega^nK(\R^n\oplus V)\rightarrow \Omega^nL(\R^n\oplus V)$
is a weak equivalence for all $n$ and $V$. Hence in each simplicial degree $q$ the map
\[
Q_q(K)(V)\rightarrow Q_q(L)(V)
\]
is a weak equivalence. The result follows since both $Q_{\bullet}K$ and $Q_{\bullet}L$ are good.
\end{proof}

\begin{Prop}
$Q$ commutes with sequential colimits.
\end{Prop}

\begin{proof}
Let $L_0\rightarrow L_1\rightarrow L_2\rightarrow \cdots$ be a sequence of orthogonal spectra
with colimit $L=\colim_i L_i$.
Colimits commute with geometric realization, thus it is enough to show that
\[
\colim_i Q_q L_i =Q_q L
\]
for every simplicial degree $q$.
Sequential colimits are level-wise constructions, so it is enough to check that
the colimit of $Q_qL_0(V)\rightarrow Q_qL_1(V)\rightarrow \cdots$ is $Q_qL(V)$ for each $V$.
We now inspect the definition of $Q_q-(V)$. Colimits commute with wedge, but in general an
arbitrary colimit does not commute with smash products. However, sequential colimits (of based spaces)
commute with the functor $X\mapsto A\wedge X$, where $A$ is some fixed space. Therefore it is enough to show
that 
\[
\colim_i \Omega^{n_0}L_i(\R^{n_0}\oplus V) = \Omega^{n_0}L(\R^{n_0}\oplus V)
\]
for fixed $n_0$ and $V$. This is true since the colimit is sequential and $S^{n_0}$ compact.
\end{proof}

\begin{Prop}
$Q$ preserves l-cofibrations of well-pointed orthogonal spectra.
\end{Prop}

\begin{proof}
We start with an l-cofibration $A\rightarrow L$ of well-pointed orthogonal spectra. 
By proposition~\ref{Prop:lcofexpthm} the map
\[
\Omega^n A(\R^n\oplus V)\rightarrow \Omega^n L(\R^n\oplus V)
\]
is a closed cofibration of well-pointed spaces for all $n$ and $V$. 
By the smash product theorem for well-pointed cofibrations of spaces 
(proposition~12 in~\cite{Strom:72}) it follows that
$Q_q A\rightarrow Q_qL$ is an l-cofibration of well-pointed orthogonal spectra for every $q$.
Both $Q_{\bullet}A$ and $Q_{\bullet}L$ are good, and 
by the gluing theorem for l-cofibrations and the 
filtration of the geometric realization
it now follows that
$QA\rightarrow QL$ is an l-cofibration.
\end{proof}

\begin{Lem}\label{Lem:BEwcontr}
The classifying space $B\mathscr{E}$ is contractible.
\end{Lem}

\begin{proof}
Direct sum induces a map $\oplus:B\mathscr{E}\times B\mathscr{E}\rightarrow B\mathscr{E}$. 
Since there is a natural transformation from the projection $\mathscr{E}\times \mathscr{E}\rightarrow \mathscr{E}$ 
onto the first factor to the direct sum, we get a homotopy between $\oplus$ and the projection 
$\pr_1:B\mathscr{E}\times B\mathscr{E}\rightarrow B\mathscr{E}$. Similarly we get a homotopy
$\oplus\simeq\pr_2$.

Now choose a basepoint $*$ in $B\mathscr{E}$ and consider the composition of the inclusion
\[
i:B\mathscr{E}\times *\rightarrow B\mathscr{E}\times B\mathscr{E}
\] 
with $\pr_1$, $\oplus$ and $\pr_2$. We get homotopies
\[
\id=\pr_1\circ i\simeq \oplus\circ i\simeq \pr_2\circ i=*\quad.
\]
Thus $B\mathscr{E}$ is contractible.
\end{proof}

We will now start proving that $\eta_L$ is a $\pi_*$-iso when $L$ is well-pointed. The proof is
divided into three parts: First we show the result for meta-stable well-pointed $L$. Next we filter
any cofibrant $L$ as the colimit of orthogonal spectra of the first type. At last we use
cofibrant approximation to prove the general case.

\begin{Lem}\label{lem:Qsigmastab}
If $L$ is meta-stable and well-pointed, then $\eta:L\rightarrow QL$ is a $\pi_*$-isomorphism.
\end{Lem}

\begin{proof}
Evaluating at a level $V$ we land in topological spaces. Here we
also have unbased homotopy colimits. Let $\widetilde{QL(V)}$ be the geometric realization of the simplicial
space with $q$-simplexes given by:
\[
\widetilde{QL(V)}_{q}=\coprod_{n_0,\ldots,n_q\geq 0}\mathscr{E}(\R^{n_{q-1}},\R^{n_q})\times
\cdots\times\mathscr{E}(\R^{n_{0}},\R^{n_1})
\times\Omega^{n_0} L(\R^{n_0}\oplus V)\quad.
\]
Consider the diagram
\[
\begin{CD}
&& L(V) && \\
&& @VVV &&\\
B\mathscr{E} @>>> \widetilde{QL(V)} @>>> QL(V)\\
&& @VVV &&\\
&& B\mathscr{E} &&
\end{CD}\quad.
\]
By lemma~\ref{Lem:sigmastab} and proposition~\ref{prop:hocolimqf}, there exists a constant $d$ such that
$\widetilde{QL(\R^n)}\rightarrow B\mathscr{E}$ is a $(2n+d)$-quasi fibration. But $B\mathscr{E}$ is
contractible, therefore the map $L(\R^n)\rightarrow \widetilde{QL(\R^n)}$ is $(2n+d-1)$-connected.

The horizontal part is a cofibration sequence. 
Now we use corollary~\ref{Cor:collweakcontr} and
that $B\mathscr{E}$ is contractible to see that
\[
\pi_i(\widetilde{QL(\R^n)})\rightarrow \pi_i(QL(\R^n))
\]
is an isomorphism for all $i$.

Putting things together we see that there exists a constant $c$ such that
\[
L(\R^n)\rightarrow QL(\R^n)
\]
is $(2n+c)$-connected. And it follows that $L\rightarrow QL$ is a $\pi_*$-iso.
\end{proof}

\begin{Rem}
In the proof above we compared an unbased homotopy colimit over a topological category, $\widetilde{QL(V)}$, 
with the corresponding based homotopy colimit, $QL(V)$. In the more elementary
case where the category is discrete, there is a general result due to E. Dror Farjoun
that compares unbased and based homotopy colimits, see proposition~18.8.4 in~\cite{Hirschhorn:03}.
\end{Rem}

\begin{Lem}
If $L$ is cofibrant, then the natural map $L\rightarrow QL$ is a $\pi_*$-iso.
\end{Lem}

\begin{proof}
We first prove this when $L$ is cellular. By proposition~\ref{Prop:mstabcellfiltr}, there exists a sequence
$L_0\rightarrow L_1\rightarrow \cdots$ of orthogonal spectra with colimit $L$ such that
each $L_i$ is meta-stable and well-pointed, while the maps $L_i\rightarrow L_{i+1}$ are
l-cofibrations. Applying $Q$ to this sequence we get
\[
\begin{CD}
L_0@>>> L_1 @>>> \cdots @>>> L\\
@V{\simeq}VV  @V{\simeq}VV && @VVV\\
QL_0@>>> QL_1 @>>> \cdots @>>> QL
\end{CD}\quad.
\]
The vertical maps $L_i\rightarrow QL_i$ are $\pi_*$-isomorphism by lemma~\ref{lem:Qsigmastab}.
Since both sequences consist of l-cofibrations, it follows that $L\rightarrow QL$
also is a $\pi_*$-iso.

For the general case we use that a cofibrant $L$ is a retract of some $L'$ which is cellular.
It follows that $L\rightarrow QL$ is a retract of $L'\rightarrow QL'$. Thus $L\rightarrow QL$
is a $\pi_*$-iso.
\end{proof}

\begin{Prop}\label{prop:LQLstable}
If $L$ is well-pointed, then the natural map $L\rightarrow QL$ is a $\pi_*$-iso.
\end{Prop}

\begin{proof}
Consider the diagram
\[
\begin{CD}
\Gamma L @>>> Q\Gamma L\\
@VVV @VVV\\
L @>>> QL
\end{CD}\quad.
\]
The top map is a $\pi_*$-iso by the previous lemma. The left map is a level equivalence by theorem~\ref{Thm:cofrepl}.
Since $Q$ preserves level equivalences between well-pointed orthogonal spectra, lemma~\ref{lem:Qpresle}, the right map
is also a level equivalence. It follows that the bottom map also is a $\pi_*$-iso.
\end{proof}

Let $\mathscr{E}_k$ be the full subcategory of $\mathscr{E}$ having objects $\R^n$ for $n\geq k$.
We will compare $QL$ with the homotopy colimit of $\Omega^n L(\R^n\oplus -)$ over this subcategory. Let
$Q^kL$ be the geometric realization of the simplicial orthogonal spectrum whose $q$-simplexes are:
\[
Q_q^kL(V)=\bigvee_{n_0,\ldots,n_q\geq k}\mathscr{E}(\R^{n_0},\R^{n_1})_+\wedge\cdots\wedge\mathscr{E}(\R^{n_{q-1}},\R^{n_q})_+
\wedge\Omega^{n_0} L(V\oplus \R^{n_0})\quad.
\]

\begin{Lem}\label{Lem:Qklevelequiv}
The inclusion $Q^kL\rightarrow QL$ is a level equivalence for all $k$.
\end{Lem}

\begin{proof}
We will show that the functor $(\R^k\oplus-):\mathscr{E}\rightarrow\mathscr{E}_k$ induces a 
map $f:QL\rightarrow Q^kL$ which is a homotopy inverse to $\incl$.
First consider the composition $\incl\circ f: QL\rightarrow QL$.
This map is induced by $\R^n\mapsto \R^k\oplus\R^n$, considered as an endofunctor on $\mathscr{E}$.
But we have a natural transformation $\tau:\id_{\mathscr{E}}\rightarrow (\R^k\oplus-)$
which includes $\R^n$ as the last $n$-coordinates of $\R^{k+n}$.
By proposition~\ref{prop:hocolimhomotopy}, which also holds in the based case since
the formulas for the simplicial homotopy still work, we get a homotopy between
$\id_{QL}$ and $\incl\circ f: QL\rightarrow QL$ induced by $\tau$.

The opposite composition, 
$f\circ\incl$ is also induced by $(\R^k\oplus-)$ considered as an endofunctor on
$\mathscr{E}_k$, and the same natural transformation 
gives a homotopy $f\circ\incl\simeq \id_{Q^kL}$.
\end{proof}

\begin{Lem}
The classifying space $B\mathscr{E}_k$ is contractible.
\end{Lem}

\begin{proof}
The proof is similar to that of lemma~\ref{Lem:BEwcontr}. Notice that also in this case
there are natural transformations from the projections $\mathscr{E}_k\times\mathscr{E}_k\rightarrow\mathscr{E}_k$
to the direct sum.
\end{proof}

\begin{Lem}
If $L$ is meta-stable and well-pointed, then $QL$ is an $\Omega$-spectrum.
\end{Lem}

\begin{proof}
To prove this we have to show that $QL(\R^n)\rightarrow \Omega QL(\R^{n+1})$ is
a weak equivalence for all $n$. Fix $n$. Observe that by lemma~\ref{Lem:Qklevelequiv}
it is enough to show that for any $\lambda$ there is a $k$ such that
\[
Q^kL(\R^n)\rightarrow \Omega Q^kL(\R^{n+1})
\]
is $\lambda$-connected.

As in the proof of lemma~\ref{lem:Qsigmastab} we have 
unbased homotopy limits $\widetilde{Q^kL(V)}$ when $V$ is fixed.
From the diagram
\[
\begin{CD}
&& \Omega^k L(V\oplus\R^k) && \\
&& @VVV &&\\
B\mathscr{E}_k @>>> \widetilde{Q^k L(V)} @>>> Q^k L(V)\\
&& @VVV &&\\
&& B\mathscr{E}_k &&
\end{CD}
\]
we estimate the connectivity of the map $\Omega^k L(\R^{n+k})\rightarrow Q^kL(\R^n)$.
By lemma~\ref{Lem:sigmastab} and proposition~\ref{prop:hocolimqf} and the fact that $B\mathscr{E}_k$
is contractible, there exists a constant $d$ such that
\[
\Omega^k L(\R^{n+k})\rightarrow \widetilde{Q^k L(\R^n)}
\]
is $(2n+k+d)$-connected. Using that the horizontal part is a cofiber sequence and that $B\mathscr{E}_k$
is contractible, it follows by corollary~\ref{Cor:collweakcontr}
that
\[
\widetilde{Q^k L(\R^n)} \rightarrow Q^k L(\R^n)
\]
is a weak equivalence. Thus there exists a $d$ such that for all $n$ and $k$ the map
\[
\Omega^k L(\R^{n+k})\rightarrow Q^kL(\R^n)
\]
is $(2n+k+d)$-connected.
Now inspect the diagram
\[
\begin{CD}
\Omega^k L(\R^{n+k}) @>>> \Omega^{k+1} L(\R^{n+k+1})\\
@VVV @VVV\\
Q^kL(\R^n) @>>> \Omega Q^kL(\R^{n+1})
\end{CD}\quad.
\]
There is a constant $c$ such that the map on the top and the two vertical maps are
$(2n+k+c)$-connected. Thus we can, by increasing $k$, ensure that the map on the bottom is
$\lambda$-connected for any fixed $n$.
\end{proof}

\begin{Lem}
If $L$ is cofibrant, then $QL$ is an $\Omega$-spectrum.
\end{Lem}

\begin{proof}
First assume that $L$ is cellular. By proposition~\ref{Prop:mstabcellfiltr},
there exists a sequence $L_0\rightarrow L_1\rightarrow \cdots$
with colimit $L$ such that each $L_i$ is a meta-stable well-pointed orthogonal spectrum, and each map
$L_i\rightarrow L_{i+1}$ is an l-cofibration.
By the proposition above we see that $QL$ is the colimit of a sequence $QL_0\rightarrow QL_1\rightarrow \cdots$
of $\Omega$-spectra where the maps are l-cofibrations.
Now we have
\[
QL(V)\cong \colim_i QL_i(V)\xrightarrow{\simeq} \colim_i\Omega QL_i(V\oplus\R)\cong \Omega QL(V\oplus\R)\quad.
\]
The map in the middle is a weak equivalence since each $QL_i(V)\rightarrow \Omega QL_i(V\oplus\R)$ is, and
both colimits are sequential over unbased closed cofibrations. The last map is a homeomorphism since
the colimit system is sequential and $S^1$ is compact.

General cofibrant $L$ are retracts of some cellular $L'$. Thus $QL$ is the retract of some $\Omega$-spectrum $QL'$.
But then the map $QL(V)\rightarrow \Omega QL(V\oplus \R)$ is a retract of $QL'(V)\rightarrow \Omega QL'(V\oplus \R)$,
and the former map must be a weak equivalence since the latter already is.
\end{proof}

\begin{Prop}\label{prop:QLOmega}
If $L$ is well-pointed, then $QL$ is an $\Omega$-spectrum.
\end{Prop}

\begin{proof}
Consider the level equivalence
\[
Q\Gamma L\rightarrow QL\quad.
\]
Since $\Gamma L$ is cofibrant, we know that $Q\Gamma L$ is an $\Omega$-spectrum.
It is an elementary fact that a well-pointed orthogonal spectrum level equivalent to an $\Omega$-spectrum is itself an
$\Omega$-spectrum. 
\end{proof}

\begin{Cor}
If $f:L\rightarrow K$ is a $\pi_*$-iso between well-pointed orthogonal spectra, then $Qf:QL\rightarrow QK$
is a level equivalence.
\end{Cor}

\begin{proof}
This follows since a $\pi_*$-iso between $\Omega$-spectra is a level equivalence. See~\cite{MandellMaySchwedeShipley:01} lemma~8.11.
\end{proof}

We will now describe the monoidal structure of $Q$ with respect to $\wedge$.
To do this, we first define maps of orthogonal spectra
\[
\mu:\Omega^n L(\R^n\oplus -)\wedge \Omega^m K(\R^m\oplus -)\rightarrow \Omega^{n+m} (L\wedge K)(\R^{n+m}\oplus -)\quad,
\]
natural for $L$ and $K$ in $\mathscr{IS}$ and $\R^n$ and $\R^m$ in $\mathscr{E}$.
We take the external viewpoint of the smash product, and let the map
\[
\mu:\Omega^n L(\R^n\oplus V)\wedge \Omega^m K(\R^m\oplus W)\rightarrow \Omega^{n+m} (L\wedge K)(\R^{n+m}\oplus V\oplus W)
\]
be given by sending the point represented by $\alpha:S^n\rightarrow L(\R^n\oplus V)$ and $\beta:S^m\rightarrow K(\R^m\oplus W)$
to the point represented by the composition
\footnotesize
\[
S^{n+m}\xrightarrow{\alpha\wedge\beta}L(\R^n\oplus V)\wedge K(\R^m\oplus W)\rightarrow 
(L\wedge K)(\R^n\oplus V\oplus \R^m\oplus W)\cong (L\wedge K)(\R^{n+m}\oplus V\oplus W)\quad.
\]
\normalsize
To see that $\mu$ is a map of orthogonal spectra, we should check that it commutes with left and right assembly. 
Let us inspect the case of the right assembly. Since the identity on $L\wedge K$ is a spectrum map, the external viewpoint
gives that
\[
\begin{CD}
L(V)\wedge K(W)\wedge S^U @>>> (L\wedge K)(V\oplus W)\wedge S^U\\
@V{\id\wedge\sigma}VV @VV{\sigma}V\\
L(V)\wedge K(W\oplus U) @>>> (L\wedge K)(V\oplus W\oplus U)
\end{CD}
\]
commutes. Using this fact, it is easily checked by chasing that the following diagram also commutes:
\footnotesize
\[
\begin{CD}
\Omega^n L(\R^n\oplus V)\wedge \Omega^m K(\R^m\oplus W)\wedge S^U @>>> \Omega^{n+m} (L\wedge K)(\R^{n+m}\oplus V\oplus W)\wedge S^U \\
@V{\id\wedge\sigma}VV @VV{\sigma}V\\
\Omega^n L(\R^n\oplus V)\wedge \Omega^m K(\R^m\oplus W\oplus U) @>>> \Omega^{n+m} (L\wedge K)(\R^{n+m}\oplus V\oplus W\oplus U)
\end{CD}\quad.
\]
\normalsize
The other case is done similarly. We conclude that $\mu$ is a map of orthogonal spectra.

Naturality in $\mathscr{E}\times \mathscr{E}$ means that the map from $\Omega^n L(\R^n\oplus V)\wedge \Omega^m K(\R^m\oplus W)$
induced by linear isometries $f:\R^n\rightarrow\R^{n'}$ and $g:\R^m\rightarrow \R^{m'}$ corresponds to the map from
$\Omega^{n+m} (L\wedge K)(\R^{n+m}\oplus V\oplus W)$ induced by $f\oplus g$. It is enough to check this when $g$ is the identity on $\R^m$.
Let $d$ be the codimension of $f$ and $\alpha$ and $\beta$ as above. It all boils down to the commutativity of the following diagram:
\[
\begin{CD}
S^{n'+m}\cong S^d\wedge S^{n+m} \\
@V{\id\wedge\alpha\wedge\beta}VV \\
S^d\wedge L(\R^n\oplus V)\wedge K(\R^m\oplus W) 
@>>> S^d\wedge (L\wedge K)(\R^n\oplus V \oplus \R^m\oplus W)\\
@V{\bar{\sigma}\wedge\id}VV @VV{\bar{\sigma}}V \\
L(\R^d\oplus \R^n\oplus V)\wedge K(\R^m\oplus W) 
@>>> (L\wedge K)(\R^d\oplus \R^n\oplus V \oplus \R^m\oplus W)\\
&& @V{\cong}VV\\
&& (L\wedge K)(\R^{n'+m}\oplus V \oplus W)
\end{CD}\quad.
\]

\begin{Prop}
$\mu$ induces a natural transformation $QL\wedge QK\rightarrow Q(L\wedge K)$ and $Q$ becomes a lax monoidal functor 
with respect to $\wedge$. In addition we have $\mu_{L,K}\circ(\eta_L\wedge\eta_K)=\eta_{L\wedge K}$.
\end{Prop}

\begin{proof}
Inspecting the definition of $Q_{\bullet}L$, we see that $\mu$ together with direct sum 
$\mathscr{E}(\R^n,\R^{n'})\times\mathscr{E}(\R^m,\R^{m'})\rightarrow \mathscr{E}(\R^{n+m},\R^{n'+m'})$
give a simplicial map
\[
Q_{\bullet}L\wedge Q_{\bullet}K\rightarrow Q_{\bullet}(L\wedge K)\quad.
\]
The geometric realization gives the natural transformation $QL\wedge QK\rightarrow Q(L\wedge K)$.
Associativity follows from associativity of $\oplus$ on $\mathscr{E}$ and $\mu$ on $\Omega^nL(\R^n\oplus-)$.

The natural inclusion for the sphere spectrum, $\eta_S:S\subseteq QS$, satisfies left and right unity:
To see this one observes that 
$\mu:\Omega^n L(\R^n\oplus -)\wedge \Omega^m K(\R^m\oplus -)\rightarrow \Omega^{n+m} (L\wedge K)(\R^{n+m}\oplus -)$
is equal to the right assembly if $m=0$ and $K=S$, and equal to the left assembly if $n=0$ and $L=S$.
This gives left and right unity for $\mu$. To get form here to $Q$ we use in addition the unit of $\mathscr{E}$.

To check the last formula of the proposition, we observe that $\mu=\id$ when $n=m=0$.
\end{proof}

Recall from remark~\ref{Rem:ISnotconvenient} that $Q$ cannot be symmetric.
However we have:

\begin{Prop}
$Q$ is skew-symmetric. This means that
there is a natural transformation $\iota:QL\rightarrow QL$ such that $\iota^2=\id$ and the diagram
\[
\begin{CD}
QL\wedge QK @>{\iota\wedge\iota}>> QL\wedge QK @>{\text{twist}}>> QK\wedge QL\\
@VVV && @VVV\\
Q(L\wedge K) @>{\iota}>> Q(L\wedge K) @>{\text{twist}}>> Q(K\wedge L)
\end{CD}
\]
commutes. Furthermore, $\iota\eta=\eta$ and $\iota$ is level equivalent to the identity
when $L$ is well-pointed.
\end{Prop}

\begin{proof}
Let $r_n:\R^n\rightarrow \R^n$ be the isometries which reverses the standard basis.
Conjugating an isometric embedding $\R^n\rightarrow \R^m$ with $r_n$ and $r_m$
gives a functor
\[
\operatorname{conj}:\mathscr{E}\rightarrow \mathscr{E}\quad.
\]
Let $G$ denote the functor $\mathscr{E}\rightarrow \mathscr{IS}$ given by $\R^n\mapsto \Omega^n L(\R^n\oplus -)$.
We now construct a natural transformation $\iota$ from $G$ to $G\circ \operatorname{conj}$ by
sending a point $\alpha$ in $\Omega^n L(V\oplus\R^n)$ to $\iota(\alpha)$ defined by the commutativity of
\[
\begin{CD}
S^n @>{\alpha}>> L(V\oplus\R^n)\\
@V{S^{r_n}}VV @VV{L(\id_V\oplus r_n)}V\\
S^n @>{\iota(\alpha)}>> L(V\oplus\R^n)
\end{CD}\quad\quad\quad\quad\quad.
\]
This natural transformation gives a map of simplicial spaces:
\[
Q_{\bullet}L \xrightarrow{\iota} Q_{\bullet}L\quad.
\]
Taking the geometric realization we get the natural transformation we are looking for.
It is easily seen that $\iota^2=\id$.

To get the commutativity of the main diagram of the proposition, we observe that for
maps $\alpha: S^n\rightarrow L(V\oplus\R^n)$ and $\beta: S^m\rightarrow K(W\oplus\R^m)$
we have that
\[
\text{twist}\circ\iota(\alpha\wedge \beta)=\iota(\beta)\wedge\iota(\alpha)
\]
as maps $S^{n+m}\rightarrow (K\wedge L)(V\oplus W\oplus \R^{n+m})$. Here $\text{twist}$ is the
isomorphism $L\wedge K\cong K\wedge L$.

To see that $\iota$ is level equivalent to the identity, we inspect the diagram
\[
\begin{CD}
L @= L\\
@V{\eta}VV @VV{\eta}V\\
QL @>{\iota}>> QL
\end{CD}\quad.
\]
Commutativity follows since $r_0=\id$ on $\R^0$. When $L$ is well-pointed, both vertical maps are $\pi_*$-isomorphisms, 
so we see that $\iota$ induces the identity on homotopy groups. But $QL$ is an $\Omega$-spectrum, therefore it
follows that $\iota$ also induces the identity 
\[
\pi_q QL(V)\rightarrow \pi_q QL(V)
\]
for all levels $V$ and $q\geq0$. 
\end{proof}

\begin{Prop}
There is a natural map $\alpha_{L,K}:QL\times QK\rightarrow Q(L\times K)$, this is a level equivalence if $L$ and $K$
are well-pointed.
With $\eta_*$ and $\alpha$ the functor $Q$ is lax monoidal with respect to $\times$.
\end{Prop}

\begin{proof}
$\alpha$ is defined similarly to $\mu$. Given $(\beta_1,\beta_2)\in \Omega^n L(\R^n\oplus V)\times \Omega^m K(\R^m\oplus V)$
we can suspend $\beta_1$ to a point in $\Omega^{n+m} L(\R^{n+m}\oplus V)$ using the inclusion $i_1:\R^n\rightarrow \R^n\oplus\R^m=\R^{n+m}$
and suspend $\beta_2$ to a point in $\Omega^{n+m} K(\R^{n+m}\oplus V)$ using the inclusion $i_2:\R^m\rightarrow \R^n\oplus\R^m=\R^{n+m}$.
This gives a point in $\Omega^{n+m} (L\times K)(\R^{n+m}\oplus V)$ by the canonical homeomorphism
\[
\Omega^{n+m} L(\R^{n+m}\oplus V)\times\Omega^{n+m} K(\R^{n+m}\oplus V)\cong \Omega^{n+m} (L\times K)(\R^{n+m}\oplus V)\quad.
\]
This natural transformation
\[
\Omega^n L(\R^n\oplus V)\times \Omega^m K(\R^m\oplus V)\rightarrow \Omega^{n+m} (L\times K)(\R^{n+m}\oplus V)
\]
together with direct sum on $\mathscr{E}$ induce the natural map
\[
\alpha_{L,K}:QL\times QK\rightarrow Q(L\times K)\quad.
\]
This clearly satisfies associativity and unity with respect to $\eta_*:*\rightarrow Q*$.

For all $L$ and $K$, we have a commutative diagram
\[
\begin{CD}
L\times K @= L\times K\\
@V{\eta_L\times \eta_K}VV @VV{\eta_{L\times K}}V\\
QL\times QK @>{\alpha_{L,K}}>> Q(L\times K)
\end{CD}\quad.
\]
If $L$ and $K$ are well-pointed, then the vertical maps are $\pi_*$-isomorphisms. And it follows that $\alpha_{L,K}$
also is a $\pi_*$-iso. But $QL\times QK$ and $Q(L\times K)$ are both $\Omega$-spectra, hence $\alpha_{L,K}$
is a level equivalence.
\end{proof}

\chapter{Equivariance for orthogonal spectra}\label{Chapt:equiv}

We now give an exposition of the theory of
equivariant orthogonal spectra. Much of the material presented here
can also be found in~\cite{MandellMay:02}. 
However, there are some new results. The author would like to point out three
novelties: We introduce new types of cells, induced- and orbit- cells, and provide a
lax symmetric orbit cofibrant replacement functor $\tilde{\Gamma}$, see theorem~\ref{Thm:orbitcofrepl}.
The second new result is a formula for the geometric fixed points of induced $G$-spectra,
see proposition~\ref{Prop:geomfpandindsp}. We also introduce a diagonal map
for the iterated smash products $L^{\wedge q}$ of an orthogonal spectrum $L$, see definition~\ref{Def:diagonalmap}.
When $L$ is cofibrant, the diagonal map is an isomorphism into the geometric fixed points,
see proposition~\ref{Prop:diagonalisomorphism}.

In this chapter $G$ will be a compact Lie group, but for some arguments we 
restrict to the case where $G$ is a finite, discrete group. 
Genuine equivariance means to allow any $G$-representation 
when indexing our spectrum. Naive equivariance means to allow
trivial representations only.
However, lemma~\ref{Lem:ChangeUniverse} provides change of universe functors,
and they are equivalences of categories. Hence, 
\textbf{there is only one category of orthogonal $G$-spectra}.
This category has many model structures, and these model structures do depend
on the choice of a $G$-universe. So the modifiers ``genuine'' and ``naive''
apply to notions such as weak equivalences,
cellularity, model structures, geometric fixed points etc. 

\section{Preliminaries}

Let us now introduce the relevant terminology and notation for $G$-categories.
This material can also be found in chapter~II~\S{}1 of~\cite{MandellMay:02},
but is included here for the convenience of the reader.

A \textit{topological $G$-category} is a category $\mathscr{C}_G$
such that its hom sets $\mathscr{C}_G(C,D)$ are based topological $G$-spaces,
and composition 
\[
\circ:\mathscr{C}_G(D,E)\wedge \mathscr{C}_G(C,D)\rightarrow\mathscr{C}_G(C,E)
\]
is a continuous $G$-map. We think about the elements of $\mathscr{C}_G(C,D)$
as the non-equivariant maps $C\rightarrow D$, and we call them the \textit{arrows}
of $\mathscr{C}_G$. Observe that $\id_C\in\mathscr{C}_G(C,C)$ is fixed by the $G$-action;
this is easily deduced from the fact that the composition $\circ$ is $G$-equivariant.

Given a topological $G$-category $\mathscr{C}_G$, we can form a topological category
$G\mathscr{C}$. This category has the same objects, but the hom sets are given by
taking $G$-fixed points:
\[
G\mathscr{C}(C,D)=\mathscr{C}_G(C,D)^G\quad.
\]
Hence $G\mathscr{C}(C,D)$ is a based topological space, and composition is continuous.
We think about the elements of $G\mathscr{C}(C,D)$ as the $G$-equivariant maps $C\rightarrow D$, 
and we call them \textit{$G$-maps}.

If we ever encounter a situation where the $G$-spaces of arrows, $\mathscr{C}_G(D,E)$,
are given as unbased $G$-spaces, we implicitly add disjoint $G$-fixed basepoints.

\begin{Exa}
We let $\Top_{G*}$ denote the topological $G$-category of compactly 
generated $G$-spaces (=weak Hausdorff $k$-spaces with $G$-action)
and non-equivariant maps. $G$ acts on the space $\Top_{G*}(X,Y)$ of arrows
$X\rightarrow Y$ by conjugation.
Written out explicitly, the element $g\in G$ sends a map $f:X\rightarrow Y$ to
the composition $gfg^{-1}$.
The $G$-maps from $X$ to $Y$
are the $G$-fixed points
\[
G\Top_*(X,Y)=\Top_{G*}(X,Y)^G\quad.
\]
\end{Exa}

A \textit{continuous $G$-functor} $F:\mathscr{C}_G\rightarrow\mathscr{D}_G$ between 
topological $G$-categories is a functor $F$ such that
\[
F:\mathscr{C}_G(C,D)\rightarrow\mathscr{D}_G(F(C),F(D))
\]
is a continuous $G$-equivariant map of based $G$-spaces.
It follows that $F$ induces a functor $G\mathscr{C}\rightarrow G\mathscr{D}$.

A \textit{natural $G$-transformation} $\alpha:F_1\rightarrow F_2$ between continuous
$G$-functors $\mathscr{C}_G\rightarrow\mathscr{D}_G$ consists of $G$-maps
$\alpha:F_1(C)\rightarrow F_2(C)$ for every object $C$ in $\mathscr{C}_G$
such that the diagrams
\[
\begin{CD}
F_1(C) @>{F_1(f)}>> F_1(D)\\
@V{\alpha}VV @VV{\alpha}V\\
F_2(C) @>{F_2(f)}>> F_2(D)
\end{CD}
\]
commute in $\mathscr{D}_G$ for all $f\in\mathscr{C}_G(C,D)$.

Now we begin defining orthogonal $G$-spectra. The first thing we have to do is to somehow
specify the $G$-representations we will allow for indexing. 
We prefer the notion of a good collection, but one could also talk about $G$-universes.
After that, we define the topological $G$-categories $\mathscr{I}_G$.
And $\mathscr{I}_G$-spaces will relate to orthogonal $G$-spectra, just as $\mathcal{I}$-spaces
relate to orthogonal spectra.

\begin{Def}
Let $\mathscr{V}$ be a collection of finite dimensional real
$G$-inner product spaces. We call $\mathscr{V}$
a \textit{good collection} if it contains the trivial representations and 
is closed under direct sum. $\mathscr{V}$ is called
a \textit{very good collection} if it in addition to being a good collection,
is closed under passage to subrepresentations.
\end{Def}

A \textit{$G$-universe} $U$ is a sum of countably many copies of each real $G$-inner product
spaces in some set of irreducible representations of $G$ that includes the trivial representation.
$U$ is \textit{complete} if it contains all irreducible representations. $U$ is \textit{trivial}
if it contains only trivial representations.
Observe that there is a correspondence between very good collections and universes. Given $U$ one can define
a collection $\mathscr{V}(U)$ consisting of all $G$-representations
isomorphic to a finite dimensional sub $G$-inner product space of $U$. Given a very good collection $\mathscr{V}$
one can pick one representative for every isomorphism class of irreducible representations contained in $\mathscr{V}$.
This set of irreducible representations generate a $G$-universe.

\begin{Def}
Let $\mathscr{V}$ be a good collection of $G$-representations. Define $\mathscr{I}^{\mathscr{V}}_G$ to be the
(unbased) topological $G$-category whose objects are those of $\mathscr{V}$ and whose arrows (=non-equivariant maps)
are the linear isometric isomorphisms. $G$ acts on the space $\mathscr{I}^{\mathscr{V}}_G(V,W)$ of 
arrows $V\rightarrow W$ by conjugation. Let $G\mathscr{I}^{\mathscr{V}}$ be the topological category
with the same objects, but with $G$-maps $V\rightarrow W$ as morphisms. We have
\[
G\mathscr{I}^{\mathscr{V}}(V,W)=\mathscr{I}^{\mathscr{V}}_G(V,W)^G\quad.
\]
\end{Def}

\begin{Rem}\label{Rem:omitcollection}
We will usually omit the collection $\mathscr{V}$ from the notation. Thus we
write $\mathscr{I}_G$ instead of $\mathscr{I}_G^{\mathscr{V}}$.
Mostly we will be interested in the case where $\mathscr{V}=\All$ is the collection
of all $G$-representations. The other extreme case is $\mathscr{V}=\text{triv}$, the collection 
containing only the trivial $G$-representations. We shall see below that up
to equivalence of categories the choice of $\mathscr{V}$ does not matter. However,
it plays an important role for model structures on $G\mathscr{IS}$.
\end{Rem}

\begin{Def}
An \textit{$\mathscr{I}_G$-space} is a continuous $G$-functor $X:\mathscr{I}_G\rightarrow \Top_{G*}$.
Let $\mathscr{I}_G\Top_*$ be the topological $G$-category of $\mathscr{I}_G$-spaces and arrows the non-equivariant
natural transformations $X\rightarrow Y$. $G$ acts on $\mathscr{I}_G\Top_*(X,Y)$ by conjugation.
We define $G\mathscr{I}\Top_*$ to be the topological category  with the same objects and natural $G$-maps
$X\rightarrow Y$ as morphisms. We have
\[
G\mathscr{I}\Top_*(X,Y)=\mathscr{I}_G\Top_*(X,Y)^G\quad.
\]
\end{Def}

It is not obvious how to define a topology on $\mathscr{I}_G\Top_*(X,Y)$. Here is how to do it.
First choose a skeleton $\sk \mathscr{I}_G$ of $\mathscr{I}_G$. (A \textit{skeleton} $\sk\mathscr{C}$ for a category 
$\mathscr{C}$ is defined as a full subcategory such that 
each object in $\mathscr{C}$ is isomorphic to a unique object in $\sk\mathscr{C}$.)
Observe that $\sk \mathscr{I}_G$ is small.
Given $\mathscr{I}_G$-spaces $X$ and $Y$, we can consider the product of the function
spaces $F(X(V),Y(V))$ over $V\in\sk\mathscr{I}_G$. Now $\mathscr{I}_G\Top_*(X,Y)$
lies as a subset inside $\prod_{V\in\sk\mathscr{I}_G}F(X(V),Y(V))$, and we give it 
the subspace topology.

Here is an example of an $\mathscr{I}_G$-space:

\begin{Exa}
Let $S$ be the functor $\mathscr{I}_G\rightarrow \Top_{G*}$ sending
the $G$-representation $V$ to $S^V$, the one-point-compactification of $V$.
\end{Exa}

We now define orthogonal $G$-spectra:

\begin{Def}
An \textit{orthogonal $G$-spectrum} $L$ is an $\mathscr{I}_G$-space $L$ together
with a natural $G$-map $\sigma:L(V)\wedge S^W\rightarrow L(V\oplus W)$
such that unit and associativity diagrams commute. 
These diagrams are identical to those found in definition~\ref{Def:ortsp}.
Let $\mathscr{I}_G\mathscr{S}$ denote the topological $G$-category of 
orthogonal $G$-spectra and arrows $f:L\rightarrow K$ that commute with $\sigma$.
In general $f$ is non-equivariant, and $G$ acts on
$\mathscr{I}_G\mathscr{S}(L,K)$ by conjugation. Let $G\mathscr{IS}$ be the
topological category of orthogonal $G$-spectra and $G$-maps. We have
\[
G\mathscr{IS}(L,K)=\mathscr{I}_G\mathscr{S}(L,K)^G\quad.
\]
\end{Def}

Observe that $S$ is an orthogonal $G$-spectrum. We call $S$ the \textit{sphere spectrum}.

Similar to the description of orthogonal spectra given by theorem~\ref{thm:orthspasdiagr},
we now define topological $G$-categories $\mathscr{J}_G$ such that 
$\mathscr{J}_G$-spaces are the same as orthogonal $G$-spectra. This is also done
in chapter~II~\S{}4 of~\cite{MandellMay:02}.

\begin{Def}
Let $\mathscr{V}$ be some good collection of $G$-representations.
The objects of $\mathscr{J}_G^{\mathscr{V}}$
are the same as the objects of $\mathscr{I}_G^{\mathscr{V}}$, the finite dimensional $G$-representations $V$
contained in $\mathscr{V}$. 
Let $\mathscr{E}(V,W)$ be the $G$-space of (non-equivariant) linear isometries $V\hookrightarrow W$. 
$G$ acts on $\mathscr{E}(V,W)$ by conjugation.
And define $E(V,W)$ to be the $G$-space of 
pairs $(f,w)$ where $f:V\rightarrow W$ is a linear isometry and $w\in W$ is orthogonal to 
the linear subspace $f(V)$.
$E(V,W)$ is a vector bundle over $\mathscr{E}(V,W)$, and we define the $G$-space of morphisms 
$\mathscr{J}_G^{\mathscr{V}}(V,W)$
to be the Thom space of $E(V,W)$. (First apply fiber-wise one-point compactification to $E(V,W)$,
then identify the points at $\infty$.)
The $G$-action on $\mathscr{J}_G^{\mathscr{V}}(V,W)$ is explicitly given as follows; an element $g\in G$
sends the pair $(f,w)$ to $(gfg^{-1},g(w))$.
Composition
\[
\circ:\mathscr{J}_G^{\mathscr{V}}(W,U)\wedge\mathscr{J}_G^{\mathscr{V}}(V,W)\rightarrow\mathscr{J}_G^{\mathscr{V}}(V,U)
\]
is defined by the formula $(h,u)\circ(f,w)=(h\circ f,h(w)+u)$. The identity of $V$ in $\mathscr{J}_G^{\mathscr{V}}$ 
is represented by $(\id_V,0)$. 
Direct sum gives $\mathscr{J}_G^{\mathscr{V}}$ a symmetric monoidal structure:
\[
\oplus:\mathscr{J}_G^{\mathscr{V}}(V,W)\wedge\mathscr{J}_G^{\mathscr{V}}(V',W')
\rightarrow \mathscr{J}_G^{\mathscr{V}}(V\oplus V,W\oplus W')
\]
is defined by $(f,w)\oplus(f',w')=(f\oplus f',(w,w'))$. 
\end{Def}

Observe that when $V\subseteq W$ we have the identification:
\[
\mathscr{J}_G^{\mathscr{V}}(V,W)\cong O(W)_+\wedge_{O(W-V)} S^{W-V}\quad.
\]
As usual we let $G\mathscr{J}^{\mathscr{V}}$ denote the $G$-fixed category of
$\mathscr{J}_G^{\mathscr{V}}$. The two categories has the same objects and 
\[
G\mathscr{J}^{\mathscr{V}}(V,W)=\mathscr{J}_G^{\mathscr{V}}(V,W)^G\quad.
\]
We follow the convention from remark~\ref{Rem:omitcollection}
and omit the collection $\mathscr{V}$ from the notation, thus
writing $\mathscr{J}_G$ and $G\mathscr{J}$.

\textit{$\mathscr{J}_G$-spaces} are defined as continuous $G$-functors $\mathscr{J}_G\rightarrow{\Top_G}_*$.
We have an external smash product sending a pair of $\mathscr{J}_G$-spaces to a $\mathscr{J}_G\times \mathscr{J}_G$-space.
The (internal) smash product of $\mathscr{J}_G$-spaces is given as the left Kan extension along $\oplus$
of the external smash product.
Similar to theorem~\ref{thm:orthspasdiagr} we have the following result:

\begin{Thm}
The symmetric monoidal category of orthogonal $G$-spectra is isomorphic to the symmetric monoidal category of
$\mathscr{J}_G$-spaces.
\end{Thm}

This is theorem~II.4.3 in~\cite{MandellMay:02}. For a proof mimic~\S{}23 in~\cite{MandellMaySchwedeShipley:01}.

\subsection{Shift desuspension functors}

The equivariant \textit{shift desuspension functors} $F_V:\Top_{G*}\rightarrow\mathscr{I}_G\mathscr{S}$
are defined for all $G$-representations $V$. For based $G$-spaces $A$, the orthogonal $G$-spectrum $F_VA$ is
given at level $W$ by
\[
(F_VA)(W)=\mathscr{J}_G(V,W)\wedge A\quad.
\]
Observe that $F_V$ is the left adjoint to evaluation at level $V$. We have:
\[
\mathscr{I}_G\mathscr{S}(F_VA,L)\cong\Top_{G*}(A,L(V))
\]
for all $V$, $A$ and $L$.

\section{Change of universe functors}

Reading the definition it seems that the categories $\mathscr{I}_G\mathscr{S}$ and $G\mathscr{IS}$ 
depend on the choice of a good collection $\mathscr{V}$ of $G$-representations. As we soon
shall see, up to canonical 
equivalence of categories this choice does not matter. Therefore, it is our point of view
that the notion of an orthogonal $G$-spectrum is well-defined, without any modifier determining a
choice of $\mathscr{V}$. ``Naive'' and ``genuine'' are examples of such modifiers.
However, these modifiers will later play an important role, for example 
when considering extra structure on $\mathscr{I}_G\mathscr{S}$
and when constructing associated functors.

The key lemma is: 

\begin{Lem}\label{Lem:ChangeUniverse}
Let $V\in\mathscr{V}$ be a $G$-representation and $n=\dim V$. For an orthogonal $G$-spectrum
$L$ the evaluation $G$-map
\[
\mathscr{I}_G(\R^n,V)_+\wedge L(\R^n)\rightarrow L(V) 
\]
induces a $G$-homeomorphism
\[
\alpha:\mathscr{I}_G(\R^n,V)_+\wedge_{O(n)} L(\R^n)\rightarrow L(V) \quad.
\]
\end{Lem} 

\begin{proof}
This is lemma~V.1.1 in~\cite{MandellMay:02}. The proof is illustrative, so
we include it here:
The evaluation is a $G$-map since $L$ is a $G$-functor.
The $O(n)$-action on $\mathscr{I}_G(\R^n,V)_+$ and $L(\R^n)$ commutes with the $G$-action, hence
the $\alpha$ is a $G$-map.
Choose any linear isomorphism $f:\R^n\rightarrow V$. We get the inverse to $\alpha$ by sending 
$y\in L(V)$ to the point represented by $(f,L(f^{-1})(y))$.
\end{proof}

\begin{Def}
Let $\mathscr{V}$ and $\mathscr{V}'$ be good collections of $G$-representations.
We define the $G$-functor 
$I^{\text{triv}}_{\mathscr{V}}:\mathscr{I}^{\text{triv}}_G\mathscr{S}\rightarrow \mathscr{I}^{\mathscr{V}}_G\mathscr{S}$
by letting
\[
(I^{\text{triv}}_{\mathscr{V}}L)(V)=\mathscr{I}_G(\R^n,V)_+\wedge_{O(n)} L(\R^n)\quad,\text{ where $n=\dim V$.}
\]
In addition there are forgetful functors 
$I_{\text{triv}}^{\mathscr{V}'}:\mathscr{I}^{\mathscr{V}'}_G\mathscr{S}\rightarrow \mathscr{I}^{\text{triv}}_G\mathscr{S}$.
And we define $I_{\mathscr{V}}^{\mathscr{V}'}$ as the composite
\[
\mathscr{I}^{\mathscr{V}'}_G\mathscr{S}\xrightarrow{I_{\text{triv}}^{\mathscr{V}'}} \mathscr{I}^{\text{triv}}_G\mathscr{S}
\xrightarrow{I^{\text{triv}}_{\mathscr{V}}} \mathscr{I}^{\mathscr{V}}_G\mathscr{S}\quad.
\]
\end{Def}

\begin{Thm}
$I_{\mathscr{V}}^{\mathscr{V}'}$ is an isomorphism between the categories
$\mathscr{I}^{\mathscr{V}'}_G\mathscr{S}$ and $\mathscr{I}^{\mathscr{V}}_G\mathscr{S}$.
\end{Thm}

\begin{proof}
It is enough to check that $I^{\text{triv}}_{\mathscr{V}}$ is an isomorphism of categories.
Its inverse is the forgetful functor $I_{\text{triv}}^{\mathscr{V}}$. And by
the definition of the former it is easily seen that $I^{\text{triv}}_{\mathscr{V}}\circ I_{\text{triv}}^{\mathscr{V}}$
is naturally isomorphic to the identity functor on $\mathscr{I}^{\text{triv}}_G\mathscr{S}$.
Lemma~\ref{Lem:ChangeUniverse} above provides a natural isomorphism
\[
I_{\text{triv}}^{\mathscr{V}}\circ I^{\text{triv}}_{\mathscr{V}}\cong \id\quad.
\]
This finishes the proof.
\end{proof}

\section{Notions of equivalence}

Unlike the category of orthogonal spectra, $\mathscr{IS}$, where we considered just two
notions of equivalence, namely level equivalences and $\pi_*$-isomorphisms, we will in the equivariant
case define many different classes of weak equivalences. The reason for this phenomenon is three
choices influencing our definition. These are
\begin{itemize}
\item[] the choice of level-wise versus stable,
\item[] the choice of a family $\mathscr{F}$ of subgroups of $G$, and
\item[] the choice of a good collection $\mathscr{V}$ of $G$-representations.
\end{itemize}
Of course, not all possible sets of choices are equally interesting. We shall
point out a few interesting examples. Often one has a particular application
of the theory in mind when considering a notion of equivalence. For example,
we will define cyclotomic $\pi_*$-isomorphisms below, in order to study $THH$ and $TC$
of orthogonal ring spectra (with involution) in chapter~\ref{Chapt:trace}.

A \textit{family of subgroups} of $G$ is defined as a collection $\mathscr{F}$
of $H\subseteq G$ closed under passage to conjugates and subgroups.

Let us begin by considering the level-wise cases.

\begin{Def}
Let $f:L\rightarrow K$ be a map of orthogonal $G$-spectra.
We say that $f$ is a \textit{level-wise $(\mathscr{F},\mathscr{V})$-equivalence} if
for every $H\in\mathscr{F}$ and $V\in\mathscr{V}$ the map
\[
L(V)^H\rightarrow K(V)^H
\]
is a weak equivalence of topological spaces.
\end{Def}

We now define homotopy groups. 
In non-equivariant homotopy theory we have one homotopy group for every $q\in\Z$,
whereas in equivariant homotopy theory we index our homotopy groups
by an integer $q$ and a subgroup $H$ of $G$.
For orthogonal $G$-spectra the homotopy groups also depend on the choice of a 
good collection $\mathscr{V}$. 

\begin{Def}
Let $U$ be a $G$-universe associated to $\mathscr{V}$.
The homotopy groups are defined by
\[
\pi^H_q L=\begin{cases}
\colim_{V\subset U}\pi_q(\Omega^V L(V))^H&\text{ if $q\geq0$, and}\\
\colim_{\R^q\subseteq V\subset U}\pi_0(\Omega^{V-\R^q} L(V))^H&\text{ if $q\leq 0$.}\\
\end{cases}
\]
\end{Def}

\begin{Rem}
Observe that for fixed $H$ the homotopy group $\pi_q^H L$ does not really depend on
which $G$-representations $\mathscr{V}$ contains, but rather on the $H$-representations
appearing as the restriction of some $V$ in $\mathscr{V}$. To see this, assume that 
$\phi:V\rightarrow W$ is an $H$-linear isometric isometry. The diagram
\[
\begin{CD}
S^V @>{f}>> L(V)\\
@V{\phi_*}V{\cong}V @V{\cong}V{L(\phi)}V\\
S^W @>{g}>> L(W)
\end{CD}
\]
gives a homeomorphism between the (non-equivariant) spaces $\left(\Omega^VL(V)\right)^H$
and $\left(\Omega^VL(V)\right)^H$, defined by sending $f$ to the unique $g$ making
the diagram commute.
\end{Rem}

Stable notions of equivalence are now defined as follows:

\begin{Def}\label{Def:stableFVpiiso}
Let $f:L\rightarrow K$ be a map of orthogonal $G$-spectra.
We say that $f$ is an \textit{$(\mathscr{F},\mathscr{V})$-$\pi_*$-isomorphism} if
for every $H\in\mathscr{F}$ and $q\in \Z$ the map
\[
\pi_q^H L\rightarrow \pi_q^H K
\]
is an isomorphism.
\end{Def}

We now have overwhelmingly many notions of equivalence for orthogonal $G$-spectra.
Let us now give names to the extreme cases and some other interesting examples.

\begin{Def}
The \textit{naive level-equivalences} are the level-wise $(\All,\text{triv})$-equivalences,
the \textit{genuine level-equivalences} are the level-wise $(\All,\All)$-equivalences,
the \textit{naive $\pi_*$-isomorphisms} are the $(\All,\text{triv})$-$\pi_*$-isomorphism, and
the \textit{genuine $\pi_*$-isomorphisms} are the $(\All,\All)$-$\pi_*$-isomorphism.
\end{Def}

Here the family $\All$ is the collection of all subgroups of $G$, 
$\text{triv}$ is the collection of the trivial $G$-representations, and the good collection
$\All$ is the collection of all $G$-representations.

\begin{Rem}\label{Rem:comparewequiv}
Clearly every genuine level-equivalence is a naive level-equivalence. Furthermore 
genuine $\pi_*$-isomorphisms are naive $\pi_*$-isomorphisms by theorem~V.1.7 in~\cite{MandellMay:02}.
Because of lemma~\ref{Lem:ChangeUniverse} it is tempting to think that the converse statements also
must be true. However this is \textbf{not} the case. The reason in the level-case, is that the 
$H$-fixed points of $\mathscr{I}_G(\R^n,V)_+\wedge_{O(n)} L(\R^n)$ is generally not equal to
$\mathscr{I}_G(\R^n,V)_+\wedge_{O(n)} L(\R^n)^H$.
\end{Rem}

If the Lie group $G$ is $S^1$ or $O(2)$ we define:

\begin{Def}\label{Def:cyclpiiso}
Let $\mathscr{F}$ be the family of finite cyclic subgroups of $S^1$.
The \textit{cyclotomic $\pi_*$-isomorphisms} are the $(\mathscr{F},\All)$-$\pi_*$-isomorphisms.
\end{Def}

By the inclusion $S^1\subset O(2)$, the family $\mathscr{F}$ is also a family of subgroups of $O(2)$.

$\Omega$-spectra in the equivariant setting are given as follows:

\begin{Def}\label{Def:OmegaGsp}
An orthogonal $G$-spectrum $L$ is an $\Omega$-$G$-spectrum if all maps
\[
L(V)^H\rightarrow \left(\Omega^WL(V\oplus W)\right)^H
\]
are weak equivalences of spaces. Here $V$ and $W$ are $G$-representations, and $H$ a closed subgroup of $G$. 
\end{Def}

\section{Cellular orthogonal $G$-spectra}

In this section we recall the notion of an $F^{\mathscr{V}}I_G$-cell from~\cite{MandellMay:02},
we also introduce new types of cells, namely induced $G$-cells and orbit $G$-cells.

\begin{Def}
Given a good collection $\mathscr{V}$ of $G$-representations. Choose a skeleton $\sk\mathscr{I}^{\mathscr{V}}_G$ for
$\mathscr{I}^{\mathscr{V}}_G$. The \textit{$F^{\mathscr{V}}I_G$-cells} is the set of all maps
\[
F_V\left(S^{n-1}\times G/H\right)_+\rightarrow F_V\left(D^{n-1}\times G/H\right)_+\quad,
\]
where $n\geq0$, $V\in\sk\mathscr{I}_G$ and $H$ is a closed subgroup of $G$.
\end{Def}

The two extreme cases are:

\begin{Def}
The \textit{naive $FI_G$-cells} are the $F^{\text{triv}}I_G$-cells; we allow only trivial $G$-representations $V$.
The \textit{genuine $FI_G$-cells} are the $F^{\All}I_G$-cells; we allow all $G$-representations.
\end{Def}

\begin{Rem}\label{Rem:equivdepofrep}
To what degree does the equivariant structure of the orthogonal $G$-spectrum $F_V(D^n\times G/H)_+$
depend on the representation $V$? 

Assume that $H$ is trivial. Choose a (non-equivariant) isometric isomorphism $\phi:\R^m\rightarrow V$.
Evaluating $F_V(D^n\times G)_+$ at some level $W$ we get
\[
\mathscr{J}_G(V,W)\wedge(D^n\times G)_+\quad.
\]
Represent a point by a triple $(f,x,g)$. We have a map into $\mathscr{J}_G(\R^m,W)\wedge(D^n\times G)_+$
given by sending $(f,x,g)$ to $(fg\phi,x,g)$. And it is easily checked that this map is a $G$-map.
Hence we have a $G$-isomorphism of orthogonal $G$-spectra
\[
F_V(D^n\times G)_+\cong F_{\R^m}(D^n\times G)_+\quad.
\]

Now assume that $H$ is non-trivial.
Whenever we have an $H$-linear isometric isomorphism $\phi:U\rightarrow V$, then
$F_V(D^n\times G)_+$ and $F_{U}(D^n\times G)_+$ are $G$-isomorphic. The isomorphism
is defined by sending $(f,x,[g])$ to $(fg\phi g^{-1},x,[g])$. The $H$-linearity of $\phi$
ensures that this map does not depend on choice of $g$ representing the class $[g]\in G/H$.

Up to $G$-isomorphism, the orthogonal $G$-spectrum $F_V(D^n\times G/H)_+$ depends only on 
$V$ as an $H$-representation.
\end{Rem}

In our applications we will need even more general cells. Therefore we define:

\begin{Def}\label{Def:IndGFIcells}
For any closed subgroup $H$ of $G$, let $\All_H$ be the collection of all finite dimensional 
$H$-representations. Choose skeletons $\sk\mathscr{I}^{\All_H}_H$ for all $\mathscr{I}^{\All_H}_H$.
Define the set of \textit{induced $G$-cells} to be the set $\Ind_GFI$ of all maps
\[
\left(F_V S^{n-1}_+\right)\wedge_H G_+\rightarrow \left(F_V D^{n}_+\right)\wedge_H G_+\quad,
\]
where $n\geq0$, $H$ is a closed subgroup of $G$, and $V$ is a finite dimensional $H$-representation 
in $\sk\mathscr{I}^{\mathscr{V}_H}_H$.
\end{Def}

The functor $-\wedge_H G_+$ assigns to an $H$-space its induced $G$-space. It also takes orthogonal $H$-spectra
to orthogonal $G$-spectra. For precise definitions see section~\ref{sect:indspect} below.

\begin{Rem}
If $V$ is the restriction of a $G$-representation, then 
$\left(F_V D^{n}_+\right)\wedge_H G_+\cong F_V(D^n_+\wedge G/H_+)$.
This shows that all genuine $FI_G$-cells are induced $G$-cells.
\end{Rem}

In order to get a symmetric cofibrant replacement functor we will consider an even nastier kind of cells. 
In this case we restrict ourselves to finite $G$. The cells are defined as follows:

\begin{Def}\label{Def:OrbGFIcells}
Let $H$ be any finite group, and let $\mathscr{V}_H$ be the collection of all finite dimensional 
$H$-representations. Choose skeletons $\sk\mathscr{I}^{\mathscr{V}_H}_H$ for all $\mathscr{I}^{\mathscr{V}_H}_H$.
Define the set of \textit{orbit $G$-cells} to be the set $\Orb_GFI$ of all maps
\[
\left( F_V S^{n-1}_+\wedge G_+\right)/H\rightarrow \left( F_V D^{n}_+\wedge G_+\right)/H\quad,
\]
where $n\geq0$, $V$ is a finite dimensional $H$-representation 
in $\sk\mathscr{I}^{\mathscr{V}_H}_H$, $H$ acts trivially on $S^{n-1}$ and $D^n$, 
and the action of $H$ on $G$ is given via a group homomorphism $H\rightarrow G$.
\end{Def}

\begin{Rem}
If $H$ is a subgroup of $G$ acting on $G$ via the inclusion, then
we have the following identification:
\[
\left( F_V D^{n}_+\wedge G_+\right)/H \cong \left(F_V D^{n}_+\right)\wedge_H G_+\quad. 
\]
Thus we see that any induced $G$-cell is an orbit $G$-cell.
\end{Rem}

In order to bring things more down to earth, we will now write out explicitly
what an orbit $G$-cell looks like at some level $W$.

\begin{Exa}
Let $\phi:H\rightarrow G$ denote the group homomorphism.
Let $W$ be a $G$-representation. We are now going to evaluate $F_V D^{n}_+\wedge G_+$
at level $W$ and specify the $G$ and $H$ actions.
By definition of the shift desuspension we have
\[
\left(F_V D^{n}_+\wedge G_+\right)(W)=\mathscr{J}(V,W)\wedge D^n_+\wedge G_+\quad.
\]
Recall that a non-basepoint in $\mathscr{J}(V,W)$ consists of an isometric embedding $f:V\rightarrow W$
and a $w\in W$ orthogonal to $f$. Thus a tuple $(f,w,x,g)$ represents a point in the space above.
An $h\in H$ gives an isometry $h:V\rightarrow V$, and $h$ sends $f$ to the composition $fh^{-1}$.
Via $\phi$ the element $h$ acts on $g$ by sending $g$ to $g\phi(h)$. Dividing out by $H$ we see that
\[
(fh^{-1},w,x,g)\quad\text{and}\quad (f,w,x,g\phi(h))
\]
are identified in $\left(F_V D^{n}_+\wedge G_+\right)/H(W)$.
Furthermore, on this space we have an action of $G$. Let $\gamma$ be an element in $G$.
It acts from the left by sending $(f,w,x,g)$ to $(\gamma f,\gamma(w),x,\gamma g)$.
\end{Exa}

This example also shows the reason for the name ``orbit $G$-cell''. Contrary to all
other types of cells, the action of $O(W)$ on the $W$'th level of the orbit cell is not necessarily free.
And non-free actions have more than one type of orbits.

\begin{Rem}\label{Rem:orbitcellredundancy}
Notice the following redundancy in the definition of the orbit $G$-cell $\left( F_V D^{n}_+\wedge G_+\right)/H$:
If there is a kernel $K$ of the map $H\rightarrow G\times O(V)$, then
we have the identification
\[
\left( F_V D^{n}_+\wedge G_+\right)/H\cong \left( F_V D^{n}_+\wedge G_+\right)/J\quad,
\]
where $J$ is the quotient $H/K$.
\end{Rem}

\begin{Exa}
Consider the orbit $G$-cells when $G$ is the trivial group. They have the form
\[
\left(F_V S^{n-1}_+\right)/H\rightarrow \left(F_V D^{n}_+\right)/H\quad.
\]
If $V$ is a non-trivial $H$-representation, then this orbit cell is not isomorphic
to any $FI$-cell.
\end{Exa}

Now we are ready to define the various types of relative $G$-equivariant cellular maps. 

\begin{Def}\label{def:FIGcell}
Let $K$ be either the set of $F^{\mathscr{V}}I_G$-cells for some good collection,
the set of induced $G$-cells, or the set of orbit $G$-cells.
A map $i:A\rightarrow L$ of orthogonal $G$-spectra is \textit{relative $K$-cellular} if:
\begin{itemize}
\item[] $i(A)$ is a sub-$G$-spectrum of $L$.
\item[] There is a set $C$ of sub-$G$-spectra $L_{\alpha}$ such that each 
$L_{\alpha}$ contains $i(A)$ and $\bigcup_{\alpha\in C}L_{\alpha}=L$.
\item[] $C$ is partially ordered by inclusion. We write $\beta\leq\alpha$ 
if $L_{\beta}\subseteq L_{\alpha}$. And for all $\alpha$ the
set $P_{\alpha}=\{ \beta\in C\;|\; \beta<\alpha\}$ is finite.
\item[] For every $\alpha\in C$ there is pushout diagram with $G$-equivariant maps:
\[
\begin{CD}
\partial E @>{i}>> E\\
@VVV @VVV\\
\bigcup_{\beta<\alpha} L_{\beta} @>>> L_{\alpha}
\end{CD}\quad,
\]
where $\partial E\xrightarrow{i} E$ is a cell in $K$.
\end{itemize}
\end{Def}

Corresponding to the different kinds of relative cellular maps we have q-cofibrations.
The \textit{naive q-cofibrations} are the retracts of relative $FI^{\text{triv}}_G$-cellular maps,
the \textit{genuine q-cofibrations} are the retracts of relative $FI^{\All}_G$-cellular maps, 
the \textit{induced q-cofibrations} are the retracts of relative $\Ind_GFI$-cellular maps, and
the \textit{orbit q-cofibrations} are the retracts of relative $\Orb_GFI$-cellular maps.

\begin{Rem}
Observe that naive q-cofibrations are genuine q-cofibrations, that genuine q-cofibrations
are induced q-cofibrations, and that induced q-cofibrations are orbit q-cofibrations.
All these kinds of q-cofibrations are both l- and h-cofibrations.
\end{Rem}

\section{Model structures on $G\mathscr{IS}$}\label{sect:ms}

A justified question is what combinations of cofibrations and equivalences give
model structures on $G\mathscr{IS}$. In this subsection we shall briefly 
comment this by recalling some results from~\cite{MandellMay:02}.
But first we recall the notion of a compactly generated model category,
see definition~5.9 in~\cite{MandellMaySchwedeShipley:01}:

Roughly speaking 
a model category $\mathscr{A}$ is \textit{compactly generated} if there exist sets of
maps $I$ and $J$ in $\mathscr{A}$, which can be used in the small object argument, and 
such that the fibrations in $\mathscr{A}$ are the maps which satisfy the right lifting
property with respect to all maps in $J$ and the acyclic fibrations are the maps which
satisfy the right lifting property with respect to all maps in $I$.
The maps in $I$ are called the \textit{generating cofibrations} and the maps in
$J$ are called the \textit{generating acyclic cofibrations}.

Here are our model categories:
\begin{itemize}
\item[] The \textit{naive level-wise model structure} on $G\mathscr{IS}$ has naive
level-equivalences as weak equivalences and naive q-cofibrations as cofibrations. See
theorem~III.2.4 in~\cite{MandellMay:02}.
\item[] The \textit{genuine level-wise model structure} on $G\mathscr{IS}$ has genuine
level-equivalences as weak equivalences and genuine q-cofibrations as cofibrations. Again see
theorem~III.2.4 in~\cite{MandellMay:02}.
\item[] The \textit{naive stable model structure} on $G\mathscr{IS}$ has naive
$\pi_*$-isomorphisms as weak equivalences and naive q-cofibrations as cofibrations. See
theorem~III.4.2 in~\cite{MandellMay:02}.
\item[] The \textit{genuine stable model structure} on $G\mathscr{IS}$ has genuine
$\pi_*$-isomorphisms as weak equivalences and genuine q-cofibrations as cofibrations. 
The fibrant objects are the $\Omega$-$G$-spectra. Again see
theorem~III.4.2 in~\cite{MandellMay:02}.
\item[] Let $\mathscr{F}$ be a family of subgroups of $G$. The \textit{stable $\mathscr{F}$-model structure}
has $(\mathscr{F},\All)$-$\pi_*$-isomorphisms as weak equivalences. The generating cofibrations are 
those genuine $G$-cells $F_V (D^n\times G/H)_+$ where $H\in\mathscr{F}$. See 
theorem~IV.6.5 in~\cite{MandellMay:02}.
\item[] In particular there is a \textit{stable cyclotomic model structure} on
orthogonal $S^1$- and $O(2)$-spectra. The weak equivalences are the cyclotomic
$\pi_*$-isomorphisms, and the cofibrations are constructed allowing only the cells
$F_V (D^n\times G/H)_+$ where $H\subset S^1$ ($\subset O(2)$) is finite cyclic.
\end{itemize}
All these model structures are compactly generated.

\begin{Rem}
Induced $G$-cells and orbit $G$-cells are introduced in this thesis. And it
is beyond the scope of this work to find model structures on $G\mathscr{IS}$
where the cofibrations are generated by these classes of cells.
However, if such model structures exist, they probably have better properties.
For the case with orbit $G$-cells, theorem~\ref{Thm:orbitcofrepl} gives a hint about this.
\end{Rem}

\section{Categorical fixed points}

Let $H$ be a closed subgroup of $G$. In this section we will describe how
to take the $H$-fixed points of an orthogonal $G$-spectrum. The basic
properties of this construction is given in proposition~\ref{Prop:catfp}.
We will also define the notion of a free orthogonal $G$-spectrum, and state some
elementary observations.

\begin{Def}
Assume that $L$ is an orthogonal $G$-spectrum and $N$ a closed normal subgroup of $G$.
Let $J$ be the quotient $G/N$.
The \textit{categorical $N$-fixed point spectrum} $L^N$ is given by
\[
L^N(\R^n)=L(\R^n)^N\quad,
\]
for trivial representations $\R^n$. Clearly $L^N$ is an orthogonal spectrum.
Notice that $L^N$ has a $J$-action. 
By the appropriate change of universe functor we define $L^N(V)$
for any $J$-representation $V$. Thus $L^N$ is an orthogonal $J$-spectrum.
\end{Def}

More generally we could define $L^H$ for any closed subgroup $H$ of $G$
by first restricting $L$ to an orthogonal $N_H$-spectrum. Here $N_H$ denotes the
normalizer of $H$ in $G$. Taking the $H$-fixed point spectrum
$L^H$ we get an orthogonal $W_H=N_H/H$-spectrum.

Again let $N$ be a closed normal subgroup of $G$, and let
$\epsilon:G\rightarrow J=G/N$ be the quotient map. We shall now
specify a right adjoint $\epsilon^*$ to the functor $(-)^N$. 
Given an orthogonal $J$-spectrum $K$, we define the orthogonal $G$-spectrum 
$\epsilon^*K$ by giving $K(\R^n)$ the $N$-trivial $G$-action. We have
\[
(\epsilon^*K)(\R^n)=\epsilon^*(K(\R^n))\quad.
\]
We extend $(\epsilon^*K)(V)$ to all $G$-representations $V$
by the appropriate change of universe functor.

The main properties of the categorical fixed points can now be summarized
in the following proposition:

\begin{Prop}\label{Prop:catfp}
Let $L$ be an orthogonal $G$-spectrum and $K$ an orthogonal $J$-spectrum. There
is a natural isomorphism
\[
G\mathscr{IS}(\epsilon^*K,L)\cong J\mathscr{IS}(K,L^N)\quad.
\]
Furthermore, if $V$ is a $G$-representation and $A$ a based $G$-space, then we have that
\[
\left(F_V A\right)^N=\begin{cases}
F_V (A^N) & \text{if $V$ is an $N$-trivial $G$-representation, and}\\
* & \text{otherwise.}
\end{cases}
\]
The functor $(-)^N$ preserves naive and genuine q-cofibrations, it also preserves
acyclic naive q-cofibrations, but not acyclic genuine q-cofibrations.
\end{Prop}

For a proof see the propositions~V.3.5 and~V.3.10 in~\cite{MandellMay:02}.

We now define what a free orthogonal $G$-spectrum is:

\begin{Def}
Assume that $L$ is an orthogonal $G$-spectrum. We say that $L$ is \textit{free}
if $L(V)^H=*$ for all non-trivial closed subgroups $H$ of $G$ and all $G$-representations $V$. 
\end{Def}

Observe that if $L$ is free, then the categorical $H$-fixed point spectrum $L^H$ is trivial for
all $H\neq\{1\}$. However, the converse of this statement is not true. It is possible for an
orthogonal $G$-spectrum to have $L(\R^m)^H=*$ for all $\R^m$ and $H\neq\{1\}$, but $L(V)^H\neq*$
for some non-trivial $G$-representation $V$ and $H\neq\{1\}$. For examples of such $L$, 
consider the orbit $G$-cells.

\begin{Rem}\label{Rem:cellsoffreeGspect}
We now make some elementary observations: Let $L$ be an orthogonal $G$-spectrum.
\begin{itemize}
\item[] Assume that $L$ is both free and naive $FI_G$-cellular. Then all cells occurring are of the form
$F_{\R^m}(D^n\times G)_+$. 
\item[] Assume that $L$ is both free and genuine $FI_G$-cellular. The all cells occurring are of the
form $F_V(D^n\times G)_+$, where $V$ is some $G$-representation, but
due to remark~\ref{Rem:equivdepofrep} it can always be assumed that $V$ is trivial. Hence, $L$
is actually naive $FI_G$-cellular.
\item[] Assume that $L$ is both free and $\Ind_GFI$-cellular. An induced cell $\left(F_V D^n_+\right)\wedge_H G_+$
have non-trivial $H$ fixed points, so unless $H$ is trivial $L$ cannot be free. 
Recall that $V$ was an $H$-representation. Thus $V=\R^m$ for some $m$ since $H=\{1\}$.
Hence, all cells of $L$ are of the form $F_{\R^m}(D^n\times G)_+$. This shows that $L$
actually is naive $FI_G$-cellular.
\item[] Assume that $L$ is both free and $\Orb_G FI$-cellular. Then group homomorphism
$H\rightarrow G$ of a cell in $L$ must be trivial, and all cells can thus be written as
\[
\left( F_V D^{n}_+\wedge G_+\right)/H\cong \left( F_V D^{n}_+\right)/H \wedge G_+
\]
for some $H$-representation $V$.
\end{itemize}
\end{Rem}

We end this section by a simple, but important observation.

\begin{Prop}\label{Prop:freenaiveisgenuine}
If $f:L\rightarrow K$ is a naive level-equivalence between free orthogonal $G$-spectra, then
$f$ is also a genuine level-equivalence.
\end{Prop}

\begin{proof}
Let $V$ be a $G$-representation and $H$ a subgroup of $G$. 
We must check that $f^H:L(V)^H\rightarrow K(V)^H$ is a weak equivalence of spaces.
There are two cases to consider. Assume first that $H$ is the trivial group.
Since $f$ is a naive level-equivalence, 
the map
\[
L(V)^H=L(V)\cong L(\R^m)\xrightarrow{f} K(\R^m)\cong K(V)=K(V)^H
\]
is a weak equivalence of spaces. Here $m=\dim V$, and 
the identifications $L(V)\cong L(\R^m)$ and $K(V)\cong K(\R^m)$
are non-equivariant.

The other case is when $H$ is non-trivial. Then by freeness, both $L(V)^H$ and $K(V)^H$ are equal to the
trivial orthogonal spectrum $*$. Thus it is a tautology that
\[
f^H:L(V)^H\rightarrow K(V)^H
\]
is a weak equivalence of spaces.
\end{proof}

\section{Geometrical fixed points}

We now define the geometric fixed point functor, and give the relevant results.
This functor certainly depends on the choice of a collection $\mathscr{V}$ of
$G$-representations. In fact, if $\mathscr{V}=\text{triv}$, then the geometric
fixed points are equal to the categorical fixed points. However, we are going to
use the convention that for the purpose of taking geometric fixed points
if not otherwise specified, then $\mathscr{V}$ is understood to
be the collection $\All$ of all $G$-representations.

We follow the exposition given in chapter~V~\S{}4 of~\cite{MandellMay:02} closely,
and begin with some categorical preliminaries.

Let $E$ denote the short exact sequence,
\[
0\rightarrow N\rightarrow G\xrightarrow{\epsilon} J\rightarrow 0\quad,
\]
of Lie groups. Here $N$ is a closed normal subgroup of $G$.
We now define a topological $J$-category $\mathscr{J}_E$ as follows:
The objects are the $G$-representations $V$ contained in
our collection $\mathscr{V}$. The morphisms from $V$ to $W$
are the $N$-fixed points of $\mathscr{J}_G(V,W)$. This means that
a non-basepoint arrow of $\mathscr{J}_E(V,W)$ can be represented by
a pair $(f,w)$, where $f$ is an $N$-linear isometry $V\rightarrow W$ and
$w$ is a point in $W^N$ orthogonal to $f(V^N)$. Observe that 
$\mathscr{J}_E=G\mathscr{J}$ when $N=G$, $\mathscr{J}_E=\mathscr{J}_G$
when $N$ is trivial, and $\mathscr{J}_E=\mathscr{J}_J^{\text{triv}}$
when $\mathscr{V}=\text{triv}$.

Let $\phi:\mathscr{J}_E\rightarrow\mathscr{J}_J$ denote the $J$-functor
which sends the $G$-representation $V$ to the $J$-representation $V^N$,
and the arrow $(f,w)\in\mathscr{J}_E(V,W)$ to $(f^{N},w)\in\mathscr{J}_J(V^N,W^N)$.
We think about $\phi$ as the $N$-fixed point functor.

Let $\nu:\mathscr{J}_J\rightarrow\mathscr{J}_E$ be the $J$-functor
which sends the $J$-representation $V$ to the $G$-representation
$\epsilon^*V$, and the arrow $(f,w)\in\mathscr{J}_J(V,W)$ maps to
$(f,w)\in\mathscr{J}_E(\epsilon^*V,\epsilon^*W)$. We think about $\nu$
as the pullback along $\epsilon:G\rightarrow J$.

Observe that $\phi\circ\nu=\id:\mathscr{J}_J\rightarrow \mathscr{J}_J$.

\begin{Def}\label{Def:forgetandprolong}
Let $\mathscr{J}_E\Top_*$ denote the category of $\mathscr{J}_E$-spaces,
namely continuous $J$-functors $\mathscr{J}_E\rightarrow{\Top_J}_*$. The functors
$\phi$ and $\nu$ induce forgetful functors denoted by
\[
\U_{\phi}:\mathscr{J}_J\Top_*\rightarrow \mathscr{J}_E\Top_*
\quad\text{and}\quad
\U_{\nu}:\mathscr{J}_E\Top_*\rightarrow \mathscr{J}_J\Top_*
\]
respectively. Left Kan extension along $\phi$ and $\nu$ gives prolongation
functors
\[
\bbP_{\phi}:\mathscr{J}_E\Top_*\rightarrow \mathscr{J}_J\Top_*
\quad\text{and}\quad
\bbP_{\nu}:\mathscr{J}_J\Top_*\rightarrow \mathscr{J}_E\Top_*
\]
left adjoint to $\U_{\phi}$ and $\U_{\nu}$ respectively.
We have $\U_{\nu}\circ \U_{\phi}=\id$ and $\bbP_{\phi}\circ \bbP_{\nu}\cong\id$.
\end{Def}

\begin{Def}\label{Def:geomfpfunct}
Let $\Fix^N:\mathscr{J}_G\Top_*\rightarrow\mathscr{J}_E\Top_*$ be the functor
sending an orthogonal $G$-spectrum $L$ to the $\mathscr{J}_E$-space $\Fix^N L$
given by
\[
(\Fix^N L)(V)=L(V)^N
\]
and with evaluation $J$-maps
\[
(\Fix^N L)(V)\wedge\mathscr{J}_E(V,W)=
L(V)^N\wedge \mathscr{J}_G(V,W)^N\xrightarrow{\ev^N} L(W)^N=(\Fix^N L)(W)\quad.
\]
Here the maps $\ev:L(V)\wedge \mathscr{J}_G(V,W)\rightarrow L(W)$ are the evaluation $G$-map
of $L$. Define the \textit{geometric fixed point functor} 
$\Phi^N:\mathscr{J}_G\Top_*\rightarrow\mathscr{J}_J\Top_*$
to be the composition $\bbP_{\phi}\circ\Fix^N$.
\end{Def}

\begin{Constr}\label{Constr:cattogeomfp}
There is a natural $J$-map from the categorical fixed points
to the geometrical fixed points. This map $L^N\rightarrow \Phi^N L$
is defined as follows: Observe that $L^N=\U_{\nu}\Fix^N L$.
The adjunction between $\U_{\phi}$ and $\bbP_{\phi}$ has a unit
$\eta:\id\rightarrow\U_{\phi}\circ \bbP_{\phi}$. And the natural $J$-map
above is constructed as the composition
\[
L^N=\U_{\nu}\Fix^N L\xrightarrow{\U_{\nu}\eta}\U_{\nu}\U_{\phi}\bbP_{\phi}\Fix^N L=\bbP_{\phi}\Fix^N L=\Phi^N L\quad.
\]
\end{Constr}

In order to compute with the geometric fixed points we have the following proposition:

\begin{Prop}\label{Prop:fundresultofgeomfp}
For a $G$-representation $V$ in $\mathscr{V}$ and a based $G$-space $A$, we have
\[
\Phi^N(F_V A)\cong F_{V^N} A^N\quad.
\]
Furthermore, if $K$ is the pushout of $B\leftarrow A\xrightarrow{i} L$ in the category of
orthogonal $G$-spectra and $i$ is a closed inclusion, then $\Phi^NK$ is the
pushout of $\Phi^NB\leftarrow \Phi^NA\rightarrow \Phi^NL$. The functor $\Phi^N$
preserves q-cofibrations and acyclic q-cofibrations.
\end{Prop}

For a proof see proposition~V.4.5 in~\cite{MandellMay:02}.

\section{Induced orthogonal $G$-spectra}\label{sect:indspect}

In this section we will define the notion of an induced $G$-spectrum,
see definition~\ref{Def:indGsp}. The main topic is to study the geometric fixed points 
of induced $G$-spectra. Our main result is proposition~\ref{Prop:geomfpandindsp}.
There is an annoying condition in this proposition, but so far no
counterexample has been found.

\subsection{Basic facts about induced $G$-spaces}

Let $H$ be a subgroup of $G$, and $N$ a normal subgroup of $G$ which is contained in $H$.
All subgroups are closed, and $i:H\rightarrow G$ denotes the inclusion.
We will first define orthogonal $G$-spectra induced from 
orthogonal $H$-spectra. To do this we need some basic facts about equivariant spaces.

The setup for the groups occurring in this subsection can be expressed by
two short exact sequences. We have
\[
\begin{CD}
0 @>>> N @>>> H @>>> J_0 @>>> 0\\
&& @V{=}VV @V{i}VV @VV{i_1}V\\
0 @>>> N @>>> G @>>> J @>>> 0
\end{CD}\quad.
\]
We denote the first sequence by $E_0$ and the second sequence by $E$.
Let $j:E_0\rightarrow E$ be the map of sequences given by the diagram above.

Assume that $X$ is a based $G$-space and $Y$ a based $H$-space. By forgetting
part of the $G$-action on $X$ we get the $H$-space $i^*X$. Observe that $i^*$ is a functor
from based $G$-space to based $H$-spaces. This functor has both a left and a right adjoint.
The right adjoint of $i^*$ is the \textit{coinduced $G$-space}, and is given by
the formula $F_H(G_+,Y)$, where $F_H$ denotes the space of based $H$-maps. 
However, we will not meet coinduced $G$-spaces
in this thesis. Therefore we do not write out this side of the theory.
The left adjoint is defined by sending $Y$ to the based $G$-space
\[
Y\wedge_H G_+\quad.
\]
This is the quotient space of $Y\wedge G_+$ where $(hy,g)$ is identified with $(y,gh)$
for all $h\in H$. An element $\gamma\in G$ sends $(y,g)$ to $(y,\gamma g)$.
We call $Y\wedge_H G_+$ the \textit{induced $G$-space}.

A basic lemma for induced spaces is:

\begin{Lem}\label{Lem:basicGhomeoI}
Giving smash products the diagonal action, there is a natural 
$G$-homeomorphism 
\[
X\wedge (Y\wedge_H G_+)\cong (i^*X\wedge Y)\wedge_H G_+\quad.
\]
In particular, $X\wedge (G/H)_+\cong (i^*X)\wedge_H G_+$.
\end{Lem}

\begin{proof}
A point $(x;y,g)$ in $X\wedge (Y\wedge_H G_+)$ is sent to $(g^{-1}x,y;g)$ in $(i^*X\wedge Y)\wedge_H G_+$.
\end{proof}

We will now see how the $N$-fixed point functor commutes with the induced space functor.
Recall that $J=G/N$ and $J_0=H/N$.
Note that $X^N$ is a $J$-space and $Y^N$ is a $J_0$-space. 
We have:

\begin{Lem}\label{Lem:fixedpointsofinducedspaces}
There is a natural $J$-homeomorphism
\[
(Y\wedge_H G_+)^N\cong Y^N\wedge_{J_0}J_+\quad.
\]
\end{Lem}

\begin{proof}
We construct map both ways and end the proof by observing that they are inverses to each other.
Assume that $(y,g)$ is a point in $(Y\wedge_H G_+)^N$. For all $n\in N$ we have
\[
(y,g)=(y,ng)=(y,gn')=(n'y,g)\quad,
\]
where $n'=g^{-1}ng\in N$. Thus $y=n'y$.
Since conjugating with $g$ is an automorphism of $N$, we see that $y$ is an $N$-fixed point.
We therefore send $(y,g)\in (Y\wedge_H G_+)^N$ to $(y,[g])\in Y^N\wedge_{J_0}J_+$.

Given a point $(y,[g])$ in $Y^N\wedge_{J_0}J_+$, we choose an element $g\in G$ representing
the class $[g]\in J=G/N$. We claim that the map sending $(y,[g])$ to $(y,g)\in (Y\wedge_H G_+)^N$
is well-defined. Suppose that $g'\in G$ was another choice of element representing $[g]$, then
\[
(y,g')=(y,gg^{-1}g')=(g^{-1}g'y,g)=(y,g)\quad.
\]
The last equality follows since $g^{-1}g'\in N$ and $y\in Y^N$.
Obviously, the two maps are inverses to each other.
\end{proof}

\subsection{``Change'' functors for equivariant orthogonal spectra and $\mathscr{J}_E$-spaces}
 
Now consider the case of equivariant orthogonal spectra. First we define the change of group functor,
then we look at its adjoints and at last we determine how the left adjoint, the induced spectra,
interact with the geometric fixed point functors. 
The inclusion $i:H\rightarrow G$ induces a change of group functor
from orthogonal $G$-spectra to orthogonal $H$-spectra. Assume that $K$ is an orthogonal $G$-spectrum, and
$L$ an orthogonal $H$-spectrum.

\begin{Def}
The \textit{change of group functor} is defined by letting $i^*K$ be the orthogonal $H$-spectrum
given by $(i^*K)(i^*V)=i^*K(V)$ for $G$-representations $V$. We extend $i^*K$ to $H$-representations
$W$ not of the form $i^*V$, by the change of universe functor.
\end{Def}

The right adjoint of $i^*$, the \textit{coinduced orthogonal $G$-spectrum} functor,
is defined at level $V$ as $F_H(G_+,L)(V)=F_H(G_+,L(i^*V))$, where $F_H$ is the space of based $H$-maps.
For further details, see proposition~V.2.4 in~\cite{MandellMay:02}.

We now look at induced orthogonal $G$-spectra. The definition can be found in proposition~V.2.3
in~\cite{MandellMay:02}, but we include it here for the convenience of the reader.

\begin{Def}\label{Def:indGsp}
The \textit{induced orthogonal $G$-spectrum}, $L\wedge_H G_+$, is given by
$(L\wedge_H G_+)(V)=L(i^*V)\wedge_H G_+$ for $G$-representations $V$.
\end{Def}

How to define the evaluation $H$-maps will be discussed below in the more general setting
of induced $\mathscr{J}_E$ spaces. See remark~\ref{Rem:DefindorthGspeval}.

Let us now build the theory for ``change'' functors for $\mathscr{J}_E$-spaces.
Unlike the case of equivariant orthogonal spectra, the change of universe functors
will not always be equivalences of categories. This subtle difference forces us
to be extremely careful regarding universes. We begin with an example:

\begin{Exa}
A key ingredient in the proof of lemma~\ref{Lem:ChangeUniverse}, was that for any
good collection $\mathscr{V}$ of $G$-representations and any $V\in\mathscr{V}$,
there is a trivial $G$-representation $\R^n$ such that $\R^n$ and $V$
were isomorphic in the category $\mathscr{J}_G^{\mathscr{V}}$. The analogous
statement for $\mathscr{J}_E$ is in general not true:

Let $E$ be the sequence $0\rightarrow C_2\rightarrow S^1\rightarrow S^1/C_2\rightarrow 0$,
and let $\mathscr{V}$ be the collection of all $S^1$-representations.
Identify $S^1$ with the unit circle in $\C$ and consider the representation $\C$, where $S^1$
acts by multiplication. Any $C_2$-linear isometry from a trivial $S^1$-representation $\R^n$
into $\C$ must map $C_2$-fixed points onto $C_2$-fixed points. Since $\C^{C_2}=0$, there
is no isomorphism in $\mathscr{J}_E^{\mathscr{V}}$ between $\C$ and some trivial representation.

Therefore, the forgetful functor from $\mathscr{J}_E^{\mathscr{V}}$-spaces to
$\mathscr{J}_E^{\text{triv}}$-spaces cannot be an isomorphism of categories.
\end{Exa}

\begin{Def}
Let $E$ be a short exact sequence of compact Lie groups.
Assume that $\mathscr{V}\subseteq\mathscr{V}'$ are two good collections of $G$-representations. 
Then we have a forgetful functor
$\U:\mathscr{J}_E^{\mathscr{V}'}\Top_*\rightarrow\mathscr{J}_E^{\mathscr{V}}\Top_*$. 
By left Kan extension, we define a prolongation functor
\[
\bbP:\mathscr{J}_E^{\mathscr{V}}\Top_*\rightarrow\mathscr{J}_E^{\mathscr{V}'}\Top_*
\]
left adjoint to $\U$. These are the \textit{change of universe functors} for
$\mathscr{J}_{E}$-spaces.
\end{Def}

Let us now look at the change of sequence. As above we consider a diagram
\[
\begin{CD}
0 @>>> N @>>> H @>>> J_0 @>>> 0\\
&& @V{=}VV @V{i}VV @VV{i_1}V\\
0 @>>> N @>>> G @>>> J @>>> 0
\end{CD}\quad.
\]
This diagram is a map $j:E_0\rightarrow E$ of short exact sequences.
Let $\mathscr{V}$ be a good collection of $G$-representations.
Now let $i^*\mathscr{V}$ be defined as the collection of $H$-representations $i^*V$,
which are the restriction of some $V\in\mathscr{V}$. And we define:

\begin{Def}
The \textit{change of sequence functor}
\[
j^*:\mathscr{J}_E^{\mathscr{V}}\Top_*\rightarrow\mathscr{J}_{E_0}^{i^*\mathscr{V}}\Top_*
\]
is given by sending the $\mathscr{J}_E^{\mathscr{V}}$-space $X$ to
the $\mathscr{J}_{E_0}^{i^*\mathscr{V}}$-space given by
\[
(j^*X)(i^*V)=i_1^*X(V)
\]
for $H$-representations $i^*V$ being the restriction of a $G$-representation $V$ in $\mathscr{V}$.
\end{Def}

The change of group functor for equivariant orthogonal spectra had both left and right adjoints.
These were the induced and coinduced spectra. Analogously we now define induced and coinduced 
$\mathscr{J}_{E}$-spaces. Since the coinduced $\mathscr{J}_{E}$-spaces will not
be used later in the thesis, we only sketch the definition: Given a $\mathscr{J}_{E_0}^{i^*\mathscr{V}}$-space 
$Y$ we let the \textit{coinduced $\mathscr{J}_E^{\mathscr{V}}$-space} $F_{J_0}(J_+,Y)$
be given by
\[
F_{J_0}(J_+,Y)(V)=F_{J_0}(J_+,Y(i^*V))\quad.
\]
The induced $\mathscr{J}_{E}$-spaces are defined as follows:

\begin{Def}
Let $Y$ be a $\mathscr{J}_{E_0}^{i^*\mathscr{V}}$-space, then the 
\textit{induced $\mathscr{J}_E^{\mathscr{V}}$-space}
$Y\wedge_{J_0}J_+$ is given by
\[
(Y\wedge_{J_0}J_+)(V)=Y(i^*V)\wedge_{J_0}J_+
\]
for $G$-representations $V$ in $\mathscr{V}$. The evaluation $J$-maps are given by
\begin{align*}
\mathscr{J}_E^{\mathscr{V}}(V,W)\wedge (Y\wedge_{J_0}J_+)(V)
&= \mathscr{J}_G(V,W)^N\wedge(Y(i^*V)\wedge_{J_0}J_+)\\
&\cong \left( i_1^*\mathscr{J}_G(V,W)^N\wedge Y(i^*V)\right)\wedge_{J_0}J_+\\
&\cong \left( \mathscr{J}_H(i^*V,i^*W)^N\wedge Y(i^*V)\right)\wedge_{J_0}J_+\\
& \rightarrow Y(i^*W)\wedge_{J_0} J_+= (Y\wedge_{J_0} J_+)(W)
\end{align*}
for $G$-representations $V$ and $W$ in $\mathscr{V}$.
\end{Def}

Here we have used the $G$-homeomorphism $X\wedge (Y\wedge_H G_+)\cong (i^*X\wedge Y)\wedge_H G_+$
from lemma~\ref{Lem:basicGhomeoI}, and the evaluation
$\mathscr{J}_{E_0}(i^*V,i^*W)\wedge Y(i^*V)\rightarrow Y(i^*W)$ for the 
$\mathscr{J}_{E_0}^{i^*\mathscr{V}}$-space $Y$.

\begin{Rem}\label{Rem:DefindorthGspeval}
If $N$ is trivial, then $\mathscr{J}_E=\mathscr{J}_G$. Hence, the last part of the definition above
tells us how to define the evaluation $G$-maps for the induced orthogonal $G$-spectrum functor.
\end{Rem}

\subsection{Geometric fixed points and induced spectra}

We will now consider how the geometric $N$-fixed point functors interact with induced orthogonal $G$-spectra.
In order to prove the result we have to make an assumption relating representations of $N$ and $G$.
The author has not found any counterexample to the condition.

As before we consider the diagram
\[
\begin{CD}
0 @>>> N @>>> H @>>> J_0 @>>> 0\\
&& @V{=}VV @V{i}VV @VV{i_1}V\\
0 @>>> N @>>> G @>>> J @>>> 0
\end{CD}\quad.
\]
$E_0$ is the first short exact sequence and $E$ the second.
Let $j:E_0\rightarrow E$ be the map of sequences given by the diagram.

\begin{Prop}\label{Prop:geomfpandindsp}
Let $j:E_0\rightarrow E$ be as above.
Suppose that for any $N$-representation $W$ there exists a $G$-representation $V$
and an $N$-linear isometric embedding $W\rightarrow V$ such that $W^N=V^N$.
Then for orthogonal $H$-spectra $L$ there is a natural isomorphism
\[
(\Phi^NL)\wedge_{J_0}J_+\cong \Phi^N(L\wedge_H G_+)\quad.
\]
\end{Prop}

\begin{Rem}\label{Rem:abtNGconition}
It is enough to check the assumption for non-trivial irreducible $N$-representations $W$.
The condition certainly holds whenever the quotient group $J=G/N$ is finite.
Because in this case one can take $V$ to be the induced $G$-representation of $W$.

If one fixes the compact Lie group $G$, one can try to list normal subgroups and their representations
and then check explicitly if the condition holds. To check that the
assumption holds in the cases $G=S^1$ and $G=O(2)$, is an easy exercise left to the reader.
\end{Rem}

The proof of the proposition is to check commutativity of the following diagram:
\[
\begin{CD}
\mathscr{J}_H\Top_* 
@>{\Fix^N}>> \mathscr{J}_{E_0}^{\mathscr{W}}\Top_*
@>{\bbP_{\phi_0}}>> \mathscr{J}_{J_0}\Top_*\\
@V{=}VV @V{\U}VV @VV{=}V\\
\mathscr{J}_H\Top_* 
&& \mathscr{J}_{E_0}^{i^*\mathscr{V}}\Top_*
@>{\bbP_{\phi_0}^{i^*\mathscr{V}}}>> \mathscr{J}_{J_0}\Top_*\\
@V{-\wedge_H G_+}VV @V{-\wedge_{J_0} J_+}VV @VV{-\wedge_{J_0} J_+}V\\
\mathscr{J}_G\Top_* 
@>{\Fix^N}>> \mathscr{J}_{E}^{\mathscr{V}}\Top_*
@>{\bbP_{\phi}}>> \mathscr{J}_{J}\Top_*
\end{CD}\quad.
\]
Here $\mathscr{V}$ denotes the collection of all $G$-representations, while 
$\mathscr{W}$ are all $H$-representations. We have not specified universes for
the categories of orthogonal $H$-, $G$-, $J_0$- and $J$-spectra.
This is since the change of universe is an isomorphism for these cases,
and we can change universe whenever needed.
The functors $\Fix^N$ are defined with respect to genuine equivariance.
Similarly, the prolongation functors $\bbP_{\phi_0}$ and $\bbP_{\phi}$,
given in definition~\ref{Def:forgetandprolong}, use all representations.
However, the prolongation functor $\bbP_{\phi_0}^{i^*\mathscr{V}}$
is defined using only those representations being restrictions from $J$ and $G$.

We prove three lemmas, each checking commutativity of one of the squares.
Let us start with the hardest:

\begin{Lem}
Suppose that for any $N$-representation $W$ there exists a $G$-representation $V$
and an $N$-linear isometric embedding $W\rightarrow V$ such that $W^N=V^N$.
Let $Y$ be a $\mathscr{J}_{E_0}^{\mathscr{W}}$-space. Then we have a natural isomorphism
\[
\bbP_{\phi_0}^{i^*\mathscr{V}}\U Y \cong \bbP_{\phi_0}Y\quad.
\]
\end{Lem}

\begin{proof}
Since $\bbP_{\phi_0}^{i^*\mathscr{V}}$ and $\bbP_{\phi_0}$ are left adjoints, while
$\U$ is a right adjoint, the proof of this lemma cannot be abstract category theory.
By change of universe for orthogonal $J_0$-spectra it is enough to check
that we have a natural isomorphism when evaluating at the trivial $J_0$-representations $\R^n$.
We begin by writing out both sides evaluated at $\R^n$ explicitly. We have
\[
\left(\bbP_{\phi_0}^{i^*\mathscr{V}}\U Y\right)(\R^n)
=\int^{i^*V\in i^*\mathscr{V}}\mathscr{J}_{J_0}((i^*V)^N,\R^n)\wedge Y(i^*V)
\]
and
\[
\left(\bbP_{\phi_0} Y\right)(\R^n)
=\int^{W\in \mathscr{W}}\mathscr{J}_{J_0}(W^N,\R^n)\wedge Y(W)\quad.
\]
Recall that $i^*\mathscr{V}$ are the $H$-representations of the form $i^*V$ for some $G$-representation $V$,
while $\mathscr{W}$ are all $H$-representations. 
Since $i^*\mathscr{V}\subset\mathscr{W}$, we clearly have a natural map
\[
\bbP_{\phi_0}^{i^*\mathscr{V}}\U Y\rightarrow \bbP_{\phi_0} Y\quad.
\]
Now we show that this map is surjective. Pick a point in $(\bbP_{\phi_0} Y)(\R^n)$.
It is represented by a pair $((f,u),y)$ in $\mathscr{J}_{J_0}(W^N,\R^n)\wedge Y(W)$ for
some $H$-representation $W$. By the assumption, we can choose a $G$-representation $V$
together with an $N$-linear isometric embedding $g:W\rightarrow V$ such that $W^N\cong V^N$.
What we have is an arrow $(g,0)$ in $\mathscr{J}_{E_0}(W,i^*V)$. This arrow is a relation
in $(\bbP_{\phi_0} Y)(\R^n)$ between the point we picked and the point
$((f(g^N)^{-1},u),(g,0)(y))$ in $\mathscr{J}_{J_0}((i^*V)^N,\R^n)\wedge Y(i^*V)$.
This new point is in the image of $(\bbP_{\phi_0}^{i^*\mathscr{V}}\U Y)(\R^n)$.

To show injectivity we consider a generating relation in $(\bbP_{\phi_0} Y)(\R^n)$
between $((f_1,u_1),y_1)$ in $\mathscr{J}_{J_0}(W_1^N,\R^n)\wedge Y(W_1)$ and 
$((f_2,u_2),y_2)$ in $\mathscr{J}_{J_0}(W_2^N,\R^n)\wedge Y(W_2)$.
By the argument above we can choose liftings of these points to $(\bbP_{\phi_0}^{i^*\mathscr{V}}\U Y)(\R^n)$.
And the proof is completed by showing that the lifted points are related in
$(\bbP_{\phi_0}^{i^*\mathscr{V}}\U Y)(\R^n)$.

A generating relation in $(\bbP_{\phi_0} Y)(\R^n)$ is specified by
a triple
\[
\left( (f,u), (h,w), y\right) \in \mathscr{J}_{J_0}(W_2^N,\R^n)\wedge 
\mathscr{J}_{E_0}(W_1,W_2)\wedge Y(W_1)\quad.
\]
Here $f:W_2^N\rightarrow \R^n$ is an isometric embedding, $u$ a point in $\R^n$ orthogonal
to $f(W_2^N)$, $h:W_1\rightarrow W_2$ an $N$-linear isometric embedding, $w$ a point in $W_2^N$
orthogonal to $h(W_1^N)$ and $y$ a point in $Y(W_1)$.
This generating relation identifies the point
\[
((f_1,u_1),y_1)=\left( (fh^N,u+f(w)),y\right)
\quad\text{in}\quad\mathscr{J}_{J_0}(W_1^N,\R^n)\wedge Y(W_1)
\]
and
\[
((f_2,u_2),y_2)=\left( (f,u),(h,w)(y)\right)
\quad\text{in}\quad\mathscr{J}_{J_0}(W_2^N,\R^n)\wedge Y(W_2)\quad.
\]
Liftings of these points to $(\bbP_{\phi_0}^{i^*\mathscr{V}}\U Y)(\R^n)$ are 
given by $G$-representations $V_1$ and $V_2$, and $N$-linear
isometric embeddings $g_1:W_1\rightarrow V_1$ and $g_2:W_2\rightarrow V_2$
such that $W_1^N=V_1^N$ and $W_2^N=V_2^N$. The lifting of the first point is given by
\[
\left( (f_1(g_1^N)^{-1},u_1),(g_1,0)(y_1)\right)
\quad\text{in}\quad\mathscr{J}_{J_0}((i^*V_1)^N,\R^n)\wedge Y(i^*V_1)
\]
and similarly for the second point. Unfortunately, we cannot automatically complete the diagram
\[
\begin{CD}
W_1 @>{h}>> W_2\\
@V{g_1}VV @VV{g_2}V\\
V_1 && V_2
\end{CD}
\]
by an arrow at the bottom. However, we can form the $H$-representation
$i^*V_1\oplus_{W_1}i^*V_2$, and by the assumption there exists a $G$-representation
$V$ together with an $N$-linear isometric embedding $i^*V_1\oplus_{W_1}i^*V_2\rightarrow V$,
such that their $N$-fixed points agree. Now we have $N$-linear isometric embeddings
$h_1:V_1\rightarrow V$ and $h_2:V_2\rightarrow V$. Putting these maps into the diagram
above we get a commutative pentagon. Also observe that $V^N=V_2^N=W_2^N$.
Thus we have the arrow $(h_2,0)$ in $\mathscr{J}_{E_0}(i^*V_2,i^*V)$
and the arrow $(h_1,w)$ in $\mathscr{J}_{E_0}(i^*V_1,i^*V)$.

By $(h_2,0)$ we have a relation in $(\bbP_{\phi_0}^{i^*\mathscr{V}}\U Y)(\R^n)$ between
the second lifting,
\[
\left( (f_2(g_2^N)^{-1},u_2),(g_2,0)(y_2)\right)
\quad\text{in}\quad\mathscr{J}_{J_0}((i^*V_2)^N,\R^n)\wedge Y(i^*V_2)\quad,
\]
and
\[
\left( (f_2(g_2^N)^{-1}(h_2^N)^{-1},u_2),(h_2,0)(g_2,0)(y_2)\right)
\quad\text{in}\quad\mathscr{J}_{J_0}((i^*V)^N,\R^n)\wedge Y(i^*V)\quad.
\]
Now look at the generating relation in $(\bbP_{\phi_0}^{i^*\mathscr{V}}\U Y)(\R^n)$
specified by
\begin{multline*}
\left( (f_2(g_2^N)^{-1}(h_2^N)^{-1},u_2),(h_1,w),(g_1,0)(y_1)\right)\\
\quad\text{in}\quad\mathscr{J}_{J_0}((i^*V)^N,\R^n)\wedge \mathscr{J}_{E_0}(i^*V_1,i^*V)\wedge Y(i^*V_1)\quad.
\end{multline*}
This relation identifies
\begin{multline*}
\left( (f_2(g_2^N)^{-1}(h_2^N)^{-1}h_1^N,u_2+f_2(g_2^N)^{-1}(h_2^N)^{-1}(w)),(g_1,0)(y_1)\right)\\
= \left( (f_1(g_1^N)^{-1},u_1),(g_1,0)(y_1)\right)
\end{multline*}
in $\mathscr{J}_{J_0}((i^*V_1)^N,\R^n)\wedge Y(i^*V_1)$ with 
\begin{multline*}
\left( (f_2(g_2^N)^{-1}(h_2^N)^{-1},u_2),(h_1,w)(g_1,0)(y_1)\right)\\
= \left( (f_2(g_2^N)^{-1}(h_2^N)^{-1},u_2),(h_2,0)(g_2,0)(y_2)\right)
\end{multline*}
in $\mathscr{J}_{J_0}((i^*V)^N,\R^n)\wedge Y(i^*V)$.
Thus we have a relation in $(\bbP_{\phi_0}^{i^*\mathscr{V}}\U Y)(\R^n)$
between the two lifted points. This completes the proof of the lemma.
\end{proof}

The two other lemmas are easy:

\begin{Lem}
Let $L$ be an orthogonal $H$-spectrum. We have a natural isomorphism
\[
(\U\Fix^NL)\wedge_{J_0}J_+\cong \Fix^N(L\wedge_H G_+)
\]
of $\mathscr{J}_{E}^{\mathscr{V}}$-spaces.
\end{Lem}

\begin{proof}
Let $V$ be a $G$-representation. We evaluate both sides of the natural isomorphism
at $V$. The left side at level $V$ becomes:
\begin{align*}
\left((\U\Fix^NL)\wedge_{J_0}J_+\right)(V) &= (\Fix^NL)(i^*V)\wedge_{J_0} J_+\\
&= L(i^*V)^N\wedge_{J_0} J_+\quad.
\end{align*}
And the right side at level $V$ becomes:
\begin{align*}
\Fix^N(L\wedge_H G_+)(V)&=\left((L\wedge_H G_+)(V)\right)^N\\
&= \left(L(i^*V)\wedge_H G_+\right)^N\quad.
\end{align*}
From lemma~\ref{Lem:fixedpointsofinducedspaces} we recall
the natural homeomorphism $(Y\wedge_H G_+)^N\cong Y^N\wedge_{J_0}J_+$
for $H$-spaces $Y$. Setting $Y=L(i^*V)$ we get 
\[
\left((\U\Fix^NL)\wedge_{J_0}J_+\right)(V) \cong\Fix^N(L\wedge_H G_+)(V)\quad.
\]
To check the fact that the evaluation $J$-maps 
$\mathscr{J}_E(V,W)\wedge\left((\U\Fix^NL)\wedge_{J_0}J_+\right)(V)\rightarrow\left((\U\Fix^NL)\wedge_{J_0}J_+\right)(W)$
and $\mathscr{J}_E(V,W)\wedge\Fix^N(L\wedge_H G_+)(V)\rightarrow\Fix^N(L\wedge_H G_+)(W)$
agree is left as an exercise to the reader.
\end{proof}

\begin{Lem}
The diagram
\[
\begin{CD}
\mathscr{J}_{E_0}^{i^*\mathscr{V}}\Top_*
@>{\bbP_{\phi_0}^{i^*\mathscr{V}}}>> \mathscr{J}_{J_0}\Top_*\\
@V{-\wedge_{J_0} J_+}VV @VV{-\wedge_{J_0} J_+}V\\
\mathscr{J}_{E}^{\mathscr{V}}\Top_*
@>{\bbP_{\phi}}>> \mathscr{J}_{J}\Top_*
\end{CD}
\]
commutes.
\end{Lem}

\begin{proof}
It is enough to check commutativity of the corresponding diagram of
right adjoints:
\[
\begin{CD}
\mathscr{J}_{E_0}^{i^*\mathscr{V}}\Top_*
@<{\U_{\phi_0}^{i^*\mathscr{V}}}<< \mathscr{J}_{J_0}\Top_*\\
@A{j^*}AA @AA{i_1^*}A\\
\mathscr{J}_{E}^{\mathscr{V}}\Top_*
@<{\U_{\phi}}<< \mathscr{J}_{J}\Top_*
\end{CD}\quad.
\]
Let $L$ be an orthogonal $J$-spectrum and $V$ a $G$-representation. 
Now compare $\U_{\phi_0}^{i^*\mathscr{V}}i_1^*L$ and $j^*\U_{\phi}L$ at the level $i^*V$.
We have:
\begin{align*}
(\U_{\phi_0}^{i^*\mathscr{V}}i_1^*L)(i^*V) &= (i_1^*L)\left( (i^*V)^N\right)\\
&= (i_1^*L)\left( i_1^*(V^N)\right)\\
&= i_1^*\left(L(V^N)\right)\\
&= i_1^*\left((\U_{\phi}L)(V)\right)\\
&= (j^*\U_{\phi}L)(i^*V)\quad.
\end{align*}
Again we leave to the reader to check that the evaluation $J_0$-maps agree.
\end{proof}

Together the three lemmas above prove proposition~\ref{Prop:geomfpandindsp}.

\section{A symmetric cofibrant replacement functor}

In this section we mainly work with non-equivariant orthogonal spectra.
Recall the definition of the cofibrant replacement functor $\Gamma$, see theorem~\ref{Thm:cofrepl}.
Due to lemma~\ref{Lem:twistcells} there are two obstructions to
symmetry of $\Gamma$. The obstructions are that the twists of the disks, 
$D^{n_1}\times D^{n_2}\cong D^{n_2}\times D^{n_1}$,
and the twists of the indexing spaces, $\R^{m_1}\oplus\R^{m_2}\cong \R^{m_2}\oplus\R^{m_1}$,
are not the identity maps. If we are content with getting orbit cells instead of 
naive cells, then we can divide out by these twists when performing the small object argument.
We will explain this in detail below, thus constructing a symmetric functor $\tilde{\Gamma}$.
The section ends with a subsection containing various results concerning this functor and equivariance.

Our theorem says:

\begin{Thm}\label{Thm:orbitcofrepl}
There is an endofunctor $\tilde{\Gamma}$ on orthogonal spectra having the following properties:
\begin{itemize}
\item[] $\tilde{\Gamma}L$ is orbit cofibrant for all $L$.
\item[] If $K\rightarrow L$ is the inclusion of a subspectrum, then $\tilde{\Gamma}K\rightarrow\tilde{\Gamma}L$
is an orbit q-cofibration.
\item[] $\tilde{\Gamma}$ comes with a natural level-wise acyclic fibration $\tilde{\Gamma}L\rightarrow L$.
\item[] There is a natural quotient map $\Gamma L\rightarrow \tilde{\Gamma}L$.
\item[] There is a symmetric natural transformation 
$\phi:\tilde{\Gamma}L\wedge\tilde{\Gamma}K\rightarrow\tilde{\Gamma}(L\wedge K)$
and a canonical map $S\rightarrow\tilde{\Gamma}S$.
\end{itemize}
\end{Thm}

To prove this, we begin with the construction of $\tilde{\Gamma}$:

\begin{Constr}
We modify the small object argument.
The basic step is to introduce a new gluing construction.

Suppose that $p:A\rightarrow L$ is a map of orthogonal spectra.
Let $C_{n,m}$ be the set of all
diagrams
\[
\begin{CD}
F_{\R^{m}}S^{n-1}_+ @>{F_{\R^{m}}i_n}>> F_{\R^{m}}D^{n}_+\\
@V{f}VV @VV{g}V\\
A @>{p}>> L
\end{CD}\quad.
\]
Recall that the orthogonal spectrum $G(p)$ was defined as the pushout of
\[
A \xleftarrow{f}\bigvee F_{\R^{m}}S^{n-1}_+ \xrightarrow{F_{\R^{m}}i_n} \bigvee F_{\R^{m}}D^{n}_+\quad,
\]
where the wedge runs over all diagrams in $C_{n,m}$ with $n,m\geq0$. 

Recall from remark~\ref{Rem:cellsymmetries} that $\Sigma_{n}$ and $\Sigma_{m}$
both act on a cell $F_{\R^{m}}S^{n-1}_+ \rightarrow F_{\R^{m}}D^{n}_+$ in $C_{n,m}$.
A permutation $\sigma\in \Sigma_n$ gives a map $\sigma:D^n\rightarrow D^n$, and
$\sigma$ acts on diagrams $\alpha$ of $C_{n,m}$ by composing
$f_{\alpha}$ and $g_{\alpha}$ with $F(\sigma)$. 
A permutation $\rho\in \Sigma_m$ gives an isometry $\rho:\R^m\rightarrow\R^m$,
and $\rho$ acts on diagrams $\alpha$ of $C_{n,m}$ by composing 
$f_{\alpha}$ and $g_{\alpha}$ with $F_{\rho}$.
Altogether we get an action of $\Sigma_n\times\Sigma_m$ on $G_{n,m}(p)$.
Divide out by this action and define $\tilde{G}(p)$ as the union over all quotients:
\[
\tilde{G}(p)=\bigcup_{n,m} G_{n,m}(p)/ \Sigma_n\times\Sigma_m\quad.
\]
\end{Constr}

Now we proceed as before and define $\tilde{\Gamma}L$ by iterating
the gluing construction $\tilde{G}(p)$. Start with $*\xrightarrow{p_0}L$,
define $\tilde{G}^1(L)=\tilde{G}(p_0)$, and let $p_1$ be the canonical map $\tilde{G}^1(L)\rightarrow L$.
Inductively we get $G^{i-1}(L)\subseteq\tilde{G}^i=\tilde{G}(p_{i-1})\xrightarrow{p_i}L$.

\begin{Def}
The \textit{orbit cofibrant replacement functor} $\tilde{\Gamma}$ is defined for $L\in\mathscr{IS}$
as the colimit of the $\tilde{G}^i(L)$'s. 
\end{Def}

We now begin to prove the statements in theorem~\ref{Thm:orbitcofrepl}.
As a first result we justify the name ``\textbf{orbit cofibrant} replacement functor''
by showing:

\begin{Prop}
$\tilde{\Gamma}L$ is orbit cofibrant for any orthogonal spectrum $L$.
\end{Prop}

\begin{proof}
It is enough to consider the natural inclusion $j$ of gluing construction $A\xrightarrow{j} \tilde{G}(p)$,
and show that this map is a relative orbit q-cofibration.

Let $\alpha$ be a cell in $C_{n,m}$. There is a subgroup $H_{\alpha}$ of $\Sigma_n\times\Sigma_m$
of symmetries fixing the diagram
\[
\begin{CD}
F_{\R^{m}}S^{n-1}_+ @>{F_{\R^{m}}i_n}>> F_{\R^{m}}D^{n}_+\\
@V{f}VV @VV{g}V\\
A @>{p}>> L
\end{CD}
\]
corresponding to $\alpha$. Observe that $\tilde{G}(p)$ also can be described as the pushout of
\[
A \xleftarrow{f}\bigvee_{\alpha} \left(F_{\R^{m}}S^{n-1}_+\right)/H_{\alpha} \rightarrow 
\bigvee_{\alpha}\left(F_{\R^{m}}D^{n}_+\right)/H_{\alpha}\quad,
\]
where $\alpha$ runs through one representative for every $\Sigma_n\times\Sigma_m$ orbit of $C_{n,m}$,
$n,m\geq 0$. By Illman's theorem~\cite{Illman:83} we may triangulate 
$F_{\R^{m}}D^{n}_+$ evaluated at $\R^m$ as a finite $O(m)\times H_{\alpha}$-CW-complex.
All $O(m)$-orbits are free. Divide out by the $H_{\alpha}$-action.  
This describes $\left(F_{\R^{m}}S^{n-1}_+\right)/H_{\alpha} \rightarrow 
\bigvee_{\alpha}\left(F_{\R^{m}}D^{n}_+\right)/H_{\alpha}$ as a relative orbit cellular map.
Hence $A\xrightarrow{j} \tilde{G}(p)$ is also relative orbit cellular.
\end{proof}

Since $\tilde{G}(p)$ is a quotient of $G(p)$ it follows that $\tilde{\Gamma}L$ is a quotient of $\Gamma L$.
The canonical map $S\rightarrow\tilde{\Gamma}S$ is constructed just as before.
To prove the next two statements of the theorem, namely that $\tilde{\Gamma}$ takes inclusions to
orbit q-cofibrations, and that the natural map $\tilde{\Gamma}L\rightarrow L$ is a level-wise acyclic
fibration, we just copy the proofs of the propositions~\ref{Prop:Gammaincltoqcof} 
and~\ref{Prop:Gammalevelacyclfibr} respectively. 

What now remains is to define the natural map 
$\phi:\tilde{\Gamma}L\wedge\tilde{\Gamma}K\rightarrow\tilde{\Gamma}(L\wedge K)$, and
show that it is symmetric. 

\begin{Constr}
The construction is similar to construction~\ref{Constr:GammaMult}. We inspect the $\tilde{G}^i(L)$'s and the
$\tilde{G}^j(K)$'s and define inductively maps
\[
\phi_{i,j}:\tilde{G}^i(L)\wedge \tilde{G}^j(K)\rightarrow \tilde{G}^{i+j-1}(L\wedge K)
\]
such that diagrams 
similar to those in construction~\ref{Constr:GammaMult} commute.
By taking the colimit as both $i$ and $j$ tend to infinity, we get our natural transformation
$\tilde{\Gamma} L\wedge \tilde{\Gamma} K\rightarrow \tilde{\Gamma}(L\wedge K)$.

We proceed by induction on $i+j$.
Assume that $\alpha$ and $\beta$ are the cells given by the diagrams
\[
\begin{CD}
F_{\R^{m}}S^{n-1}_+ @>>> F_{\R^{m}}D^{n}_+\\
@V{f}VV @VV{g}V\\
G^{i-1}(L) @>>> L
\end{CD}\quad\text{and}\quad
\begin{CD}
F_{\R^{m'}}S^{n'-1}_+ @>>> F_{\R^{m'}}D^{n'}_+\\
@V{f'}VV @VV{g'}V\\
G^{j-1}(K) @>>> K
\end{CD}\quad.
\]
By the construction of $\tilde{G}^i(L)$ and $\tilde{G}^j(K)$ there are unique lifts of 
$\alpha$ and $\beta$ to diagrams
$\bar{\alpha}$ and $\bar{\beta}$:
\[
\begin{CD}
F_{\R^{m}}S^{n-1}_+ @>>> F_{\R^{m}}D^{n}_+\\
@V{f}VV @VV{\bar{g}}V\\
\tilde{G}^{i-1}(L) @>>> \tilde{G}^i(L)
\end{CD}\quad\text{and}\quad
\begin{CD}
F_{\R^{m'}}S^{n'-1}_+ @>>> F_{\R^{m'}}D^{n'}_+\\
@V{f'}VV @VV{\bar{g}'}V\\
\tilde{G}^{j-1}(K) @>>> \tilde{G}^j(K)
\end{CD}\quad.
\]
As in the old construction, we see that $\bar{\alpha}\square\bar{\beta}$ together with the map 
\[
\tilde{G}^{i-1}(L)\wedge \tilde{G}^j(K) \cup \tilde{G}^i(L)\wedge \tilde{G}^{j-1}(K)
\xrightarrow{\phi_{i-1,j}\cup\phi_{i,j-1}}
\tilde{G}^{i+j-2}(L\wedge K)
\]
determines a cell $\delta$ in $G^{i+j-1}(L\wedge K)$.
Let $H_{\alpha}$ be the subgroup of $\Sigma_n\times\Sigma_m$ that preserves $\alpha$,
$H_{\beta}$ the analogous subgroup for $\beta$,
and $H_{\delta}$ the subgroup of $\Sigma_{n+n'}\times\Sigma_{m+m'}$ that preserves $\delta$.
Now observe that $H_{\alpha}\amalg H_{\beta}$ then must be contained in $H_{\delta}$.
This ensures that all relations in $\tilde{G}^i(L)\wedge \tilde{G}^j(K)$
also give relations in $\tilde{G}^{i+j-1}(L\wedge K)$. Hence, $\phi_{i,j}$ is well-defined.
\end{Constr}

We finish the proof of theorem~\ref{Thm:orbitcofrepl} by showing:

\begin{Prop}
The natural transformation 
$\phi:\tilde{\Gamma}L\wedge\tilde{\Gamma}K\rightarrow\tilde{\Gamma}(L\wedge K)$
is symmetric.
\end{Prop}

\begin{proof}
Comparing $\Gamma$ with $\tilde{\Gamma}$ we see that there is a commutative diagram
\[
\begin{CD}
\Gamma L\wedge \Gamma K @>{\phi}>> \Gamma(L\wedge K)\\
@VVV @VVV\\
\tilde{\Gamma} L\wedge \tilde{\Gamma} K @>{\phi}>> \tilde{\Gamma}(L\wedge K)
\end{CD}\quad,
\]
where the vertical maps are quotient maps.
Proposition~\ref{Prop:Gammaskew} says that $\Gamma$ is skew-symmetric, with $\iota:\Gamma L\rightarrow \Gamma L$
measuring the failure of symmetry. Recall that $\iota$ was defined by flipping both the disks $D^n$
and the indexing spaces $\R^m$ of the cells in $\Gamma L$. Similarly we can define $\iota$ for $\tilde{\Gamma}L$.
And the proof of proposition~\ref{Prop:Gammaskew} works also in this case, and yields that $\tilde{\Gamma}$
is skew-symmetric. But when construction $\tilde{\Gamma}$ we divided out by all permutations of coordinates on both
$D^n$ and $\R^{m}$. This shows that $\iota=\id$ for $\tilde{\Gamma}L$. Hence, symmetry follows.
\end{proof}

\subsection{Equivariant features of $\Gamma$ and $\tilde{\Gamma}$}

For orthogonal $G$-spectra we have defined many different notions of cofibrancy.
In each case one could ask for a cofibrant replacement functor. 
For the category of orthogonal spectra we constructed such a functor $\Gamma$ in
theorem~\ref{Thm:cofrepl}. We also have the orbit cofibrant replacement functor $\tilde{\Gamma}$ 
defined above.
It seems unlikely that these functors can be used in the equivariant setting, but
by some miracle $\Gamma L$ and $\tilde{\Gamma}L$ are naive cofibrant and orbit cofibrant respectively,
and the natural maps $\Gamma L\rightarrow L$ and $\tilde{\Gamma}L\rightarrow L$
are naive level-equivalences. 
The author thinks it is unlikely that these maps are genuine level-equivalences.

In order to apply the functors $\Gamma$ and $\tilde{\Gamma}$ to orthogonal $G$-spectra,
we assume that $G$ is a finite and discrete group. By indexing
our orthogonal $G$-spectra $L$ by trivial representations only, we see that
such $L$ are the same as orthogonal spectra with $G$-action, i.e. functors $G\rightarrow\mathscr{IS}$.
Hence applying $\Gamma$ or $\tilde{\Gamma}$ to $L$ yields a new orthogonal spectrum
with $G$-action. We can ask whether or not $\Gamma L$ and $\tilde{\Gamma}L$ are cofibrant
(for one of the $G$-equivariant notions of cofibrancy).

We discuss the functor $\Gamma$ first. And we have:

\begin{Prop}\label{Prop:Gammaequivariantly}
Assume that $G$ is a finite group.
Let $L$ be an orthogonal $G$-spectrum, then $\Gamma L$ is naive $G$-cellular.
\end{Prop}

\begin{proof}
By induction it is enough to consider the gluing construction. Suppose that $p:A\rightarrow L$
is a $G$-equivariant map between orthogonal $G$-spectra. 
We must construct a relative $FI_G$-cellular structure on $A\rightarrow {G}(p)$.

Recall that $C$ denotes the set of all diagrams 
\[
\begin{CD}
F_{\R^m}S^{n-1}_+ @>>> F_{\R^m}D^n_+\\
@VVV @VVV\\
A @>{p}>> L
\end{CD}\quad,
\]
where $n,m\geq0$. Let $\alpha$ be a diagram in $C$, and let $H_{\alpha}$ be the subgroup
of $G$ consisting of those elements $g$ which preserve the cell $\alpha$.
We can now describe the gluing construction $G(p)$ equivariantly by the pushout diagram
\[
\begin{CD}
\bigvee F_{\R^{m}}S^{n-1}_+\wedge (G/H_{\alpha})_+ @>>> \bigvee F_{\R^{m}}D^{n}_+\wedge (G/H_{\alpha})_+\\
@VVV @VVV\\
A @>>> {G}(p)
\end{CD}\quad,
\]
where the wedge runs through one representative $\alpha$ for each $G$-orbit in $C$.
This implies that $A\rightarrow {G}(p)$ is a
relative naive $FI_G$-cellular map. 
\end{proof}

\begin{Lem}\label{Lem:fixedpointsandGamma}
Let $L$ be an orthogonal $G$-spectrum and $H$ a subgroup of $G$, then 
\[
(\Gamma L)^H(\R^m)=\Gamma (L^H)(\R^m)\quad,
\]
for all $\R^m$.
\end{Lem}

\begin{proof}
To see this we inspect the gluing construction for a map $p:A\rightarrow L$.
Let $\alpha$ be a diagram in $C$.
The result follows from the observation that the diagram
\[
\begin{CD}
F_{\R^m} S^{n-1}_+ @>>> F_{\R^m} D^{n}_+\\
@V{f}VV @VV{g}V\\
A @>{p}>> L
\end{CD}
\]
is fixed by $H$ if and only if $f$ and $g$ map into $A^H$ and $L^H$ respectively.
\end{proof}

\begin{Cor}\label{Cor:Gammaalsoequiv}
Let $L$ be an orthogonal $G$-spectrum.
The map $\Gamma L\rightarrow L$ induces a weak equivalence
\[
(\Gamma L)^H(\R^m)\rightarrow L^H(\R^m)
\]
for all $\R^m$ and $H$.
Consequently, the map $\Gamma L\rightarrow L$ is a naive $G$-equivariant level-equivalence.
\end{Cor}

This corollary says that $\Gamma$ is a cofibrant replacement functor for the naive model structure
on $G\mathscr{IS}$. But in general a naive level-equivalence is not a genuine level-equivalence,
see remark~\ref{Rem:comparewequiv}. 

We now turn to the case of the orbit cofibrant replacement functor $\tilde{\Gamma}$, and
show corresponding results:

\begin{Prop}\label{Prop:tildeGammaequivariantly}
Assume that $G$ is a finite group.
Let $L$ be an orthogonal $G$-spectrum, then $\tilde{\Gamma} L$ is orbit $G$-cellular.
\end{Prop}

\begin{proof}
It is enough to consider the gluing construction $\tilde{G}(p)$ for 
a $G$-equivariant map $p:A\rightarrow L$ between orthogonal $G$-spectra, and
we must construct a relative $\Orb_GFI$-cellular structure on $A\rightarrow \tilde{G}(p)$.

Recall that for fixed $n$ and $m$, the set $C_{n,m}$ consists of all non-trivial diagrams 
\[
\begin{CD}
F_{\R^m}S^{n-1}_+ @>>> F_{\R^m}D^n_+\\
@VVV @VVV\\
A @>{p}>> L
\end{CD}\quad.
\]
Observe that the group $G\times\Sigma_n\times \Sigma_m$ acts on $C'_{n,m}$.
Let $\alpha$ be a diagram in $C'_{n,m}$, and let $H_{\alpha}$ be the subgroup
of $G\times\Sigma_n\times \Sigma_m$ consisting of those elements $(g,\sigma,\rho)$ 
which preserve the cell $\alpha$.
We can now describe the gluing construction $\tilde{G}(p)$ equivariantly by the pushout diagram
\[
\begin{CD}
\bigvee \left(F_{\R^{m}}S^{n-1}_+\wedge G_+\right)/H_{\alpha} @>>> \bigvee \left(F_{\R^{m}}D^{n}_+\wedge G_+\right)/H_{\alpha}\\
@VVV @VVV\\
A @>>> \tilde{G}(p)
\end{CD}\quad,
\]
where the wedge runs through all $n$ and $m$ and
one representative $\alpha$ for each $G\times\Sigma_n\times \Sigma_m$-orbit in $C'_{n,m}$.

We now appeal to Illman's theorem~\cite{Illman:83} in order to show that each map
\[
\left(F_{\R^{m}}S^{n-1}_+\wedge G_+\right)/H_{\alpha} \rightarrow \left(F_{\R^{m}}D^{n}_+\wedge G_+\right)/H_{\alpha}
\]
is relative $\Orb_GFI$-cellular. To see this consider $\left(F_{\R^{m}}D^{n}_+\wedge G_+\right)(\R^m)$
as an $(O(m)\times H_{\alpha})$-space and triangulate to describe it as an $(O(m)\times H_{\alpha})$-CW-complex.
Dividing by $H_{\alpha}$ we get that the map above is orbit $G$-cellular.
This implies that $A\rightarrow \tilde{G}(p)$ is a
relative $\Orb_GFI$-cellular map. And the result follows.
\end{proof}

We have:

\begin{Prop}
For all orthogonal $G$-spectra $L$ the map $\tilde{\Gamma}L\rightarrow L$ is a naive level-equivalence.
\end{Prop}

\begin{proof}
It is enough to show that for all $G$-equivariant diagrams
\[
\begin{CD}
(S^{n-1}\times G/H)_+ @>>> \tilde{\Gamma}L(\R^m)\\
@VVV @VVV\\
(D^n\times G/H)_+ @>>> L(\R^m)
\end{CD}
\]
there is a lift $(D^n\times G/H)_+\rightarrow\tilde{\Gamma}L(\R^m)$.
Since $S^{n-1}$ is compact and $G$ finite, 
the map on the top factors through some $\tilde{G}^i(L)(\R^m)$.
Restricting the left side to some element of $G/H$ we now get the diagram
\[
\begin{CD}
S^{n-1}_+ @>>> \tilde{G}^i(L)(\R^m)\\
@VVV @VVV\\
D^n_+ @>>> L(\R^m)
\end{CD}\quad.
\]
Observe that the action of $H$ preserves this diagram.
And by construction of $\tilde{G}^{i+1}(L)$ we certainly have a lift $D^n_+\rightarrow \tilde{G}^{i+1}(L)(\R^m)$.
This lift is $H$-equivariant. Now a basic adjunction gives a $G$-equivariant lift
\[
(D^n\times G/H)_+\rightarrow \tilde{G}^{i+1}(L)(\R^m)\subset \tilde{\Gamma}L(\R^m)
\]
solving the first lifting problem.
\end{proof}

This proposition shows that $\tilde{\Gamma}$ is an orbit cofibrant replacement functor for
the naive $G$-equivariant model structure on $G\mathscr{IS}$.

\section{The diagonal map}

We consider the diagonal map in two cases, without and with involution.
For an $FI$-cellular orthogonal spectra $L$ (without involution), 
we describe a $C_q$-equivariant cell structure 
on the iterated smash product $L^{\wedge q}$. 
This uses the induced cells, see definition~\ref{Def:IndGFIcells}.
After this we construct the diagonal map $L^{\wedge q}\rightarrow \Phi^{C_r}(L^{\wedge rq})$.
It is defined for arbitrary $L$, but if $L$ is cofibrant, then the diagonal map is an isomorphism.
The proof uses the equivariant cell structure. 

In the second subsection, we repeat this for orthogonal spectra $L$ 
with involution. Recall that $D_{2q}$ denotes the dihedral group of order $2q$.
In the involutive case the diagonal map $L^{\wedge q}\rightarrow \Phi^{C_r}(L^{\wedge rq})$
is $D_{2q}$-equivariant, and an isomorphism when $L$ is genuine $\Z/2$-equivariantly cofibrant. 

\subsection{Without involution}

Let $L$ be an orthogonal spectrum. Consider the iterated smash product $L^{\wedge q}$.
Since $\wedge$ is symmetric, we get a $\Sigma_q$ action on $L^{\wedge q}$ by permuting factors.
In our applications we will only consider actions of cyclic groups, so for simplicity we restrict
our attention to the action of $C_q$ on $L^{\wedge q}$. 

Now assume that $L$ is $FI$-cellular. The following result describes an induced
$C_q$-cellular structure on $L^{\wedge q}$:

\begin{Prop}\label{Prop:indCqcellstroniteratedprod}
Let $C$ be the partially ordered set of cells for an $FI$-cellular orthogonal spectrum $L$.
The $q$-fold product $C^{\times q}$ has a $C_q$ action, and let $D$ be the
set of $C_q$-orbits. The $C_q$-equivariant structure on $L^{\wedge q}$ is given as follows:
\begin{itemize}
\item[] For each $[\boldsymbol{\alpha}]$ in $D$ there 
is a subspectrum $(L^{\wedge q})_{[\boldsymbol{\alpha}]}$ of
$L^{\wedge q}$, and 
$\bigcup_{[\boldsymbol{\alpha}]\in D}(L^{\wedge q})_{[\boldsymbol{\alpha}]}=L^{\wedge q}$.
\item[] $D$ is partially ordered by inclusion. 
We write $[\boldsymbol{\beta}]\leq[\boldsymbol{\alpha}]$ 
if $(L^{\wedge q})_{[\boldsymbol{\beta}]}\subseteq (L^{\wedge q})_{[\boldsymbol{\alpha}]}$. 
And for all $[\boldsymbol{\alpha}]$ the set 
$P_{[\boldsymbol{\alpha}]}=\{ [\boldsymbol{\beta}]\in D\;|\; [\boldsymbol{\beta}]<[\boldsymbol{\alpha}]\}$
is finite.
\item[] For every $[\boldsymbol{\alpha}]\in D$ there is 
pushout diagram with $C_q$-equivariant maps:
\[
\begin{CD}
\left(F_{\R^{ms}}S^{ns-1}_+\right)\wedge_{C_s}(C_q)_+
@>>> \left(F_{\R^{ms}}D^{ns}_+\right)\wedge_{C_s}(C_q)_+\\
@VVV @VVV\\
\bigcup_{[\boldsymbol{\beta}]<[\boldsymbol{\alpha}]} (L^{\wedge q})_{[\boldsymbol{\beta}]} 
@>>> (L^{\wedge q})_{[\boldsymbol{\alpha}]}
\end{CD}\quad,
\]
where $C_s$ is a subgroup of $C_q$, and $C_s$ acts 
on $\R^{ms}$, $S^{ns-1}$ and $D^{ns}$ by permuting the coordinates.
\end{itemize}
Moreover, each $\left(F_{\R^{ms}}S^{ns-1}_+\right)\wedge_{C_s}(C_q)_+\rightarrow 
\left(F_{\R^{ms}}D^{ns}_+\right)\wedge_{C_s}(C_q)_+$
can be subdivided as a relative $\Ind{C_q}FI$-cellular map. Hence, $L^{\wedge q}$ is
an induced cofibrant orthogonal $C_q$-spectrum.
\end{Prop}

\begin{proof}
Recall from proposition~\ref{Prop:prodcell} that 
$L^{\wedge q}$ has an $FI$-cellular structure with $C^{\times q}$ as its set of cells.
An $\boldsymbol{\alpha}=(\alpha_1,\ldots,\alpha_q)\in C^{\times q}$ represents an orbit 
$[\boldsymbol{\alpha}]$ in $D$ and we 
define $(L^{\wedge q})_{[\boldsymbol{\alpha}]}$ to be
\[
\bigcup_{\rho\in C_q}L_{\alpha_{\rho(1)}}\wedge\cdots\wedge L_{\alpha_{\rho(q)}}\quad.
\]
Clearly, $\bigcup_{[\boldsymbol{\alpha}]\in D}(L^{\wedge q})_{[\boldsymbol{\alpha}]}=L^{\wedge q}$.

Assume that $(L^{\wedge q})_{[\boldsymbol{\beta}]}\subseteq(L^{\wedge q})_{[\boldsymbol{\alpha}]}$
for some $\boldsymbol{\alpha}=(\alpha_1,\ldots,\alpha_q)$
and $\boldsymbol{\beta}=(\beta_1,\ldots,\beta_q)$.
Then each $L_{\beta_i}$ must be contained in some $L_{\alpha_j}$.
Since there are only finitely many $\delta\in C$ such that $L_{\delta}\subset L_{\alpha_j}$,
it follows that each $P_{[\boldsymbol{\alpha}]}$ is finite.

Consider some $\boldsymbol{\alpha}=(\alpha_1,\ldots,\alpha_q)$ in $C^{\times q}$, and
let $C_s$ be the subgroup of $C_q$ acting trivially on $\boldsymbol{\alpha}$.
To be more explicit:
Let $s$ be the greatest integer dividing $q$, where $t=\frac{q}{s}$, such that 
$\alpha_i=\alpha_j$ whenever $i\equiv j$ modulo $t$. We now see that $C_s$
acts trivially on the $q$-tuple $(\alpha_1,\ldots,\alpha_q)$.
By the $FI$-cellular structure on $L$ there is a pushout diagram
\[
\begin{CD}
F_{\R^{m_i}}S^{n_i-1}_+@>>> F_{\R^{m_i}}D^{n_i}_+\\
@VVV @VVV\\
\bigcup_{\beta<\alpha_i} L_{\beta} @>>> L_{\alpha_i}
\end{CD}
\]
for every $\alpha_i$. Smashing these diagrams together, as in 
lemma~\ref{Lem:wedgepushout}, we get a pushout diagram
\[
\begin{CD}
F_{\R^{m_1+\cdots+m_q}}S^{n_1+\cdots+n_q-1}_+@>>> F_{\R^{m_1+\cdots+m_q}}D^{n_1+\cdots+n_q}_+\\
@VVV @VVV\\
\bigcup_{\boldsymbol{\beta}<\boldsymbol{\alpha}} L_{\beta_1}\wedge\cdots\wedge L_{\beta_q} 
@>>> L_{\alpha_1}\wedge\cdots\wedge L_{\alpha_q}
\end{CD}\quad,
\]
where the maps are $C_s$-equivariant. $C_q$ acts on such diagrams, and taking the union over
all diagrams in the $C_q$-orbit of our $\boldsymbol{\alpha}$, we get a pushout diagram
\[
\begin{CD}
\left(F_{\R^{m_1+\cdots+m_q}}S^{n_1+\cdots+n_q-1}_+\right)\wedge_{C_s}(C_q)_+
@>>> \left(F_{\R^{m_1+\cdots+m_q}}D^{n_1+\cdots+n_q}_+\right)\wedge_{C_s}(C_q)_+\\
@VVV @VVV\\
\bigcup_{[\boldsymbol{\beta}]<[\boldsymbol{\alpha}]} (L^{\wedge q})_{[\boldsymbol{\beta}]} 
@>>> (L^{\wedge q})_{[\boldsymbol{\alpha}]}
\end{CD}\quad.
\]
Since $m_i=m_j$ and $n_i=n_j$ whenever $i\equiv j$ modulo $t$, we see that 
$C_s$ acts on $\R^{m_1+\cdots+m_q}$, $S^{n_1+\cdots+n_q-1}$ and $D^{n_1+\cdots+n_q}$ by
permuting coordinates. Now put $m=m_1+\cdots+m_t$ and $n=n_1+\cdots+n_t$.

The last statement of the proposition follows by applying Illman's
equivariant triangulation theorem for finite groups, see~\cite{Illman:78},
to produce a $C_s$-triangulation of $D^{ns}$.
\end{proof}

\begin{Rem}
Assume that $S\rightarrow L$ is a relative $FI$-cellular map. Recall the notation
$sL^{\wedge q-1}$ for the subspectrum 
\[
S\wedge L^{\wedge q-1}\cup L\wedge S\wedge L^{\wedge q-2}\cup \cdots \cup L^{\wedge q-1}\wedge S\subseteq L^{\wedge q}\quad.
\]
The proposition above immediately gives a $C_q$-equivariant description of $sL^{\wedge q-1}$.
It is the sub-$C_q$-spectrum of $L^{\wedge q}$ built using only those $[\boldsymbol{\alpha}]$ in $D$
with at least one $\alpha_i$ equal to the cell of $S$.
\end{Rem}

Next we construct the diagonal map.
Let $L$ be any orthogonal spectrum. As above the $q$-fold smash product $L^{\wedge q}$
is an orthogonal $C_q$-spectrum via the action that cyclically permutes the factors.
Similarly $L^{\wedge rq}$ is an orthogonal $C_{rq}$-spectrum. 
The diagonal will be a map $L^{\wedge q}\rightarrow\Phi^{C_r}(L^{\wedge rq})$. This exists for all positive numbers
$r$ and $q$.

\begin{Constr}
Let $E$ be the short exact sequence $0\rightarrow C_r\rightarrow C_{rq}\rightarrow C_q\rightarrow 0$.
Recall the definition of the category $\mathscr{J}_E$. 
Let $i:\mathscr{J}_E^{\text{reg}}\rightarrow \mathscr{J}_E$ be the
full subcategory whose objects are direct sums of regular $C_{rq}$-representations.
For an inner product space $U$, we have that $U^{\oplus rq}$ is an object of $\mathscr{J}_E^{\text{reg}}$.
First we construct a $C_q$-map
\[
L^{\wedge q}(U^{\oplus q})\rightarrow \Fix^{C_r}(L^{\wedge rq})(U^{\oplus rq})\quad.
\]
Using coends we can express the smash product $L^{\wedge q}(U^{\oplus q})$ as
\[
\int^{d_1,\ldots,d_q}\mathscr{J}(\R^{d_1}\oplus\cdots\oplus\R^{d_q},U^{\oplus q})
\wedge L(\R^{d_1})\wedge\cdots\wedge L(\R^{d_q})\quad.
\]
A point in this space is given by $(f,u;x_1,\ldots,x_q)$, where $f:\R^{d_1}\oplus\cdots\R^{d_q}\rightarrow U^{\oplus q}$
is an isometric embedding, $u$ is a point in $U^{\oplus q}$ orthogonal to $f$, and $x_i$ lies in $L(\R^{d_i})$.
To describe the $C_q$-action, we let $a:U^{\oplus q}\rightarrow U^{\oplus q}$ be the map which sends $(u_1,\ldots,u_q)$
to $(u_q,u_1,\ldots,u_{q-1})$, and we define 
$b:\R^{d_1}\oplus\cdots\oplus\R^{d_q}\rightarrow\R^{d_q}\oplus\R^{d_1}\oplus\cdots\oplus\R^{d_{q-1}}$
similarly. The preferred generator of $C_q$ then acts by sending
\[
(f,u;x_1,\ldots,x_q)\quad\text{to}\quad (afb^{-1},au;x_q,x_1,\ldots,x_{q-1})\quad.
\]
Analogously, we have that $L^{\wedge rq}(U^{\oplus rq})$ is equal to the coend
\[
\int^{d_1,\ldots,d_{rq}}\mathscr{J}(\R^{d_1}\oplus\cdots\oplus\R^{d_{rq}},U^{\oplus rq})
\wedge L(\R^{d_1})\wedge\cdots\wedge L(\R^{d_{rq}})\quad.
\]
To define the diagonal map, take a point
\[
(f,u;x_1,\ldots,x_q)\quad\text{in}\quad\mathscr{J}(\R^{d_1}\oplus\cdots\oplus\R^{d_q},U^{\oplus q})
\wedge L(\R^{d_1})\wedge\cdots\wedge L(\R^{d_q})
\]
and map it to
\[
\left(f^{\oplus r},(u,\ldots,u);(x_1,\ldots,x_q),(x_1,\ldots,x_q),\ldots,(x_1,\ldots,x_q)\right)
\]
in
\scriptsize
\[
\mathscr{J}((\R^{d_1}\oplus\cdots\oplus\R^{d_q})
\oplus\cdots\oplus(\R^{d_1}\oplus\cdots\oplus\R^{d_q}),U^{\oplus rq})
\wedge ( L(\R^{d_1})\wedge\cdots\wedge L(\R^{d_q}) )
\wedge \cdots\wedge ( L(\R^{d_1})\wedge\cdots\wedge L(\R^{d_q}))\quad.
\]
\normalsize
We easily see that the image is a $C_r$-fixed point, and that the diagonal map is $C_q$-equivariant.
The $C_q$-map, as constructed above, can be rewritten as a natural transformation
\[
\U_i\U_{\phi}L^{\wedge q}\rightarrow \U_i\Fix^{C_r}(L^{\wedge rq})
\]
of $\mathscr{J}^{\text{reg}}_E$-spaces, where $\U_i$ is the forgetful functor from
$\mathscr{J}^{\text{reg}}_E$-spaces to $\mathscr{J}_E$-spaces, and $\U_{\phi}$
denotes the functor given in the definition of
geometrical fixed points, see definition~\ref{Def:forgetandprolong}. 
By left Kan extension, we have a left adjoint $\bbP_i$ to $\U_i$, and the counit of this adjunction
is a natural transformation $\bbP_i\U_i\rightarrow\id$. Recall that $\bbP_{\phi}$
denotes the left adjoint to $\U_{\phi}$.
\end{Constr}

\begin{Def}\label{Def:diagonalmap}
The \textit{diagonal map} $L^{\wedge q}\rightarrow\Phi^{C_r}(L^{\wedge rq})$
is defined as the composition
\[
L^{\wedge q}\rightarrow \bbP_{\phi}\bbP_i\U_i\Fix^{C_r}(L^{\wedge rq})\rightarrow
\bbP_{\phi}\Fix^{C_r}(L^{\wedge rq})=\Phi^{C_r}(L^{\wedge rq})\quad,
\]
where the first map comes from the construction above, and the second map
is induced by the counit $\bbP_i\U_i\rightarrow\id$.
\end{Def}

\begin{Rem}\label{Rem:degeniteratedprod}
Assume that $S\rightarrow L$ is the inclusion of a subspectrum. Considering the restriction
of the diagonal map to $sL^{\wedge q-1}$, we get a diagram
\[
\begin{CD}
sL^{\wedge q-1} @>>> \Phi^{C^r}(s^rL^{\wedge rq-r})\\
@V{\subseteq}VV @VV{\subseteq}V\\
L^{\wedge q} @>>> \Phi^{C_r}(L^{\wedge rq})
\end{CD}\quad,
\]
where $s^rL^{\wedge rq-r}$ is the subspectrum of $L^{\wedge rq}$ given as
\[
s^rL^{\wedge rq-r}=\bigcup_i L^{\wedge i-1}\wedge S \wedge L^{\wedge q-1} \wedge S \wedge  L^{\wedge q-1} \wedge S \wedge \cdots
\wedge S \wedge L^{\wedge q-i}\quad.
\]
The existence of this diagram follows from inspection of the construction.
\end{Rem}

\begin{Exa}
Let us inspect the diagonal map in the case $L=F_V A$.
Then
\[
L^{\wedge q}= F_{V^{\oplus q}} (A^{\wedge q})\quad\text{and}\quad L^{\wedge rq}= F_{V^{\oplus rq}} (A^{\wedge rq})\quad.
\]
Computing the geometric $C_r$-fixed points of the last orthogonal spectrum, we get by
the formula in proposition~\ref{Prop:fundresultofgeomfp} that
\[
\Phi^{C_r}(L^{\wedge rq})=F_{(V^{\oplus rq})^{C_r}} (A^{\wedge rq})^{C_r}=F_{V^{\oplus q}} (A^{\wedge q})=L^{\wedge q}\quad.
\]
Inspecting definitions, we see that the diagonal map is an isomorphism in this case.
\end{Exa}

More generally we have:

\begin{Prop}\label{Prop:diagonalisomorphism}
If $L$ is a cofibrant orthogonal spectrum, then the diagonal map
$L^{\wedge q}\rightarrow \Phi^{C_r}(L^{\wedge rq})$ is a $C_q$-equivariant isomorphism.
\end{Prop}

It is enough to prove the result in the case where $L$ is $FI$-cellular.
The main idea is to use the equivariant description of the iterated products, given in 
proposition~\ref{Prop:indCqcellstroniteratedprod}.
But in order to apply this description we have to compute the 
geometric fixed points of induced cells. 

\begin{Lem}\label{Lem:gomefpofindcells}
Let $H$ be a subgroup and $N$ a normal subgroup of a compact Lie group $G$. 
Suppose that for any $N$-representation $W$ there exists a $G$-representation $U$
and an $N$-linear isometric embedding $W\rightarrow U$ such that $W^N=U^N$.
Assume that
$V$ is an $H$-representation and $A$ a based $H$-space. Then we have:
\[
\Phi^N\left( (F_V A)\wedge_H G_+\right)\cong\begin{cases}
(F_{V^N}A^N)\wedge_{J_0} J_+& \text{\begin{tabular}{l}if $N$ is a subgroup of $H$\\ $\quad$ with $J_0=H/N$, $J=G/N$, and\end{tabular}}\\
*&\text{\begin{tabular}{l}otherwise.\end{tabular}}
\end{cases}
\]
\end{Lem}

\begin{proof}
Assume first that $N$ is contained in $H$. Then proposition~\ref{Prop:geomfpandindsp} applies, and we have
\[
\Phi^N\left( (F_V A)\wedge_H G_+\right)\cong \left(\Phi^N (F_V A)\right)\wedge_{J_0} J_+\quad.
\]
Furthermore, proposition~\ref{Prop:fundresultofgeomfp} yields that
\[
\left(\Phi^N (F_V A)\right)\wedge_{J_0} J_+ \cong(F_{V^N}A^N)\wedge_{J_0} J_+\quad.
\]

Next assume that $N$ is not contained in $H$. 
The following is an elementary fact:
If $X$ is a based $H$-space, then the induced
$G$-space, $X\wedge_H G_+$,
has no non-trivial $N$-fixed points. It follows that 
$\Fix^N\left( (F_V A)\wedge_H G_+\right)=*$.
Consequently, also
$\Phi^N\left( (F_V A)\wedge_H G_+\right)=*$.
\end{proof}

By remark~\ref{Rem:abtNGconition} the condition is always true for finite $G$.
We are now ready to prove the proposition:

\begin{proof}
Let $L$ be an $FI$-cellular orthogonal spectrum with cells $C$.
We now apply proposition~\ref{Prop:indCqcellstroniteratedprod}.
Let $D$ be the $C_q$-orbits of $C^{\times q}$ and let $D'$ be the $C_{rq}$-orbits of $C^{\times rq}$.
Let $\epsilon^*$ be the map $C^{\times q}\rightarrow C^{\times rq}$ given by
\[
\epsilon^*(\alpha_1,\alpha_2,\ldots,\alpha_q)=
(\alpha_1,\alpha_2,\ldots,\alpha_q,\alpha_1,\alpha_2,\ldots,\alpha_q,\ldots,\alpha_1,\alpha_2,\ldots,\alpha_q)\quad.
\]
The $q$-tuple $(\alpha_1,\alpha_2,\ldots,\alpha_q)$ is repeated $r$ times. Passing to orbits
we get a map $\epsilon^*:D\rightarrow D'$.

Now we claim that:
\[
\Phi^{C_r}(L^{\wedge rq})=\bigcup_{[\boldsymbol{\beta}]\in D} \Phi^{C_r} (L^{\wedge rq})_{[\epsilon^*\boldsymbol{\beta}]}\quad.
\]
To prove this claim, we show by induction on the number of elements in $P_{[\boldsymbol{\alpha}]}$ that
\[
\Phi^{C_r}(L^{\wedge rq})_{[\boldsymbol{\alpha}]}=
\bigcup_{[\epsilon^*\boldsymbol{\beta}]<[\boldsymbol{\alpha}]} \Phi^{C_r} (L^{\wedge rq})_{[\epsilon^*\boldsymbol{\beta}]}\quad.
\]
There are two cases to consider when proving the induction step:
If $[\boldsymbol{\alpha}]=[\epsilon^*\boldsymbol{\beta}]$ for some $[\boldsymbol{\beta}]$ in $D$, then the induction
hypothesis is trivially true. The other case is when $[\boldsymbol{\alpha}]$ is not of the form $[\epsilon^*\boldsymbol{\beta}]$
for any $[\boldsymbol{\beta}]$ in $D$. We consider the diagram
\footnotesize
\[
\begin{CD}
\underset{[\boldsymbol{\delta}]<[\boldsymbol{\alpha}]}{\bigcup} \Phi^{C_r}(L^{\wedge rq})_{[\boldsymbol{\delta}]} 
&\leftarrow& \Phi^{C_r}\left((F_{\R^{ms}}S^{ns-1}_+)\wedge_{C_s}(C_{rq})_+\right)
&\rightarrow& \Phi^{C_r}\left((F_{\R^{ms}}D^{ns}_+)\wedge_{C_s}(C_{rq})_+\right)\\
@V{=}VV @VV{=}V @VV{=}V\\
\underset{[\epsilon^*\boldsymbol{\beta}]<[\boldsymbol{\alpha}]}{\bigcup} \Phi^{C_r}(L^{\wedge rq})_{[\epsilon^*\boldsymbol{\beta}]}
&\leftarrow& * &\rightarrow& * 
\end{CD}\quad.
\]
\normalsize
The left vertical map is an equality by induction.
Since $[\boldsymbol{\alpha}]$ is not of the form  $[\epsilon^*\boldsymbol{\beta}]$,
if follows that $C_r$ is not a subgroup of $C_s$, and hence lemma~\ref{Lem:gomefpofindcells}
implies that the two other vertical maps also are equalities.
By proposition~\ref{Prop:indCqcellstroniteratedprod}, the pushout of the top row
is $\Phi^{C_r}(L^{\wedge rq})_{[\boldsymbol{\alpha}]}$. This finishes the proof of the claim.

Our next claim is that the restriction of the diagonal map to $(L^{\wedge q})_{[\boldsymbol{\alpha}]}$ 
gives an isomorphism
\[
(L^{\wedge q})_{[\boldsymbol{\alpha}]}\xrightarrow{\cong} \Phi^{C_r}(L^{\wedge rq})_{[\epsilon^*\boldsymbol{\alpha}]}
\]
for all $[\boldsymbol{\alpha}]$ in $D$. We prove the claim by induction on the number of
elements in $P_{[\boldsymbol{\alpha}]}$. Consider the diagram
\footnotesize
\[
\begin{CD}
\underset{[\boldsymbol{\beta}]<[\boldsymbol{\alpha}]}{\bigcup} (L^{\wedge q})_{[\boldsymbol{\beta}]}
&\leftarrow& \left(F_{\R^{ms}}S^{ns-1}_+\right)\wedge_{C_s}(C_q)_+
&\rightarrow& \left(F_{\R^{ms}}D^{ns}_+\right)\wedge_{C_s}(C_q)_+\\
@VVV @VVV @VVV\\
\underset{[\boldsymbol{\beta}]<[\boldsymbol{\alpha}]}{\bigcup} \Phi^{C_r}(L^{\wedge rq})_{[\epsilon^*\boldsymbol{\beta}]}
&\leftarrow& \Phi^{C_r} \left(  (F_{\R^{mrs}}S^{nrs-1}_+)\wedge_{C_{rs}}(C_{rq})_+\right)
&\rightarrow& \Phi^{C_r} \left(  (F_{\R^{mrs}}D^{nrs}_+)\wedge_{C_{rs}}(C_{rq})_+\right)
\end{CD}\quad.
\]
\normalsize
Here the vertical maps are instances of the diagonal map. 
The row on the top comes from the $C_q$-equivariant description of $(L^{\wedge q})_{[\boldsymbol{\alpha}]}$,
while the bottom row is $\Phi^{C_r}$ applied to the $C_{rq}$-equivariant description
of $(L^{\wedge rq})_{[\epsilon^*\boldsymbol{\alpha}]}$.
By induction 
and the previous claim, the left vertical map is an isomorphism.
Since $C_r$ is a subgroup of $C_{rs}$, lemma~\ref{Lem:gomefpofindcells}
implies that the two other vertical maps also are isomorphisms.
Taking the row-wise pushouts proves our claim, see proposition~\ref{Prop:indCqcellstroniteratedprod}.

Clearly, the two claims together prove the proposition.
\end{proof}

\begin{Rem}
By the proof above, we see that the diagonal map restricted to $sL^{\wedge q-1}$ is an isomorphism
\[
sL^{\wedge q-1} \cong \Phi^{C^r}(s^rL^{\wedge rq-r})\quad.
\]
\end{Rem}

\subsection{With involution}

The dihedral group $D_{2q}$ has generators $x$ and $y$, and relations $x^{q}=y^2=1$ and $xy=yx^{-1}$.
Let $L$ be an orthogonal $\Z/2$-spectrum, the involution
is determined by a map $\iota:L\rightarrow L$. Now we can extend the $C_q$-action on $L^{\wedge q}$
to a $D_{2q}$-action by letting $y$ act by
\[
L\wedge L\wedge \cdots \wedge L\xrightarrow{\text{reverse}}L\wedge \cdots \wedge L\wedge L
\xrightarrow{\iota\wedge\cdots \wedge \iota\wedge \iota} L\wedge \cdots \wedge L\wedge L\quad.
\]

Assume that $L$ is genuine $FI_{\Z/2}$-cellular with $C$ as its set of $\Z/2$-cells.
Then the dihedral group $D_{2q}$ acts on $C^{\times q}$, $x$ acts by permuting 
factors cyclically, while $y$ sends a $q$-tuple $(\alpha_1,\alpha_2,\ldots,\alpha_q)$ to
$(\alpha_q,\ldots,\alpha_2,\alpha_1)$.
Now $L^{\wedge q}$ gets an induced $D_{2q}$-cellular structure by an argument similar
to proposition~\ref{Prop:indCqcellstroniteratedprod}.

Next we construct the diagonal map for $L$ with involution.
The diagonal will be a $D_{2q}$-map $L^{\wedge q}\rightarrow\Phi^{C_r}(L^{\wedge rq})$. 

\begin{Constr}
Let $E$ be the short exact sequence $0\rightarrow C_r\rightarrow D_{2rq}\rightarrow D_{2q}\rightarrow 0$.
Let $i:\mathscr{J}_E^{\text{reg}}\rightarrow \mathscr{J}_E$ be the
full subcategory whose objects are direct sums of regular $D_{2rq}$-representations.
For each inner product space $U$, we get an object $U^{\oplus 2rq}$ of $\mathscr{J}_E^{\text{reg}}$.
Let $\U_i$ be the forgetful functor from
$\mathscr{J}^{\text{reg}}_E$-spaces to $\mathscr{J}_E$-spaces, and by left Kan extension, we 
define $\bbP_i$.

As before we construct a $D_{2q}$-map
\[
\U_i\U_{\phi}L^{\wedge q}\rightarrow \U_i\Fix^{C_r}(L^{\wedge rq})
\]
of $\mathscr{J}^{\text{reg}}_E$-spaces. This map is given on level $U^{\oplus 2rq}$ as
a natural $D_{2q}$-equivariant transformation
\[
L^{\wedge q}(U^{\oplus 2q})\rightarrow \Fix^{C_r}(L^{\wedge 2rq})(U^{\oplus 2rq})\quad.
\]
To define this, we explicitly write out both sides using coends. Picking a point on the left side,
\[
(f,u;x_1,\ldots,x_q)\quad\text{in}\quad\mathscr{J}(\R^{d_1}\oplus\cdots\oplus\R^{d_q},U^{\oplus q})
\wedge L(\R^{d_1})\wedge\cdots\wedge L(\R^{d_q})
\]
we map it to
\[
\left(f^{\oplus r},(u,\ldots,u);(x_1,\ldots,x_q),(x_1,\ldots,x_q),\ldots,(x_1,\ldots,x_q)\right)
\]
in
\scriptsize
\[
\mathscr{J}((\R^{d_1}\oplus\cdots\oplus\R^{d_q})
\oplus\cdots\oplus(\R^{d_1}\oplus\cdots\oplus\R^{d_q}),U^{\oplus rq})
\wedge ( L(\R^{d_1})\wedge\cdots\wedge L(\R^{d_q}) )
\wedge \cdots\wedge ( L(\R^{d_1})\wedge\cdots\wedge L(\R^{d_q}))\quad.
\]
\normalsize
Again it is easily seen that the image is a $C_r$-fixed point, 
and that the diagonal map is $D_{2q}$-equivariant.
\end{Constr}

\begin{Def}\label{Def:diagonalmapinvolutive}
The \textit{diagonal map} for $L$ with involution, $L^{\wedge q}\rightarrow\Phi^{C_r}(L^{\wedge rq})$,
is defined as the composition
\[
L^{\wedge q}\rightarrow \bbP_{\phi}\bbP_i\U_i\Fix^{C_r}(L^{\wedge rq})\rightarrow
\bbP_{\phi}\Fix^{C_r}(L^{\wedge rq})=\Phi^{C_r}(L^{\wedge rq})\quad,
\]
where the first map comes from the construction above, and the second map
is induced by the counit $\bbP_i\U_i\rightarrow\id$.
\end{Def}

Assume that $L$ is genuine $FI_{\Z/2}$-cellular.
Using explicit induced $D_{2q}$- and $D_{2rq}$-cellular structures on $L^{\wedge q}$ and
$L^{\wedge rq}$ respectively, we prove the following result an argument similar to
the proof of proposition~\ref{Prop:diagonalisomorphism}:

\begin{Prop}\label{Prop:diagonalisomorphisminvolutive}
If $L$ is a cofibrant orthogonal $\Z/2$-spectrum, then the diagonal map 
$L^{\wedge q}\rightarrow \Phi^{C_r}(L^{\wedge rq})$ is a $D_{2q}$-equivariant isomorphism.
\end{Prop}

Further details are omitted.

\section{Miscellaneous results}

In this section we list or show various results which will be used later in the thesis.

\subsection{About orbit q-cofibrations of orthogonal spectra}

\begin{Prop}\label{Prop:orbitsquareproduct}
Assume that $i:A\rightarrow L$ and $j:B\rightarrow K$ are orbit q-cofibrations
of (non-equivariant) orthogonal spectra, then also
\[
i\square j:L\wedge B\cup A\wedge K\rightarrow L\wedge K
\]
is an orbit q-cofibration.
\end{Prop}

\begin{proof}
We can assume that $i$ and $j$ are relative $\Orb FI$-cellular maps.
We will apply the method used when proving proposition~\ref{Prop:prodcell}.
In that proof we use only one non-formal property about the set of cells, namely
that the $\square$ product of two cells yields a new cell. Once this property has been checked for
orbit cells, the present result follows.

Consider two (non-equivariant) orbit cells:
\[
\left(F_{V_1} S^{n_1-1}_+\right)/H_1\rightarrow\left(F_{V_1} D^{n_1}_+\right)/H_1\quad\text{and}\quad
\left(F_{V_2} S^{n_2-1}_+\right)/H_2\rightarrow\left(F_{V_2} D^{n_2}_+\right)/H_2\quad.
\]
Here $V_1$ and $V_2$ are $H_1$- and $H_2$-representations respectively.
Let $V=V_1\oplus V_2$ be the $H_1\times H_2$-representation defined by letting $H_1$ act trivially on $V_2$,
and $H_2$ act trivially on $V_1$. And $\square$ of the cells above can now be written as
\[
\left( F_{V} (D^{n_1}\times S^{n_2-1}\cup S^{n_1-1}\times D^{n_2} )_+ \right)/(H_1\times H_2)
\rightarrow \left( F_{V} (D^{n_1}\times D^{n_2} )_+ \right)/(H_1\times H_2)\quad.
\]
This is again an orbit cell.
\end{proof}

\begin{Prop}\label{Prop:orbitpreservele}
Assume that $L$ is an orbit cofibrant orthogonal spectrum and $X\rightarrow Y$ a level-equivalence
of orthogonal spectra, then also
\[
L\wedge X\rightarrow L\wedge Y
\]
is a level-equivalence.
\end{Prop}

\begin{proof}
There is no loss of generality by assuming that $L$ is $\Orb FI$-cellular.
By analogy with proposition~\ref{Prop:seqcelldef}, there exists a sequence 
$*=L_0\rightarrow L_1\rightarrow L_2\rightarrow\cdots$ such that $L$ is its colimit and
each $L_{i+i}$ is the pushout of a diagram
\[
L_i\leftarrow\bigvee \left(F_{V} S^{n-1}_+\right)/H\rightarrow\bigvee\left(F_{V} D^{n}_+\right)/H\quad.
\]
Hence, each $L_i\rightarrow L_{i+1}$ is an h-cofibration.
Smashing the sequence with $X\rightarrow Y$, we get
\[
\begin{CD}
L_0 \wedge X @>>> L_1 \wedge X @>>> L_2 \wedge X @>>> \cdots\\
@VVV    @VVV     @VVV\\
L_0 \wedge Y @>>> L_1 \wedge Y @>>> L_2 \wedge Y @>>> \cdots
\end{CD}\quad.
\]
The horizontal maps are again h-cofibrations, by lemma~12.2 in~\cite{MandellMaySchwedeShipley:01}.
Hence, it is enough to show that each $L_i\wedge X\rightarrow L_i\wedge Y$ is a level-equivalence.
Proceeding inductively, what we have to show is that
\[
\left(F_{V} A\right)/H\wedge X\rightarrow\left(F_{V} A\right)/H\wedge Y
\]
is a level-equivalence when $A$ is a sphere or a disk and $H$ and $V$ arbitrary.
By remark~\ref{Rem:orbitcellredundancy}, we can always assume that $H\rightarrow O(V)$ is injective.
Then the statement above is a consequence of lemma~\ref{Lem:htpyoforbitcell} below. 
And we are done modulo proving the lemma.
\end{proof}

Inspired by lemma~15.5 in~\cite{MandellMaySchwedeShipley:01} we prove:

\begin{Lem}\label{Lem:htpyoforbitcell}
Let $A$ be a based CW-complex, $H$ a finite group, $V$ an $H$-representation and $X$ any
orthogonal spectrum. Assume that $H\rightarrow O(V)$ is injective.
Then the quotient map
\[
\left( EH_+\wedge_H F_V A\right)\wedge X\rightarrow \left(F_{V} A\right)/H\wedge X
\]
is a level-equivalence. Consequently, the functor $\left(F_{V} A\right)/H\wedge -$
preserves level-equivalences.
\end{Lem}

\begin{proof}
We evaluate both sides of the quotient map at some level $\R^m$. In order to
write things out we choose a linear isometry $V\rightarrow\R^m$. Then we have
\[
\left((F_{V} A)\wedge X\right)(\R^m)\cong
O(m)_+\wedge_{O(\R^m-V)}(A\wedge X(\R^m-V))\quad,
\]
and
\[
\left( (EH_+\wedge F_V A)\wedge X\right)(\R^m)\cong
O(m)_+\wedge_{O(\R^m-V)}(EH_+\wedge A\wedge X(\R^m-V))\quad.
\]
The quotient $EH\times O(m)\rightarrow O(m)$ is an $(H\times O(\R^m-V))$-equivariant homotopy equivalence
since $O(m)$ is a free $(H\times O(\R^m-V))$-space that can be triangulated as a finite
$(H\times O(\R^m-V))$-CW-complex by~\cite{Illman:83}.

We compare the description above via the quotient map $EH\times O(m)\rightarrow O(m)$.
Dividing out by the $H$-action, we get a weak equivalence
\begin{multline*}
\left(O(m)_+\wedge_{O(\R^m-V)}(EH_+\wedge A\wedge X(\R^m-V))\right)/H\\
\rightarrow\left(O(m)_+\wedge_{O(\R^m-V)}(A\wedge X(\R^m-V))\right)/H\quad.
\end{multline*}
And the result follows.
\end{proof}

Analogous to proposition~\ref{Prop:qcofincolim}
we have the following proposition for orbit q-cofibrations:

\begin{Prop}\label{Prop:orbitqcofincolim}
Assume that we have a map between two sequences of orthogonal spectra:
\[
\begin{CD}
K_0 @>>> K_1 @>>> K_2 @>>> \cdots\\
@VVV    @VVV     @VVV\\
L_0 @>>> L_1 @>>> L_2 @>>> \cdots
\end{CD}\quad.
\]
If $K_0\rightarrow L_0$ is an orbit q-cofibration, and
$K_i\cup_{K_{i-1}}L_{i-1}\rightarrow L_i$ is an orbit q-cofibration for every $i\geq0$, then
\[
\colim_i K_i\rightarrow \colim_i L_i
\]
is also an orbit q-cofibration. 
\end{Prop}

\begin{proof}
If orbit q-cofibrations were the cofibrations of a model structure on $\mathscr{IS}$, then
the formal proof given for proposition~\ref{Prop:qcofincolim} would apply.
Lacking such a model structure we give a direct proof.

Observe that the following diagram is pushout for all $n$:
\[
\begin{CD}
L_n\cup_{K_n}K_{n+1} @>>> L_n\cup_{K_n}\colim K_i\\
@VVV @VVV\\
L_{n+1} @>>> L_{n+1}\cup_{K_{n+1}} \colim K_i
\end{CD}\quad.
\]
The left vertical map is an orbit q-cofibration by assumption, hence the right vertical map
is also an orbit q-cofibration.
We thus get a sequence
\[
\colim K_i\rightarrow L_0\cup_{K_0}\colim K_i \rightarrow L_1\cup_{K_1}\colim K_i \rightarrow \cdots \rightarrow \colim L_i
\]
of orbit q-cofibrations. 
An elementary argument, similar to the last part of the proof of proposition~\ref{Prop:seqcelldef},
shows that given a sequence $X_0\rightarrow X_1\rightarrow X_2\rightarrow \cdots$ of
orbit q-cofibrations maps, the induced map $X_0\rightarrow \colim X_i$ is also an orbit q-cofibration.
Applied to our situation we see that $\colim_i K_i\rightarrow \colim_i L_i$ is an orbit q-cofibration. 
\end{proof}

\subsection{About the small object argument}

The small object argument, see~\cite{DwyerSpalinski:95} or~\cite{Hirschhorn:03}, is the most common way to produce both
cofibrant and fibrant replacement functors in a model category. We have
already used this construction when defining our functor $\Gamma$. 
In general we have proposition~10.5.16 in~\cite{Hirschhorn:03} which says:

\begin{Prop}
If $\mathscr{C}$ is a cocomplete category and $I$ a set of maps
in $\mathscr{C}$ that permits the small object argument, then there is a functorial
factorization of every map in $\mathscr{C}$ into a relative $I$-cell complex followed
by a map having the right lifting property with respect to any map in $I$.
\end{Prop}

So given $I$, we get a functor taking a map $f:X\rightarrow Y$ in $\mathscr{C}$ to a factorization
\[
X\xrightarrow{j_I} Q_I(f)\xrightarrow{p_I} Y\quad.
\]
And $Q_I$ is defined as the colimit of a sequence $X=Q^0_I(f)\rightarrow Q^1_I(f)\rightarrow Q^2_I(f)\rightarrow\cdots$,
where each step is a gluing construction.

We now compare the small object arguments in two different categories:

\begin{Lem}\label{Lem:smallobjandfunctors}
Assume that $F:\mathscr{C}\rightarrow\mathscr{D}$ is a functor between cocomplete categories.
Let $I$ and $J$ be sets of maps in $\mathscr{C}$ and $\mathscr{D}$ respectively, and suppose
that they permit the small object argument. If $F$ takes $I$ into $J$, then there is a functorial 
diagram
\[
\begin{CD}
F(X)@>{F(j_I)}>> F(Q_I(f)) @>{F(p_I)}>> F(Y)\\
@V{=}VV @VVV @VV{=}V\\
F(X)@>{j_J}>> Q_J(F(f)) @>{p_J}>> F(Y)
\end{CD}
\]
for every map $f:X\rightarrow Y$ in $\mathscr{C}$.
\end{Lem}

\begin{proof}
Assume by induction that we have a sequence
\[
F(X)\rightarrow F(Q^i_I(f))\rightarrow Q^i_J(F(f))\rightarrow F(Y)\quad.
\]
Now consider one $I$-cell in $Q^{i+1}_I(f)$ relative to $Q^i_I(f)$. Such cells are determined by
a diagram
\[
\begin{CD}
A @>>> B\\
@VVV @VVV\\
Q^i_I(f) @>>> Y
\end{CD}\quad,
\]
where $A\rightarrow B$ is a map in $I$. Now apply $F$ to this diagram, and use the induction hypothesis to
form
\[
\begin{CD}
F(A) @>>> F(B)\\
@VVV @VVV\\
F(Q^i_I(f)) @>>> F(Y)\\
@VVV @VV{=}V\\
Q^i_J(F(f)) @>>> F(Y) 
\end{CD}\quad.
\]
Comparing the pushout of the upper square with the pushout of the outer square, we get the sequence
\[
F(X)\rightarrow F(Q^{i+1}_I(f))\rightarrow Q^{i+1}_J(F(f))\rightarrow F(Y)\quad.
\]
\end{proof}

\begin{Exa}\label{Exa:fibreplandgeomfp}
A fibrant replacement functor $Q^G$ for the category of orthogonal $G$-spectra can
be constructed by applying the small object argument to a set of maps $K$ called
the generating acyclic q-cofibrations. See definition~III.4.6 in~\cite{MandellMay:02}.
By proposition~\ref{Prop:fundresultofgeomfp} $\Phi^N$ takes $K$ to the generating
acyclic q-cofibrations of orthogonal $G/N$-spectra. Hence, there is a natural transformation
\[
\Phi^N Q^G L\rightarrow Q^{G/N} \Phi^N L
\]
for orthogonal $G$-spectra $L$.
\end{Exa}

\subsection{About geometric fixed points}

We recall corollary~V.4.6 and proposition~V.4.7 in~\cite{MandellMay:02}, and enhance the last result
to also cover a new case:

\begin{Prop}\label{Prop:geomfpofsuspandprod}
For based $G$-spaces $A$, the geometric $N$-fixed points of the suspension spectrum are given by
\[
\Phi^N F_0 A\cong F_0 A^N\quad.
\]
For orthogonal $G$-spectra $K$ and $L$, there is a natural $J$-map
\[
\alpha:\Phi^N K\wedge \Phi^N L\rightarrow \Phi^N(K\wedge L)
\]
of orthogonal $J$-spectra, and $\alpha$ is an isomorphism if $K$ and $L$ are cofibrant.
Furthermore, $\alpha$ is also an isomorphism if $K$ is a suspension spectrum and $L$ is arbitrary.
\end{Prop}

\begin{proof}
The first parts are cited from the reference. To show the last part, 
we do an explicit calculation. Assume that  $K=F_0A$, where $A$ is a based $G$-space.
We then have
\[
\Fix^N(F_0 A\wedge L)\cong \Fix^N(A\wedge L)\cong A^N\wedge (\Fix^NL)  \cong F_0A^N\wedge (\Fix^NL)\quad.
\]
Here $F_0$ on the left is the suspension from $G$-spaces to orthogonal $G$-spectra, while $F_0$ on
the right is the suspension from $G$-spaces to $\mathscr{J}_E$-spaces.
Since $\bbP_{\phi}$ is strong symmetric monoidal, see~I.2.14 in~\cite{MandellMay:02}, we get that
\[
\Phi^N(F_0 A\wedge L)=\bbP_{\phi}\Fix^N(F_0 A\wedge L)\cong\bbP_{\phi}F_0A^N\wedge \bbP_{\phi}(\Fix^NL)= F_0A^N\wedge\Phi^NL\quad.
\]
This completes the proof.
\end{proof}

Geometric fixed points can be used in order to recognize $G$-equivariant $\pi_*$-isomorphisms:

\begin{Prop}\label{Prop:geomfpdetectspiiso}
Assume that $\mathscr{F}$ is a family of normal subgroups of $G$. Let $f:K\rightarrow L$
be a map of orthogonal $G$-spectra. The following statements are equivalent:
\begin{itemize}
\item[] $f$ is an $(\mathscr{F},\All)$-$\pi_*$-isomorphism.
\item[] $f$ induces isomorphisms $\pi^N_q K\rightarrow \pi^N_q L$ for all $q$ and $N\in\mathscr{F}$.
\item[] $\Phi^Nf:\Phi^N K\rightarrow \Phi^N L$ is non-equivariantly a $\pi_*$-isomorphism for all $N\in\mathscr{F}$.
\end{itemize}
\end{Prop}

Compare this result to theorem~4.7 in~\cite{GreenleesMay:95b}.

\begin{proof}
The first two statements are equivalent by definition. Since compact Lie groups have the descending chain property,
we can do induction on the size of $\mathscr{F}$. 

Let us show that the last statement implies the second.
Let $\mathscr{F}$ be some family of normal subgroups of $G$.
Assume that we are given a map $f:K\rightarrow L$ such that 
$\Phi^Hf:\Phi^H K\rightarrow \Phi^H L$ is non-equivariantly a $\pi_*$-isomorphism for all $H\in\mathscr{F}$.
Let $N$ be any normal subgroup of $\mathscr{F}$. What we want to check is that
$f$ induces isomorphisms $\pi^N_* K\rightarrow \pi^N_* L$.

Let $\mathscr{F}[N]$ be the family of subgroups of $G$ that do not contain $N$. Observe that the intersection
$\mathscr{F}\cap\mathscr{F}[N]$ is a family properly contained in $\mathscr{F}$. By the induction hypothesis,
we get that $f$ induces isomorphisms $\pi^H_* K\rightarrow \pi^H_* L$ for all
$H\in \mathscr{F}\cap\mathscr{F}[N]$. Observe that the family $\mathscr{F}\cap\mathscr{F}[N]$
consists of all proper subgroups of $N$.

Recall the notion of a universal $\mathscr{F}[N]$-space. It is a $G$-CW-complex $E\mathscr{F}[N]$ such that
$E\mathscr{F}[N]^H\simeq*$ for $H\in\mathscr{F}[N]$ and $E\mathscr{F}[N]^H=\emptyset$ for $H\not\in\mathscr{F}[N]$.
Furthermore, $\tilde{E}\mathscr{F}[N]$ is defined as the cofiber of the map $E\mathscr{F}[N]_+\rightarrow S^0$.

Let us now restrict $G$-actions to $N$-actions. Observe that the restriction of $E\mathscr{F}[N]$ is
a universal $\mathscr{F}\cap\mathscr{F}[N]$-space. We already know that $f:K\rightarrow L$
is an $\mathscr{F}\cap\mathscr{F}[N]$-equivalence. And proposition~IV.6.7 in~\cite{MandellMay:02}
now implies that
\[
f\wedge\id:K\wedge E\mathscr{F}[N]_+\rightarrow L\wedge E\mathscr{F}[N]_+
\]
is a $\pi_*$-isomorphism of orthogonal $N$-spectra.

Let $Q$ denote a fibrant replacement functor for the genuine model structure on orthogonal $G$-spectra.
Consider the diagram
\[
\begin{CD}
Q(K\wedge E\mathscr{F}[N]_+)^N @>>> Q(K)^N @>>> Q(K\wedge\tilde{E}\mathscr{F}[N])^N\\
@VVV @VVV @VVV\\
Q(L\wedge E\mathscr{F}[N]_+)^N @>>> Q(L)^N @>>> Q(L\wedge\tilde{E}\mathscr{F}[N])^N
\end{CD}
\]
of non-equivariant orthogonal spectra.
We have just shown that the left vertical map is a $\pi_*$-isomorphism.
By proposition~V.4.17 in~\cite{MandellMay:02}, the orthogonal spectrum
$Q(K\wedge\tilde{E}\mathscr{F}[N])^N$ is naturally $\pi_*$-isomorphic to $\Phi^N K$.
By statement three, we get that
\[
\Phi^N K\simeq Q(K\wedge\tilde{E}\mathscr{F}[N])^N\rightarrow Q(L\wedge\tilde{E}\mathscr{F}[N])^N\simeq \Phi^N L
\]
is a $\pi_*$-isomorphism. Now it follows by long exact sequences of homotopy groups, that
\[
\pi_*^N(K)\cong \pi_*Q(K)^N\rightarrow \pi_*Q(L)^N\cong \pi^N_* L
\]
is an isomorphism.

To show that the second statement implies the last statement, one uses the above argument backward.
\end{proof}

\begin{Rem}
The reason for assuming that $\mathscr{F}$ consists of \textbf{normal} subgroups
is that geometric $N$-fixed points of orthogonal $G$-spectra, have been defined only
for normal $N$. Hence, proposition~V.4.17 in~\cite{MandellMay:02} supplies the homotopy equivalence
\[
Q(K\wedge\tilde{E}\mathscr{F}[N])^N\cong \Phi^N K
\]
in this case only.
\end{Rem}

\section{Cyclic and dihedral orthogonal spectra}

In this short section we will define 
comment on how 
the geometric realization of involutive simplicial, cyclic and dihedral orthogonal spectra
becomes equivariant orthogonal spectra.

\begin{Exa}
First consider an involutive simplicial orthogonal spectrum.
This is a functor $L_{\bullet}:\catDeltaT^{\op}\rightarrow \mathscr{IS}$.
Taking the geometric realization level-wise, we get an orthogonal spectrum
$|L_{\bullet}|$ with $\Z/2$-action. By change of universe, see lemma~\ref{Lem:ChangeUniverse},
we can evaluate
this spectrum at any $\Z/2$-representation $V$.
\end{Exa}

\begin{Exa}
For an $r$-cyclic orthogonal spectrum $L_{\bullet}$ the geometric realization
$|L_{\bullet}|$ has an $S^1$-action. Using the change of universe functor,
see lemma~\ref{Lem:ChangeUniverse}, we can evaluate at any $S^1$-representation
$V$ getting an $S^1$-space
\[
|L_{\bullet}|(V)\quad.
\]
\end{Exa}

\begin{Exa}
Similarly, if $L_{\bullet}$ is an $r$-dihedral orthogonal spectrum, then using
lemma~\ref{Lem:ChangeUniverse}, we see that 
\[
|L_{\bullet}|(V)
\]
is well defined for any $O(2)$-representation $V$.
\end{Exa}

\chapter[Operads in $\mathscr{IS}$ and involution]{Operads in orthogonal spectra and involution}\label{Cha:invoper}

We begin this chapter by studying operads and orthogonal spectra. Traditionally an
operad consists of topological spaces together with composition operations.
In the book~\cite{MarklShniderStasheff:02} the definition of an operad
is extended by replacing topological spaces by objects in a symmetric monoidal category. 
We recall this definition below in full generality, but our main focus will be
operads in orthogonal spectra.

One usually designs an operad in order to study its algebras. Classically, we have the
operads (in topological spaces) $\mathcal{N}$ and $\mathcal{M}$, whose algebras are 
the commutative and the associative monoids respectively. We will here introduce an operad $\mathcal{H}$
having associative monoids with involution as its algebras. Via suspension
these results extend to orthogonal spectra. In particular an $\mathcal{H}$-algebra
in orthogonal spectra is an orthogonal ring spectrum with involution.

In section~\ref{sect:operadsbarconstr} we develop, along the lines of~\cite{May:72},
the theory of the two-sided bar construction. 
Under the hypothesis that $\mathcal{P}$ and $\mathcal{Q}$ are sufficiently ``equal up to homotopy'',
we can use this construction to replace a $\mathcal{P}$-algebra by a 
weakly equivalent $\mathcal{Q}$-algebra. This is made precise in remark~\ref{Rem:replproc}.

In section~\ref{sect:involution} we study a geometrically interesting
involution $\iota$ on the homotopy groups $\pi_*S[\Omega M]$. 
The main result of the thesis, theorem~\ref{Thm:main}, says that there exists an
orthogonal ring spectrum $R$ with involution which represents $\iota$ on $\pi_*S[\Omega M]$.
The proof uses the machinery of operads. We design an operad $\mathcal{D}_n$ in
orthogonal spectra which has $S[\Omega M]$ as an algebra. The result follows 
by showing that $\mathcal{D}_n$ is sufficiently equal to $\mathcal{H}$.

\section{Operads in $\mathscr{IS}$}

We will now begin looking at operads in orthogonal spectra. May's original
definition~\cite{May:77} describes the composition as a many-variable operation.
This can be replaced by a collection of two-variable compositions. Using this
viewpoint, one greatly reduces the complexity of the description of associativity.
This description is due to Gerstenhaber and Markl, see~\cite{MarklShniderStasheff:02}.

Let $\Sigma$ be the category with objects the finite sets $\mathbf{n}=\{1,2,\ldots,n\}$ 
for every non-negative integer $n$ and
bijections as morphisms. Here $\mathbf{0}$ is the empty set.
Thus there is no morphism $\mathbf{n}\rightarrow\mathbf{m}$
for $n\neq m$, while the endomorphisms of $\mathbf{n}$ can be identified with the
symmetric group $\Sigma_n$. Therefore, we call $\Sigma$ the symmetric groupoid.

There are composition operations
\[
\circ_i:\Sigma_m\times\Sigma_n\rightarrow \Sigma_{m+n-1}
\]
for $n\geq0$ and $1\leq i\leq m$. These are defined in the appendix. Let me recall
the ``box''-model here: For $i$ and permutations $\rho\in\Sigma_m$ and $\upsilon\in\Sigma_n$
we put boxes around the integers from $1$ to $m+n-1$
as follows:
\[
\boxed{1},\ldots,\boxed{i-1},\boxed{i, i+1,\ldots, i+n-1},\boxed{i+n},\ldots,\boxed{m+n-1}\quad.
\]
We now use $\rho$ to permute the boxes, while we use $\upsilon$ to permute the
elements in the $i$'th box. Removing the boxes one gets the permutation $\rho\circ_i\upsilon$.
This operation gives the symmetric groupoid the structure of a discrete operad. We
call this operad $\mathcal{M}$.

There is an alternative description of the composition operations using permutation matrices.
Recall that the permutation matrix of $\rho\in \Sigma_m$ is the unique $m\times m$-matrix $A$ 
such that
\[
A e_i=e_{\rho(i)}\quad\text{for all $i$.}
\]
Here $e_i$ is the $i$'th unit vector in $\R^m$. This embeds $\Sigma_m$ as a closed subgroup of $O(m)$.
A matrix in the image is called a permutation matrix, and these are exactly those matrices such that
every column and every row contain only $0$'s, except for one entry which has value $1$.

Let $A$ be the permutation matrix of $\rho\in\Sigma_m$ and $B$ the permutation matrix for $\nu\in\Sigma_n$.
We now want to describe the permutation matrix for $\rho\circ_i\nu$. We have a block decomposition of $A$ as
\[
\begin{pmatrix}
A_{11} & 0 & A_{13}\\
0 & 1 & 0\\
A_{31} & 0 & A_{33}
\end{pmatrix}\quad,
\]
where $1$ lies in the $\rho(i)$'th row and the $i$'th column. Now form the $(m+n-1)\times(m+n-1)$-matrix
\[
\begin{pmatrix}
A_{11} & 0 & A_{13}\\
0 & B & 0\\
A_{31} & 0 & A_{33}
\end{pmatrix}\quad.
\]
This is the permutation matrix for $\rho\circ_i\nu$.

In the case $n=0$ we interpret $B$ as the $0\times0$-matrix. Thus the operation $\circ_i$
deletes the $i$'th column and the $\rho(i)$'th row from the matrix $A$.

Now we are ready to define operads in orthogonal spectra. 

\begin{Def}
Let $(\mathscr{C},\wedge,S)$ be a symmetric monoidal category.
An \textit{operad in $\mathscr{C}$} is a functor $\mathcal{P}:\Sigma\rightarrow \mathscr{C}$
with $\mathcal{P}(0)=S$
together with composition operations
\[
\circ_i:\mathcal{P}(m)\wedge\mathcal{P}(n)\rightarrow\mathcal{P}(m+n-1)
\]
defined for integers $m$, $n$ and $i$ such that $n\geq0$ and $1\leq i\leq m$,
satisfying the following axioms:
\begin{itemize}
\item[i)] \textbf{Associativity:} For the iterated compositions of 
$\mathcal{P}(m)\wedge\mathcal{P}(n)\wedge\mathcal{P}(p)$, the following associativity holds:
\[
\circ_i(\circ_j\wedge\id)=\begin{cases}
\circ_{j+p-1}(\circ_i\wedge\id)(\id\wedge\pi)&\text{for $1\leq i < j$,}\\
\circ_j(\id\wedge\circ_{i-j+1})&\text{for $j\leq i <j+n$, and}\\
\circ_j(\circ_{i-n+1}\wedge\id)(\id\wedge\pi)&\text{for $j+n\leq i$.}
\end{cases}
\]
Here $\pi:\mathcal{P}(n)\wedge\mathcal{P}(p)\rightarrow\mathcal{P}(p)\wedge\mathcal{P}(n)$
is the symmetry transposition for $\wedge$.
\item[ii)] \textbf{Equivariance:} Since $\mathcal{P}$ is a functor from $\Sigma$, each $\mathcal{P}(m)$
has an action of $\Sigma_m$. We write this action on the right, and
for $\rho\in\Sigma_m$ and $\upsilon\in\Sigma_n$ the following diagram commutes:
\[
\begin{CD}
\mathcal{P}(m)\wedge \mathcal{P}(n) @>{\circ_{\rho(i)}}>> \mathcal{P}(m+n-1)\\
@V{(-.\rho)\wedge(-.\upsilon)}VV @VV{-.(\rho\circ_i\upsilon)}V\\
\mathcal{P}(m)\wedge \mathcal{P}(n) @>{\circ_{i}}>> \mathcal{P}(m+n-1)
\end{CD}\quad.
\]
\item[iii)] \textbf{Unity:}
There is a map $1:\mathcal{P}(0)\rightarrow \mathcal{P}(1)$ 
such that the following diagrams commute
for all $1\leq i\leq m$
\[
\begin{CD}
\mathcal{P}(m)\wedge \mathcal{P}(0) @>{=}>> \mathcal{P}(m)\wedge S\\
@V{\id\wedge 1}VV @VV{\cong}V\\
\mathcal{P}(m)\wedge\mathcal{P}(1) @>{\circ_i}>> \mathcal{P}(m)
\end{CD}
\quad\text{and}\quad
\begin{CD}
S\wedge \mathcal{P}(m)@>{\cong}>> \mathcal{P}(m)\\
@V{1\wedge \id}VV @VV{=}V\\
\mathcal{P}(1)\wedge\mathcal{P}(m) @>{\circ_1}>> \mathcal{P}(m)
\end{CD}
\]
\end{itemize}
\end{Def}

The following types of operads are relevant for our applications:

\begin{Def}
\begin{itemize}
\item[] We get \textit{discrete operads} by putting 
the symmetric monoidal category of sets with cross product and unit $*$, into the definition above.
\item[] We get \textit{operads in topological spaces} by putting the symmetric monoidal category of spaces with cross product and unit $*$,
into the definition above.
\item[] We get \textit{operads in orthogonal spectra} by putting $(\mathscr{IS},\wedge,S)$, the symmetric monoidal category
of orthogonal spectra, into the definition above.
\end{itemize}
\end{Def}

\begin{Rem}\label{Rem:multioperation}
The definition above is equivalent to May's ``multi-operation'' definition
of an operad. We construct the multi-operation 
\[
\gamma:\mathcal{P}(k)\wedge\mathcal{P}(j_1)\wedge\cdots\wedge\mathcal{P}(j_k)\rightarrow\mathcal{P}(j)
\quad\text{, where $j=j_1+\cdots+j_k$,}
\]
as the composition 
\[
\gamma=\circ_{j_{k-1}+\cdots+j_1+1}(\circ_{j_{k-2}+\cdots+j_1+1}\wedge\id)
\cdots(\circ_{j_1+1}\wedge\id\wedge\cdots\wedge\id)
(\circ_1\wedge\id\wedge\cdots\wedge\id)\quad.
\]
To go the other way one uses the unit and defines $\circ_i=\gamma(-;1,\ldots,1,-,1,\ldots,1)$.
\end{Rem}

\begin{Prop}\label{Prop:FPisoperad}
If $F$ is a lax symmetric monoidal functor and $\mathcal{P}$ an operad, then $F\mathcal{P}$ is also an operad.
\end{Prop}

\begin{proof}
$F\mathcal{P}$ is clearly a functor defined on the symmetric groupoid, but we redefine $F\mathcal{P}(0)$ to be $S$.
We define the composition operations for
$F\mathcal{P}$ as the maps
\[
F\mathcal{P}(m)\wedge F\mathcal{P}(n)\rightarrow F(\mathcal{P}(m)\wedge \mathcal{P}(n))\xrightarrow{F(\circ_i)} 
F\mathcal{P}(m+n-1)\quad.
\]
Proving associativity for $F\mathcal{P}$ uses associativity and symmetry for $F$, and equivariance holds since the
multiplication for $F$ is a natural transformation. The map $1:F\mathcal{P}(0)\rightarrow F\mathcal{P}(1)$
is defined as the composition
\[
F\mathcal{P}(0)=S\rightarrow FS\xrightarrow{F(1)}F\mathcal{P}(1)\quad.
\]
And unity for $F\mathcal{P}$ follows from the unity of $F$ and $\mathcal{P}$.
\end{proof}

\begin{Exa}
The functor embedding sets in $\Top$ as the discrete spaces, is lax symmetric monoidal, hence we may
consider every discrete operad as an operad in topological spaces.
\end{Exa}

\begin{Exa}[Suspension operads]\label{Exa:suspoperads}
The functor sending a space $X$ to the orthogonal spectrum 
\[
F_0(X_+)\quad\text{given at level $V$ by }F_0(X_+)(V)=X_+\wedge S^V
\]
is symmetric monoidal. Here $F_0$ is the shift desuspension functor.
If $\mathcal{C}$ is an operad in topological spaces, then we may
construct a suspension operad in orthogonal spectra by 
composing with $X\mapsto F_0(X_+)$.
Usually we will denote this operad in orthogonal spectra simply by $\mathcal{C}$, instead of
$F_0(\mathcal{C}_+)$. And we will call an operad in orthogonal spectra discrete if
it is the suspension of a discrete operad in topological spaces. 
\end{Exa}

Let us give names to two operads in orthogonal spectra.
\begin{itemize}
\item[] We define $\mathcal{N}(m)(V)=S^V$ for all $m$ and $V$. The operations $\circ_i$
are in this case just the canonical map $S\wedge S\rightarrow S$.
\item[] Define $\mathcal{M}$ to be the suspension of the discrete operad
$m\mapsto \Sigma_m$. The operations $\circ_i$ are then equal to the suspension of the operations
$\Sigma_m\times\Sigma_n\rightarrow \Sigma_{m+n-1}$ defined above.
\end{itemize}

We will now define an operad geared toward anti-commutative involutions.
The hyperoctahedral group $\mathcal{H}(n)$ 
is the group of rigid symmetries of a cube in $\R^n$. 
If we let the cube be $[-1,1]^n$, we can identify
$\mathcal{H}(n)$ as a closed subgroup of $O(n)$.
The matrices in the image are those such that every row and every column
have $0$'s in all except one entry, and this entry is $1$ or $-1$.

We use the matrix description to describe the composition operations.
First let $T_n$ be the $n\times n$-matrix
\[
T_n=\begin{pmatrix}
0 & \cdots & 0 & -1 \\
0 & \cdots & -1 & 0 \\
\vdots & \mbox{}_{\cdotp}\cdotp^{\cdotp} & \vdots & \vdots\\
-1 & \cdots & 0 & 0
\end{pmatrix}\quad.
\]
Now define $\circ_i:\mathcal{H}(m)\times \mathcal{H}(n)\rightarrow\mathcal{H}(m+n-1)$ as follows.
Let $A$ be an $m\times m$-matrix describing an element in $\mathcal{H}(m)$, and let $B$ be a matrix in $\mathcal{H}(n)$.
There is a block decomposition of $A$ as
\[
\begin{pmatrix}
A_{11} & 0 & A_{13}\\
0 & a & 0\\
A_{31} & 0 & A_{33}
\end{pmatrix}\quad,
\]
where $a$ is $1$ or $-1$ and lies in the $i$'th column of $A$. Now define the $(m+n-1)\times (m+n-1)$-matrix $C$ to be
\[
\begin{pmatrix}
A_{11} & 0 & A_{13}\\
0 & B & 0\\
A_{31} & 0 & A_{33}
\end{pmatrix}\quad\text{if $a=1$,}
\]
and
\[
\begin{pmatrix}
A_{11} & 0 & A_{13}\\
0 & T_n B & 0\\
A_{31} & 0 & A_{33}
\end{pmatrix}\quad\text{if $a=-1$.}
\]
Then $C$ is the matrix of $A\circ_i B$ in $\mathcal{H}(m+n-1)$.

\begin{Prop}
$\mathcal{H}$ is a discrete operad.
\end{Prop}

\begin{proof}
The right action of $\sigma\in\Sigma_n$ on $\mathcal{H}(n)$ is given by
multiplication of matrices. If $A$ is a matrix in $\mathcal{H}$ and $B$ is the
permutation matrix of $\sigma$, then $\sigma$ sends $A$ to $AB$.
And the unit in $\mathcal{H}(1)$ is the identity matrix in $O(1)$.

Using the ``matrix''-model it is easy to verify all three axioms.
\end{proof}

\begin{Rem}\label{Rem:AltHdescr}
In this remark we describe $\mathcal{H}(n)$ as $(\Z/2)^n\rtimes \Sigma_n$,
and give a formula for $\circ_i$.

Write $(\Z/2)$ multiplicatively, denote 
elements of $(\Z/2)^n$ by $\mb{x}=(x_1,x_2,\ldots,x_n)$,
and embed $(\Z/2)^n$ in $O(n)$ as the $n\times n$-matrices
having $1$ or $-1$ on the diagonal and $0$ elsewhere.
Identifying $\Sigma_n$ with the permutation matrices we see that $\mathcal{H}(n)$
is actually the product of $(\Z/2)^n$ and $\Sigma^n$ inside $O(n)$.
$(\Z/2)^n$ is normal in $\mathcal{H}(n)$ and $(\Z/2)^n\cap\Sigma^n=\{I\}$,
thus $\mathcal{H}(n)$ is a semi-direct product $(\Z/2)^n\rtimes \Sigma_n$.
Here $\Sigma_n$ acts by permutation of factors on $(\Z/2)^n$.

It is possible to treat $n\mapsto (\Z/2)^n$ as a non-equivariant operad.
The composition operations are given by
\small
\[
(x_1,\ldots,x_m)\circ_i(y_1,\ldots,y_n)=\begin{cases}
(x_1,\ldots,x_{i-1},y_1,\ldots,y_n,x_{i+1},\ldots,x_{m}) &\text{ if $x_i=1$, and}\\
(x_1,\ldots,x_{i-1},-y_n,\ldots,-y_1,x_{i+1},\ldots,x_{m}) &\text{ if $x_i=-1$.}
\end{cases}
\]
\normalsize
We now introduce the following convention: $\tau_n$ without an argument denotes the
order reversing permutation in $\Sigma_n$, while $\tau_n$ with an argument 
denotes the group homomorphism $\Z/2\rightarrow\Sigma_n$ sending $-1$ to the order reversing permutation.
Hence $\tau_n(-1)=\tau_n$, while $\tau_n(1)=\id$.
Now we can give a formula for the composition operations of $\mathcal{H}$ in terms
of the $\circ$'s of $n\mapsto (\Z/2)^n$ and $(n\mapsto \Sigma_n)=\mathcal{M}$.
Inspecting the ``matrix''-model we get
\[
(\mb{x},\rho)\circ_i(\mb{y},\upsilon)=(\mb{x}\circ_{\rho(i)} \mb{y},\rho\circ_i (\tau_n({x_{\rho(i)}})\upsilon))
\]
for $\mb{x}\in(\Z/2)^m$, $\rho\in\Sigma_m$, $\mb{y}\in(\Z/2)^n$ and $\upsilon\in\Sigma_n$.
\end{Rem}

\begin{Exa}\label{Exa:Hinos}
Via suspension the operad $\mathcal{H}$ induces an operad in orthogonal spectra. This induced operad 
will also be called $\mathcal{H}$.
\end{Exa}

Here is another example of a lax symmetric monoidal functor applied to operads in orthogonal spectra.

\begin{Exa}\label{Exa:GammaM}
Recall the cofibrant replacement functor $\tilde{\Gamma}$ from theorem~\ref{Thm:orbitcofrepl}. We have shown that
$\tilde{\Gamma}$ is lax symmetric monoidal. Hence, for any operad $\mathcal{P}$ in orthogonal spectra we get
an operad $\tilde{\Gamma}\mathcal{P}$. The natural map $\tilde{\Gamma} L\rightarrow L$ induces a map of operads
$\tilde{\Gamma}\mathcal{P}\rightarrow\mathcal{P}$. The nice thing about this new operad is that each 
$\tilde{\Gamma}\mathcal{P}(m)$ is $\Sigma_m$-equivariantly naive orbit cofibrant.

However, it is hard to say anything about the genuine homotopy of $\tilde{\Gamma}\mathcal{P}$, even in the case
where $\mathcal{P}$ is $\Sigma$-free.
\end{Exa}

The idea behind an operad is that $\mathcal{P}(n)$ can parametrize $n$-fold multiplications
on an object $L$. If we have such a parametrization, we call $L$ a $\mathcal{P}$-algebra.
The precise definition is:

\begin{Def}\label{Def:Palg}
Let $\mathcal{P}$ be an operad in the symmetric monoidal category $(\mathscr{C},\wedge,S)$.
A $\mathcal{P}$-algebra is an object $L$ in $\mathscr{C}$ together with operations 
\[
\theta_m:\mathcal{P}(m)\wedge L^{\wedge m}\rightarrow L
\]
such that the following axioms holds:
\begin{itemize}
\item[i)] \textbf{$\theta$ acts:} For all $n\geq0$ and $1\leq i\leq m$ the following
diagram commutes:
\[
\begin{CD}
\mathcal{P}(m)\wedge\mathcal{P}(n)\wedge L^{\wedge (m+n-1)} @= \mathcal{P}(m)\wedge\mathcal{P}(n)\wedge L^{\wedge (m+n-1)}\\
@V{\text{shuffle}}VV @VV{\circ_i}V\\
\mathcal{P}(m)\wedge L^{\wedge (i-1)}\wedge\mathcal{P}(n)\wedge L^{\wedge n}\wedge L^{\wedge (m-i)}
&& \mathcal{P}(m+n-1)\wedge L^{\wedge (m+n-1)} \\
@V{\id\wedge\theta_n\wedge \id}VV @VV{\theta_{m+n-1}}V\\ 
\mathcal{P}(m)\wedge L^{\wedge m}
@>{\theta_m}>> L
\end{CD}\quad.
\]
\item[ii)] \textbf{Triviality of the unit:} The diagram
\[
\begin{CD}
\mathcal{P}(0)\wedge L @= S\wedge L\\
@V{1\wedge\id}VV @VV{\cong}V\\
\mathcal{P}(1)\wedge L @>{\theta_1}>> L
\end{CD}
\]
commutes.
\item[iii)] \textbf{Equivariance:} The group $\Sigma_m$ acts from the left on $\mathcal{P}(m)$ and
acts from the right on $L^{\wedge m}$ by permutation of the factors. $\theta_m$ is equivariant in the sense
that the diagram
\[
\begin{CD}
\mathcal{P}(m)\wedge L^{\wedge m} @>{\rho\wedge \id}>> \mathcal{P}(m)\wedge L^{\wedge m}\\
@V{\id\wedge\rho}VV @VV{\theta_m}V\\
\mathcal{P}(m)\wedge L^{\wedge m} @>{\theta_m}>> L
\end{CD}
\]
commutes for every $\rho\in\Sigma_m$.
\end{itemize}
\end{Def}

\begin{Rem}
$\theta_0$ is a map $\mathcal{P}(0)=S\rightarrow L$,
and we call this map the unit of $L$.
\end{Rem}

Let us look at an example:

\begin{Exa}\label{Exa:Halginsets}
In this example we consider the discrete operad $\mathcal{H}$ in the category of sets.
When considering the category of sets, remember that ``$\wedge$'' in the definition above
is the cross product, and $S$ is the set $\{1\}$.
We now want to recognize the class of $\mathcal{H}$-algebras as a more familiar type of mathematical objects.

Assume that $X$ is an $\mathcal{H}$-algebra.
Let $1$ in $X$ denote the image of $\theta_0:\{1\}\rightarrow X$.
Now define $\mu:X\times X\rightarrow X$ and $\iota:X\rightarrow X$ by
\begin{align*}
\mu(x,y)&=\theta_2\left(\begin{pmatrix} 1& 0 \\ 0& 1\end{pmatrix},x,y\right)\quad,\text{ and}\\
\iota(x)&=\theta_1\left(\begin{pmatrix} -1\end{pmatrix},x\right)\quad.
\end{align*}
Let us now do some calculations. First we have that
\begin{align*}
\mu(\mu(x,y),z)&=\theta_2\left(\begin{pmatrix} 1& 0 \\ 0& 1\end{pmatrix},
\theta_2\left(\begin{pmatrix} 1& 0 \\ 0& 1\end{pmatrix},x,y\right),z\right)\\
&= \theta_3\left( \begin{pmatrix} 1& 0 \\ 0& 1\end{pmatrix}\circ_1 \begin{pmatrix} 1& 0 \\ 0& 1\end{pmatrix},x,y,z\right)\\
&= \theta_3\left( \begin{pmatrix} 1& 0 & 0 \\ 0& 1 & 0\\ 0& 0& 1\end{pmatrix},x,y,z\right)\quad.
\end{align*}
Similarly, we can show that $\mu(x,\mu(y,z))=\theta_3\left( \begin{pmatrix} 1& 0 & 0 \\ 0& 1 & 0\\ 0& 0& 1\end{pmatrix},x,y,z\right)$.
Hence, $\mu$ is an associative operation on $X$.
Furthermore, we have:
\begin{align*}
\mu(1,x)&=\theta_2\left(\begin{pmatrix} 1& 0 \\ 0& 1\end{pmatrix},\theta_0(1),x\right)\\
&=\theta_1\left( \begin{pmatrix} 1& 0 \\ 0& 1\end{pmatrix}\circ_1\begin{pmatrix} \mbox{} \end{pmatrix},x\right)\\
&=\theta_1\left( 1, x\right)\\
&= x\quad,
\end{align*}
and by the same methods one also calculates that $\mu(x,1)=x$. This shows that $1$ is a two-sided unit for $\mu$.
Hence, $X$ is a monoid with unit. Let us now look at the operation $\iota$. We have:
\begin{align*}
\iota(\iota(x))&= \theta_1\left(\begin{pmatrix} -1\end{pmatrix},\theta_1\left(\begin{pmatrix} -1\end{pmatrix},x\right)\right)\\
&=\theta_1\left( \begin{pmatrix} -1\end{pmatrix}\circ_1 \begin{pmatrix} -1\end{pmatrix},x\right)\\
&=\theta_1\left( 1, x\right)\\
&= x\quad.
\end{align*}
The interaction of $\mu$ and $\iota$ can be computed as follows: First we have
\begin{align*}
\iota(\mu(x,y))&= \theta_1\left(\begin{pmatrix} -1\end{pmatrix},
\theta_2\left(\begin{pmatrix} 1 & 0\\ 0& 1 \end{pmatrix},x,y\right)\right)\\
&=\theta_2\left( \begin{pmatrix} -1\end{pmatrix}\circ_1 \begin{pmatrix} 1 & 0\\ 0& 1\end{pmatrix},x,y\right)\\
&=\theta_2\left( \begin{pmatrix} 0 & -1\\ -1 & 0\end{pmatrix},x,y\right)\quad,
\end{align*}
and secondly we calculate that
\begin{align*}
\mu(\iota(y),\iota(x))&= \theta_2\left(\begin{pmatrix} 1 & 0\\ 0& 1 \end{pmatrix}, 
\theta_1\left(\begin{pmatrix} -1\end{pmatrix},y\right),
\theta_1\left(\begin{pmatrix} -1\end{pmatrix},x\right)\right)\\
&= \theta_2\left(\begin{pmatrix} 1 & 0\\ 0& 1 \end{pmatrix}\circ_1 \begin{pmatrix} -1\end{pmatrix},y,
\theta_1\left(\begin{pmatrix} -1\end{pmatrix},x\right)\right)\\
&= \theta_2\left(\begin{pmatrix} -1 & 0\\ 0& 1 \end{pmatrix},y,
\theta_1\left(\begin{pmatrix} -1\end{pmatrix},x\right)\right)\\
&= \theta_2\left(\begin{pmatrix} -1 & 0\\ 0& 1 \end{pmatrix}\circ_2 \begin{pmatrix} -1\end{pmatrix},y,x\right)\\
&= \theta_2\left(\begin{pmatrix} -1 & 0\\ 0& -1 \end{pmatrix},y,x\right)\\
&= \theta_2\left(\begin{pmatrix} 0 & -1\\ -1& 0 \end{pmatrix},x,y\right)\quad.
\end{align*}
In the last step we used the equivariance axiom for $\mathcal{H}$-algebras. What we have seen is that
$\iota^2(x)=x$ and $\iota(\mu(x,y))= \mu(\iota(y),\iota(x))$. We say that $\iota$ is an involution on the monoid
$X$ which anti-commutes with the multiplication.

It can be shown that there are no more relations for a general $\mathcal{H}$-algebra. Hence, we recognize
$X$ as a monoid with unit and anti-commutative involution. 
\end{Exa}

\begin{Prop}
If $F$ is a lax symmetric monoidal functor and $L$ a $\mathcal{P}$-algebra, then $FL$ is an $F\mathcal{P}$-algebra.
\end{Prop}

\begin{proof}
We define the operations $\theta_m'$ for $FL$. For $m\geq 1$ we define $\theta_m'$ as the composition
\[
F\mathcal{P}(m)\wedge (FL)^{\wedge m}\rightarrow F(\mathcal{P}(m)\wedge L^{\wedge m})\xrightarrow{F\theta_m}FL\quad.
\]
For $m=0$ we use the unit of $F$ to define $\theta_0'$:
\[
F\mathcal{P}(0)=S\rightarrow FS\xrightarrow{F\theta_0} FL\quad.
\]
It is an exercise to check the $F\mathcal{P}$-algebra axioms for $FL$. Notice in particular that we need symmetry of $F$
to prove both associativity and equivariance.
\end{proof}

We conclude this section by the following important observation:

\begin{Prop}\label{Prop:interpretoperadalgebras}
There are 1-1 correspondences between
\begin{itemize}
\item[] $\mathcal{M}$-algebras in $\mathscr{IS}$ and orthogonal ring spectra,
\item[] $\mathcal{N}$-algebras in $\mathscr{IS}$ and commutative orthogonal ring spectra, and
\item[] $\mathcal{H}$-algebras in $\mathscr{IS}$ and orthogonal ring spectra with involution.
\end{itemize}
\end{Prop}

\begin{proof}
When $\mathcal{P}$ is the suspension of a discrete operad, we may identify $\mathcal{P}(n)\wedge L^{\wedge n}$
with a wedge sum 
\[
\bigvee L^{\wedge n}
\]
indexed over the non-base points in $\mathcal{P}(n)(0)$.

Given an $\mathcal{M}$-algebra $L$, the map $\theta_0:S\rightarrow L$ is the unit. To get the multiplication
we restrict $\theta_2$ to the wedge summand
corresponding to $\id\in\Sigma_2$. The map on the other summand corresponds to $\mu\circ\pi$, where $\pi$
exchanges the factors of $L\wedge L$. Associativity comes from comparing $\mu\circ(\id\wedge\mu)$ and
$\mu\circ(\mu\wedge\id)$ to $\theta_3$ restricted the summand corresponding to $\id\in\Sigma_3$.

Conversely, given an orthogonal ring spectrum $L$, it becomes an $\mathcal{M}$-algebra by
defining $\theta_n$ on the wedge summand corresponding to $\id\in\Sigma_n$ to be
the multiplication map $L^{\wedge n}\rightarrow L$, and extend to the other summands by equivariance.

If $L$ is an $\mathcal{N}$-algebra, we get that $\mathcal{N}(2)\wedge L\wedge L=L\wedge L$,
and set $\mu=\theta_2$. Commutativity follows from the $\Sigma_2$-equivariance of $\theta_2$.

Conversely, given a commutative orthogonal ring spectrum $L$, it becomes an $\mathcal{N}$-algebra
by defining $\theta_n$ on $\mathcal{N}(n)\wedge L^{\wedge n}=L^{\wedge n}$ by multiplication.
It is well defined because of commutativity.

If $L$ is an $\mathcal{H}$-algebra, we let multiplication be $\theta_2$ restricted 
to the summand determined by the matrix $\begin{pmatrix}1 & 0 \\ 0 & 1 \end{pmatrix}$,
and the involution $\iota:L\rightarrow L$ is defined to be $\theta_1$ restricted 
to the summand corresponding to the matrix $\begin{pmatrix} -1 \end{pmatrix}$.
The calculation that 
$\begin{pmatrix} -1 \end{pmatrix}\circ_1\begin{pmatrix} -1 \end{pmatrix}=\begin{pmatrix} 1 \end{pmatrix}$ 
implies that $\iota^2=\id$,
and
\[
\left( \begin{pmatrix} 1 & 0\\ 0& 1 \end{pmatrix}\circ_1 
\begin{pmatrix} -1 \end{pmatrix}\right)\circ_2 \begin{pmatrix} -1 \end{pmatrix}
=\begin{pmatrix} -1 & 0\\ 0& -1 \end{pmatrix}
=\begin{pmatrix} -1 \end{pmatrix}\circ_1 \begin{pmatrix} 0 & 1\\ 1& 0 \end{pmatrix}
\]
implies that $\iota$ is an anti-homomorphism. This is analogous to the calculation in example~\ref{Exa:Halginsets}.

Conversely, given an orthogonal ring spectrum $L$ with involution, it becomes an $\mathcal{H}$-algebra
as follows: Define $\theta_n$ on the wedge summand corresponding to the matrix
\[
\begin{pmatrix} x_1 & 0 & \cdots & 0\\ 
0 & x_2 & \cdots & 0\\
\vdots &\vdots & \ddots & \vdots \\
0 & 0 & \cdots & x_n\end{pmatrix}
\]
as the composition
\[
L^{\wedge n}= L\wedge L\wedge \cdots\wedge L\xrightarrow{\iota({x_1})\wedge \iota({x_2})\wedge\cdots\wedge\iota({x_n})}
L\wedge L\wedge \cdots\wedge L\xrightarrow{\text{multiplication}}L\quad.
\]
Here $\iota(x)$ denotes the involution $\iota$ if $x=-1$, while $\iota(1)=\id$, the identity of $L$.
We extend to all of $\mathcal{H}(n)$ by $\Sigma_n$-equivariance.
\end{proof}

\section{Operads and the two sided bar construction}\label{sect:operadsbarconstr}

\subsection{Operads and monads}

We now follow the theory as presented in May's book~\cite{May:72}.
Our goal is to check that the basic results also hold for 
operads in a symmetric monoidal category $(\mathscr{C},\wedge,S)$ provided that
$\mathscr{C}$ has all small colimits. In particular the theory of this subsection applies to 
orthogonal spectra.

Assume that $\mathcal{P}$ is an operad in $\mathscr{C}$ and that $L$ is an
object in $\mathscr{C}$ which comes with a chosen map $S\rightarrow L$. 
Let $I$ denote B\"okstedts category. 
This category has objects the finite sets $\mathbf{n}$
for $n\geq 0$. The morphisms are the injective functions. 
Notice that $\Sigma$ is a subcategory of $I$.
For any morphism $\rho:\mathbf{n'}\rightarrow\mathbf{n}$ in $I$
we have maps
\[
\id\wedge \rho_*:\mathcal{P}(n)\wedge L^{\wedge n'}
\rightarrow \mathcal{P}(n)\wedge L^{\wedge n}
\]
and
\[
\rho^*\wedge\id:\mathcal{P}(n)\wedge L^{\wedge n'}
\rightarrow \mathcal{P}(n')\wedge L^{\wedge n'}\quad.
\]
The first map comes from shuffling factors according to
$\rho$, and inserting $S\rightarrow L$
for those factors in $L^{\wedge n}$ corresponding to 
points in $\mathbf{n}$, not in the image of $\rho$.
The second map comes from the identification of $S$ 
with $\mathcal{P}(0)$, $n-n'$ times, and then using
the appropriate composition operations to reduce from 
$\mathcal{P}(n)\wedge\mathcal{P}(0)^{\wedge n-n'}$
to $\mathcal{P}(n')$.

\begin{Def}
We define $PL$ to be the coequalizer of
\[
\bigvee_{\rho:\mathbf{n'}\rightarrow\mathbf{n}} \mathcal{P}(n)\wedge L^{\wedge n'}
\rightrightarrows \bigvee_{\mathbf{n}}\mathcal{P}(n)\wedge L^{\wedge n}\quad.
\]
\end{Def}

\begin{Rem}
Observe that $PL$ also can be described as the coend
\[
\int^{\mb{n}\in I}\mathcal{P}(n)\wedge L^{\wedge n}\quad.
\]
\end{Rem}

By $S\downarrow \mathscr{C}$ we mean the category of objects in $\mathscr{C}$ under $S$.
We want to show that $P$ is a monad in this category. Recall that a monad $M$ in a 
category $\mathscr{C}$ consists of a functor $M:\mathscr{C}\rightarrow\mathscr{C}$
together with natural transformations $\mu:M^2\rightarrow M$ and $\eta:\id\rightarrow M$
such that $\eta$ is a left and right unit for $\mu$, and $\mu$ is associative.

\begin{Prop}
For any operad $\mathcal{P}$ in $\mathscr{C}$,
$P$ is a monad in $S\downarrow\mathscr{C}$.
\end{Prop}

\begin{proof}
Using the unit of $\mathcal{P}$ we have a map 
\[
L\cong S\wedge L\rightarrow \mathcal{P}(1)\wedge L
\subset \bigvee_{\mathbf{m}}\mathcal{P}(m)\wedge L^{\wedge m}\rightarrow PL\quad.
\]
This is the natural transformation $\eta$. Since $L$ is under $S$, $PL$ is also
an object in $S\downarrow\mathscr{C}$ via the composition
$S\rightarrow L\xrightarrow{\eta}PL$.

To construct the multiplication $\mu:PPL\rightarrow PL$ we will use the 
composition operations of $\mathcal{P}$. Recall the definition of May's
multioperation $\gamma:\mathcal{P}(m)\wedge \mathcal{P}(n_1)\wedge \cdots \wedge\mathcal{P}(n_m)\rightarrow\mathcal{P}(n_1+\cdots+n_m)$ 
as the composition of $\circ$'s, see remark~\ref{Rem:multioperation}.
We now define $\tilde{\mu}$ as the composition
\begin{multline*}
\mathcal{P}(m)\wedge (\mathcal{P}(n_1)\wedge L^{\wedge n'_1})\wedge\cdots 
\wedge (\mathcal{P}(n_m)\wedge L^{\wedge n'_m})\\
\xrightarrow{\text{shuffle}}\mathcal{P}(m)\wedge \mathcal{P}(n_1)\wedge \cdots \wedge\mathcal{P}(n_m)
\wedge L^{\wedge n'_1}\wedge\cdots \wedge L^{\wedge n'_m}\\
\xrightarrow{\gamma\wedge\id}\mathcal{P}(n_1+n_{2}+\cdots+n_m)\wedge L^{\wedge (n'_1+n'_{2}+\cdots+ n'_m)}\quad.
\end{multline*}
Notice that $\tilde{\mu}$ is natural for $\mb{n_j}\in I^{\op}$ and for $\mb{n'_j}\in I$. Thus we have an induced map
\[
\mathcal{P}(m)\wedge (PL)^{\wedge m}\rightarrow PL\quad.
\]
Given $\rho:\mb{m'}\rightarrow\mb{m}$, it is an exercise to check that the diagram
\[
\begin{CD}
\mathcal{P}(m)\wedge (PL)^{\wedge m'} @>{\id\wedge\rho_*}>> \mathcal{P}(m)\wedge (PL)^{\wedge m}\\
@V{\rho^*\wedge\id}VV @VVV\\
\mathcal{P}(m')\wedge (PL)^{\wedge m'} @>>> PL
\end{CD}
\]
commutes. And we get our monad multiplication
\[
\mu:PPL\rightarrow PL\quad.
\]
Clearly $\eta$ is a left and right unit. Associativity of $\mu$ follows from associativity rules for the composition operations
$\circ_j$.
\end{proof}

If $M$ is a monad in a category $\mathscr{C}$, then we recall that an $M$-algebra is an object $L$ of $\mathscr{C}$
together with a map $\theta:ML\rightarrow L$, such that $\theta\mu=\theta M(\theta)$ and $\theta\eta_L=\id_L$.
See chapter~VI in~\cite{MacLane:98}.

\begin{Prop}
There is a natural one-to-one correspondence between $\mathcal{P}$-algebras and
$P$-algebras.
\end{Prop}

\begin{proof}
Given a $P$-algebra $L$, we define the $\mathcal{P}$-algebra maps $\theta_m$ as
the compositions
\[
\mathcal{P}(m)\wedge L^{\wedge m}\rightarrow PL
\xrightarrow{\theta}L\quad.
\]

Conversely, if $L$ is a $\mathcal{P}$-algebra, we check that the following diagram commutes
\[
\begin{CD}
\mathcal{P}(m)\wedge L^{\wedge m'} @>{\id\wedge\rho_*}>> \mathcal{P}(m)\wedge L^{\wedge m}\\
@V{\rho^*\wedge\id}VV @VV{\theta_m}V\\
\mathcal{P}(m')\wedge L^{\wedge m'} @>{\theta_{m'}}>> PL
\end{CD}
\]
for all $\rho:\mb{m'}\rightarrow\mb{m}$ in $I$. Thus we have an induced map $PL\rightarrow L$,
and we take this as a definition of $\theta$.
\end{proof}

\begin{Cor}\label{Lem:PLisPalg}
$PL$ is a $\mathcal{P}$-algebra, and for any map $f:K\rightarrow  L$ in $\mathscr{C}$,
the induced map $Pf:PK\rightarrow PL$ is a map of $\mathcal{P}$-algebras.
\end{Cor}

\begin{proof}
The multiplication $\mu:PPL\rightarrow PL$ gives $PL$ a $\mathcal{P}$-algebra structure.
And naturality implies that $Pf:PK\rightarrow PL$ is a $\mathcal{P}$-algebra morphism.
\end{proof}

Let $M$ be a monad in $\mathscr{C}$. Recall from May~\cite{May:72}
that an $M$-functor is a functor $F$ with the same source as $M$ together with a natural
transformation $\lambda:FM\rightarrow F$, such that $\lambda F\eta$ is the identity and
$\lambda F\mu=\lambda\lambda$.

\begin{Prop}
If $\alpha:\mathcal{P}\rightarrow\mathcal{Q}$ is a map of operads, then $Q$ is a $P$-functor.
\end{Prop}

\begin{proof}
By functorality of the construction of $P$ form $\mathcal{P}$ it is clear that
$\alpha$ induces a morphism of monads $P\rightarrow Q$. 

Now let $\alpha_L$ denote the natural transformation $PL\rightarrow QL$, and
define $\lambda$ to be the composition
\[
QPL\xrightarrow{Q\alpha_{L}} QQL\xrightarrow{\mu'} QL\quad.
\] 
Here $\mu'$ is the multiplication for $Q$. The properties
of a $P$-functor are easily verified.
\end{proof}

\subsection{Homotopy theory of operads and their algebras}\label{subsect:htyptheoryofoperads}

Having treated the categorical theory of operads and monads we now turn toward
homotopy theory. 
Berger and Moerdijk,~\cite{BergerMoerdijk:03}, define model structure on operads in monoidal model categories. 
Their approach requires a symmetric monoidal fibrant replacement functor. See their theorem~3.1. 
We are interested in orthogonal spectra, but
this category does not possess such a functor, see our remark~\ref{Rem:ISnotconvenient}. Also
see example~4.6.4 in~\cite{BergerMoerdijk:03}.
However, we do not need a model structure on operads in orthogonal spectra. Direct methods are sufficient.

Suppose given notions of cofibration and weak equivalence for orthogonal spectra.
Let $\alpha:\mathcal{P}\rightarrow \mathcal{Q}$ be a map of operads in orthogonal spectra, 
and $f:K\rightarrow L$ a map of orthogonal spectra under $S$. We ask:
\begin{itemize}
\item[] When is $\eta: L\rightarrow PL$ a cofibration?
\item[] When is $Pf:PK\rightarrow PL$ a cofibration?
\item[] If each $\alpha:\mathcal{P}(j)\rightarrow \mathcal{Q}(j)$ is a weak equivalence
when is also $PL\rightarrow QL$ a weak equivalence? 
\item[] If $f$ is a weak equivalence, when is $Pf:PK\rightarrow PL$ also a weak equivalence? 
\end{itemize}

\begin{Rem}\label{Rem:alloperadstoodifficult}
The author originally wanted to address these questions for arbitrary operads $\mathcal{P}$.
In this setting, the functor $\tilde{\Gamma}$ of theorem~\ref{Thm:orbitcofrepl} should be
the cofibrant replacement functor for operads. Proposition~\ref{Prop:tildeGammaequivariantly}
would ensure that $\tilde{\Gamma}\mathcal{P}$ had $\Sigma$-equivariance.
However, as lemma~\ref{Lem:operadfiltration} below shows, one must be able to analyze
the smash product over $\Sigma_j$. But getting results about $X\wedge_{\Sigma_j} Y$
when $X$ and $Y$ are $\Sigma_j$-equivariantly orbit cofibrant, turned out to be too
difficult. Therefore the author had to impose a very restrictive condition on the
operads. Further details about this condition can be found below. Still the cases
of main interest, the operads $\mathcal{D}_n$, fits into this restrictive framework
because of theorem~\ref{Thm:ocreplforDn}.

A milder hypothesis that could work is to assume $\Sigma_j$-freeness of $\mathcal{P}(j)$.
But for our applications the stricter condition is sufficient.
\end{Rem}

We will answer the questions by giving
sufficient criteria in the propositions below. 
Central in all arguments
is a filtration of $PL$. It is given by defining $F_jPL$
to be the coequalizer of
\[
\bigvee_{
\begin{matrix}\scriptstyle 
\rho:\mathbf{n'}\rightarrow\mathbf{n}\\
\scriptstyle n\leq j
\end{matrix}}
\mathcal{P}(n)\wedge L^{\wedge n'}
\rightrightarrows \bigvee_{n=0}^j\mathcal{P}(n)\wedge L^{\wedge n}\quad.
\]
By $F_0PL$ we will understand the sphere spectrum $S$, and $S=F_0PL\rightarrow PL$
is the unit of $PL$. As coends we can write
\[
F_jPL=\int^{n\leq j}\mathcal{P}(n)\wedge L^{\wedge n}\quad.
\]

\begin{Lem}\label{Lem:operadfiltration}
$F_1PL=\mathcal{P}(1)\wedge L$,
for $j\geq 2$ there are pushout squares
\[
\begin{CD}
\mathcal{P}(j)\wedge_{\Sigma_j}sL^{\wedge j-1} @>>> F_{j-1}PL\\
@VVV @VVV\\
\mathcal{P}(j)\wedge_{\Sigma_j}L^{\wedge j} @>>> F_{j}PL
\end{CD}
\]
and $\colim_j F_jPL=PL$.
\end{Lem}

Here $sL^{\wedge j-1}$ is an abbreviation for 
$(S\wedge L^{\wedge j-1})\cup (L\wedge S\wedge L^{\wedge j-2})\cup \cdots \cup(L^{\wedge j-1}\wedge S)$.
The proof that follows is categorical, so this lemma holds for any symmetric monoidal category which
has small colimits.

\begin{proof}
To see that $PL=\colim_j F_jPL$ we consider the diagram
\[
\begin{CD}
\cdots @>>> \bigvee_{
\begin{matrix}\scriptstyle 
\rho:\mathbf{n}\rightarrow\mathbf{m}\\
\scriptstyle m\leq j-1
\end{matrix}}
\mathcal{P}(m)\wedge L^{\wedge n} @>>>
\bigvee_{
\begin{matrix}\scriptstyle 
\rho:\mathbf{n}\rightarrow\mathbf{m}\\
\scriptstyle m\leq j
\end{matrix}}
\mathcal{P}(m)\wedge L^{\wedge n} @>>> \cdots\\
&& \downdownarrows && \downdownarrows \\
\cdots @>>> \bigvee_{m=0}^{j-1}\mathcal{P}(m)\wedge L^{\wedge m}
@>>> \bigvee_{m=0}^{j}\mathcal{P}(m)\wedge L^{\wedge m} @>>> \cdots
\end{CD}
\]
and use that taking colimits and taking coequalizers commute.

To see that the diagram is pushout we will use a trick involving coends. We now fix $j$.
Suppose that $\mb{m}$ and $\mb{n}'$ are objects in B\"okstedts category $I$, we then
have a pushout square of based sets
\small
\[
\begin{CD}
\{\theta:\mb{m}\rightarrow \mb{n}'\quad\text{, where $m<n'=j$}\}_+
&\rightarrow& \{\theta:\mb{m}\rightarrow\mb{n}'\quad\text{, where $m<j$ and $n'\leq j$}\}_+\\
@VVV @VVV\\
\{\theta:\mb{m}\rightarrow \mb{n}'\quad\text{, where $m\leq n'=j$}\}_+
&\rightarrow& \{\theta:\mb{m}\rightarrow\mb{n}'\quad\text{, where $m\leq j$ and $n'\leq j$}\}_+
\end{CD}\quad.
\]
\normalsize
Here the $\theta$'s are injective maps. We interpret the sets on the left to be $*$ if $n'\neq j$,
and all sets are $*$ if $n'>j$. The subscript $+$ means that we have added an extra basepoint.
Varying $\mb{m}\in I^{\op}$ and $\mb{n}'\in I$ we see that the collection of these sets is a 
functor $I^{\op}\times I\rightarrow \Ens_*$.

Now smash the diagram above with $\mathcal{P}(n)$ on the left and $L^{\wedge m'}$ on the right.
We get a pushout diagram
\small
\[
\begin{CD}
\mathcal{P}(n)\wedge\{\mb{m}\xrightarrow{\theta} \mb{n}'\;|\;m<n'=j\}_+\wedge L^{\wedge m'}
&\rightarrow& \mathcal{P}(n)\wedge\{\mb{m}\xrightarrow{\theta}\mb{n}'\;|\;m<j\text{, }n'\leq j\}_+\wedge L^{\wedge m'}\\
@VVV @VVV\\
\mathcal{P}(n)\wedge\{\mb{m}\xrightarrow{\theta} \mb{n}'\;|\;m\leq n'=j\}_+\wedge L^{\wedge m'}
&\rightarrow& \mathcal{P}(n)\wedge\{\mb{m}\xrightarrow{\theta}\mb{n}'\;|\;m\text{, }n'\leq j\}_+\wedge L^{\wedge m'}
\end{CD}
\]
\normalsize
of functors $I^{\op}\times I^{\op}\times I\times I\rightarrow\mathscr{IS}$.
Pushouts and coends commute, so applying the iterated coend to the diagram yields a pushout.
We calculate the corners of the resulting diagram by first taking the coend over $\mb{n}\in I$,
then over $\mb{m}\in I$.
We have:
\begin{align*}
\int^{\mb{m}}\int^{\mb{n}}
&\mathcal{P}(n)\wedge\{\theta:\mb{m}\rightarrow \mb{n}\;|\;m<n=j\}_+\wedge L^{\wedge m}\\
&= \int^{\mb{m}} 
\mathcal{P}(j)\wedge\{\theta:\mb{m}\rightarrow \mb{j}\;|\;m<j\}_+\wedge L^{\wedge m}\\
&= \mathcal{P}(j)\wedge_{\Sigma_j}sL^{\wedge j-1}\quad,\\ &\\
\int^{\mb{m}}\int^{\mb{n}}
&\mathcal{P}(n)\wedge\{\theta:\mb{m}\rightarrow \mb{n}\;|\;m\leq n=j\}_+\wedge L^{\wedge m}\\
&= \int^{\mb{m}} 
\mathcal{P}(j)\wedge\{\theta:\mb{m}\rightarrow \mb{j}\;|\;m\leq j\}_+\wedge L^{\wedge m}\\
&= \mathcal{P}(j)\wedge_{\Sigma_j}L^{\wedge j}\quad,\\ &\\
\int^{\mb{m}}\int^{\mb{n}}
&\mathcal{P}(n)\wedge\{\theta:\mb{m}\rightarrow\mb{n}\;|\;m<j\text{ and }n\leq j\}_+\wedge L^{\wedge m}\\
&= \int^{\mb{m}}\mathcal{P}(m)\wedge T(m<j)\wedge L^{\wedge m}\\
&= \int^{n<j}\mathcal{P}(n)\wedge L^{\wedge n}\\
&= F_{j-1}PL\quad\quad\quad\text{and}\\ &\\
\end{align*}
\begin{align*}
\int^{\mb{m}}\int^{\mb{n}}
&\mathcal{P}(n)\wedge\{\theta:\mb{m}\rightarrow\mb{n}\;|\;m\leq j\text{ and }n\leq j\}_+\wedge L^{\wedge m}\\
&= \int^{\mb{m}}\mathcal{P}(m)\wedge T(m\leq j)\wedge L^{\wedge m}\\
&= \int^{n\leq j}\mathcal{P}(n)\wedge L^{\wedge n}\\
&= F_{j}PL\quad.
\end{align*}
Here $T(m<j)$ is the functor sending $\mb{m}$ to $S^0$ if $m<j$ and to $*$ if $m\geq j$. And $T(m\leq j)$
is defined similarly. 

Thus the calculations show that the resulting diagram, which must be a pushout by construction, is equal to
\[
\begin{CD}
\mathcal{P}(j)\wedge_{\Sigma_j}sL^{\wedge j-1} @>>> F_{j-1}PL\\
@VVV @VVV\\
\mathcal{P}(j)\wedge_{\Sigma_j}L^{\wedge j} @>>> F_{j}PL
\end{CD}\quad.
\]
This concludes the proof.
\end{proof}

From now on we will only consider operads in orthogonal spectra. In order to prove the propositions
we will only work with those operads $\mathcal{P}$ such that each $\mathcal{P}(j)$ can be written
equivariantly as $X\wedge (\Sigma_j)_+$ for some (non-equivariant) orthogonal spectrum $X$.

\begin{Prop}\label{Prop:Pconserveqcof}
If $S\rightarrow\mathcal{P}(1)$ and
$S\rightarrow L$ are orbit q-cofibrations, and each $\mathcal{P}(j)$ can be written equivariantly
as a product $X\wedge (\Sigma_j)_+$ with $X$ being an orbit cofibrant (non-equivariant) orthogonal spectrum,
then the unit $\eta:L\rightarrow PL$
is an orbit q-cofibration.
\end{Prop}

\begin{proof}
First fix $j\geq 2$.
By proposition~\ref{Prop:orbitsquareproduct} $sL^{\wedge (j-1)}\rightarrow L^{\wedge j}$
is an orbit q-cofibration. Applying $\mathcal{P}(j)\wedge_{\Sigma_j}-$ we get
\[
\mathcal{P}(j)\wedge_{\Sigma_j}sL^{\wedge (j-1)}=X\wedge sL^{\wedge (j-1)}\rightarrow X\wedge L^{\wedge j}=
\mathcal{P}(j)\wedge_{\Sigma_j} L^{\wedge j}\quad,
\]
for some orbit cofibrant $X$ depending on $j$. Using proposition~\ref{Prop:orbitsquareproduct} again,
this map is an orbit q-cofibration.

Orbit q-cofibrations are stable under cobase change.
In the filtration for $PL$ we now have that each $F_{j-1}PL\rightarrow F_jPL$ is an orbit q-cofibration
for $j\geq2$. Observe that also $L\rightarrow \mathcal{P}(1)\wedge L =F_1PL$
is an orbit q-cofibration.

We now have a sequence of orbit q-cofibrations and by proposition~\ref{Prop:orbitqcofincolim}
we have that
\[
\eta:L\rightarrow \colim_j F_jPL=PL
\]
is an orbit q-cofibration.
\end{proof}

\begin{Prop}\label{Prop:Ppreservecof}
Let $f:K\rightarrow L$ be a map under $S$. 
If $f$ and $S\rightarrow K$ are orbit q-cofibrations and each $\mathcal{P}(j)$ can be written equivariantly
as a product $X\wedge (\Sigma_j)_+$ with $X$ being an orbit cofibrant (non-equivariant) orthogonal spectrum,
then $Pf:PK\rightarrow PL$ is also an orbit q-cofibration.
\end{Prop}

\begin{proof}
By proposition~\ref{Prop:orbitqcofincolim} we must show that 
\[
F_1PK=\mathcal{P}(1)\wedge K\rightarrow \mathcal{P}(1)\wedge L= F_1PL
\]
is an orbit q-cofibration, and that for every $j\geq2$ the map
\[
F_jPK\cup_{F_{j-1}PK}F_{j-1}PL\rightarrow F_{j}PL
\]
also is an orbit q-cofibration. The first statement follows directly from the assumptions together with
proposition~\ref{Prop:orbitsquareproduct}. The second statement is proved as follows:

Fix some $j\geq2$ and consider the diagram
\[
\begin{CD}
sK^{\wedge (j-1)} @>>> K^{\wedge j}\\
@VVV @VVV\\
sL^{\wedge (j-1)} @>>> L^{\wedge j}
\end{CD}\quad.
\]
Observe that 
\[
K^{\wedge j}\cup_{sK^{\wedge (j-1)}} sL^{\wedge (j-1)}= 
K\wedge L^{\wedge j-1}\cup L\wedge K\wedge L^{\wedge j-2}\cup\cdots\cup L^{\wedge j-1}\wedge K\quad.
\]
By proposition~\ref{Prop:orbitsquareproduct} the canonical map
\[
K^{\wedge j}\cup_{sK^{\wedge (j-1)}} sL^{\wedge (j-1)}\rightarrow L^{\wedge j} 
\]
is an orbit q-cofibration. Since $\mathcal{P}(j)=X\wedge (\Sigma_j)_+$ for some orbit cofibrant $X$,
the functor $\mathcal{P}(j)\wedge_{\Sigma_j}-$
takes orbit q-cofibrations to orbit q-cofibrations, see 
proposition~\ref{Prop:orbitsquareproduct}.
It follows that the map
\[
\mathcal{P}(j)\wedge_{\Sigma_j}K^{\wedge j}\cup_{\mathcal{P}(j)\wedge_{\Sigma_j}sK^{\wedge (j-1)}}
\mathcal{P}(j)\wedge_{\Sigma_j}sL^{\wedge (j-1)}\rightarrow \mathcal{P}(j)\wedge_{\Sigma_j}L^{\wedge j}
\]
is an orbit q-cofibration. Now take a look at the diagram
\[
\begin{CD}
F_{j-1}PK @<<< \mathcal{P}(j)\wedge_{\Sigma_j}sK^{\wedge (j-1)} @>>> \mathcal{P}(j)\wedge_{\Sigma_j}K^{\wedge j}\\
@VVV @VVV @VVV\\
F_{j-1}PL @<<< \mathcal{P}(j)\wedge_{\Sigma_j}sL^{\wedge (j-1)} @>>> \mathcal{P}(j)\wedge_{\Sigma_j}L^{\wedge j}
\end{CD}\quad.
\]
Apply lemma~\ref{Lem:polem} and use that orbit q-cofibrations are stable under cobase change to conclude that 
\[
F_jPK\cup_{F_{j-1}PK}F_{j-1}PL\rightarrow F_{j}PL
\]
is an orbit q-cofibration.
\end{proof}

\begin{Prop}\label{Prop:piisoofoperadsonL}
If for each $j$ the map $\alpha:\mathcal{P}(j)\rightarrow \mathcal{Q}(j)$ can be written equivariantly
as a product $X\wedge(\Sigma_j)_+\rightarrow Y\wedge(\Sigma_j)_+$ where $X\rightarrow Y$
is a $\pi_*$-isomorphism between orbit cofibrant 
orthogonal spectra, and $S\rightarrow L$
is an orbit q-cofibration, then $PL\rightarrow QL$ is a $\pi_*$-isomorphism.
\end{Prop}

\begin{proof}
Since each $\mathcal{P}(j)$ is a product $X\wedge (\Sigma_j)_+$ where $X$ is orbit cofibrant,
the first part of the proof of proposition~\ref{Prop:Pconserveqcof} shows that each map
\[
F_{j-1}PL\rightarrow F_jPL
\]
is an orbit q-cofibration. Similarly we have that each $F_{j-1}QL\rightarrow F_jQL$ also is an orbit q-cofibration.
Observe that any orbit q-cofibration is an l-cofibration. 
Hence by proposition~\ref{Prop:lcofseqpiiso} it is enough to show that each $F_jPL\rightarrow F_jQL$
is a $\pi_*$-isomorphism.

Fix $j\geq1$. 
We now perform a little trick using the cofibrant replacement functor $\Gamma$:
Let $K$ be $sL^{\wedge (j-1)}$ or $L^{\wedge j}$. 
In the argument that follows $K$ could in fact be any
orthogonal $\Sigma_j$-spectrum. The diagram
\[
\begin{CD}
\mathcal{P}(j)\wedge_{\Sigma_j} \Gamma K @>{\simeq}>> \mathcal{Q}(j)\wedge_{\Sigma_j} \Gamma K\\
@V{\simeq}VV @VV{\simeq}V\\
\mathcal{P}(j)\wedge_{\Sigma_j} K @>>> \mathcal{Q}(j)\wedge_{\Sigma_j} K
\end{CD}
\]
can be written as
\[
\begin{CD}
X\wedge \Gamma K @>{\simeq}>> Y\wedge \Gamma K\\
@V{\simeq}VV @VV{\simeq}V\\
X\wedge K @>>> Y\wedge K
\end{CD}
\]
for some orbit cofibrant $X$ and $Y$. Since $\Gamma K\rightarrow K$ is a level-equivalence,
proposition~\ref{Prop:orbitpreservele} implies that the vertical maps are level-equivalences.
The map at the top is a $\pi_*$-iso by proposition~\ref{Prop:os:wedge}.
Thus the two maps
\[
\mathcal{P}(j)\wedge_{\Sigma_j} sL^{\wedge (j-1)}\rightarrow \mathcal{Q}(j)\wedge_{\Sigma_j} sL^{\wedge (j-1)}\quad\text{and}\quad
\mathcal{P}(j)\wedge_{\Sigma_j} L^{\wedge j}\rightarrow \mathcal{Q}(j)\wedge_{\Sigma_j}L^{\wedge j}
\]
are $\pi_*$-isomorphisms.

Now we prove by induction that $F_jPL\rightarrow F_jQL$ is a $\pi_*$-iso. For $j=1$ this follows directly
from the argument above since $F_1PL=\mathcal{P}(1)\wedge L$ and $F_1QL=\mathcal{Q}(1)\wedge L$. For the induction step
we consider the diagram
\[
\begin{CD}
F_{j-1}PL @<<< \mathcal{P}(j)\wedge_{\Sigma_j}sL^{\wedge (j-1)} @>{i}>> \mathcal{P}(j)\wedge_{\Sigma_j}L^{\wedge j}\\
@VVV @VVV @VVV\\
F_{j-1}QL @<<< \mathcal{Q}(j)\wedge_{\Sigma_j}sL^{\wedge (j-1)} @>{i'}>> \mathcal{Q}(j)\wedge_{\Sigma_j}L^{\wedge j}
\end{CD}\quad.
\]
The vertical maps are $\pi_*$-isos and the maps marked $i$ and $i'$ are l-cofibrations. 
By proposition~\ref{Prop:os:gluepiisolcof}
we get that the row-wise pushout, $F_jPL\rightarrow F_jQL$, is again a $\pi_*$-iso.
\end{proof}

\begin{Prop}\label{Prop:pifreePpreservepiiso}
If each $\mathcal{P}(j)$ can be written equivariantly
as a product $X\wedge (\Sigma_j)_+$ with $X$ being an orbit cofibrant (non-equivariant) orthogonal spectrum,
and the maps $S\rightarrow L$ and $S\rightarrow K$ are
orbit q-cofibrations, and $f:L\rightarrow K$ a $\pi_*$-isomorphism under $S$, then the map $PL\rightarrow PK$
is also a $\pi_*$-isomorphism.
\end{Prop}

\begin{proof}
It is enough to show that each $F_jPL\rightarrow F_jPK$
is a $\pi_*$-iso. This follows from proposition~\ref{Prop:lcofseqpiiso}.

We have $F_1PL=\mathcal{P}(1)\wedge L$ and $F_1PK=\mathcal{P}(1)\wedge K$, so 
applying the trick of the previous proof to the diagram
\[
\begin{CD}
\mathcal{P}(1)\wedge \Gamma L @>>> \mathcal{P}(1)\wedge \Gamma K\\
@VVV @VVV\\
\mathcal{P}(1)\wedge L @>>> \mathcal{P}(1)\wedge K
\end{CD}
\]
we see that the natural map 
$F_1PL\rightarrow F_1PK$ is a $\pi_*$-iso.

Let $cL^{\wedge j}$ denote the cofiber of $sL^{\wedge j-1}\rightarrow L^{\wedge j}$. 
This map is an orbit cofibration by proposition~\ref{Prop:orbitsquareproduct}, hence also an
l-cofibration. This implies that 
$cL^{\wedge j}$ has the homotopy type of the homotopy cofiber. Similarly we can define $cK^{\wedge j}$.

Observe that proposition~\ref{Prop:squarepiiso} also holds for orbit q-cofibrations.
To see this notice that its proof is formal. In the orbit cofibrant case, we can use proposition~\ref{Prop:orbitpreservele}
instead of proposition~\ref{Prop:os:wedge}, and proposition~\ref{Prop:orbitsquareproduct} instead 
of proposition~\ref{Prop:os:qcofsquare}.
By induction on $j$ and using that the conclusion of proposition~\ref{Prop:squarepiiso} for the induction step,
we prove that the map
\[
cL^{\wedge j}\rightarrow cK^{\wedge j}
\]
is a $\pi_*$-isomorphism.

We now claim that the map
$\mathcal{P}(j)\wedge_{\Sigma_j} cL^{\wedge j}\rightarrow\mathcal{P}(j)\wedge_{\Sigma_j} cK^{\wedge j}$
is a $\pi_*$-iso. 
To check the claim,
recall that $\mathcal{P}(j)=X\wedge (\Sigma_j)_+$, where $X$ is orbit cofibrant.
Now inspect the diagram
\[
\begin{CD}
&& \Gamma X \wedge \Gamma cL^{\wedge j} @>{\simeq}>> \Gamma X \wedge \Gamma cK^{\wedge j} &&\\
&& @V{\simeq}VV @VV{\simeq}V &&\\
&& X \wedge \Gamma cL^{\wedge j} @>>> X \wedge \Gamma cK^{\wedge j} &&\\
&& @V{\simeq}VV @VV{\simeq}V &&\\
\mathcal{P}(j)\wedge_{\Sigma_j} cL^{\wedge j} @= X \wedge cL^{\wedge j} @>>> 
X \wedge \Gamma cK^{\wedge j} @= \mathcal{P}(j)\wedge_{\Sigma_j} cK^{\wedge j}
\end{CD}\quad.
\]
The propositions~\ref{Prop:orbitpreservele} and~\ref{Prop:os:wedge} show that the maps marked with $\simeq$
are $\pi_*$-isos. The claim follows.

Now inspect the map between cofiber sequences given by the filtration:
\[
\begin{CD}
F_{j-1}PL @>>> F_jPL @>>> \mathcal{P}(j)\wedge_{\Sigma_j} cL^{\wedge j}\\
@VVV @VVV @VVV\\
F_{j-1}PK @>>> F_jPK @>>> \mathcal{P}(j)\wedge_{\Sigma_j} cK^{\wedge j}
\end{CD}\quad.
\]
The first map is a $\pi_*$-iso by induction, while the last map is a $\pi_*$-iso by
the argument above. It follows that the middle vertical map 
also is a $\pi_*$-iso.
\end{proof}

\subsection{Operads and simplicial objects}

We are now going to discuss how a monad $P$ induced by an operad 
extends to simplicial orthogonal spectra.
Recall that a simplicial orthogonal spectrum is a functor
\[
L_{\bullet}:\catDelta^{\op}\rightarrow \mathscr{IS}\quad.
\]
But $P$ has domain orthogonal spectra under the sphere spectrum, so we cannot form
the composition $PL_{\bullet}$ unless we demand that $L_{\bullet}$ has a chosen lifting
to a simplicial object in $S\downarrow\mathscr{IS}$.
Luckily this is not a restrictive condition; a functor $\catDelta^{\op}\rightarrow (S\downarrow \mathscr{IS})$
is equivalent to a simplicial orthogonal spectrum $L_{\bullet}$ together with
a chosen map $S\rightarrow L_0$.

Geometric realization of a simplicial object in $S\downarrow\mathscr{IS}$
yields an orthogonal spectrum under $S$. And we ask if $P$ and $|-|$ commutes:

\begin{Prop}
Let $L_{\bullet}$ be a simplicial object in $S\downarrow\mathscr{IS}$
and $\mathcal{P}$ any operad in orthogonal spectra,
then there is a natural isomorphism
\[
\nu:|PL_{\bullet}|\rightarrow P|L_{\bullet}|\quad,
\]
such that the following two diagrams commute:
\[
\begin{CD}
|L_{\bullet}|@>{|\eta|}>>|PL_{\bullet}|\\
@V{=}VV @VV{\nu}V\\
|L_{\bullet}|@>{\eta}>> P|L_{\bullet}|
\end{CD}
\quad\text{and}\quad
\begin{CD}
|PPL_{\bullet}|@>{P\nu\;\nu}>>PP|L_{\bullet}|\\
@V{|\mu|}VV @VV{\mu}V\\
|PL_{\bullet}|@>{\nu}>>P|L_{\bullet}|
\end{CD}\quad.
\]
\end{Prop}

\begin{proof}
We use the filtration of $P$ and construct $\nu:|F_jPL_{\bullet}|\rightarrow F_jP|L_{\bullet}|$ by induction. 
For $j=1$ we let
\[
|F_1PL_{\bullet}|=|\mathcal{P}(1)\wedge L_{\bullet}|
\rightarrow \mathcal{P}(1)\wedge |L_{\bullet}|=F_1P|L_{\bullet}|
\]
be the natural isomorphism given by lemma~\ref{Lem:sosrealizprod}. 
For the inductive step we consider the diagram
\[
\begin{CD}
|\mathcal{P}(j)\wedge_{\Sigma_j}L^{\wedge j}_{\bullet}| @<<<
|\mathcal{P}(j)\wedge_{\Sigma_j}sL^{\wedge j-1}_{\bullet}| @>>> |F_{j-1}PL_{\bullet}|\\
@VVV @VVV @VVV\\
\mathcal{P}(j)\wedge_{\Sigma_j}|L_{\bullet}|^{\wedge j} @<<<
\mathcal{P}(j)\wedge_{\Sigma_j}s|L_{\bullet}|^{\wedge j-1} @>>> F_{j-1}P|L_{\bullet}|
\end{CD}\quad.
\]
The first two vertical maps are isomorphisms by lemma~\ref{Lem:sosrealizprod}, while the last is an
isomorphism by induction. It follows that the map between the pushouts also is an isomorphism.

Since $\nu:|PL_{\bullet}|\rightarrow P|L_{\bullet}|$ 
is a colimit of isomorphisms, it is itself an isomorphism.

To see that $\nu$ is unital we just observe from the construction above that the following
diagram commutes:
\[
\begin{CD}
|L_{\bullet}|@>{|\eta|}>> |F_1PL_{\bullet}|@>>>|PL_{\bullet}|\\
@V{=}VV @VVV @VV{\nu}V\\
|L_{\bullet}|@>{\eta}>> F_1P|L_{\bullet}| @>>> P|L_{\bullet}|
\end{CD}\quad.
\]

At last we want to check that $\nu\circ |\mu|$ is equal to $\mu\circ P\nu\circ\nu$.
Since the natural map
$\bigvee_j\mathcal{P}(j)\wedge K^{\wedge j}\rightarrow PK$ is surjective for any orthogonal spectrum 
$K$ under $S$, it is enough to check that the following diagram commutes:
\[
\begin{CD}
\left| \bigvee_j\mathcal{P}(j)\wedge\left(\bigvee_k\mathcal{P}(k)\wedge L_{\bullet}^{\wedge k}\right)^{\wedge j}\right| @>>>
\bigvee_j\mathcal{P}(j)\wedge\left(\bigvee_k\mathcal{P}(k)\wedge| L_{\bullet}|^{\wedge k}\right)^{\wedge j}\\
@VVV @VVV\\
\left| \bigvee_j\mathcal{P}(j)\wedge L_{\bullet}^{\wedge j}\right| @>>>
\bigvee_j\mathcal{P}(j)\wedge| L_{\bullet}|^{\wedge j}
\end{CD}\quad.
\]
To do this recall the stepwise definition of $\mu$, use that wedge and geometric realization commute
and lemma~\ref{Lem:sosrealizprod}. This finishes the proof.
\end{proof}

\begin{Cor}\label{cor:realsimpPalg}
Let $\mathcal{P}$ be an operad in orthogonal spectra.
The geometric realization is a functor from simplicial $P$-algebras to $P$-algebras.
\end{Cor}

\begin{proof}
A simplicial $P$-algebra is a functor 
\[
L_{\bullet}:\catDelta^{\op}\rightarrow\{ P\text{-algebras}\}\quad.
\]
Its geometric realization is defined by realizing the underlying simplicial orthogonal spectrum.
We define $\theta$ for $|L_{\bullet}|$ as the composition
\[
P|L_{\bullet}|\underset{\cong}{\xleftarrow{\nu}}|PL_{\bullet}|\rightarrow|L_{\bullet}|\quad.
\]
Here the last map uses the $P$-algebra structure of each $L_q$.
Using the two diagrams in the proposition above, we easily see that this
is a $P$-algebra.

Functorality follows from naturality of $\nu$.
\end{proof}

\subsection{The bar construction}

We are now going to recall the definition of May's two-sided bar construction.
In~\cite{May:72} May uses this bar construction in relation with operads in topological spaces.
In this subsection we will prove his results for operads in orthogonal spectra.
Due to remark~\ref{Rem:alloperadstoodifficult} we will not consider arbitrary operads,
but only those $\mathcal{P}$ such that each $\mathcal{P}(j)$ can be written equivariantly
as a product $X\wedge (\sigma_j)_+$ for some orbit cofibrant $X$.
We are particularly interested in improving a $\mathcal{P}$-algebra to a homotopy equivalent 
$\mathcal{Q}$-algebra, when given a map of operads $\mathcal{P}\rightarrow\mathcal{Q}$ such that each
$\mathcal{P}(j)\rightarrow\mathcal{Q}(j)$ can be written equivariantly as the
product of a $\pi_*$-iso $X\rightarrow Y$ and $(\Sigma_j)_+$.

\begin{Def}
Let $C$ be a monad, $F$ a $C$-functor and $X$ a $C$-algebra, then we define
$B_{q}(F,C,X)=FC^qX$. We have face and degeneracy operators given by
\begin{align*}
d_0=\lambda &,\quad \lambda:FC^qX\rightarrow FC^{q-1}X\quad,\\
d_i=FC^{i-1}\mu &,\quad \mu:C^{q-i+1}X\rightarrow C^{q-i}X \quad\text{for $0<i<q$},\\
d_q=FC^{q-1}\theta &, \quad \theta:CX\rightarrow X\quad\text{and}\\
s_i=FC^{i}\eta &,\quad \eta:C^{q-i}X\rightarrow C^{q-i+1}X\quad.
\end{align*}
And we define $B(F,C,X)$ as the geometric realization of $B_{\bullet}(F,C,X)$.
\end{Def}

We now specify the situation we are interested in. 
Let $\mathcal{P}$ be an operad in orthogonal spectra, and $L$ a $\mathcal{P}$-algebra. 
In addition let $\alpha:\mathcal{P}\rightarrow\mathcal{Q}$ be a map of operads.
Now assume that:
\begin{itemize}
\item[] the unit $S\rightarrow L$ is an orbit q-cofibration,
\item[] the unit $S\rightarrow \mathcal{P}(1)$ is an orbit q-cofibration and each 
$\mathcal{P}(j)$ can be written equivariantly
as a product $X\wedge (\Sigma_j)_+$ with $X$ being an orbit cofibrant (non-equivariant) orthogonal spectrum,
\item[] the unit $S\rightarrow \mathcal{Q}(1)$ is an orbit q-cofibration and 
each $\mathcal{Q}(j)$ can be written equivariantly
as a product $Y\wedge (\Sigma_j)_+$ with $Y$ being an orbit cofibrant (non-equivariant) orthogonal spectrum, and
\item[] each $\alpha:\mathcal{P}(j)\rightarrow\mathcal{Q}(j)$ can be written equivariantly
as a product $X\wedge(\Sigma_j)_+\rightarrow Y\wedge(\Sigma_j)_+$ where $X\rightarrow Y$
is a $\pi_*$-isomorphism.
\end{itemize}
These are standing assumptions for the rest of this subsection.

In order to do calculations with the bar constructions we
should know that $[q]\mapsto QP^qL$ is good.

\begin{Lem}\label{Lem:Bisgood}
Under the assumptions above $[q]\mapsto QP^qL$ is good.
\end{Lem}

\begin{proof}
The $i$'th degeneracy operator is defined as
\[
s_i= QP^i\eta\quad,\text{where }\eta:P^{q-i}L\rightarrow P^{q-i+1}L\quad,
\]
for $0\leq i\leq q$. We will check that $s_i$ is an orbit q-cofibration.

We use proposition~\ref{Prop:Pconserveqcof} and proposition~\ref{Prop:Ppreservecof}. By simultaneous
induction we prove that the unit $S\rightarrow P^jL$ and $\eta:P^{j}L\rightarrow P^{j+1}L$
are orbit q-cofibrations. Since both $P$ and $Q$ preserves orbit q-cofibrations 
this implies that $s_i= QP^i\eta$ is an orbit q-cofibration.
\end{proof}

\begin{Lem}
$B(Q,P,L)$ is a $Q$-algebra.
\end{Lem}

\begin{proof}
By corollary~\ref{cor:realsimpPalg} it is enough to show that
$[q]\mapsto QP^qL$ is a simplicial $Q$-algebra.

For a given $q$ the $Q$-algebra structure map $\theta:QQP^qL\rightarrow QP^qL$ is defined as
the multiplication $\mu'$ of $Q$.
By naturality of $\mu'$ we easily see that the $s_i$'s for $0\leq i\leq q$ and
the $d_i$'s for $1\leq i\leq q$ are $Q$-algebra morphisms. They are all defined as
$Qf$ for suitable $f$'s. Only $d_0$ requires more checking. Recall that $d_0$ is defined
to be $\lambda:QP^qL\rightarrow QP^{q-1}L$, but $\lambda=\mu'Q\alpha$ and associativity of
$\mu'$ implies that $d_0$ is a $Q$-algebra morphism.
\end{proof}

\begin{Prop}\label{Prop:firstrepl}
Since $L$ is $P$-algebra, the evaluation $B(P,P,L)\rightarrow L$ is a map of 
$P$-algebras and a $\pi_*$-isomorphism.
\end{Prop}

\begin{proof}
Define 
\[
f_q :PP^qL\rightarrow L
\]
to be the iterated composition 
\[
PP^qL\xrightarrow{d_q}P^qL\xrightarrow{d_{q-1}}P^{q-1}L\rightarrow\cdots
\rightarrow PL\xrightarrow{\theta}L\quad.
\]
This defines a simplicial map from $[q]\mapsto PP^{q}L$ to the constant simplicial orthogonal
spectrum $[q]\mapsto L$.
Its realization is the evaluation $B(P,P,L)\rightarrow L$. Since $L$ is a $P$-algebra the following diagram
commutes for all $q$:
\[
\begin{CD}
PPP^qL@>{P f_q}>> PL\\
@V{\mu}VV @VV{\theta}V\\
PP^qL @>{f_q}>> L
\end{CD}\quad,
\]
thus the collection $f_q$ is a map of simplicial $P$-algebras, 
and by corollary~\ref{cor:realsimpPalg}
it follows that $B(P,P,L)\rightarrow L$ is also a $P$-algebra map.

Using the unit $\eta:L\rightarrow PL$ we define a coretraction for the evaluation map.
On the level of $q$-simplices it is defined by
\[
s_0^q\eta: L\rightarrow PP^qL\quad.
\]
Be warned that $L\rightarrow B(P,P,L)$ is not a $P$-algebra map.

The composition $L\rightarrow B(P,P,L)\rightarrow L$ is clearly the identity.
Composing in the opposite order we get $B(P,P,L)\rightarrow L\rightarrow B(P,P,L)$
and the resulting map is the realization of
\[
s_0^q\eta f_q:PP^qL\rightarrow PP^qL\quad.
\]
Now it is easy to see that the maps $h_i:PP^qL\rightarrow PP^{q+1}L$ for $0\leq i\leq q$,
defined by 
\[
h_i=s_0^i\eta d_0^i:PP^qL\rightarrow PP^{q+1}L\quad,
\]
give a simplicial homotopy between $s_0^q\eta f_q$ and the identity. Thus
$B(P,P,L)\rightarrow L$ is a $\pi_*$-isomorphism.
\end{proof}

\begin{Prop}\label{Prop:secrepl}
Under the assumptions above
the map $B(P,P,L)\rightarrow B(Q,P,L)$ 
induced by $\alpha$ is both a $\pi_*$-isomorphism and a map of $P$-algebras.
\end{Prop}

\begin{proof}
By inductive use of proposition~\ref{Prop:Pconserveqcof} 
the map $S\rightarrow P^qL$ is an orbit q-cofibration for every $q$.
By proposition~\ref{Prop:piisoofoperadsonL} it follows that $PP^qL\rightarrow QP^qL$ is  
a $\pi_*$-isomorphism for all $q$. The simplicial
orthogonal spectra $[q]\mapsto PP^qL$ and $[q]\mapsto QP^qL$ both are good by lemma~\ref{Lem:Bisgood},
hence $B(P,P,L)\rightarrow B(Q,P,L)$ is a $\pi_*$-iso. 

The last claim is easily checked: $\alpha$ induces a $P$-algebra structure on $QP^qL$, and the map
\[
\big([q]\mapsto PP^qL\big)\rightarrow \big([q]\mapsto QP^qL\big)
\]
is a map of simplicial $P$-algebras. We now appeal to corollary~\ref{cor:realsimpPalg}.
\end{proof}

\begin{Rem}\label{Rem:replproc}
Together the two propositions~\ref{Prop:firstrepl} and~\ref{Prop:secrepl} give a procedure for replacing
a $P$-algebra by a $P$-equivalent $Q$-algebra. We have
\[
B(Q,P,L)\leftarrow B(P,P,L)\rightarrow L\quad.
\]
Here both maps are $P$-algebra maps and $\pi_*$-isomorphisms. In addition the first $P$-algebra is also a $Q$-algebra.
Moreover, the following proposition shows that $B(Q,P,L)$ is unique up to $\pi_*$-isomorphism of $Q$-algebras.
\end{Rem}

\begin{Prop}\label{Prop:uniquerepl}
In addition to the assumptions above assume that   
$S\rightarrow L'$ and $S\rightarrow A$ are orbit q-cofibrations. If 
$A\leftarrow L'\rightarrow L$ are $\pi_*$-isomorphisms of $P$-algebras, with 
$A$ a $Q$-algebra, then there are $\pi_*$-isomorphisms
\[
A\leftarrow B(Q,P,L')\rightarrow B(Q,P,L)
\]
of $Q$-algebras.
\end{Prop}

\begin{proof}
By functorality of the two-sided bar construction we have maps 
\[
B(Q,P,L) \leftarrow B(Q,P,L')\rightarrow B(Q,P,A)\rightarrow B(Q,Q,A)\rightarrow A\quad.
\]
All these maps are easily seen to be maps of $Q$-algebras.
By proposition~\ref{Prop:firstrepl} the last map is a $\pi_*$-iso.

By proposition~\ref{Prop:pifreePpreservepiiso} and corollary~\ref{cor:realsimpPalg} the first two maps are
$\pi_*$-isomorphisms.

Combining the propositions~\ref{Prop:piisoofoperadsonL} and~\ref{Prop:pifreePpreservepiiso}, we prove by induction that
each 
\[
QP^qA\rightarrow Q^{q+1}A
\]
is a $\pi_*$-iso. Hence by corollary~\ref{cor:realsimpPalg} also the map
$B(Q,P,A)\rightarrow B(Q,Q,A)$ is a $\pi_*$-iso.
\end{proof}

\section{Involution operads on $S[\Omega M]$}\label{sect:involution}

The main concern of the previous sections has been to set up a theory
for operads in orthogonal spectra.
Let $M$ be a compact manifold and $\xi$ a vector bundle over $M$.
In this section we will apply this theory 
to construct involutions, depending on $\xi$, on orthogonal ring spectra $R$,
which are $\pi_*$-isomorphic to $S[\Omega M]$. 
Thus $R$ and $S[\Omega M]$ have identical homotopy groups, and we want
the involution on $R$ to coincide with the involution $\iota$ on 
$\pi_*S[\Omega M]$, where $\iota$ is given as follows:

\begin{Def}\label{Def:iota}
Assume that $\xi$ is an $n$-vector bundle.
A class in $\pi_q S[\Omega M]$ is represented by a map
\[
\alpha:S^{q+k}\rightarrow \Omega M_+\wedge S^k\quad.
\]
Parallel transportation in $\xi$ along loops in $M$ gives a homomorphism
\[
P:\Omega M\rightarrow GL(\R^n)\quad.
\]
Using $P$, we define a map $\bar{P}:S^n\wedge \Omega M_+\rightarrow \Omega M_+\wedge S^n$ by sending $(v,\gamma)$
to $(\bar{\gamma},P(\gamma)(v))$. We have transported $v$ along $\gamma$ and reversed the loop.
The \textit{involution $\iota$} is now defined by sending the class of $\alpha$ to the class of the composition
\[
S^{n}\wedge S^{q+k}\xrightarrow{\id\wedge \alpha} S^n\wedge \Omega M_+\wedge S^k
\xrightarrow{\bar{P}\wedge\id}\Omega M_+\wedge S^{k+n}\quad.
\]
\end{Def}

The strategy now is to construct operads $\mathcal{D}_n$ in orthogonal spectra, 
subsection~\ref{subsect:theconstruction}, for positive integers $n$. When $n$ is the
fiber dimension of $\xi$, we show in subsection~\ref{subsect:action} 
that $S[\Omega M]$ is a $\mathcal{D}_n$-algebra.
The next step, subsection~\ref{subsect:homotopydiscreteness}, 
is to show that $\mathcal{D}_n$ and $\mathcal{H}$ are 
``homotopy equivalent'' operads. In subsection~\ref{subsect:maintheorem}
we bring everything together and state and prove the main theorem.

\subsection{The construction}\label{subsect:theconstruction}

Here we will design operads $\mathcal{D}_n$ in orthogonal spectra for every positive integer $n$.
Their purpose is to encode the involution on $S[\Omega M]$ given by an $n$-dimensional vector bundle $\xi$
over a manifold $M$. In subsection~\ref{subsect:action} below we will see that 
$S[\Omega M]$ is a $\mathcal{D}_n$-algebra. The main result of this subsection is:

\begin{Thm}\label{Thm:Dnisoperad}
$\mathcal{D}_n$, as constructed below, is an operad in orthogonal spectra. And there is a map of operads
$\mathcal{D}_n\rightarrow\mathcal{H}$.
\end{Thm}

The proof will be given at the end of this subsection, before that we have to construct the operad.
We will also provide an orbit cofibrant replacement for $\mathcal{D}_n$, see theorem~\ref{Thm:ocreplforDn} below.

Our first aim is to define for each $j$ orthogonal $\Sigma_j$-spectra $\mathcal{D}_n(j)$.
But before reaching this aim we have to introduce
the topological groups $D_n(j;V)$. Here $V$ is a finite dimensional real inner product space.
We write the group operation of $D_n(j;V)$ multiplicatively, $1$ is the unit,
and we define this group by specifying generators and relations. 
The generators have the form
$(\phi,r)$, where $\phi:\R^n\hookrightarrow V$ is an isometric embedding
and $1\leq r\leq j$ an integer. Notice that the generators form a topological space.
There are two classes of relations. These are:
\begin{itemize}
\item[i)] Cancellation of repeated pairs: For all $\phi$ and $r$ we set
\[
(\phi,r)(\phi,r)=1
\]
\item[ii)] Orthogonal pairs commute: Whenever $\phi\perp\psi$ and $r\neq r'$, we set
\[
(\phi,r)(\psi,r')=(\psi,r')(\phi,r)\quad.
\]
\end{itemize}
There is an increasing filtration of $D_n(j;V)$. Let $F_m D_n(j;V)$ be the subset of all elements
represented by words with $m$ or fewer letters. As a topological space
$F_m D_n(j;V)$ is a quotient of $\coprod_{i=0}^{m} X^{\times i}$, where $X$ is the space of generators.
And the group $D_n(j;V)$ has the topology of the union (=colimit topology).
Notice that an isometric embedding $V\rightarrow V'$ induces a homomorphism $D_n(j;V)\rightarrow D_n(j;V')$.
We are now ready to define $\mathcal{D}_n(j)$.

\begin{Def}
Let $\mathcal{D}_n(j)$ be the orthogonal $\Sigma_j$-spectrum given by the formula
\[
\mathcal{D}_n(j)(V)=(D_n(j;V)\times \Sigma_j)_+\wedge S^V\quad.
\]
\end{Def}

The right assembly, $\sigma:\mathcal{D}_n(j)(V)\wedge S^W\rightarrow \mathcal{D}_n(j)(V\oplus W)$,
is induced by the natural isometric embedding $V\rightarrow V\oplus W$. 
Let $(x,\rho,v)$ be a point in $(D_n(j;V)\times \Sigma_j)_+\wedge S^V$ and $\nu\in\Sigma_j$.
Then the right $\Sigma_j$-action $\mathcal{D}_n(j)(V)\times \Sigma_j\rightarrow \mathcal{D}_n(j)(V)$
is simply given by 
\[
((x,\rho,v),\nu)\mapsto (x,\rho\nu,v)\quad.
\]

\begin{Lem}
Each $\mathcal{D}_n(j)$ can be written equivariantly
as a product $X\wedge (\Sigma_j)_+$, where $X$ is a (non-equivariant) orthogonal spectrum.
\end{Lem}

\begin{proof}
Let $X(V)=D_n(j;V)_+\wedge S^V$.
\end{proof}

Now recall the definition of the operad $\mathcal{H}$. Here we use the description given in remark~\ref{Rem:AltHdescr}.
And we consider $\mathcal{H}$ to be an operad in orthogonal spectra via suspension, see example~\ref{Exa:suspoperads}.
We will now define for each $j$ a map
\[
\mathcal{D}_n(j)\rightarrow\mathcal{H}(j)\quad.
\]
After defining the operad structure on $\mathcal{D}_n$ we will see that
the collection of these maps defines a map of operads. We begin the construction by defining 
group homomorphisms $p:D_n(j;V)\rightarrow (\Z/2)^j$. 
Recall that we write $\Z/2$ multiplicatively.
The maps $p$ are given on generators by
\[
p(\phi,r)=(1,\ldots,1,-1,1,\ldots,1)\quad,
\]
where the $-1$ lies in the $r$'th factor. The map above is now defined as
\[
\begin{CD}
\mathcal{D}_n(j)(V)&=& (D_n(j;V)\times \Sigma_j)_+\wedge S^V\\
&& @V{(p\times \id)_+\wedge\id}VV\\
\mathcal{H}(j)(V) &=& ((\Z/2)^j \times \Sigma_j)_+\wedge S^V
\end{CD}\quad.
\]

The construction of composition operations on $\mathcal{D}_n$ is quite abstract. 
Let me therefore suggest to the reader to
take a look at how the action of $\mathcal{D}_n$ on $S[\Omega M]$ is defined, see
the first part of subsection~\ref{subsect:action}, before proceeding
with the details given here. The action has been the author's guideline when defining
the composition operation. How to define the operad structure on $\mathcal{D}_n$ is forced
by the formulas for the action. 
So keeping the main geometrical idea behind the action in mind, 
will probably help the reader to understand this subsection. 

Assume for a moment that we can define non-$\Sigma$-equivariant composition 
operations $\circ_i:D_n(j;V)\times D_n(k;W)\rightarrow D_n(j+k-1;V\oplus W)$.
Using these our next goal is to define composition operations on the product $D_n(j;V)\times \Sigma_j$.
There will be a twist by $p$ in the definition of $\circ_i$ on the product.
Write $p_i(x)$ for the $i$'th
factor of $p(x)$. Let $(x,\rho)$ and $(y,\upsilon)$
be elements in $D_n(j;V)\times\Sigma_j$ and $D_n(k;W)\times \Sigma_k$ respectively, 
then the formula is
\[
(x,\rho)\circ_i(y,\upsilon)=
(x\circ_{\rho(i)}y,\rho\circ_i(\tau_k({p_{\rho(i)}(x)})\upsilon))\quad.
\]
Here $\tau_k$ denotes the group homomorphism $\Z/2\rightarrow \Sigma_k$ which sends $-1$ to the
order reversing permutation. Compare this formula with
the formula defining $\mathcal{H}$ in remark~\ref{Rem:AltHdescr}.

\paragraph{Details:}
We cannot postpone the details any more.
To define the composition operations 
\[
\circ_i:D_n(j;V)\times D_n(k;W)\rightarrow D_n(j+k-1;V\oplus W)
\]
we first define a homomorphism $c_i:D_n(k;W)\rightarrow D_n(j+k-1;V\oplus W)$.
Next we define a left action, $\vdash_i$, depending on $i$, of $D_n(j;V)$ on $D_n(j+k-1;V\oplus W)$.
Then $\circ_i$ is given by the formula
\[
x\circ_i y=x\vdash_i c_i(y)\quad.
\]
Let $i_1:V\rightarrow V\oplus W$ and $i_2:W\rightarrow V\oplus W$
be the natural inclusions. 
$c_i$ is defined on generators of $D_n(k;W)$ by
\[
c_i(\phi,r)=(i_2\phi,r+i-1)\quad.
\]
To define $\vdash_i$, we first introduce an automorphism, $z\mapsto\bar{z}$ of $D_n(j+k-1;V\oplus W)$.
This automorphism depends on $i$ and is 
defined on generators by
\[
(\phi,r)\mapsto\begin{cases}
(\phi,r) &\text{for $r<i$,}\\
(\phi,k+2i-r-1)&\text{for $i\leq r< k+i$, and}\\
(\phi,r) &\text{for $k+i\leq r$.}
\end{cases}
\]
For a generator $(\phi,r)$ in $D_n(j;V)$
and an element $z\in D_n(j+k-1;V\oplus W)$ we now define $\vdash_i$ by
\[
(\phi,r)\vdash_i z=\begin{cases}
(i_1\phi,r)z&\text{for $r<i$,}\\
(i_1\phi,i+k-1)\cdots(i_1\phi,i)\bar{z}&\text{for $r=i$, and}\\
(i_1\phi,r+k-1)z&\text{for $r>i$.}
\end{cases}
\]

To prove formulas containing the composition operators we often do induction on the 
length of a word in the $D_n(-;-)$ groups. We also need some basic formulas.
These are:

\begin{Lem}\label{Lem:circformulas}
Let $x$ be a word and $(\phi,r)$ a generator of $D_n(j;V)$, let $y$ be a word and $(\psi,r')$
a generator of $D_n(k;W)$, and let $z$ be a word in $D_n(j+k-1;V\oplus W)$. We have:
\begin{itemize}
\item[i)] $1\circ_i 1=1$,
\item[ii)] $x\vdash_i z= (x\circ_i 1)z$ if $p_i(x)=1$,
\item[iii)] $x\vdash_i z= (x\circ_i 1)\bar{z}$ if $p_i(x)=-1$,
\item[iv)] $(x(\phi,r))\circ_i 1=(x\circ_i 1)(\phi,r)$ if $r<i$,
\item[v)] $(x(\phi,r))\circ_i 1=(x\circ_i 1)(\phi,i+k-1)\cdots(\phi,i)$ if $r=i$ and $p_i(x)=1$, 
\item[vi)] $(x(\phi,r))\circ_i 1=(x\circ_i 1)(\phi,i)\cdots(\phi,i+k-1)$ if $r=i$ and $p_i(x)=-1$,
\item[vii)] $(x(\phi,r))\circ_i 1=(x\circ_i 1)(\phi,r+k-1)$ if $r>i$,
\item[viii)] $x\circ_i (y(\psi,r'))=(x\circ_i y)(\psi,r'+i-1)$ if $p_i(x)=1$,
\item[ix)] $x\circ_i (y(\psi,r'))=(x\circ_i y)(\psi,k+i-r')$ if $p_i(x)=-1$,
\item[x)] $(x(\phi,r))\circ_i y=(x\circ_i y)(\phi,r)$ if $r<i$, and
\item[xi)] $(x(\phi,r))\circ_i y=(x\circ_i y)(\phi,r+k-1)$ if $r>i$.
\end{itemize}
Observe that we have omitted the canonical inclusions $i_1$ and $i_2$ from the notation.
\end{Lem}

\begin{proof}
All formulas are checked by inspecting the definition of $\circ_i$. 
To illustrate the techniques involved we write out the proofs for vi) and x).

\paragraph{vi):} We have
\begin{align*}
(x(\phi,r))\circ_i 1 &= x\vdash_i ((\phi,i)\vdash_i 1)\\
&=x\vdash_i \big( (\phi,i+k-1)\cdots (\phi,i)\big)\\
&=(x\circ_i 1)\overline{(\phi,i+k-1)\cdots (\phi,i) }\\
&=(x\circ_i 1)(\phi,i)\cdots (\phi,i+k-1)\quad.
\end{align*}
Here we have used iii) and that $\vdash_i$ is a group action.

\paragraph{x):} Observe that $(\phi,r)c_i(y)=c_i(y)(\phi,r)$ by the ``orthogonal pairs commute''
relation in $D_n(j+k-1;V\oplus W)$. Assume that $p_i(x)=1$, so that ii) applies. Now we have:
\begin{align*}
(x(\phi,r))\circ_i y &= x\vdash_i \big( (\phi,r)\vdash_i c_i(y)\big)\\
&= \big(x\circ_i 1\big) \big( (\phi,r) c_i(y)\big)\\
&= (x\circ_i 1) c_i(y)(\phi,r) \\
&= \big(x\vdash_i c_i(y))(\phi,r) \\
&= (x\circ_i y)(\phi,r)\quad.
\end{align*}
If $p_i(x)=-1$ we use iii) instead of ii) in the above calculation.
\end{proof}

We have the following proposition telling us how $p$ and $\circ$ interact:

\begin{Prop}\label{Prop:pandcirc}
Let $x\in D_n(j;V)$ and $y\in D_n(k;W)$.
Then
\[
p_h(x\circ_i y)=\begin{cases}
p_h(x)&\text{for $h<i$,}\\
p_{h-i+1}(y)&\text{for $i\leq h<i+k$ and $p_i(x)=1$,}\\
-p_{i+k-h}(y)&\text{for $i\leq h<i+k$ and $p_i(x)=-1$, and}\\
p_{h-k+1}(x)&\text{for $i+k\leq h$.}
\end{cases}
\]
\end{Prop}

In the case $i\leq h<i+k$ we can rewrite the formula as
$p_h(x\circ_i y)=p_i(x)\cdot p_{\tau_k({p_i(x)})(h-i+1)}(y)$.

\begin{proof}
The cases $h<i$ and $i+k\leq h$
follow immediately from the definitions. Assume that $i\leq h<i+k$. 
The proof proceeds by induction on the length of $x$.
If $x=1$, then $p_i(x)=1$ and we have
\[
p_h(x\circ_i y)=p_h(1\circ_i y)=p_h(c_i(y))=p_{h-i+1}(y)\quad.
\]

Now let $x=(\phi,r)x'$ and assume that the formula is true for $x'$.
There are six cases to consider. For $r$ we have three possibilities $r<i$, $r=i$ or $r>i$,
and $p_i(x')$ can be $1$ or $-1$. We check two cases carefully, and leave the four others to the reader.

If $r<i$ and $p_i(x')=-1$, then also $p_i(x)=-1$. The left side of the formula becomes
\[
p_h(((\phi,r)x')\circ_i y)=p_h((\phi,r)(x'\circ_i y))
=p_h((\phi,r))\cdot p_h(x'\circ_i y)=p_h(x'\circ_i y)\quad,
\]
while the right side is $-p_{i+k-h}(y)$.
Thus the formula holds by the induction hypothesis.

If $r=1$ and $p_i(x')=1$, then $p_i(x)=p_i((\phi,i)x')=p_i(\phi,i)\cdot p_i(x')=(-1)\cdot 1=-1$.
Observe that by the definition of the involution
$z\mapsto \bar{z}$ on $D_n(j+k-1;V\oplus W)$ we have $p_h(\bar{z})=p_{k+2i-h-1}(z)$.
Calculating the left side of the formula we get
\[
p_h(((\phi,i)x')\circ_i y)=p_h((\phi,i+k-1)\cdots(\phi,i)\overline{(x'\circ_i y)})
=-p_{k+2i-h-1}(x'\circ_i y)\quad,
\]
and
the right side is $-p_{i+k-h}(y)$.
The induction hypothesis says that $p_i(x')=1$ and thus $p_{k+2i-h-1}(x'\circ_i y)=p_{i+k-h}(y)$
since $i\leq k+2i-h-1< i+k$. Therefore the formula is true for $x$.
\end{proof}

We now deduce the arithmetic rules for the composition operators on $D_n(-;-)$.

\begin{Prop}\label{Prop:circonDn}
Let $x\in D_n(j;V)$, $y\in D_n(k;W)$ and $z\in D_n(l;U)$.
The following associativity formula holds:
\[
(x\circ_i y)\circ_h z=\begin{cases}
(x\circ_h z)\circ_{i+l-1}y&\text{ for $h<i$,}\\
x\circ_i(y\circ_{h-i+1}z)&\text{for $i\leq h<i+k$ and $p_i(x)=1$,}\\
x\circ_i(y\circ_{i+k-h}z)&\text{for $i\leq h<i+k$ and $p_i(x)=-1$, and}\\
(x\circ_{h-k+1}z)\circ_i y&\text{for $i+k\leq h$.}
\end{cases}
\]
\end{Prop}

\begin{proof}
To complete this proof we have to do induction three times. Luckily we can reduce the number of cases
using the following observation:
Suppose that we already have proved the first case of the formula,
$(x\circ_i y)\circ_h z=(x\circ_h z)\circ_{i+l-1}y $ for $h<i$, then by
inserting $x=x'$, $y=z'$, $z=y'$, $l=k'$, $h=i'$ and $i=h'-k'+1$ we immediately get that
\[
(x'\circ_{h'-k'+1}z')\circ_{i'} y'=(x'\circ_{i'}y')\circ_{h'} z'\quad\text{for $i'<h'-k'+1$.}
\]
And this is the last case. Hence we need only to prove the first three cases. 

\paragraph{First induction:} We begin the proof by showing that
\[
(x\circ_i 1)\circ_h 1=\begin{cases}
(x\circ_h 1)\circ_{i+l-1}1 & \text{for $h<i$,}\\
x\circ_i(1\circ_{h-i+1}1) & \text{for $i\leq h<i+k$ and $p_i(x)=1$, and}\\
x\circ_i(1\circ_{i+k-h}1) & \text{for $i\leq h<i+k$ and $p_i(x)=-1$.}\\
\end{cases}
\]
This is proved by induction on the length of $x$.
Assume that $x=x'(\phi,r)$. To complete the induction step 
we need to check the following cases individually:
\begin{itemize}
\item[]\begin{tabular}{l}
$h<i$ and $r<h$,\\
$h<i$, $r=h$ and $p_h(x')=1$,\\
$h<i$, $r=h$ and $p_h(x')=-1$,\\
$h<i$ and $h<r<i$,\\
$h<i$, $r=i$ and $p_i(x')=1$,\\
$h<i$, $r=i$ and $p_i(x')=-1$,\\
$h<i$ and $i<r$,\\
$i\leq h<i+k$, $r<i$ and $p_i(x')=1$,\\
$i\leq h<i+k$, $r<i$ and $p_i(x')=-1$,\\
$i\leq h<i+k$, $r=i$ and $p_i(x')=1$,\\
$i\leq h<i+k$, $r=i$ and $p_i(x')=-1$,\\
$i\leq h<i+k$, $r>i$ and $p_i(x')=1$, and\\
$i\leq h<i+k$, $r>i$ and $p_i(x')=-1$.\\
\end{tabular}
\end{itemize}
All cases are straight forward to check using various formulas from lemma~\ref{Lem:circformulas}.
As an illustration we verify two of the cases.

For example if $r=h<i$ and $p_h(x')=-1$ then 
\begin{align*}
((x'(\phi,r))\circ_i 1)\circ_h 1 
&= ((x'\circ_i 1)(\phi,h)) \circ_h 1\\
&= ((x'\circ_i 1)\circ_h 1) (\phi,h)\cdots (\phi, h+l-1)\\
&= ((x'\circ_h 1)\circ_{i+l-1}1) (\phi,h)\cdots (\phi, h+l-1)\\
&= ((x'\circ_h 1)(\phi,h)\cdots (\phi, h+l-1) )\circ_{i+l-1} 1\\
&= ((x'(\phi,r))\circ_h 1)\circ_{i+l-1} 1\quad.
\end{align*}
Here we have used the formulas~iv) and vi) from lemma~\ref{Lem:circformulas} and the
induction hypothesis for the middle step.

In the case $i\leq h<i+k$, $r>i$ and $p_i(x')=1$, we have
\begin{align*}
((x'(\phi,r))\circ_i 1)\circ_h 1 
&= ((x'\circ_i 1)(\phi,r+k-1)) \circ_h 1\\
&= ((x'\circ_i 1)\circ_h 1)(\phi,r+k+l-2)\\
&= (x'\circ_i(1\circ_{h-i+1}1))(\phi,r+k+l-2)\\
&= (x'(\phi,r))\circ_i(1\circ_{h-i+1}1)\quad.
\end{align*}
Here we have used the induction hypothesis and formula~vii) from lemma~\ref{Lem:circformulas}.
The other cases are left as exercises.

\paragraph{Second induction:} Next we do induction on the length of $y$ to show that
\[
(x\circ_i y)\circ_h 1= \begin{cases}
(x\circ_h 1)\circ_{i+l-1}y & \text{for $h<i$,}\\
x\circ_i(y\circ_{h-i+1}1) & \text{for $i\leq h<i+k$ and $p_i(x)=1$, and}\\
x\circ_i(y\circ_{i+k-h}1) & \text{for $i\leq h<i+k$ and $p_i(x)=-1$.}
\end{cases}
\]
Observe that for $y=1$ we have the formula we proved above. 
Assume that $y=y'(\psi,r')$. Again there are several cases to consider:
\begin{itemize}
\item[]\begin{tabular}{l}
$h<i$ and $p_i(x)=1$,\\
$h<i$ and $p_i(x)=-1$,\\
$i\leq h<i+k$, $p_i(x)=1$ and $r'<h-i+1$,\\
$i\leq h<i+k$, $p_i(x)=1$, $r'=h-i+1$ and $p_{h-i+1}(y')=1$,\\
$i\leq h<i+k$, $p_i(x)=1$, $r'=h-i+1$ and $p_{h-i+1}(y')=1$,\\
$i\leq h<i+k$, $p_i(x)=1$ and $r'>h-i+1$,\\
$i\leq h<i+k$, $p_i(x)=-1$ and $r'<i+k-h$,\\
$i\leq h<i+k$, $p_i(x)=-1$, $r'=i+k-h$ and $p_{i+k-h}(y')=1$,\\
$i\leq h<i+k$, $p_i(x)=-1$, $r'=i+k-h$ and $p_{i+k-h}(y')=-1$, and\\
$i\leq h<i+k$, $p_i(x)=-1$ and $r'>i+k-h$.\\
\end{tabular}
\end{itemize}

For example if $h<i$ and $p_i(x)=1$, then
\begin{align*}
(x\circ_i (y'(\psi,r')))\circ_h 1 &=
((x\circ_i y')(\psi,k+i-r'))\circ_h 1\\
&= ((x\circ_i y')\circ_h 1)(\psi,k+i-r')\\
&= ((x\circ_h 1)\circ_{i+l-1}y')(\psi,k+i-r')\\
&= (x\circ_h 1)\circ_{i+l-1}(y'(\psi,r'))\quad.
\end{align*}
We have used that $p_{i+l-1}(x\circ_h 1)=p_i(x)=1$, the formulas~viii) and~vi) from the lemma and the induction hypothesis.

Let us check one more case. If $i\leq h<i+k$, $r'=i+k-h$,  
$p_i(x)=-1$ and $p_{i+k-h}(y')=1$, then
\begin{align*}
(x\circ_i (y'(\psi,r')))\circ_h 1 &=
((x\circ_i y')(\psi,i+k-r'))\circ_h 1\\
&=((x\circ_i y')(\psi,h))\circ_h 1\\
&= ((x\circ_i y')\circ_h 1)(\psi,h)\cdots(\psi,h+l-1)\\
&= (x\circ_i(y'\circ_{i+k-h}1))(\psi,h)\cdots(\psi,h+l-1)\\
&= x\circ_i((y'\circ_{i+k-h}1)(\psi,i+k-h+l-1)\cdots(\psi,i+k-h))\\
&= x\circ_i((y'(\psi,i+k-h))\circ_{i+k-h} 1)\\
&= x\circ_i((y'(\psi,r'))\circ_{i+k-h} 1)\quad.
\end{align*}
We have used the formulas~ix), vi) and~v) from the lemma, the induction hypothesis and that 
$p_h(x\circ_i y')=-1$.

\paragraph{Third induction:} At last we use
induction on the length of $z$ to prove that
\[
(x\circ_i y)\circ_h z=\begin{cases}
(x\circ_h z)\circ_{i+l-1}y&\text{ for $h<i$,}\\
x\circ_i(y\circ_{h-i+1}z)&\text{for $i\leq h<i+k$ and $p_i(x)=1$, and}\\
x\circ_i(y\circ_{i+k-h}z)&\text{for $i\leq h<i+k$ and $p_i(x)=-1$.}
\end{cases}
\]
Observe that the previous induction proves this formula in the case $z=1$.
Now assume that $z=z'(\theta,r'')$. These are the cases to consider:
\begin{itemize}
\item[]\begin{tabular}{l}
$h<i$ and $p_h(x)=1$,\\
$h<i$ and $p_h(x)=-1$,\\
$i\leq h<i+k$, $p_i(x)=1$ and $p_{h-i+1}(y)=1$,\\
$i\leq h<i+k$, $p_i(x)=1$ and $p_{h-i+1}(y)=-1$,\\
$i\leq h<i+k$, $p_i(x)=-1$ and $p_{i+k-h}(y)=1$, and\\
$i\leq h<i+k$, $p_i(x)=-1$ and $p_{i+k-h}(y)=-1$.\\
\end{tabular}
\end{itemize}

We write out two of the cases: If $h<i$ and $p_h(x)=1$, then
\begin{align*}
(x\circ_i y)\circ_h (z'(\theta,r'')) 
&=((x\circ_i y)\circ_h z')(\theta,r''+h-1)\\
&=((x\circ_h z')\circ_{i+l-1}y)(\theta,r''+h-1)\\
&=((x\circ_h z')(\theta,r''+h-1))\circ_{i+l-1}y\\
&=(x\circ_h (z'(\theta,r'')))\circ_{i+l-1}y\quad.
\end{align*}
We have here used that $p_h(x\circ_i y)=p_h(x)=1$, the formulas~viii) and~x) from the lemma and
the induction hypothesis.

In the case $i\leq h<i+k$, $p_i(x)=1$, $p_{h-i+1}(y)=-1$ we have
\begin{align*}
(x\circ_i y)\circ_h (z'(\theta,r'')) 
&= ((x\circ_i y)\circ_h z')(\theta,l+h-r'')\\
&= (x\circ_i(y\circ_{h-i+1}z'))(\theta,l+h-r'')\\
&= x\circ_i((y\circ_{h-i+1}z')(\theta,l+h-i+1-r''))\\
&= x\circ_i(y\circ_{h-i+1}(z'(\theta,r'')))\quad.
\end{align*}
We have used that $p_h(x\circ_i y)=p_{h-i+1}(y)=-1$, the formulas~viii) and~ix) of the lemma and
the induction hypothesis.

The checking of all remaining cases is left to the reader.
\end{proof}

Recall that our aim is to define composition operations for the products $D_n(j;V)\times \Sigma_j$.
The formula is
\[
(x,\rho)\circ_i(y,\upsilon)=
(x\circ_{\rho(i)}y,\rho\circ_i(\tau_k({p_{\rho(i)}(x)})\upsilon))\quad.
\]
And we want to check associativity, equivariance and unity. We begin with equivariance.

\begin{Lem}\label{Lem:circequivariance}
Let $(x,\rho)\in D_n(j;V)\times \Sigma_j$, $\rho'\in \Sigma_j$, $(y,\upsilon)\in D_n(k;W)\times \Sigma_k$
and $\upsilon'\in\Sigma_k$.
There is a right action of $\Sigma_j$ on $D_n(j;V)\times \Sigma_j$ defined by
\[
(x,\rho).\rho'=(x,\rho\rho')
\]
and we have
\[
((x,\rho).\rho')\circ_i((y,\upsilon).\upsilon')
=((x,\rho)\circ_{\rho'(i)}(y,\upsilon)).(\rho'\circ_i\upsilon')\quad.
\]
\end{Lem}

\begin{proof}
This proof is easy. We have:
\begin{align*}
((x,\rho).\rho')\circ_i((y,\upsilon).\upsilon')
&=(x,\rho\rho')\circ_i(y,\upsilon\upsilon')\\
&=(x\circ_{\rho\rho'(i)}y, (\rho\rho')\circ_i(\tau_k({p_{\rho\rho'(i)}(x)})\upsilon\upsilon'))\\
&=(x\circ_{\rho\rho'(i)}y, \rho\circ_{\rho'(i)}(\tau_k({p_{\rho\rho'(i)}(x)})\upsilon)(\rho'\circ_i\upsilon'))\\
&=(x\circ_{\rho\rho'(i)}y, \rho\circ_{\rho'(i)}(\tau_k({p_{\rho\rho'(i)}(x)})\upsilon)).(\rho'\circ_i\upsilon')\\
&=((x,\rho)\circ_{\rho'(i)}(y,\upsilon)).(\rho'\circ_i\upsilon')\quad.
\end{align*}
\end{proof}

\begin{Lem}\label{Lem:circassociativity}
Let $(x,\rho)\in D_n(j;V)\times \Sigma_j$, $(y,\upsilon)\in D_n(k;W)\times \Sigma_k$ and
$(z,\mu)\in D_n(l;U)\times \Sigma_l$.
The following associativity holds for $\circ$:
\[
((x,\rho)\circ_i (y,\upsilon))\circ_h (z,\mu)=\begin{cases}
((x,\rho)\circ_h (z,\mu))\circ_{i+l-1}(y,\upsilon)&\text{ for $h<i$,}\\
(x,\rho)\circ_i((y,\upsilon)\circ_{h-i+1}(z,\mu))&\text{for $i\leq h<i+k$, and}\\
((x,\rho)\circ_{h-k+1}(z,\mu))\circ_i (y,\upsilon)&\text{for $i+k\leq h$.}
\end{cases}
\]
\end{Lem}

\begin{proof}
As in the proof of proposition~\ref{Prop:circonDn}, we observe that the first case
of this formula implies the last case. Hence it is enough to check the first two cases.

Recall the formula defining the composition operators on $D_n(-;-)\times \Sigma_-$.
We may rewrite the formula as
\[
(x,\rho)\circ_i(y,\upsilon)=\begin{cases}
(x\circ_{\rho(i)}y,\rho\circ_i\upsilon)&\text{if $p_{\rho(i)}(x)=1$, and}\\
(x\circ_{\rho(i)}y,\rho\circ_i(\tau_k\upsilon))&\text{if $p_{\rho(i)}(x)=-1$.}
\end{cases}
\]
Here we have used the convention that $\tau_k$ without an argument denotes the
order reversing permutation in $\Sigma_k$, while $\tau_k$ with an argument 
denotes the group homomorphism $\Z/2\rightarrow\Sigma_k$ sending $-1$ to the order reversing permutation.
This formula will be applied many times throughout this proof, both forward and backward.

\paragraph{A special case for $h<i$:} 
First we assume that the permutations $\rho$, $\upsilon$ and $\mu$ are the identity. 
\begin{align*}
\big( (x,\id_j)\circ_i(y,\id_k)\big)\circ_h(z,\id_l)
&=(x\circ_i y,\id_j\circ_i\tau_k({p_i(x)}))\circ_h(z,\id_l)\\
&=((x\circ_i y)\circ_h z,(\id_j\circ_i\tau_k({p_i(x)}))\circ_h\tau_l({p_h(x)}))\\
&=((x\circ_h z)\circ_{i+l-1} y, (\id_j\circ_h\tau_l({p_h(x)}))\circ_{i+l-1}\tau_k({p_i(x)}))\\
&=(x\circ_h z,\id_j\circ_h\tau_l({p_h(x)}))\circ_{i+l-1}(y,\id_k)\\
&=\big( (x,\id_j)\circ_h (z,\id_l)\big)\circ_h (y,\id_k)\quad.
\end{align*}
In addition to the formula we have used 
the associativity for the composition operators, proposition~\ref{Prop:circonDn}
for $D_n(-;-)$ and lemma~\ref{Lem:circonM} for $\circ$ of the permutations, 
and the calculations:
$\id_j\circ_i\tau_k({p_i(x)})(h)=h$, $p_h(x\circ_i y)=p_h(x)$,
$\id_j\circ_h\tau_l({p_h(x)})(i+l-1)=i+l-1$ and $p_{i+l-1}(x\circ_h z)=p_i(x)$. 

\paragraph{A special case for $i\leq h<i+k$:} Also here we let the permutations
$\rho$, $\upsilon$ and $\mu$ be the identity.
If $p_i(x)=1$ and $p_{h-i+1}(y)=1$, then
\begin{align*}
\big( (x,\id_j)\circ_i(y,\id_k)\big)\circ_h(z,\id_l)
&= (x\circ_i y,\id_j\circ_i\id_k)\circ_h(z,\id_l)\\
&= ((x\circ_i y)\circ_h z,(\id_j\circ_i\id_k)\circ_h\id_l)\\
&= (x\circ_i(y\circ_{h-i+1}z),\id_j\circ_i(\id_k\circ_{h-i+1}\id_l))\\
&= (x,\id_j)\circ_i(y\circ_{h-i+1}z, \id_k\circ_{h-i+1}\id_l)\\
&= (x,\id_j)\circ_i\big( (y,\id_k)\circ_{h-i+1}(z,\id_l)\big)\quad.
\end{align*}
If $p_i(x)=1$ and $p_{h-i+1}(y)=-1$, then
\begin{align*}
\big( (x,\id_j)\circ_i(y,\id_k)\big)\circ_h(z,\id_l)
&= (x\circ_i y,\id_j\circ_i\id_k)\circ_h(z,\id_l)\\
&= ((x\circ_i y)\circ_h z,(\id_j\circ_i\id_k)\circ_h\tau_l)\\
&= (x\circ_i(y\circ_{h-i+1}z),\id_j\circ_i(\id_k\circ_{h-i+1}\tau_l))\\
&= (x,\id_j)\circ_i(y\circ_{h-i+1}z, \id_k\circ_{h-i+1}\tau_l)\\
&= (x,\id_j)\circ_i\big( (y,\id_k)\circ_{h-i+1}(z,\id_l)\big)\quad.
\end{align*}
If $p_i(x)=-1$ and $p_{h-i+1}(y)=1$, then
\begin{align*}
\big( (x,\id_j)\circ_i(y,\id_k)\big)\circ_h(z,\id_l)
&= (x\circ_i y,\id_j\circ_i\tau_k)\circ_h(z,\id_l)\\
&= ((x\circ_i y)\circ_{2i+k-h-1} z,(\id_j\circ_i\tau_k)\circ_h\tau_l)\\
&= (x\circ_i(y\circ_{h-i+1}z), \id_j\circ_i(\tau_k\circ_{h-i+1}\tau_l))\\
&= (x\circ_i(y\circ_{h-i+1}z), \id_j\circ_i(\tau_{k+l-1}(\id_k\circ_{h-i+1}\id_l)))\\
&= (x,\id_j)\circ_i(y\circ_{h-i+1}z, \id_k\circ_{h-i+1}\id_l)\\
&= (x,\id_j)\circ_i\big( (y,\id_k)\circ_{h-i+1}(z,\id_l)\big)\quad.
\end{align*}
If $p_i(x)=-1$ and $p_{h-i+1}(y)=-1$, then
\begin{align*}
\big( (x,\id_j)\circ_i(y,\id_k)\big)\circ_h(z,\id_l)
&= (x\circ_i y,\id_j\circ_i\tau_k)\circ_h(z,\id_l)\\
&= ((x\circ_i y)\circ_{2i+k-h-1} z,(\id_j\circ_i\tau_k)\circ_h\id_l)\\
&= (x\circ_i(y\circ_{h-i+1}z),\id_j\circ_i(\tau_k\circ_{h-i+1}\id_l))\\
&= (x\circ_i(y\circ_{h-i+1}z),\id_j\circ_i(\tau_{k+l-1}(\id_k\circ_{h-i+1}\tau_l))))\\
&= (x,\id_j)\circ_i(y\circ_{h-i+1}z, \id_k\circ_{h-i+1}\tau_l)\\
&= (x,\id_j)\circ_i\big( (y,\id_k)\circ_{h-i+1}(z,\id_l)\big)\quad.
\end{align*}
In addition to the formula we have used proposition~\ref{Prop:circonDn}
and lemma~\ref{Lem:circonM}, and some small calculations. Most notably that
$\tau_k\circ_{h-i+1}\id_l=\tau_{k+l-1}(\id_k\circ_{h-i+1}\tau_l)$.

\paragraph{The general case:}
Now we will use equivariance to get associativity in the case
where the permutations $\rho$, $\upsilon$ and $\mu$ are arbitrary. We have:
\footnotesize
\begin{align*}
((x,\rho)\circ_i (y,\upsilon))\circ_h (z,\mu)
&=\big(((x,\id_j)\circ_{\rho(i)}(y,\id_k))\circ_{(\rho\circ_i\upsilon)(h)}(z,\id_l)\big)
.\big((\rho\circ_i\upsilon)\circ_h\mu\big)\quad,\\
((x,\rho)\circ_h (z,\mu))\circ_{i+l-1}(y,\upsilon)
&=\big(((x,\id_j)\circ_{\rho(h)}(z,\id_l))\circ_{(\rho\circ_h\mu)(i+l-1)}(y,\id_k)\big)
.\big((\rho\circ_h\mu)\circ_{i+l-1}\upsilon\big)
\quad,\\
(x,\rho)\circ_i((y,\upsilon)\circ_{h-i+1}(z,\mu))
&=\big((x,\id_j)\circ_{\rho(i)}((y,\id_k)\circ_{\upsilon(h-i+1)}(z,\id_l))\big)
.\big(\rho\circ_i(\upsilon\circ_{h-i+1}\mu)\big)
\quad\text{and}\\
((x,\rho)\circ_{h-k+1}(z,\mu))\circ_i (y,\upsilon)
&=\big(((x,\id_j)\circ_{\rho(h-k+1)}(z,\id_l))\circ_{(\rho\circ_h\mu)(i)}(y,\id_k)\big)
.\big((\rho\circ_{h-k+1}\mu)\circ_{i}\upsilon\big)\quad.
\end{align*}
\normalsize
Using that $\circ$ is associative for permutations, and that associativity holds
when the permutations are the identity (the special cases), we now get the general result.
\end{proof}

\begin{Lem}
$(1,\id)$ is the unit for $\circ$ on $D_n(j;V)\times \Sigma_j$.
\end{Lem}

\begin{proof}
This is obvious. Simple calculations show that 
for all $(x,\rho)\in D_n(j;W)\times \Sigma_j$ we have
\[
(x,\rho)\circ_i(1,\id_1)=(x,\rho)
\]
and
\[
(1,\id_1)\circ_1(x,\rho)=(x,\rho)\quad.
\]
\end{proof}

We now complete the proof of theorem~\ref{Thm:Dnisoperad}:

\begin{proof}
The main issue is to construct the composition operations, and then to verify
the axioms. We already have defined
\[
\circ_i:(D_n(j;V)\times \Sigma_j)_+\wedge (D_n(k;W)\times \Sigma_k)_+
\rightarrow (D_n(j+k-1;V\oplus W)\times \Sigma_{j+k-1})_+
\]
by the formula
\[
(x,\rho)\circ_i(y,\upsilon)=
(x\circ_{\rho(i)}y,\rho\circ_i(\tau_k({p_{\rho(i)}(x)})\upsilon))\quad.
\]
Recall that $\mathcal{D}_n(j)(V)$ is defined as $(D_n(j;V)\times \Sigma_j)_+\wedge S^V$. Let 
$\tilde{\circ}_i:\mathcal{D}_n(j)(V)\wedge\mathcal{D}_n(k)(W)\rightarrow\mathcal{D}_n(j+k-1)(V\oplus W)$ 
be induced by these. We think of $\tilde{\circ}_i$ as exterior composition operations.
And we check that they are natural transformations that coequalizers
\[
\mathcal{D}_n(j)\tilde{\wedge} S\tilde{\wedge}\mathcal{D}_n(k)
\rightrightarrows \mathcal{D}_n(j)\tilde{\wedge}\mathcal{D}_n(k)\quad.
\]
Therefore we have induced maps
\[
\circ_i:\mathcal{D}_n(j)\wedge\mathcal{D}_n(k)\rightarrow\mathcal{D}_n(j+k-1)\quad.
\]

Associativity for $\circ$ on $\mathcal{D}_n$ follows from associativity
for $\circ$ on $j\mapsto D_n(j;-)\times\Sigma_j$, likewise for equivariance
when we let the right action of $\Sigma_j$ on $\mathcal{D}_n(j)$ be induced
from the right action of $\Sigma_j$ on $D_n(j;-)\times \Sigma_j$. The unity 
axiom follows similarly.

To see that
the collection of maps $\mathcal{D}_n(j)\rightarrow \mathcal{H}(j)$ defines a map of operads
one compares the formula above and the formula given in remark~\ref{Rem:AltHdescr}.
Recall that $p$ denotes the group homomorphism $D_n(j;V)\rightarrow (\Z/2)^{j}$,
and notice that $p(x\circ_i y)=p(x)\circ_i p(y)$ by proposition~\ref{Prop:pandcirc}.
\end{proof}

We round up this subsection by providing an orbit cofibrant replacement for $\mathcal{D}_n$:

\begin{Thm}\label{Thm:ocreplforDn}
There exists an operad in orthogonal spectra, which we denote by $\hat{\Gamma}\mathcal{D}_n$,
such that
\begin{itemize}
\item[] there is a map of operads $\hat{\Gamma}\mathcal{D}_n\rightarrow\mathcal{D}_n$, 
\item[] for each $j$ the map $\hat{\Gamma}\mathcal{D}_n(j)\rightarrow\mathcal{D}_n(j)$ is a level-equivalence, and
\item[] each $\hat{\Gamma}\mathcal{D}_n(j)$ can be described $\Sigma_j$-equivariantly as
a product $X\wedge (\Sigma_j)_+$, where $X$ is non-equivariant and orbit cofibrant. 
\end{itemize}
\end{Thm}

\begin{proof}
We have a non-$\Sigma$ version $\mathcal{D}'_n$ of the operad $\mathcal{D}_n$ given by
\[
\mathcal{D}'_n(j)(V)=D_n(j;V)_+\wedge S^V\quad.
\]
In analogy to what we did in the proof above, 
we have non-$\Sigma$ composition operations $\circ_i$ for $\mathcal{D}'_n$ induced by the $\circ_i$'s on $D_n(-;V)_+$.
Observe that $\mathcal{D}_n(j)\cong\mathcal{D}'_n(j)\wedge (\Sigma_j)_+$.
The idea is now to apply the orbit cofibrant replacement functor $\tilde{\Gamma}$
to $\mathcal{D}'_n$. For the definition and the properties of $\tilde{\Gamma}$, see theorem~\ref{Thm:orbitcofrepl}.

We now define the $\circ_i$'s for $\tilde{\Gamma}\mathcal{D}'_n$ as compositions
\[
\tilde{\Gamma}\mathcal{D}'_n(j)\wedge\tilde{\Gamma}\mathcal{D}'_n(k)\xrightarrow{\phi}
\tilde{\Gamma}\left(\mathcal{D}'_n(j)\wedge \mathcal{D}'_n(k)\right)\xrightarrow{\tilde{\Gamma}(\circ_i)}
\tilde{\Gamma}\mathcal{D}'_n(j+k-1)\quad.
\]
Since $\tilde{\Gamma}$ is symmetric, these $\circ_i$'s for $\tilde{\Gamma}\mathcal{D}'_n$
will satisfy associativity relations analogous to those given in proposition~\ref{Prop:circonDn}.

Now define $\hat{\Gamma}\mathcal{D}_n$ by
\[
\hat{\Gamma}\mathcal{D}_n(j)=\tilde{\Gamma}\mathcal{D}'_n(j)\wedge (\Sigma_j)_+\quad.
\]
We have right $\Sigma_j$ actions as usual. And the $\circ_i$'s on $\hat{\Gamma}\mathcal{D}_n$
are defined by the same formula as before. Lemma~\ref{Lem:circequivariance} is formal and the argument yields
that $\circ_i$'s on $\hat{\Gamma}\mathcal{D}_n$ are equivariant. Furthermore, the argument of 
lemma~\ref{Lem:circassociativity} is also formal, thus $\circ_i$'s on $\hat{\Gamma}\mathcal{D}_n$
are also associative. Hence $\hat{\Gamma}\mathcal{D}_n$ is an operad in orthogonal spectra.

The natural level-equivalence $\tilde{\Gamma}L\rightarrow L$ from theorem~\ref{Thm:orbitcofrepl},
induces the map of operads $\hat{\Gamma}\mathcal{D}_n\rightarrow\mathcal{D}_n$.
Clearly for each $j$ the map $\hat{\Gamma}\mathcal{D}_n(j)\rightarrow\mathcal{D}_n(j)$ is a level-equivalence.
And since $\tilde{\Gamma}L$ is orbit cofibrant for any $L$ the last statement follows.
\end{proof}

\begin{Rem}\label{Rem:MtGDn}
By the construction of $\hat{\Gamma}\mathcal{D}_n$ it is easily seen that there exists a map
of operads $f:\mathcal{M}\rightarrow \hat{\Gamma}\mathcal{D}_n$
such that the composition 
\[
\mathcal{M}\xrightarrow{f} \hat{\Gamma}\mathcal{D}_n\rightarrow \mathcal{D}_n\rightarrow\mathcal{H}
\]
is the standard inclusion.
The identity element in $D_n(j;V)$ gives an inclusion
\[
S\rightarrow \mathcal{D}'_n(j)
\]
for every $j$. Recall that $\tilde{\Gamma}$ comes with a unit map $S\rightarrow\tilde{\Gamma}S$.
So we get a map
\[
S\rightarrow \tilde{\Gamma}\mathcal{D}'_n(j)\quad.
\]
Smashing with $(\Sigma_j)_+$ yields 
\[
f:\mathcal{M}(j)\rightarrow \hat{\Gamma}\mathcal{D}_n(j)\quad.
\]
\end{Rem}

\subsection{The action of $\mathcal{D}_n$ on $S[\Omega M]$}\label{subsect:action}

In this subsection we will construct a $\mathcal{D}_n$-algebra structure on the 
orthogonal spectrum $S[\Omega M]$. This structure depends on an $n$-vector bundle $\xi$ over $M$.
Recall that by definition of $S[\Omega M]$, the $V$'th space is $(\Omega M)_+\wedge S^V$.
In this subsection we will prove the following theorem:

\begin{Thm}\label{Thm:DnalgstronSOmegaM}
Let $M$ be a compact smooth manifold and $\xi$ an $n$-vector bundle
over $M$, then the orthogonal spectrum $S[\Omega M]$
has a $\mathcal{D}_n$-algebra structure which depends on $\xi$.
\end{Thm}

Before proving this theorem there are some preliminary considerations and constructions.
First we should agree on a suitable model for the loop space $\Omega M$.
See~\cite{AdamsHilton:56} or subsection~5.1 in~\cite{CarlssonMilgram:95}
for the definition of ``Moore loops''. We modify this definition slightly to get piecewise
smooth ``Moore loops''.

Let $m_0$ be a base point in $M$.
For technical reasons it is important to have an associative multiplication (composition of loops)
and that every loop is piecewise smooth.
Here we define such a space $\Omega M$ as the geometrical realization
of a simplicial monoid. A $q$-simplex is a piecewise smooth map
\[
\gamma:\Delta^q\times I\rightarrow M
\]
together with a piecewise affine map 
\[
l:\Delta^q\rightarrow \left[ 0, \infty\right)
\]
such that $\gamma(\mb{t},0)=\gamma(\mb{t},1)=m_0$ for all $\mb{t}\in\Delta^q$ and whenever
$l(\mb{t})=0$ then $\gamma(\mb{t},s)=m_0$ for all $s\in I$. 
Here $l(\mb{t})$ is thought of as the ``length'' of the loop $s\mapsto \gamma(\mb{t},s)$. 

If we have two $q$-simplices $(\gamma_1,l_1)$ and $(\gamma_2,l_2)$ we multiply (compose) these
as follows: Let $l=l_1+l_2$ and define $\gamma$ by
\[
\gamma(\mb{t},s)=\begin{cases}
\gamma_1(\mb{t},\frac{s(l_1+l_2)}{l_1})&\text{if $s(l_1+l_2)<l_1$,}\\
m_0&\text{if $s(l_1+l_2)=l_1$, and}\\
\gamma_2(\mb{t},\frac{s(l_1+l_2)-l_1}{l_2})&\text{if $s(l_1+l_2)>l_1$.}
\end{cases}
\]
Here we have divided the interval $I$ into two pieces, the ratio between their lengths being $l_1(\mb{t})$
to $l_2(\mb{t})$. On the first piece we use $\gamma_1$ and on the second we use $\gamma_2$.
Associativity of the composition follows. We will use the notation $\centerdot$ for this operation.

Notice that $\Omega M$ has the correct homotopy type. Let $\Omega' M$ be the geometrical realization
of the simplicial set having $q$-simplices the piecewise smooth maps $\gamma:\Delta^q\times I\rightarrow M$
such that $\gamma(\mb{t},0)=\gamma(\mb{t},1)=m_0$. Then we can compare this space with $\Omega M$.
There is an inclusion 
\[
i:\Omega' M\rightarrow \Omega M
\]
defined by setting $l$ constant equal to $1$. And we have a retraction
\[
r:\Omega M\rightarrow \Omega' M
\]
by forgetting $l$. Clearly $ri=\id$. And it is possible to construct a simplicial homotopy $ir\simeq\id$.
Therefore $\Omega M\simeq \Omega'M$. Furthermore, $\Omega' M$ is homotopic to the space of continuous maps
$(I,\{0,1\})\rightarrow (M,m_0)$, see chapter~17 in~\cite{Milnor:63}.

There is an involution on $\Omega M$. We write $(\gamma,l)\mapsto \overline{(\gamma,l)}$, and it is defined by
sending $\gamma$ to the reversed loop, 
\[
\overline{\gamma}(\mb{t},s)=\gamma(\mb{t},1-s)\quad,
\]
while leaving $l$ unchanged. Notice that the involution is an anti-homomorphism. This means that
\[
\overline{ (\gamma_1,l_1)\centerdot (\gamma_2,l_2)}=\overline{ (\gamma_2,l_2)}\centerdot\overline{ (\gamma_1,l_1)}\quad.
\]

We will often simplify the notation for a loop in $\Omega M$, and leave the length $l$
out of the notation.

The construction of the $\mathcal{D}_n$-algebra structure on $S[\Omega M]$ will use a connection $\nabla$
on $\xi$. We need to take parallel transportation along piecewise smooth loops in $M$.
However, the choice of connection will not carry any information up to homotopy.

Choose a connection $\nabla$ on $\xi$ and an isomorphism $\R^n\cong \xi_{m_0}$.
Then parallel transportation yields a continuous map
\[
P:\Omega M\rightarrow GL(\R^n)
\]
such that
\begin{itemize}
\item[] $P(\gamma_1\centerdot\gamma_2)=P(\gamma_2)P(\gamma_1)$ for all piecewise smooth loops $\gamma_1$ and $\gamma_2$, and
\item[] $P(\overline{\gamma})=P(\gamma)^{-1}$ for all piecewise smooth loops $\gamma$.
\end{itemize}
For more about parallel transportation see remark~17.4 in~\cite{MadsenTornehave:97}.

Given a finite dimensional real inner product space $V$ and an isometric embedding $\phi:\R^n\rightarrow V$,
we write $V$ as the sum $\phi(\R^n) + V^{\perp}$, where $V^{\perp}$ is the orthogonal complement of
$\phi(\R^n)$ in $V$. Given a piecewise smooth loop $\gamma$ in $M$, we define the map $\phi^*(\gamma):V\rightarrow V$
by using $P(\gamma)$ on $\phi(\R^n)$ while leaving $V^{\perp}$ unchanged. For $v=\phi(u)+w$, $u\in \R^n$ and $w\in V^{\perp}$
we have
\[
\phi^*(\gamma)(v)=\phi(P(\gamma)(u))+w\quad.
\]
Let $\phi$ and $\psi$ be isometric embeddings of $\R^n$ in $V$, and $\gamma$, $\gamma_1$ and $\gamma_2$
be piecewise smooth loops in $M$, then:
\begin{itemize}
\item[] $\phi^*(\gamma_1\centerdot\gamma_2)=\phi^*(\gamma_2)\phi^*(\gamma_1)$,
\item[] $\phi^*(\overline{\gamma})=\phi^*(\gamma)^{-1}$, and
\item[] $\phi^*(\gamma_1)\psi^*(\gamma_2)=\psi^*(\gamma_2)\phi^*(\gamma_1)$ if $\phi\perp\psi$.
\end{itemize}
Notice that the map $(\phi,\gamma)\mapsto \phi^*(\gamma)$ is continuous when $\phi$ lies in the
space of isometric embeddings, $\gamma$ lies in $\Omega M$ and the image lies in $GL(V)$.

Next we define an action of the group $D_n(j;V)$ on the space
$F(S^V,S^V)\wedge (\Omega M_+)^{\wedge j}$. Recall that $F(X,Y)$ denotes the space of based maps
$X\rightarrow Y$. Let $(\phi,r)$ be a generator of $D_n(j;V)$ and $(f;\gamma_1,\ldots,\gamma_j)$
a point in $F(S^V,S^V)\wedge(\Omega M_+)^{\wedge j}$.
Then we define the action by the formula
\[
(\phi,r).(f;\gamma_1,\ldots,\gamma_j)=
(\phi^*(\gamma_r)\circ f;\gamma_1,\ldots,\gamma_{r-1},\overline{\gamma_r},\gamma_{r+1},\ldots,\gamma_j)\quad.
\]
It is easily seen that
\[
(\phi,r).\big((\phi,r).(f;\gamma_1,\ldots,\gamma_j)\big)=(f;\gamma_1,\ldots,\gamma_j)
\]
and if $\phi\perp\psi$ and $r\neq r'$, then
\[
(\phi,r).\big((\psi,r').(f;\gamma_1,\ldots,\gamma_j)\big)
=(\psi,r').\big((\phi,r).(f;\gamma_1,\ldots,\gamma_j)\big)\quad.
\]
Thus we have a well defined group action.

Also $\Sigma_j$ acts from the left on the space $F(S^V,S^V)\wedge(\Omega M_+)^{\wedge j}$.
This action is by permutation of the loops. For $\rho\in\Sigma_j$ we have
\[
\rho.(f;\gamma_1,\ldots,\gamma_j)=(f;\gamma_{\rho^{-1}(1)},\ldots,\gamma_{\rho^{-1}(j)})\quad.
\]

Now define $\bar{\theta}_j:\mathcal{D}_n(j)\wedge(\Omega M_+)^{\wedge j}\rightarrow S[\Omega M]$
by commutativity of the following diagram:
\[
\begin{CD}
S^V\wedge\left(D_n(j;V)\times\Sigma_j \right)_+\wedge F(S^V,S^V)\wedge(\Omega M_+)^{\wedge j}
&\rightarrow& S^V\wedge F(S^V,S^V)\wedge(\Omega M_+)^{\wedge j}\\
@AAA @VVV\\
\left(D_n(j;V)\times \Sigma_j \right)_+\wedge S^V \wedge(\Omega M_+)^{\wedge j}
&\rightarrow& \Omega M_+\wedge S^V
\end{CD}\quad.
\]
The top map combines the group actions, first apply the $\Sigma_j$-action, then the $D_n(j;V)$-action.
The left map is the inclusion at $\id\in F(S^V,S^V)$ and the right map evaluates $f\in F(S^V,S^V)$ on $S^V$
and multiplies (composes) the loops. The map at the bottom is $\bar{\theta}_j$ evaluated at $V$.
Clearly $\bar{\theta}_j$ commutes with assemblies for $\mathcal{D}_n(j)\wedge(\Omega M_+)^{\wedge j}$
and $S[\Omega M]$, and is thus a well defined map of orthogonal spectra.

Via a series of adjunctions there is for orthogonal spectra $L$ and $K$ and a based space $A$,
a one-to-one correspondence between maps $L\wedge A\rightarrow K$ and maps
$L\wedge F_0A\rightarrow K$. The adjunctions are:
\begin{multline*}
\mathscr{IS}(L\wedge A,Y)\cong\Top_*(A,\mathscr{IS}(L,K))\\
\cong\Top_*(A,F(L,K)(0))\cong\mathscr{IS}(F_0A,F(L,K))\cong\mathscr{IS}(L\wedge F_0A,K)\quad.
\end{multline*}
Applied to $\bar{\theta}_j$ we get our map
\[
\theta_j:\mathcal{D}_n(j)\wedge S[\Omega M]^{\wedge j}\rightarrow S[\Omega M]\quad.
\]

Alternatively, it is possible to give a more explicit description of $\bar{\theta}_j$. 
Let $x=(\phi_1,r_1)\ldots(\phi_s,r_s)$ be a point in $D_n(j;V)$, $\rho\in\Sigma_j$, $v\in V$
and $(\gamma_1,\ldots,\gamma_j)$ loops in $\Omega M^{j}$. We now want to give a formula
for $\bar{\theta}_j(x,\rho,v;\gamma_1,\ldots,\gamma_j)$.
The permutation $\rho$ permutes the loops, and each $(\phi_t,r_t)$ reverses the loop at the $r_t$'th position.
Therefore we define $\delta_i$ to be the loop given by
\[
\delta_i=\begin{cases}
\gamma_{\rho^{-1}(i)} &\text{if $p_i(x)=1$, and}\\
\overline{\gamma_{\rho^{-1}(i)}} &\text{if $p_i(x)=-1$.}
\end{cases}
\]
Each $(\phi_t,r_t)$ also changes the vector $v$ in $V$ by parallel transportation along the loop at
the $r_t$'th position. But notice that the direction along the loop
in which one should perform the parallel transportation, depends on the number of occurrences of the number $r_t$
among $r_{t+1},\ldots,r_s$. Therefore we define the sign $\epsilon_t$ by
\[
\epsilon_t= p_{r_t}((\phi_{t+1},r_{t+1})\ldots(\phi_s,r_s))\quad.
\]
Calculating, we get that
\[
\bar{\theta}_j(x,\rho,v;\gamma_1,\ldots,\gamma_j)=(\delta_1\centerdot\cdots\centerdot\delta_j,
\phi_1^*(\gamma_{\rho^{-1}(r_1)})^{\epsilon_1}\cdots\phi_s^*(\gamma_{\rho^{-1}(r_s)})^{\epsilon_s}(v))\quad.
\]

\begin{Rem}\label{Rem:identifyinvolution}
To identify the involution we should pay special attention to the case where $j=1$. 
Inspect the map $\bar{\theta}_1$ at level $V=\R^n$, and at the
points in $\mathcal{D}_n(1)$ given by $x=(\id_{\R^n},1)$ and $\rho=\id$.
We send $(v,\gamma)$ to
\[
\bar{\theta}_1(x,\rho,v;\gamma)=(\bar{\gamma},(\id_{\R^n})^*(\gamma)(v))=(\bar{\gamma},P(\gamma)(v))\quad.
\]
Recall that $P(\gamma)$ is the parallel transport in $\xi$ along $\gamma$.
The resulting map
\[
S^n\wedge \Omega M_+\rightarrow \Omega M_+\wedge S^n
\]
is precisely the map $\bar{P}$, which we use to define the involution $\iota$
on $\pi_*S[\Omega M]$. See definition~\ref{Def:iota}.
\end{Rem}

We complete the proof of theorem~\ref{Thm:DnalgstronSOmegaM} by showing that the maps $\theta_j$
for $j\geq 0$, is a $\mathcal{D}_n$-algebra structure on $S[\Omega M]$.

\begin{proof}
We have to check the axioms given in definition~\ref{Def:Palg}. Triviality of the unit and equivariance
are easily seen to hold. It remains to show that $\theta$ acts.

Let $Z$ be the space $F(S^{V\oplus W},S^{V\oplus W})\wedge (\Omega M_+)^{\wedge j+k-1}$. A point $z$ in $Z$
can be written as $(f;\gamma_1,\ldots,\gamma_{j+k-1})$, where $f$ is an endomorphism of $S^{V\oplus W}$ and
the $\gamma$'s are loops in $M$. The main ingredient of the proof will be to define several group actions 
on $Z$, understand how these interact with each other and how the actions relate to $\theta$ and $\circ$.

By $\mathbb{A}$ we will denote the action of $D_n(j+k-1;V\oplus W)$ on $Z$ defined above. Recall that the
formula on generators is
\[
\mathbb{A}\left( (\phi,r),(f;\gamma_1,\ldots,\gamma_{j+k-1})\right)=
(\phi^*(\gamma_r)\circ f;\gamma_1,\ldots,\gamma_{r-1},\overline{\gamma_r},\gamma_{r+1},\ldots,\gamma_{j+k-1})\quad.
\]

The action of $\Sigma_{j+k-1}$ on $Z$ will in this proof be denoted by $\mathbb{B}$ and is given by
\[
\mathbb{B}\left(\rho,(f;\gamma_1,\ldots,\gamma_{j+k-1})\right)=
(f;\gamma_{\rho^{-1}(1)},\ldots,\gamma_{\rho^{-1}(j+k-1)})\quad.
\]

Depending on $i$ there are actions $A_i$ of $D_n(j;V)$ on $Z$. We define the action $A_i$ of the generator $(\phi,r)$
on $z=(f;\gamma_1,\ldots,\gamma_{j+k-1})$ by the formula
\footnotesize
\[
A_i((\phi,r),z)=\begin{cases}
(\phi^*(\gamma_{r})\circ f;\gamma_{1},\ldots,\gamma_{r-1},\overline{\gamma_{r}},\gamma_{r+1},\ldots,\gamma_{j+k-1})&
\text{if $r<i$,}\\
(\phi^*(\delta)\circ f;
\gamma_{1},\ldots,\gamma_{i-1},\overline{\gamma_{i+k-1}},\ldots,\overline{\gamma_{i}},\gamma_{i+k},\ldots,\gamma_{j+k-1})&
\text{if $r=i$, and}\\
(\phi^*(\gamma_{r+k-1})\circ f;\gamma_{1},\ldots,\gamma_{r+k-2},\overline{\gamma_{r+k-1}},\gamma_{r+k},\ldots,\gamma_{j+k-1})&
\text{if $r>i$.}
\end{cases}
\]
\normalsize
Here $\delta$ is the composition $\gamma_{i}\centerdot\cdots\centerdot\gamma_{i+k-1}$.
We have implicitly changed the target of the isometric embedding $\phi$ to be $V\oplus W$ 
via the canonical map $V\rightarrow V\oplus W$.

Also depending on $i$ there are actions $B_i$ of $\Sigma_j$ on $Z$. Let $z=(f;\gamma_1,\ldots,\gamma_{j+k-1})$
be a point in $Z$.
The action $B_i$ of $\rho\in\Sigma_j$ is given by putting boxes around the $\gamma$'s as follows:
\[
\boxed{\gamma_1},\ldots,\boxed{\gamma_{i-1}},\boxed{\gamma_i, \gamma_{i+1},\ldots, \gamma_{i+k-1}},\boxed{\gamma_{i+n}},
\ldots,\boxed{\gamma_{j+k-1}}\quad.
\]
And we use $\rho$ to permute the boxes. The action leaves $f$ unchanged. The result is called $B_i(\rho,z)$.

Define the action $\alpha_i$ of $D_n(k;W)$ on $Z$ by the formula:
\[
\alpha_i((\phi,r),z)=(\phi^*(\gamma_{r+i-1})\circ f;\gamma_{1},\ldots,\gamma_{r+i-2},\overline{\gamma_{r+i-1}},
\gamma_{r+i},\ldots,\gamma_{j+k-1})\quad.
\]
Here we understand the target of $\phi$ to be $V\oplus W$ via the canonical map $W\rightarrow V\oplus W$.

The action $\beta_i$ of $\Sigma_k$ on $Z$ is given by permuting the loops $\gamma_{i},\ldots,\gamma_{i+k-1}$.
For $\upsilon\in\Sigma_k$ we have
\[
\beta_i(\upsilon,z)=(f;\gamma_1,\ldots,\gamma_{i-1},\gamma_{\upsilon^{-1}(1)+i-1},
\ldots,\gamma_{\upsilon^{-1}(k)+i-1},\gamma_{i+k},\ldots,\gamma_{j+k-1})\quad.
\]

Recall also the definition of $\circ_i:D_n(j;V)\times D_n(k;W)\rightarrow D_n(j+k-1;V\oplus W)$.
We had homomorphisms $c_i:D_n(k;W)\rightarrow D_n(j+k-1;V\oplus W)$ given by $c_i(\phi,r)=(\phi,r+i-1)$,
and actions $\vdash_i$ of $D_n(j;V)$ on $D_n(j+k-1;V\oplus W)$.
For a generator $(\phi,r)$ in $D_n(j;V)$ and an element $y\in D_n(j+k-1;V\oplus W)$, $\vdash_i$ is given by
\[
(\phi,r)\vdash_i y=\begin{cases}
(\phi,r)y&\text{for $r<i$,}\\
(\phi,i+k-1)\cdots(\phi,i)\bar{y}&\text{for $r=i$, and}\\
(\phi,r+k-1)y&\text{for $r>i$,}
\end{cases}
\]
where $y\mapsto\bar{y}$ is an automorphism of $D_n(j+k-1;V\oplus W)$
defined on generators by
\[
(\phi,r)\mapsto\begin{cases}
(\phi,r) &\text{for $r<i$,}\\
(\phi,k+2i-r-1)&\text{for $i\leq r< k+i$, and}\\
(\phi,r) &\text{for $k+i\leq r$.}
\end{cases}
\]
Now if $x\in D_n(j;V)$ and $y\in D_n(k;W)$, then
\[
x\circ_i y= x\vdash_i c_i(y)\quad.
\]
Recall that $\circ_i$ on $\left(D_n(j;V)\times\Sigma_j\right)\times\left(D_n(k;W)\times \Sigma_k\right)$
is defined by the formula
\[
(x,\rho)\circ_i(y,\upsilon)=
(x\circ_{\rho(i)}y,\rho\circ_i(\tau_k({p_{\rho(i)}(x)})\upsilon))\quad.
\]

By $\mathbb{E}$ we will denote the map
\[
S^{V\oplus W}\wedge F(S^{V\oplus W},S^{V\oplus W})\wedge(\Omega M_+)^{\wedge(j+k-1)}\rightarrow \Omega M_+\wedge S^{V\oplus W}
\]
given by evaluating and composing.

Let $\tilde{Z}$ be the space $F(S^V,S^V)\wedge F(S^W,S^W)\wedge (\Omega M_+)^{\wedge (j+k-1)}$. 
There is a natural map $\tilde{Z}\rightarrow Z$ given by taking the smash product of $f_1:S^V\rightarrow S^V$
and $f_2:S^W\rightarrow S^W$. Observe that
the actions $A_i$, $\alpha_i$, $B_i$ and $\beta_i$ lift to actions $\tilde{A}_i$, $\tilde{\alpha}_i$, $\tilde{B}_i$ 
and $\tilde{\beta}_i$ on $\tilde{Z}$.

We have designed $\tilde{\alpha_i}$ and $\tilde{\beta_i}$ such that they 
correspond to smashing the
actions on $F(S^W,S^W)\wedge(\Omega M_+)^{\wedge k}$ in the definition of $\bar{\theta}_k$
with the trivial actions on $F(S^V,S^V)\wedge(\Omega M_+)^{\wedge k-1}$.
Up to shuffling the factor $F(S^W,S^W)$, we have that
\[
(\id_{S^V};\gamma_1,\ldots,\gamma_{i-1},y.(\id_{S^W};\gamma_i,\ldots,\gamma_{i+k-1}),\gamma_{i+k},\ldots,\gamma_{j+k-1})
\]
is equal to
\[
\tilde{\alpha}_i\big(y,\left(\id_{S^V},\id_{S^W};\gamma_1,\ldots,\gamma_{j+k-1}\right)\big)\quad,
\]
and similar for the $\Sigma_k$-action $\beta_i$.

The actions $\tilde{A}_i$ and $\tilde{B}_i$ see the loops $\gamma_i,\ldots,\gamma_{i+k}$ as one composed
loop. If we by 
\[
e_i:F(S^V,S^V)\wedge F(S^W,S^W)\wedge (\Omega M_+)^{\wedge (j+k-1)}\rightarrow F(S^V,S^V)\wedge (\Omega M_+)^{\wedge j}\wedge S^W
\]
denote the map given by evaluating $f_2:S^W\rightarrow S^W$ on $S^W$ and composing the loops $\gamma_{i},\ldots,\gamma_{i+k-1}$,
then we observe that
\begin{itemize}
\item[] $e_i$ is $D_n(j;V)$-equivariant (the action on the target is given by smashing the $D_n(j;V)$-action in the definition of
$\bar{\theta}_j$ with the trivial action on $S^W$), and
\item[] $\rho.e_i(f_1,f_2;\gamma_1,\ldots,\gamma_{j+k-1})
=e_{\rho(i)} \tilde{B}_i\big(\rho,(f_1,f_2;\gamma_1,\ldots,\gamma_{j+k-1})\big)$ 
(as usual $\rho\in \Sigma_j$ acts on the target of $e_i$ by permuting the loops).
\end{itemize}

Let $(x,\rho,v)$ be a point in $(D_n(j;V)\times\Sigma_j)_+\wedge S^V=\mathcal{D}_n(j)(V)$,
$(y,\upsilon,w)$ a point in $(D_n(k;W)\times\Sigma_k)_+\wedge S^W=\mathcal{D}_n(k)(W)$
and $(\gamma_1,\ldots,\gamma_{j+k-1})$ loops in $(\Omega M_+)^{\wedge (j+k-1)}$.
Let $z$ in $Z$ be the point $(\id_{S^{V\oplus W}};\gamma_1,\ldots,\gamma_{j+k-1})$.

By definition of the $\bar{\theta}$'s we see that the composition
\[
\begin{CD}
\mathcal{D}_n(j)(V)\wedge \mathcal{D}_n(k)(W)\wedge (\Omega M_+)^{\wedge (j+k-1)}\\
@V{\text{shuffle}}VV\\
\mathcal{D}_n(j)(V)\wedge (\Omega M_+)^{\wedge (i-1)}\wedge 
\mathcal{D}_n(k)(W)\wedge (\Omega M_+)^{\wedge j}\wedge (\Omega M_+)^{\wedge (j-i)}\\
@V{\id\wedge\bar{\theta}_k\wedge\id}VV\\
\mathcal{D}_n(j)(V)\wedge (\Omega M_+)^{\wedge j}\wedge S^W \\
@V{\bar{\theta}_j\wedge \id_{S^W}}VV \\
\Omega M_+\wedge S^{V\oplus W}
\end{CD}
\]
evaluated at this point is equal to
\[
\mathbb{E}\left( (v,w), A_{\rho(i)}(x, B_i(\rho,\alpha_i(y,\beta_i(\upsilon,z))))\right)\quad,
\]
and by the formula for $\circ_i$ and definition of $\bar{\theta}_{j+k-1}$ the composition
\[
\begin{CD}
\mathcal{D}_n(j)(V)\wedge \mathcal{D}_n(k)(W)\wedge (\Omega M_+)^{\wedge (j+k-1)}\\
@VV{\circ_i}V\\
\mathcal{D}_n(j+k-1)(V\oplus W)\wedge (\Omega M_+)^{\wedge (j+k-1)}\\
@VV{\bar{\theta}_{j+k-1}}V\\
\Omega M_+\wedge S^{V\oplus W}
\end{CD}
\]
evaluated at the same point is
\[
\mathbb{E}\left((v,w),\mathbb{A}(x\circ_{\rho(i)}y,\mathbb{B}(\rho\circ_i(\tau_k({p_{\rho(i)}(x)})\upsilon),z))\right)\quad.
\]
To finish the proof it is enough to show that
\[
A_{\rho(i)}(x, B_i(\rho,\alpha_i(y,\beta_i(\upsilon,z))))=
\mathbb{A}(x\circ_{\rho(i)}y,\mathbb{B}(\rho\circ_i(\tau_k({p_{\rho(i)}(x)})\upsilon),z))\quad.
\]
And if we can check the following formulas, then we are done:
\begin{itemize}
\item[i)] $B_i(\rho,\alpha_i(y,z))=\alpha_{\rho(i)}(y,B_i(\rho,z))$,
\item[ii)] $\mathbb{B}(\rho\circ_i\upsilon,z)=B_i(\rho,\beta_i(\upsilon,z))$,
\item[iii)] $\mathbb{A}(x\circ_i y,\mathbb{B}(\id_j\circ_i\tau_k({p_i(x)}),z))=A_i(x,\alpha_i(y,z))$, and
\item[iv)] $\rho\circ_i(\tau_k\upsilon)=(\id_j\circ_{\rho(i)}\tau_k)(\rho\circ_i\upsilon)$.
\end{itemize}
Here iv) is a special case of lemma~\ref{Lem:Mformulas}v), while i) and ii) follow directly from the definitions.
We prove formula~iii) by induction on the length of $x$. Assume first that $x=1$. In this case we have to show that
\[
\mathbb{A}(1\circ_i y,z)=\alpha_i(y,z)\quad.
\]
Using that $1\circ_i y=c_i(y)$ and checking the definitions we see that this formula holds.

Assume that $x=(\phi,r)x'$. We consider the case when $r\neq i$.
Let $d$ be a point in $D_n(j+k-1;V\oplus W)$ and $z'\in Z$. 
By the definitions of $\vdash_i$, $\mathbb{A}$ and $A_i$ we have
\[
\mathbb{A}\big( (\phi,r)\vdash_i d,z'\big)=A_i\big( (\phi,r),\mathbb{A}(d,z')\big)\quad.
\]
Setting $d=x'\circ_i y$ and $z'=\mathbb{B}(\id_j\circ_i\tau_k({p_i(x)}),z))$ we get
\[
\mathbb{A}\big( x\circ_i y ,\mathbb{B}(\id_j\circ_i\tau_k({p_i(x)}),z)) \big)=
A_i\big( (\phi,r),\mathbb{A}(x'\circ_i y,\mathbb{B}(\id_j\circ_i\tau_k({p_i(x)}),z)))\big)\quad.
\]
Notice that $p_i(x)=p_i(x')$. By induction we have that 
$\mathbb{A}(x'\circ_i y,\mathbb{B}(\id_j\circ_i\tau_k({p_i(x')}),z))=A_i(x',\alpha_i(y,z))$.
This implies
\[
\mathbb{A}(x\circ_i y,\mathbb{B}(\id_j\circ_i\tau_k({p_i(x)}),z))=A_i(x,\alpha_i(y,z))\quad.
\]

At last we consider the case when $x=(\phi,r)x'$ and $r=i$.
It is sufficient to show that
\[
\mathbb{A}\big( (\phi,i)\vdash_i d,\mathbb{B}(\id_j\circ_i\tau_k,z')\big)=A_i\big( (\phi,i),\mathbb{A}(d,z')\big)\quad,
\]
because setting $d=x'\circ_i y$ and $z'=\mathbb{B}(\id_j\circ_i\tau_k({p_i(x')}),z))$ and using induction
yields the formula for $x$.
By definition of $\vdash_i$ we have
\[
(\phi,i)\vdash_i d=(\phi,i+k-1)\cdots(\phi,i)\bar{d}\quad.
\]
Furthermore, we check the following formulas directly:
\begin{align*}
\mathbb{A}(\bar{d},\mathbb{B}(\id_j\circ_i\tau_k,z')) &= \mathbb{B}(\id_j\circ_i\tau_k,\mathbb{A}(d,z'))\quad,\text{ and}\\
\mathbb{A}\big( (\phi,i+k-1)\cdots(\phi,i),\mathbb{B}(\id_j\circ_i\tau_k,z'')\big) &= A_i((\phi,i),z'')\quad.
\end{align*}
Now only a simple calculation remains:
\begin{multline*}
\mathbb{A}\big( (\phi,i)\vdash_i d,\mathbb{B}(\id_j\circ_i\tau_k,z')\big)\\
= \mathbb{A}\big( (\phi,i+k-1)\cdots(\phi,i)\bar{d},\mathbb{B}(\id_j\circ_i\tau_k,z')\big)\\
=\mathbb{A}\big( (\phi,i+k-1)\cdots(\phi,i),\mathbb{A}(\bar{d},\mathbb{B}(\id_j\circ_i\tau_k,z'))\big)\\
=\mathbb{A}\big( (\phi,i+k-1)\cdots(\phi,i),\mathbb{B}(\id_j\circ_i\tau_k,\mathbb{A}(d,z'))\big)\\
=A_i((\phi,i),\mathbb{A}(d,z'))\quad.
\end{multline*}
\end{proof}

\subsection{Homotopy discreteness of $\mathcal{D}_n$}\label{subsect:homotopydiscreteness}

In this subsection we will compare $\mathcal{D}_n$ to $\mathcal{H}$. 
Recall that there is a map of operads $\mathcal{D}_n\rightarrow\mathcal{H}$.
This map comes from a homomorphism of groups, $p:D_n(j;V)\rightarrow (\Z/2)^j$,
whose $i$'th factor is the ``parity'' of the number of ``letters'' of the form $(-,i)$ in a word
$x\in D_n(j;V)$. We get our map of operads 
by applying the functor $(-\times \Sigma_j)_+\wedge S^V$ to $p$.

Our theorem says: 

\begin{Thm}\label{Thm:Ddiscrete}
For each $j$ the map $\mathcal{D}_n(j)\rightarrow \mathcal{H}(j)$ 
can be written equivariantly
as a product $X\wedge(\Sigma_j)_+\rightarrow Y\wedge(\Sigma_j)_+$ where $X\rightarrow Y$
is a $\pi_*$-isomorphism of (non-equivariant) orthogonal spectra.
\end{Thm}

Our first aim is to prove the theorem in the special case $j=1$.
The ``orthogonal pairs commute''-relation of $D_n(1;V)$ is void,
therefore it is not too hard to
analyze the orthogonal spectrum $\mathcal{D}_n(1)$ directly.

\begin{Lem}\label{Lem:p1piiso}
The map $p:\mathcal{D}_n(1)\rightarrow \mathcal{H}(1)$ is a $\pi_*$-isomorphism.
\end{Lem}

Recall that there is a forgetful functor $\U$ from orthogonal spectra to
prespectra. If $L$ is an orthogonal spectrum, then the $q$'th space of $\U L$
is $L(\R^q)$. By definition~\ref{Def:piiso} a map $K\rightarrow L$ of
orthogonal spectra is a $\pi_*$-iso, if the underlying map
$\U K\rightarrow \U L$ is a $\pi_*$-iso.

We say that a map $X\rightarrow Y$ of prespectra is an \textit{l-cofibration}
if for every $q$ the map at level $q$, $X_q\rightarrow Y_q$ is
an unbased closed cofibration of topological spaces.
Notice that a map $f:K\rightarrow L$ between orthogonal spectra is an
l-cofibration if and only if $\U f$ is an l-cofibration between prespectra.

Cubical diagrams of spaces and
prespectra will play a part in the proving that $p$ is a $\pi_*$-isomorphism. 
We refer to~\cite{Goodwillie:92} for the theory.
We recall the definition here: Let $T$ be a finite set, and $\mathscr{P}(T)$ the 
partially ordered set of all subsets of $T$. Let $\mathscr{C}$ be a category, usually
a category of spaces, prespectra or orthogonal spectra. A \textit{cubical diagram}
is a functor $\mathscr{X}:\mathscr{P}(T)\rightarrow\mathscr{C}$. If $T$ has $n$ elements,
then $\mathscr{X}$ is an $n$-cube.

Assume that our category $\mathscr{C}$ comes with a distinguished class of maps,
called cofibrations, and that $\mathscr{C}$ has all finite limits. 
Following Goodwillie we define $\mathscr{X}$ to be a \textit{cofibration cube}
if for every $U\subset T$ the map
\[
\colim_{V\subsetneqq U}\mathscr{X}(V)\rightarrow \colim_{V\subset U}\mathscr{X}(V)=\mathscr{X}(U)
\]
is a cofibration. The categories of spaces, prespectra and orthogonal spectra satisfy the
assumptions. We use unbased closed cofibrations for the category of spaces and l-cofibrations 
for prespectra and orthogonal spectra.
Therefore we have notions of cofibration cubes of spaces, l-cofibration cubes of
prespectra and l-cofibration cubes of orthogonal spectra.

For a given cubical diagram $\mathscr{X}$ in spaces, prespectra or orthogonal spectra,
we are often interested in the map
$\colim_{V\subsetneqq T}\mathscr{X}(V)\rightarrow \mathscr{X}(T)$
up to homotopy. However it is often easier to calculate with the homotopy colimit.
Therefore we compare these via the canonical map. We have the following result:

\begin{Prop}\label{Prop:canpiisoofcofcubes}
If $\mathscr{X}$ is a cofibration cube of spaces, then the canonical map
\[
\hocolim_{V\subsetneqq T}\mathscr{X}(V)\rightarrow \colim_{V\subsetneqq T}\mathscr{X}(V)
\]
is a weak equivalence. Furthermore, if $\mathscr{X}$ is an l-cofibration cube of
prespectra or orthogonal spectra, then the canonical map is a $\pi_*$-isomorphism.
\end{Prop}

\begin{proof}
The statement for cofibration cubes of spaces is proposition~1.16 in~\cite{Goodwillie:92}.

Assume that $\mathscr{X}$ is an l-cofibration cube of prespectra. 
Observe that $\hocolim$ and $\colim$ are level-wise constructions.
Hence the $q$'th space of $\hocolim_{V\subsetneqq T}\mathscr{X}(V)$
is $\hocolim_{V\subsetneqq T}\mathscr{X}_q(V)$, and similarly for $\colim$.
Since each $\mathscr{X}_q$ is a cofibration cube of spaces, the 
result for the case of spaces implies that the canonical map
\[
\hocolim_{V\subsetneqq T}\mathscr{X}(V)\rightarrow \colim_{V\subsetneqq T}\mathscr{X}(V)
\]
is a level-equivalence, hence also a $\pi_*$-iso. 

The result for l-cofibration cubes of orthogonal spectra is proved similarly.
\end{proof}

Using the definition to check directly if a given cube is a cofibration cube or not,
is not a very efficient method. The author has learned the following recognition criterion 
from Christian Schlichtkrull. But first some notation:
\begin{itemize}
\item[] If $V\subset U \subset T$ and $U\smallsetminus V$ contains exactly one element, 
then we call $V\subset U$ an \textit{edge of $T$}.
\item[] If $U$ and $V$ are subsets of $T$ such that $U\cap V\subset U$ 
and $U\subset U\cup V$ are edges, then we say
that $U$ and $V$ \textit{span a $2$-face of $T$}. 
\end{itemize}

\begin{Prop}\label{Prop:reccofcube}
Let $\mathscr{X}$ be a $T$-cube of spaces. If
\begin{itemize}
\item[i)] for all edges $V\subset U$ of $T$ the map $\mathscr{X}(V)\rightarrow \mathscr{X}(U)$ is an unbased
closed cofibration, and
\item[ii)] whenever $U$ and $V$ span a $2$-face of $T$ the square 
\[
\begin{CD}
\mathscr{X}(U\cap V)@>>> \mathscr{X}(U)\\
@VVV @VVV\\
\mathscr{X}(V)@>>> \mathscr{X}(U\cup V)
\end{CD}
\]
is pullback,
\end{itemize}
then $\mathscr{X}$ is a cofibration cube. The corresponding results for l-cofibration cubes of
prespectra and orthogonal spectra also hold.
\end{Prop}

\begin{proof}
Consider cubes in spaces. We assume by induction that the result holds for all $n$-cubes 
for $n<|T|$. We must show that the map
\[
b_U:\colim_{V\subsetneqq U}\mathscr{X}(V)\rightarrow \mathscr{X}(U)
\]
is an unbased closed cofibration of spaces. Observe that by the induction hypothesis it is enough
to show that this holds for $b_T$.

The case $n=1$ is trivial, and $n=2$ follows directly from Lillig's union theorem~\cite{Lillig:73}.

Choose some $t_0\in T$ and let $T'=T\smallsetminus\{t_0\}$. Define $\mathscr{Y}$ to be the 
$T'$-cube with $\mathscr{Y}(U)$ the pushout of 
\[
\mathscr{X}(U\cup\{t_0\})\leftarrow\mathscr{X}(U)\rightarrow\mathscr{X}(T')\quad.
\]
Notice that the map $b_T$ is equal to the map
\[
b'_{T'}:\colim_{V\subsetneqq T'}\mathscr{Y}(V)\rightarrow \mathscr{Y}(T')\quad.
\]
Hence it remains to show that $\mathscr{Y}$ satisfies i) and ii).

Let $V\subset U$ be an edge of $T'$. Consider the diagram
\[
\begin{CD}
\mathscr{X}(V\cup\{t_0\})@<<<\mathscr{X}(V)@>>>\mathscr{X}(T')\\
@VVV @VVV @VV{=}V\\
\mathscr{X}(U\cup\{t_0\})@<<<\mathscr{X}(U)@>>>\mathscr{X}(T')
\end{CD}\quad.
\]
The left square is pullback by ii) for $\mathscr{X}$, thus the gluing 
lemma for unbased closed cofibrations, proposition~2.5 in~\cite{Lewis:82}
applies and yields that $\mathscr{Y}(V)\rightarrow \mathscr{Y}(U)$
is a cofibration.

Next assume that $U$ and $V$ span a $2$-face of $T'$. By i) for $\mathscr{Y}$
we can assume that $\mathscr{Y}(U)$ and $\mathscr{Y}(V)$ are subspaces of
$\mathscr{Y}(U\cup V)$, and we must show that the intersection of these subspaces is
$\mathscr{Y}(U\cap V)$. Also $\mathscr{X}(U\cup\{t_0\})$, $\mathscr{X}(U\cup\{t_0\})$
and $\mathscr{X}(T')$ are subspaces of $\mathscr{Y}(U\cup V)$ and we have:
\begin{multline*}
\mathscr{Y}(U)\cap\mathscr{Y}(V) 
= \left(\mathscr{X}(U\cup\{t_0\})\cup\mathscr{X}(T')\right)\cap\left(\mathscr{X}(V\cup\{t_0\})\cup\mathscr{X}(T')\right)\\
= \left(\mathscr{X}(U\cup\{t_0\})\cap\mathscr{X}(V\cup\{t_0\})\right)\cup \mathscr{X}(T')\\
= \mathscr{X}((U\cap V)\cup\{t_0\})\cup \mathscr{X}(T')
= \mathscr{Y}(U\cap V)\quad.
\end{multline*}
Here we have used ii) for $\mathscr{X}$ on the $2$-face spanned by $U\cup\{t_0\}$ and $V\cup\{t_0\}$.

The proposition also holds for prespectra and orthogonal spectra by applying the result for spaces level-wise.
\end{proof}

The prerequisites for proving lemma~\ref{Lem:p1piiso} are now in place, and we give its proof:

\begin{proof}
It is enough to consider the map between the underlying prespectra.
The main idea of the proof is to filter $D_n(1;\R^q)$ by word length.
Let $F_mD_n(1;\R^q)$ be the set of all elements represented by words with $m$ or fewer letters.
We relate the cofiber of $F_{m-1}D_n(1;\R^q)\subset F_mD_n(1;\R^q)$ to the Stiefel manifold
$V_n(\R^q)$ of $n$-frames in $\R^q$: If $V_n(\R^q)^m$ denotes the $m$-fold cross product and 
$sV_n(\R^q)^{m-1}$ the subspace consisting of those $m$-tuples $(\phi_1,\ldots,\phi_m)$ with $\phi_r=\phi_{r+1}$
for some $r$, then the following diagram is pushout
\[
\begin{CD}
sV_n(\R^q)^{m-1} @>>> F_{m-1}D_n(1;\R^q)\\
@VVV @VVV\\
V_n(\R^q)^m @>>> F_m D_n(1;\R^q)
\end{CD}\quad.
\]
The horizontal maps send an $m$-tuple $(\phi_1,\ldots,\phi_m)$
to the word $(\phi_1,1)\cdots(\phi_m,1)$.
Notice that the diagram is natural for isometric embeddings $\R^q\rightarrow \R^{q'}$.
The filtration of $D_n(1;\R^q)$ induces a filtration of $\U\mathcal{D}_n(1)$ by letting
the $q$'th space of $F_m\U\mathcal{D}_n(1)$ be $F_mD_n(1;\R^q)\wedge S^q$.

Fix $m$. Let $T$ be the set $\{1,2,\ldots,m-1\}$. For $U\subseteq T$ define
\begin{multline*}
V_n^m(\R^q;U)=\{(\phi_1,\ldots,\phi_m)\;|\;\\ \text{for each $i$, } \phi_i \in V_n(\R^q)
\text{ and for all $r\not\in U$, we have $\phi_r=\phi_{r+1}$.}\}
\end{multline*}
This defines a $T$-cube of spaces. 
Observe that 
\[
\colim_{U\subsetneqq T} V_n^m(\R^q;U)= sV_n(\R^q)^{m-1}\quad.
\]
Define a cubical diagram, $\mathscr{V}_n^m$, of prespectra by defining the $q$'th space
of the prespectrum $\mathscr{V}_n^m(U)$ to be
\[
V_n^m(\R^q;U)_+\wedge S^q\quad.
\]
For $m\geq 2$ we have pushout diagrams of prespectra:
\[
\begin{CD}
\colim_{U\subsetneqq T}\mathscr{V}_n^m(U) @>>> F_{m-1}\U\mathcal{D}_n(1)\\
@VVV @VVV\\
\mathscr{V}_n^m(T) @>>> F_m\U\mathcal{D}_n(1)
\end{CD}\quad. 
\]

We now check that $\mathscr{V}_n^m$ is an l-cofibration cube. 
Our intention is to apply proposition~\ref{Prop:reccofcube}.
Let $U\subset U\cup\{r\}$ be an edge of $T$. Notice that
$V_n^m(\R^q;U)$ is a smooth submanifold of $V_n^m(\R^q;U\cup\{r\})$.
Hence the existence of a tubular neighborhood implies that the
inclusion $V_n^m(\R^q;U)\rightarrow V_n^m(\R^q;U\cup\{r\})$
is an unbased closed cofibration. Since the functor
$(-)_+\wedge S^q$ preserves unbased closed cofibrations
this proves that condition i) of the proposition holds for $\mathscr{V}_n^m$.
The fact that condition ii) holds follows nearly directly from the definition of
$V_n^m(\R^q,-)$.

By direct computation we now show that $\mathscr{V}_n^m(U)$ is $\pi_*$-isomorphic
to the sphere prespectrum, $S$, for all $U$. We define a map $\mathscr{V}_n^m(U)\rightarrow S$
by identifying $S^0$ with $\{1\}_+$, sending $V_n^m(\R^q;U)$ to $1$ and applying the functor 
$(-)_+\wedge S^q$ to get the map at level $q$:
\[
f_q:(\mathscr{V}_n^m(U))_q=V_n^m(\R^q;U)_+\wedge S^q\rightarrow S^q\quad.
\]
The induced map of homotopy groups is
\[
\pi_s\mathscr{V}_n^m(U)=\colim_q\pi_{q+s}(V_n^m(\R^q;U)_+\wedge S^q)
\xrightarrow{\colim_q\pi_{q+s}(f_q)}\colim_q\pi_{q+s}S^q=\pi_s S\quad.
\]
We inspect the map in the middle at some fixed $q$.
By the Hurewicz theorem the first non-trivial relative 
homotopy group of $f_q$ is isomorphic to the first non-trivial relative homology group.
And in homology suspension induces a natural isomorphism. Therefore we
inspect when
\[
\tilde{H}_s(V_n^m(\R^q;U)_+)\rightarrow\tilde{H}_s(S^0)
\]
is an isomorphism. Since $V_n^m(\R^q;U)$ is a cross product of the Stiefel manifold of
$n$-frames in $\R^q$, the range for $s$ where $\pi_{q+s}(f_q)$ is an iso 
clearly goes to $\infty$ when $q$ increases.

Let $\mathscr{S}$ be the $T$-cube with $\mathscr{S}(U)$ constant equal to the sphere prespectrum.
The computation above showed that there is a map of $T$-cubes 
\[
\mathscr{V}_n^m\rightarrow \mathscr{S}
\]
which is a $\pi_*$-iso at each $U\subseteq T$. Furthermore, both cubes are l-cofibration cubes.
Now consider the diagram
\[
\begin{CD}
\hocolim_{U\subsetneqq T}\mathscr{V}_n^m(U)@>{\simeq}>>\colim_{U\subsetneqq T}\mathscr{V}_n^m(U)\\
@V{\simeq}VV @VVV\\
\hocolim_{U\subsetneqq T}\mathscr{S}(U)@>{\simeq}>>\colim_{U\subsetneqq T}\mathscr{S}(U)
\end{CD}\quad.
\]
The left vertical map is a $\pi_*$-iso since a homotopy colimit of
$\pi_*$-iso is itself a $\pi_*$-iso. The horizontal maps are $\pi_*$-isos 
by proposition~\ref{Prop:canpiisoofcofcubes}.
But $\mathscr{S}$ is constant, so
$\colim_{U\subsetneqq T}\mathscr{S}(U)$ is equal to $S$. Hence the map
\[
\colim_{U\subsetneqq T}\mathscr{V}_n^m(U)\rightarrow S
\]
is a $\pi_*$-iso.

Recall that $\mathcal{H}(1)$ is the suspension of $(\Z/2)_+$. Thus we may identify
the prespectrum $\U\mathcal{H}(1)$ with $\bigvee_{\Z/2}S$.
We will now conclude the proof by showing, using induction, that for each $m\geq1$ the map
\[
p:F_m\U\mathcal{D}_n(1)\rightarrow \U\mathcal{H}(1)=\bigvee_{\Z/2}S=S\vee S
\]
is a $\pi_*$-iso.
Observe that $F_1 D_n(1;\R^q)$ is homeomorphic to $\{1\}\cup V_n(\R^q)$.
This implies that $S\vee\mathscr{V}_n^1(\emptyset)$ is isomorphic to $F_1\U\mathcal{D}_n(1)$,
and the induction hypothesis holds for $m=1$.

For $m\geq2$ we consider the diagram
\[
\begin{CD}
\mathscr{V}_n^m(T)@<{i}<< \colim_{U\subsetneqq T}\mathscr{V}_n^m(U) @>>> F_{m-1}\U\mathcal{D}_n(1)\\
@VVV @VVV @VV{p}V\\
S @= S @>{j}>> S\vee S
\end{CD}\quad.
\]
Depending on the parity of $m$ the map $j$ is the inclusion of the odd or even wedge summand.
All vertical maps are $\pi_*$-isos and the horizontal maps in the left square are l-cofibrations,
hence the map of the row-wise pushouts, $p:F_m\U\mathcal{D}_n(1)\rightarrow S\vee S$, is again a
$\pi_*$-iso.
\end{proof}

\begin{Lem}
$\U\mathcal{D}_n(1)$ is well-pointed.
\end{Lem}

\begin{proof}
We use the filtration $F_m\U\mathcal{D}_n(1)$ from the proof of the previous lemma.
The $q$'th space of $F_1\U\mathcal{D}_n(1)$ is
\[
S^q\vee \left(V_n(\R^q)\wedge S^q\right)\quad,
\]
hence well-pointed. Furthermore, we have seen that $\mathscr{V}_n^m$ is an l-cofibration cube.
Hence the left vertical map in the pushout diagram
\[
\begin{CD}
\colim_{U\subsetneqq T}\mathscr{V}_n^m(U) @>>> F_{m-1}\U\mathcal{D}_n(1)\\
@VVV @VVV\\
\mathscr{V}_n^m(T) @>>> F_m\U\mathcal{D}_n(1)
\end{CD} 
\]
is an l-cofibration. It follows that the map 
\[
F_{m-1}\U\mathcal{D}_n(1)\rightarrow F_m\U\mathcal{D}_n(1)
\] 
is an l-cofibration for any $m$. Therefore $\U\mathcal{D}_n(1)$ is well-pointed.
\end{proof}

The category of prespectra has a major disadvantage, it lacks a symmetric monoidal smash product.
However, for the purpose of calculating in the homotopy category we may, in several different ways,
define ``handicrafted'' or naive smash products of prespectra. 

\begin{Def}
Define the \textit{naive smash product} of prespectra $X$ and $Y$ by
\[
(X\wedge Y)_{2q}=X_q\wedge Y_q\quad\text{and}\quad(X\wedge Y)_{2q+1}=X_q\wedge Y_q\wedge S^1\quad.
\]
The structure maps are evident.
\end{Def}

And we define the $j$-fold naive smash product iteratively:
\[
(X^{1}\wedge\cdots\wedge X^{j})= X^{1}\wedge( X^2\wedge( X^3 \cdots \wedge (X^{j-1} \wedge X^{j})\cdots))\quad.
\]
\S{}11 in~\cite{MandellMaySchwedeShipley:01} explains the connection between the naive smash product of
prespectra and the smash product of orthogonal spectra. Given orthogonal spectra $L$ and $K$
there is a weak map
\[
\phi:\U L\wedge \U K\rightarrow \U(L\wedge K)\quad,
\]
and their proposition~11.9 says that $\phi$ is a $\pi_*$-iso whenever $L$ or $K$ is cofibrant.

More important for our purposes is a criterion for when the naive smash product preserves $\pi_*$-isomorphisms:

\begin{Prop}\label{Prop:preservepiisoprespectra}
Assume that $X$, $Y$ and $Z$ are well-pointed prespectra. 
If $f:Y\rightarrow Z$ is a $\pi_*$-isomorphism,
then the induced map of naive smash products 
$\id\wedge f:X\wedge Y\rightarrow X\wedge Z$ is also
a $\pi_*$-isomorphism.
\end{Prop}

\begin{proof}
We first prove the corresponding result for spaces. Let $A$, $B$ and $C$
be well-pointed spaces and $f:B\rightarrow C$ a weak equivalence. Then consider the
diagram
\[
\begin{CD}
* @<<< A\vee B @>{i}>> A\times B\\
@VVV @VVV @VVV\\
* @<<< A\vee C @>{i'}>> A\times C
\end{CD}\quad.
\]
The vertical maps are weak equivalences, and $i$ and $i'$ are 
cofibrations by Steenrod's product theorem, see theorem~6.3 in~\cite{Steenrod:67}
or theorem~6 in~\cite{Strom:68}. By proposition~\ref{Prop:gluewe}
the map $A\wedge B\rightarrow A\wedge C$ is a weak equivalence.

Let $A$ be a based CW-complex and $f:Y\rightarrow Z$ a $\pi_*$-iso of prespectra.
Theorem~7.4(i) in~\cite{MandellMaySchwedeShipley:01} says that also $\id\wedge f:A\wedge Y\rightarrow A\wedge X$
is a $\pi_*$-iso.

Assume that $B$ is a well-pointed space and $(A,*)$ a CW-approximation for $(B,*)$, see proposition~\ref{Prop:CWapprox}.
If $Y$ is a well-pointed prespectrum, then the result for spaces implies that
\[
A\wedge Y\rightarrow B\wedge Y
\]
is a level-equivalence, hence also a $\pi_*$-iso.

Now assume that $f$ is a $\pi_*$-iso between well-pointed prespectra, and $(A,*)$, $(B,*)$ as above. 
Consider the diagram
\[
\begin{CD}
A\wedge Y @>{\simeq}>> A\wedge Z\\
@V{\simeq}VV @VV{\simeq}V\\
B\wedge Y @>{\id_B\wedge f}>> B\wedge Z
\end{CD}\quad.
\]
It follows that $\id_B\wedge f$ is a $\pi_*$-iso for any well-pointed space $B$.

Next consider the homotopy groups of $X\wedge Y$. We can rewrite them as:
\[
\pi_s(X\wedge Y)=\colim_q\pi_{2q+s}(X_q\wedge Y_q)=\colim_q\pi_{q+s}(X_q\wedge Y)\quad.
\]
Since $X_q$ is well-pointed, it follows that $X_q\wedge Y\rightarrow X_q\wedge Z$ is a $\pi_*$-iso.
This proves the result.
\end{proof}

Suppose that we want to check that a map $f:X\rightarrow Y$ of prespectra is a $\pi_*$-iso. 
In order to do so, it is enough to give a weak inverse.
By a \textit{weak inverse} to $f$ we mean for each $q$
a map $g_q:Y_q\wedge S^l\rightarrow X_{q+l}$ where $l$ is some positive integer, 
such that both $g_q\circ (f_q\wedge S^l)$ and
$f_{q+l}\circ g_q$ are homotopic to suspensions.

\begin{Prop}
If $f:X\rightarrow Y$ has a weak inverse, then $f$ is a $\pi_*$-isomorphism.
\end{Prop}

\begin{proof}
First we check that $f_*:\pi_s X\rightarrow \pi_s Y$ is surjective for all $s$.
A class in $\pi_s Y$ is represented by some $\beta\in\pi_{q+s}Y_q$ and suspends
to $\beta'\in\pi_{q+s+l}(Y_q\wedge S^l)$. Let $\alpha\in\pi_{q+l+s}X_{q+l}$ be $g_q(\beta')$.
Since $f_{q+l}\circ g_q$ is homotopic to the suspension, observe that $f_{q+l}(\alpha)\simeq \beta'$. Thus
we see that the class of $\alpha$ in $\pi_sX$ maps to the class of $\beta$ in $\pi_s Y$.

To check that $f_*$ is injective, we pick an element of the kernel.
It can be represented by a $\alpha\in\pi_{q+s}X_q$ such that $f_q(\alpha)$ is null homotopic in $Y_q$.
Suspend the null homotopy by the appropriate $S^l$ and apply $g_q$. Since $g_q\circ (f_q\wedge S^l)$ 
is homotopic to the suspension, we get a null homotopy of $\alpha\wedge S^l$ in $X_{q+l}$.
\end{proof}

In contrast to proposition~\ref{Prop:preservepiisoprespectra}, we do not need that $X$
is well-pointed in order to draw the conclusion that the naive smash product functor $X\wedge-$
preserves the property of having weak inverses:

\begin{Prop}\label{Prop:whismash}
If $f:Y\rightarrow Z$ has a weak inverse and $X$ is any prespectrum, then 
$\id\wedge f:X\wedge Y\rightarrow X\wedge Z$ also has a weak inverse.
\end{Prop}

\begin{proof}
Let $g_q:Z_q\wedge S^l\rightarrow Y_{q+l}$ be the weak inverse of $f$. We will construct a
weak inverse $h$ for $\id_X\wedge f$. Observe that it is enough to define $h$ for even indexes.
And we let $h_{2q}$ be the composition
\[
(X\wedge Z)_{2q}\wedge S^{2l}\cong X_q\wedge S^l\wedge Z_q\wedge S^l
\xrightarrow{\text{suspension}\wedge g_q} X_{q+l}\wedge Y_{q+l}=(X\wedge Y)_{2(q+l)}\quad.
\]
It is easily seen that $h$ is a weak inverse as claimed.
\end{proof}

And the property of having weak inverses is closed under composition:

\begin{Prop}\label{Prop:whicompose}
If $f:X\rightarrow Y$ and $f':Y\rightarrow Z$ both have weak inverses, the composition
$f'f:X\rightarrow Z$ has also a weak inverse.
\end{Prop}

\begin{proof}
Let $g$ and $g'$ denote the respective weak inverses. We define a weak inverse $h$ 
for $f'f$ as follows: Assume $g'_q$ maps $Z_q\wedge S^{l'}$ to $Y_{q+l'}$ and
$g_{q+l'}$ maps $Y_{q+l'}\wedge S^l$ to $X_{q+l'+l}$, then we let $h_q$ be the composition
\[
Z_q\wedge S^{l'+l}\cong Z_q\wedge S^{l'}\wedge S^l\xrightarrow{g'_q\wedge S^l}
Y_{q+l'}\wedge S^l\xrightarrow{g_{q+l'}} X_{q+l'+l}\quad.
\]
\end{proof}

Now we introduce the orthogonal spectrum $\mathcal{D}'_n(j)$. Let the $V$'th space be
$D_n(j;V)_+\wedge S^V$. Notice that
\[
\mathcal{D}_n(j)\cong \mathcal{D}'_n(j)\wedge (\Sigma_j)_+\quad,
\]
and that this splitting corresponds to splitting $\mathcal{H}(j)$ as the smash product
of the suspension of $(\Z/2)^{\wedge j}_+$ and $(\Sigma_j)_+$.

\begin{Lem}\label{Lem:toppiiso}
There is a $\pi_*$-isomorphism of prespectra
$(\U\mathcal{D}_n(1))^{\wedge j}\rightarrow\U\mathcal{D}'_n(j)$.
\end{Lem}

\begin{proof}
By induction on $j$
we will construct the map $(\U\mathcal{D}_n(1))^{\wedge j}\rightarrow\U\mathcal{D}'_n(j)$
together with a weak inverse. For $j=1$ the map is the identity. 

Assume that $(\U\mathcal{D}_n(1))^{\wedge (j-1)}\rightarrow\U\mathcal{D}'_n(j-1)$
already is given. Then the map 
$(\U\mathcal{D}_n(1))^{\wedge j}\rightarrow\U\mathcal{D}_n(1)\wedge \U\mathcal{D}'_n(j-1)$
has a weak inverse by proposition~\ref{Prop:whismash}. By proposition~\ref{Prop:whicompose}
we are done once we have constructed a map 
$f:\U\mathcal{D}_n(1)\wedge \U\mathcal{D}'_n(j-1)\rightarrow\U\mathcal{D}'_n(j)$
and a weak inverse.

A map out of our naive smash product is completely determined by what it is at the $(2q)$'th spaces.
What we need is a map from $D_n(1;\R^q)_+\wedge S^q\wedge D_n(j-1;\R^q)_+\wedge S^q$
to $D_n(j;\R^2q)_+\wedge S^{2q}$. To define it we smash a suitable shuffling of $S^1$'s, 
$\operatorname{sh}:S^q\wedge S^q\cong S^{2q}$, with a
group homomorphism $\alpha:D_n(1;\R^q)\times D_n(j-1;\R^q)\rightarrow D_n(j;\R^{2q})$.

Let $i_{\text{odd}}:\R^q\rightarrow \R^{2q}$ be the inclusion of the odd coordinates, 
\[
i_{\text{odd}}(x_1,x_2,\ldots,x_q)=(x_1,0,x_2,0,\ldots,x_q,0)\quad,
\]
and $i_{\text{even}}:\R^q\rightarrow \R^{2q}$ the inclusion of the even coordinates,
\[
i_{\text{even}}(x_1,x_2,\ldots,x_q)=(0,x_1,0,x_2,\ldots,0,x_q)\quad.
\]
Our group homomorphism $\alpha$ sends $(\phi,1)$ in $D_n(1;\R^q)$ to $(i_{\text{odd}}\phi,1)$ in $D_n(j;\R^{2q})$
and $(\psi,r)$ in $D_n(j-1;\R^q)$ to $(i_{\text{even}}\psi,r+1)$ in $D_n(j;\R^{2q})$.
By the ``orthogonal pairs commute''-relation in $D_n(j;\R^{2q})$ we have
\[
(i_{\text{odd}}\phi,1)(i_{\text{even}}\psi,r+1)=(i_{\text{even}}\psi,r+1)(i_{\text{odd}}\phi,1)
\]
since $i_{\text{odd}}\phi$ is orthogonal to $i_{\text{even}}\psi$. 
This shows that the group homomorphism is well defined.

If we identify $S^q$ and $S^{2q}$ with one-point-compactifications of $\R^q$ and $\R^{2q}$
respectively, we can write the shuffling $\operatorname{sh}:S^q\wedge S^q\cong S^{2q}$ as follows:
\[
\big((v_1,v_2,\ldots,v_q),(w_1,w_2,\ldots,w_q)\big)\mapsto (v_1,w_1,v_2,w_2,\ldots,v_q,w_q)\quad.
\]
This ensures that the maps 
$D_n(1;\R^q)_+\wedge S^q\wedge D_n(j-1;\R^q)_+\wedge S^q\rightarrow D_n(j;\R^2q)_+\wedge S^{2q}$
commute strictly with the suspensions,
and thus we get our map
\[
f:(\U\mathcal{D}_n(1))^{\wedge j}\rightarrow \U\mathcal{D}'_n(j)\quad.
\]

To construct a weak inverse $g$ for $f$ we first define a group homomorphism
\[
\beta:D_n(j;\R^q)\rightarrow D_n(1;\R^q)\times D_n(j-1;\R^q)\quad.
\]
On generators $\beta$ is given by
\[
\beta(\phi,r)=\begin{cases}
\big( (\phi,1),1\big)&\text{if $r=1$, and}\\
\big( 1, (\phi, r-1)\big) &\text{if $r>1$.}
\end{cases}
\]
Now define $g_q$ as the composition
\begin{multline*}
\big(\U\mathcal{D}'_n(j)\big)_q\wedge S^q\\
=D_n(j;\R^q)_+\wedge S^{2q}
\xrightarrow{\beta\wedge \operatorname{sh}^{-1}} D_n(1;\R^q)_+\wedge S^q \wedge D_n(j-1;\R^q)_+\wedge S^q\\
=\big(\U\mathcal{D}_n(1)\wedge \U\mathcal{D}'_n(j-1)\big)_{2q}\quad.
\end{multline*}

Let $i_1:\R^q\rightarrow \R^{2q}$ be the standard inclusion (embeds $\R^q$ as the first $q$ coordinates).
The space of isometric embeddings of $\R^q$ in $\R^{2q}$ is connected, so we can choose
paths from $i_{\text{odd}}$ and $i_{\text{even}}$ to $i_1$. 
Now it is easy to see that the composition
\[
D_n(1;\R^q)\times D_n(j-1;\R^q)\xrightarrow{\alpha} D_n(j;\R^{2q})
\xrightarrow{\beta}D_n(1;\R^{2q})\times D_n(j-1;\R^{2q})
\]
is homotopic to the map of $D_n(1;-)\times D_n(j-1;-)$ induced by $i_1$.
With the opposite composition,
\[
D_n(j;\R^q)\xrightarrow{\beta} D_n(1;\R^q)\times D_n(j-1;\R^q)
\xrightarrow{\alpha} D_n(j;\R^{2q})\quad,
\]
we have to be a bit more careful. On generators this map is given by
\[
\beta\alpha(\phi,r)=\begin{cases}
(i_{\text{odd}}\phi,r)&\text{ if $r=1$, and}\\
(i_{\text{even}}\phi,r)&\text{ if $r>1$.}
\end{cases}
\]
Let $i_t$, $t\in\left[0,\frac{\pi}{2}\right]$ be the homotopy between $i_{\text{odd}}$ and $i_{\text{even}}$
given by the formula
\[
i_t(x_1,x_2,\ldots,x_q)=\big( x_1\cos t, x_1\sin t, x_2\cos t, x_2\sin t,\ldots , x_q\cos t, x_q\sin t\big)\quad.
\]
Notice that when $\mb{x}$ and $\mb{y}$ are orthogonal vectors in $\R^q$, then
\begin{itemize}
\item[] $i_t(\mb{x})$ and $i_t(\mb{y})$ are orthogonal, and
\item[] $i_0(\mb{x})$ and $i_t(\mb{y})$ are orthogonal.
\end{itemize}
Define $h_t:D_n(j;\R^q)\rightarrow D_n(j;\R^{2q})$, $t\in\left[0,\frac{\pi}{2}\right]$ on generators by
\[
h_t(\phi,r)=\begin{cases}
(i_{\text{odd}}\phi,r)&\text{ if $r=1$, and}\\
(i_{t}\phi,r)&\text{ if $r>1$.}
\end{cases}
\]
It is well defined, for $t=0$ it equals the map induced by $i_{\text{odd}}$ and for $t=\frac{\pi}{2}$
it is $\beta\alpha$. And using the path from $i_{\text{odd}}$ to $i_1$, we can extend $h_t$ to a
homotopy from $\beta\alpha$ to the map $D_n(j;\R^q)\rightarrow D_n(j;\R^{2q})$ induced by $i_1$.

Checking what happens with the spheres, we see that $g_{2q}\circ(f_q\wedge S^{2q})$ and
$f_{2q}\circ g_q$ both are homotopic to suspensions. 
\end{proof}

With the lemmas~\ref{Lem:p1piiso} and~\ref{Lem:toppiiso} 
in place it is quite easy to prove theorem~\ref{Thm:Ddiscrete}.

\begin{proof}
Recall the definition of the group homomorphism $p:D_n(j;V)\rightarrow (\Z/2)^j$.
Applying $(-)_+\wedge S^V$ we get a map of orthogonal spectra $\mathcal{D}'_n(j)\rightarrow F_0(\Z/2)_+^{\wedge j}$.
Here $F_0$ denotes the $0$'th shift desuspension functor. Notice that this map fits into a diagram of prespectra
\[
\begin{CD}
(\U\mathcal{D}_n(1))^{\wedge j}@>>>\U\mathcal{D}'_n(j)\\
@V{p^{\wedge j}}VV @VVV\\
(\U\mathcal{H}(1))^{\wedge j} @= \U F_0 (\Z/2)_+^{\wedge j}
\end{CD}\quad.
\]
The map at the top is the $\pi_*$-iso from lemma~\ref{Lem:toppiiso}, and left map
is an iterated naive smash product of the $\pi_*$-iso $p$ form lemma~\ref{Lem:p1piiso}.
Since both $\U\mathcal{D}_n(1)$ and $\U\mathcal{H}(1)$ are well-pointed, it follows
from proposition~\ref{Prop:preservepiisoprespectra} that the left map is also a $\pi_*$-iso.
This implies that $\mathcal{D}'_n(j)\rightarrow F_0(\Z/2)_+^{\wedge j}$ is a $\pi_*$-iso.
Smashing both sides with $(\Sigma_j)_+$, we get that the map 
\[
\mathcal{D}_n(j)\rightarrow\mathcal{H}(j)\quad.
\]
\end{proof}

\subsection{The main theorem}\label{subsect:maintheorem}

The following result is the main theorem of this thesis. It provides an orthogonal ring spectrum
with involution associated to a stable vector bundle over a manifold. The homotopy type of
the underlying orthogonal ring spectrum
depends only on the manifold. 

\begin{Thm}\label{Thm:main}
Let $M$ be a manifold and $\xi$ an $n$-vector bundle over $M$. There exists an orthogonal ring spectrum $R$ 
with an involution depending on $\xi$, such that $R$ is weakly homotopic in the category of 
orthogonal ring spectra to $S[\Omega M]$, and the 
involution on $R$ corresponds to $\iota$ on homotopy groups. Furthermore, up to homotopy the involution on $R$  
depends only on the stable class of $\xi$.
\end{Thm}

Here $\iota$ is the involution on $\pi_*S[\Omega M]$ given in definition~\ref{Def:iota}.

\begin{proof}
By theorem~\ref{Thm:DnalgstronSOmegaM} there is a $\mathcal{D}_n$-algebra 
structure on $S[\Omega M]$, and by theorem~\ref{Thm:Ddiscrete}
there is a map of operads $\mathcal{D}_n\rightarrow \mathcal{H}$ such that for each $j$
the map $\mathcal{D}_n(j)\rightarrow \mathcal{H}(j)$ can be written equivariantly
as a product $X\wedge(\Sigma_j)_+\rightarrow Y\wedge(\Sigma_j)_+$ where $X\rightarrow Y$
is a $\pi_*$-isomorphism.
This means that we almost have the necessary requirements for applying the replacement procedure
described in remark~\ref{Rem:replproc}. However, we do not know that
\begin{itemize}
\item[] $S\rightarrow \mathcal{D}_n(1)$ is an orbit q-cofibration, and
\item[] each $\mathcal{D}_n(j)$ can be written equivariantly
as a product $X\wedge (\Sigma_j)_+$ with $X$ being an orbit cofibrant (non-equivariant) orthogonal spectrum.
\end{itemize}
Instead of attempting to prove this, we use the orbit cofibrant replacement $\hat{\Gamma}\mathcal{D}_n$
from theorem~\ref{Thm:ocreplforDn}.
Pulling back by the map of operads $\hat{\Gamma}\mathcal{D}_n\rightarrow\mathcal{D}_n$ we see that $S[\Omega M]$
is also a $\hat{\Gamma}\mathcal{D}_n$-algebra. 
Moreover, the composition $\hat{\Gamma}\mathcal{D}_n\rightarrow\mathcal{D}_n\rightarrow\mathcal{H}$
is a map of operads, and evaluated at the $j$'th objects 
it decomposes as a product $X\wedge(\Sigma_j)_+\rightarrow Y\wedge(\Sigma_j)_+\rightarrow Z\wedge(\Sigma_j)_+$,
where $X\rightarrow Y\rightarrow Z$ are $\pi_*$-isos.
Since $S\rightarrow \mathcal{D}_n(1)$ is an inclusion, it follows
by theorem~\ref{Thm:orbitcofrepl} and the construction of $\hat{\Gamma}\mathcal{D}_n(1)$
in theorem~\ref{Thm:ocreplforDn} that $S\rightarrow \hat{\Gamma}\mathcal{D}_n(1)$ is an orbit q-cofibration. Furthermore,
each $\hat{\Gamma}\mathcal{D}_n(j)$ can be described $\Sigma_j$-equivariantly as
a product $X\wedge (\Sigma_j)_+$, where $X$ is orbit cofibrant.

Now consider the replacement procedure:
\[
B(\mathcal{H},\hat{\Gamma}\mathcal{D}_n, S[\Omega M])
\leftarrow B(\hat{\Gamma}\mathcal{D}_n,\hat{\Gamma}\mathcal{D}_n, S[\Omega M])
\rightarrow S[\Omega M]\quad.
\]
We define $R$ to be $B(\mathcal{H},\hat{\Gamma}\mathcal{D}_n, S[\Omega M])$. By the considerations above 
all maps are $\pi_*$-isos. To show that $R$ is homotopic to $S[\Omega M]$ in the category of orthogonal ring spectra,
we prove that the three maps above all are morphisms in the category of orthogonal ring spectra.

Recall that by remark~\ref{Rem:MtGDn} we have a map of operads $\mathcal{M}\rightarrow \hat{\Gamma}\mathcal{D}_n$.
Hence we have a restriction functor from the category of $\hat{\Gamma}\mathcal{D}_n$-algebras
to the category of $\mathcal{M}$-algebra. Consequently, the $\hat{\Gamma}\mathcal{D}_n$-algebra map
$B(\hat{\Gamma}\mathcal{D}_n,\hat{\Gamma}\mathcal{D}_n,S[\Omega M])\rightarrow S[\Omega M]$ is also a map of
$\mathcal{M}$-algebras, i.e. a map of orthogonal ring spectra.

Similarly we see that the map 
$B(\mathcal{H},\hat{\Gamma}\mathcal{D}_n, S[\Omega M])\leftarrow B(\hat{\Gamma}\mathcal{D}_n,\hat{\Gamma}\mathcal{D}_n,S[\Omega M])$
is a map in the category of orthogonal ring spectra, since also this map is a map of $\hat{\Gamma}\mathcal{D}_n$-algebras.

We now show that the involution only depends on the stable class of $\xi$. 
Observe that standard inclusion $i:\R^n\rightarrow \R^{n+1}$
induces group homomorphisms 
\[
D_{n+1}(j;V)\rightarrow D_n(j;V)
\] 
by sending a generator $(\phi,r)$, where $\phi:\R^{n+1}\rightarrow V$ is an isometric embedding, to
$(\phi\circ i,r)$. These group homomorphisms give rise to a map of
operads $\alpha:\mathcal{D}_{n+1}\rightarrow\mathcal{D}_n$. 
And by inspection of the construction in theorem~\ref{Thm:ocreplforDn}, we have a lifting
$\hat{\alpha}:\hat{\Gamma}\mathcal{D}_{n+1}\rightarrow\hat{\Gamma}\mathcal{D}_n$
Now notice that the pullback of the $\mathcal{D}_n$-algebra structure on $S[\Omega M]$ associated to
$\xi$ is the $\mathcal{D}_{n+1}$-algebra structure on $S[\Omega M]$ associated to $\xi\oplus\varepsilon^1$.
Here $\varepsilon^1$ denotes the trivial line bundle over $M$.
Therefore we get a $\hat{\Gamma} \mathcal{D}_{n+1}$-algebra map
$B(\hat{\Gamma}\mathcal{D}_{n+1},\hat{\Gamma}\mathcal{D}_{n+1},S[\Omega M])\rightarrow B(\mathcal{H},\hat{\Gamma}\mathcal{D}_n,S[\Omega M])$.
Feeding the diagram
\[
B(\mathcal{H},\hat{\Gamma}\mathcal{D}_n,S[\Omega M])\leftarrow 
B(\hat{\Gamma}\mathcal{D}_{n+1},\hat{\Gamma}\mathcal{D}_{n+1},\Gamma S[\Omega M])\rightarrow S[\Omega M]
\]
into proposition~\ref{Prop:uniquerepl}, we get an equivalence of $\mathcal{H}$-algebras between
\[
B(\mathcal{H},\hat{\Gamma}\mathcal{D}_n,S[\Omega M])\quad\text{and}\quad B(\mathcal{H},\hat{\Gamma}\mathcal{D}_{n+1},S[\Omega M])\quad.
\]
The first orthogonal ring spectrum has the involution associated to $\xi$, while the second has the involution associated to
$\xi\oplus\varepsilon^1$. Hence up to homotopy the involution only depends on the stable class of $\xi$.

We now check that the involution does not depend on the choice of connection. Let $\nabla_0$ and $\nabla_1$ be two
connections on $\xi$. Let $\xi\times I$ be the vector bundle over $M\times I$ induced from $\xi$ via the projection
$M\times I\rightarrow M$. And let $\nabla'_0$ and $\nabla'_1$ be the induced connections. We can define the linear combination
\[
\nabla= t\nabla'_1+(1-t)\nabla'_0\quad,
\]
where $t$ is the coordinate of $I$. We see that $\nabla$ is a connection on $\xi\times I$. And pulling $\nabla$ back over the two
inclusions $i_0,i_1:M\rightarrow M\times I$ yields $\nabla_0$ and $\nabla_1$ respectively.

The inclusion $i_0$ induces a map of $\mathcal{D}_n$-algebras
\[
S[\Omega M]\rightarrow S[\Omega(M\times I)]\quad.
\]
And we can therefore form the diagram
\[
B(\mathcal{H},\hat{\Gamma}\mathcal{D}_{n}, S[\Omega(M\times I)])\leftarrow 
B(\hat{\Gamma}\mathcal{D}_n,\hat{\Gamma}\mathcal{D}_n, S[\Omega M])\rightarrow S[\Omega M]\quad.
\]
Putting this into proposition~\ref{Prop:uniquerepl}, we get an equivalence of $\mathcal{H}$-algebras between
the orthogonal ring spectrum with involution associated to the connection $\nabla_0$ on $\xi$ and the
orthogonal ring spectrum with involution associated to the connection $\nabla$ on $\xi\times I$.
A similar consideration is also true for $i_1$ and the connection $\nabla_1$. Hence the choice of
connection is irrelevant up to homotopy of orthogonal ring spectra with involution.

To show that the involution on $R$ coincides with $\iota$ on $\pi_*S[\Omega M]$:

We first construct a commutative square
\[
\begin{CD}
F_{\R^n}S^n@>{\lambda}>> S\\
@V{f}VV @VV{i}V\\
\hat{\Gamma}\mathcal{D}_n(1) @>{p}>> \mathcal{H}(1)
\end{CD}\quad.
\]
The top map, $\lambda$, is the adjoint to the identity $S^n=S(\R^n)$. By lemma~\ref{Lem:FnSntoSpiiso}
$\lambda$ is a $\pi_*$-iso.
$p$ comes from the map of operads $\mathcal{D}_n\rightarrow \mathcal{H}$,
and the map $i$ is induced from the inclusion of the matrix $\begin{pmatrix} -1 \end{pmatrix}$
in $0$'th space of $\mathcal{H}(1)$. Thus the map $i:S\rightarrow\mathcal{H}(1)$ represents the
involution.

To construct $f$ recall 
that $\mathcal{D}_n(1)(\R^n)=D_n(1;\R^n)_+\wedge S^n$. The pair
$(\id_{\R^n},1)$ represents a point in $D_n(1;\R^n)$, and we get a map
$S^n\rightarrow \mathcal{D}_n(1)(\R^n)$ by sending $v$ to $\big( (\id_{\R^n},1), v\big)$.
By adjointness we now get a map
of orthogonal spectra $F_{\R^n}S^n\rightarrow \mathcal{D}_n(1)$,
and since $F_{\R^n}S^n$ is cofibrant, we can lift to a map
\[
f:F_{\R^n}S^n\rightarrow \hat{\Gamma}\mathcal{D}_n(1)\quad.
\]

Recall from the proof of proposition~\ref{Prop:interpretoperadalgebras}
that for an $\mathcal{H}$-algebra $L$ the involution is
the composition
\[
L\cong S\wedge L \xrightarrow{i\wedge \id} \mathcal{H}(1) \wedge L \xrightarrow{\theta_1} L\quad.
\]
Analogously, for a $\hat{\Gamma}\mathcal{D}_n$-algebra $L$ we can consider the composition
\[
F_{\R^n}S^n\wedge L\xrightarrow{f\wedge\id} \hat{\Gamma}\mathcal{D}_n(1)\wedge L\xrightarrow{\theta_1}L\quad.
\]

Now inspect the diagram
\[
\begin{CD}
S\wedge B(\mathcal{H},\hat{\Gamma}\mathcal{D}_n,S[\Omega M]) 
@>{\text{involution}}>> B(\mathcal{H},\hat{\Gamma}\mathcal{D}_n,S[\Omega M]) \\
@A{\lambda\wedge \id}AA @AA{=}A\\
F_{\R^n}S^n\wedge B(\mathcal{H},\hat{\Gamma}\mathcal{D}_n,S[\Omega M]) 
@>>> B(\mathcal{H},\hat{\Gamma}\mathcal{D}_n,S[\Omega M]) \\
@AAA @AAA\\
F_{\R^n}S^n\wedge B(\hat{\Gamma}\mathcal{D}_n,\hat{\Gamma}\mathcal{D}_n,S[\Omega M])
@>>> B(\hat{\Gamma}\mathcal{D}_n,\hat{\Gamma}\mathcal{D}_n,S[\Omega M]) \\
@VVV @VVV\\
F_{\R^n}S^n\wedge S[\Omega M]
@>>> S[\Omega M] 
\end{CD}\quad.
\]
The horizontal maps, except the first, are defined via $\hat{\Gamma}\mathcal{D}_n$-algebra 
structures on $B(\mathcal{H},\hat{\Gamma}\mathcal{D}_n,S[\Omega M])$,
$B(\hat{\Gamma}\mathcal{D}_n,\hat{\Gamma}\mathcal{D}_n,S[\Omega M])$
and $S[\Omega M]$ respectively.
Observe that all vertical maps are $\pi_*$-isomorphisms. 
The map at the top is the involution on 
$R=B(\mathcal{H},\hat{\Gamma}\mathcal{D}_n,S[\Omega M])$.
The bottom map is determined by what happens at level $\R^n$.
Evaluating at this level we get a map
\[
(F_{\R^n}S^n\wedge S[\Omega M])(\R^n)=S^n\wedge \Omega M_+\rightarrow \Omega M_+\wedge S^n= (S[\Omega M])(\R^n)\quad.
\]
By definition of $f$ above, this is the map considered in remark~\ref{Rem:identifyinvolution},
and by the remark it induces the involution $\iota$, see definition~\ref{Def:iota}, 
on the homotopy groups of $S[\Omega M]$.
\end{proof}

We end this chapter with the following conjecture:

\begin{Conj}
Suppose that $\xi_1$ and $\xi_2$ are vector bundles over $M$ with the same underlying stable spherical bundle.
Let $R_1$ and $R_2$ be the orthogonal ring spectra with involution corresponding to these vector bundles.
Then there exists an orthogonal ring spectrum $R$ with involution and maps
$R_1\xleftarrow{f_1} R\xrightarrow{f_2} R_2$ in the category of orthogonal ring spectra with involution,
such that on the underlying orthogonal spectra both $f_1$ and $f_2$ are $\pi_*$-isomorphisms.
\end{Conj}

Informally, the conjecture says that up to homotopy of $R$ the involution  
depends only on the stable class of the underlying spherical bundle of $\xi$.

The motivation for this conjecture comes from the involution on $A$-theory.
For a spherical fibration $\xi$ over $X$ Vogell defines in \S{}2 of~\cite{Vogell:85} 
an involution $\tau_{\xi}$ on $A(X)$. This involution is well defined up to homotopy. 
More recently Weiss and Williams have defined an involution on $A(X)$ via Waldhausen categories
with Spanier-Whitehead duality, see example~1.A.9 in~\cite{WeissWilliams:98}
and \S{}4.1 in~\cite{WeissWilliams:01}. Like Vogell, their involution depends
on a spherical fibration over $X$.

Morally, the $K$-theory of our orthogonal ring spectrum $R$ with involution $\iota$
should be weakly homotopy equivalent to $A(X)$ when $X=M$, $K(R)$ should
have an involution induced by $\iota$, and this involution should agree
with the involutions defined by Vogell, Weiss and Williams. If so,
the involution on $K(R)$ depends only on a stable spherical fibration,
and it is natural to believe that the to homotopy of $R$ the involution
has the same kind of dependence.

To prove the conjecture one should start with a geometric model for the space
$G_n$ of self-homotopy equivalences of $S^n$. Via a ``connection'' an $n$-spherical
bundle $\xi$ over $M$ corresponds to a map $P:\Omega M\rightarrow G_n$. But is this
map a homomorphism of monoids? Even if it is, the lack of strict inverses in $G_n$
prevents us from sending the reversed loop, $\bar{\gamma}$ to $P(\gamma)^{-1}$.
This causes trouble with the ``cancellation of repeated pairs''-relation of $\mathcal{D}_n$.
Therefore, one should blow up the operad $\mathcal{D}_n$ to handle this lack
of structure on $G_n$. After defining this huge operad, it should be possible to 
prove the conjecture in roughly the same way we have proved theorem~\ref{Thm:main}.

\chapter[$THH$ and $TC$]{$THH$ and $TC$ for orthogonal ring spectra with involution}\label{Chapt:trace}

Theorem~\ref{Thm:main} gives us an orthogonal ring spectrum with involution.
The intension behind is to define and calculate its $L$-theory, $LA$-theory,
$K$-theory, topological cyclic homology and its topological Hochschild homology.
These theories should be related via trace maps. 
Using surgery, the $LA$-theory should provide information about the homotopy
type of the automorphism space of our manifold, see~\cite{WeissWilliams:01}.
From $L$-theory there is a map $\Xi$ into the Tate construction on $K$-theory,
see~\S{}11 in~\cite{WeissWilliams:98}. Furthermore, from $K$-theory there are
trace maps into $TC$ and $THH$, see~\cite{Madsen:94}.

However, developing all of the above theory in the setting of orthogonal spectra, 
is far beyond the scope of this thesis. In this chapter we shall consider the
definition of $TC$ and $THH$ and a few basic properties. 
We follow the framework of well known theory, but there are some details
worth pointing out: In proposition~\ref{Prop:piisoiscyclpiiso} below, we
observe that it is easy to recognize cyclotomic $\pi_*$-isomorphisms between
cyclotomic spectra. Theorem~\ref{thm:THHstructure} shows that our model for $THH$ of
a cofibrant orthogonal ring spectrum is a cyclotomic spectrum in a very strong sense;
the cyclotomic structure maps $r_C:\rho_C^*\Phi^CTHH(L)\cong THH(L)$ are isomorphisms.
Due to the involution, it is important to use a model for $(n\times n)$-matrices which
is closed under transposition. Such a model is introduced in definition~\ref{Def:matrix}.

Important references for the theory of $THH$ and $TC$ in other settings 
includes~\cite{BokstedtHsiangMadsen:93}, \cite{Madsen:94}, \cite{DundasMcCarthy:96},
\cite{HesselholtMadsen:97}, \cite{Schlichtkrull:98} and~\cite{Shipley:00}.

\section{Cyclotomic orthogonal spectra and $TC$}

The purpose of this section is to define $TC$ of a cyclotomic spectrum $T$.
We will define cyclotomic spectra as certain orthogonal $S^1$-spectra
together with some extra structure. For the involutive case $T$ lies in the category of orthogonal
$O(2)$-spectra. Because orthogonal $S^1$- and $O(2)$-spectra have so many
model structures, one can easily become confused about which type of weak equivalences 
that are the correct ones to consider. Therefore, we will start this section by quickly 
listing the model structures to be used in the context of $TC$ and $THH$.
We end the section with a result, proposition~\ref{Prop:TCiso}, which tells us that our choices
of model categories were right.

\paragraph{Orthogonal spectra:}
We use the stable model structure, see definition~\ref{Def:osstablemodelstr}.
The weak equivalences $f:K\rightarrow L$ are the $\pi_*$-isomorphisms.
The fibrant orthogonal spectra are the $\Omega$-spectra.

This model structure is topological, see theorem~9.2 in~\cite{MandellMaySchwedeShipley:01}.
Via the functor $\Sing$ from topological spaces to simplicial sets, it can be shown that
every topological model category is a simplicial model category. Hence, we have homotopy
limits in $\mathscr{IS}$, and homotopy invariance holds, see theorem~18.5.3~ii) in~\cite{Hirschhorn:03}:

\begin{Prop}\label{Prop:htpyinvofholim}
Let $\mathscr{C}$ be a small category. 
If $f:K\rightarrow L$ is a map of $\mathscr{C}$-diagrams in $\mathscr{IS}$, and each $f:K(c)\rightarrow L(c)$
is a $\pi_*$-isomorphism between $\Omega$-spectra, then
\[
f_*:\holim_{\mathscr{C}} K\rightarrow\holim_{\mathscr{C}} L
\]
is also a $\pi_*$-isomorphism between $\Omega$-spectra.
\end{Prop}

\paragraph{Orthogonal $\Z/2$-spectra:}
An orthogonal $\Z/2$-spectrum has an underlying orthogonal spectrum.
We are interested in the model structure where the weak equivalences are
the $\Z/2$-maps which are $\pi_*$-isomorphisms between the underlying orthogonal spectra.
Considering an orthogonal $\Z/2$-spectrum as an orthogonal spectrum with $\Z/2$-action,
we see that functorial constructions on orthogonal spectra lift to constructions on orthogonal
$\Z/2$-spectra. For example, proposition~\ref{Prop:htpyinvofholim} holds
in this setting.

\paragraph{Orthogonal $S^1$-spectra:}
We are interested in the cyclotomic $\pi_*$-isomorphisms. They are given in definition~\ref{Def:cyclpiiso}. 
These are the maps $f:K\rightarrow L$ such that $f$ induces an isomorphism
$\pi^C_*K\rightarrow \pi^C_*L$ for all finite subgroups $C$ of $S^1$.

For some constructions we must change our orthogonal $S^1$-spectrum into an
$\Omega$-$S^1$-spectrum (=genuine fibrant orthogonal $S^1$-spectrum),
see definition~\ref{Def:OmegaGsp}. 
To achieve this, we use the fibrant replacement functor coming from the stable 
genuine model structure on orthogonal $S^1$-spectra.
This functor is constructed by the small object argument, and
we denote it by $Q^{\cy}$.

Recall the geometric fixed point functor $\Phi^C$, given in definition~\ref{Def:geomfpfunct}.
For a finite subgroup $C$ of $S^1$ it takes orthogonal $S^1$-spectra $L$ to orthogonal $S^1/C$-spectra $\Phi^CL$.
Using the group isomorphism $\rho_C:S^1\rightarrow S^1/C$, we pull back and get a new
orthogonal $S^1$-spectrum $\rho_C^*\Phi^C L$.

By proposition~\ref{Prop:fundresultofgeomfp} the functor $\rho_C^*\Phi^C$ preserves the
class of generating genuine acyclic q-cofibrations. Hence, lemma~\ref{Lem:smallobjandfunctors}
yields:

\begin{Lem}\label{Lem:Qcyandgeomfp}
There is a natural transformation $\rho_C^*\Phi^CQ^{\cy}\rightarrow Q^{\cy}\rho_C^*\Phi^C$
such that the following diagram commutes for all orthogonal $S^1$-spectra L:
\[
\begin{CD}
\rho_C^*\Phi^C L @= \rho_C^*\Phi^C L\\
@VVV @VVV\\
\rho_C^*\Phi^CQ^{\cy}L @>>> Q^{\cy}\rho_C^*\Phi^CL
\end{CD}\quad.
\]
\end{Lem}

\paragraph{Orthogonal $O(2)$-spectra:}
Again, we are interested in the cyclotomic $\pi_*$-isomorphisms, see definition~\ref{Def:cyclpiiso}. 
We have a fibrant replacement functor, $Q^{\cy}$, constructed by the small object argument in
the stable genuine model structure. Thus $Q^{\cy}L$ is an $\Omega$-$O(2)$-spectrum, for
any orthogonal $O(2)$-spectrum $L$. Let $C$ be a finite normal subgroup of $O(2)$.
Similar to the case above, we have geometric $C$-fixed point functors,
$\Phi^C$, and group isomorphisms $\rho_C:O(2)\cong O(2)/C$. We consider the composition
$\rho_C^*\Phi^C$. Also in the case of an orthogonal $O(2)$-spectrum $L$ lemma~\ref{Lem:Qcyandgeomfp} holds.

We are now ready to define cyclotomic spectra in the setting of orthogonal spectra. Compare this
definition with definition~2.2 in~\cite{HesselholtMadsen:97}. Furthermore, we introduce
the notion of an cyclotomic spectrum with involution.

\begin{Def}\label{Def:cyclotomicspect}
A \textit{cyclotomic spectrum} is an orthogonal $S^1$-spectrum $T$ together with a cyclotomic $\pi_*$-isomorphism
\[
r_C:\rho^{*}_C\Phi^C T\rightarrow T
\]
for every finite subgroup $C$ of $S^1$ such that for any pair of finite subgroups the following diagram commutes
\[
\begin{CD}
\rho^{*}_{C_r}\Phi^{C_r}\rho^{*}_{C_s}\Phi^{C_s}T @= \rho^{*}_{C_{rs}}\Phi^{C_{rs}}T\\
@V{\rho^{*}_{C_r}\Phi^{C_r}r_{C_s}}VV @VV{r_{C_{rs}}}V\\
\rho^{*}_{C_r}\Phi^{C_r} @>{r_{C_r}}>> T
\end{CD}\quad.
\]
A map of cyclotomic spectra is a map of orthogonal $S^1$-spectra which commutes with the $r_C$'s.
\end{Def}

\begin{Def}
A \textit{cyclotomic spectrum with involution} is an orthogonal $O(2)$-spectrum $T$ together with 
a cyclotomic $\pi_*$-isomorphism
\[
r_C:\rho^*_C\Phi^C T\rightarrow T
\]
for every finite subgroup $C$ of $S^1\subset O(2)$ such that the diagram in definition~\ref{Def:cyclotomicspect}
commutes for every pair of such subgroups.
A map of cyclotomic spectra with involution is a map of orthogonal $O(2)$-spectra which commutes with the $r_C$'s.
\end{Def}

Because of the maps $r_C$, it is easy to check when a map between cyclotomic spectra is a cyclotomic $\pi_*$-isomorphism:

\begin{Prop}\label{Prop:piisoiscyclpiiso}
A map $f:T_1\rightarrow T_2$ between cyclotomic spectra (with involution) is a cyclotomic $\pi_*$-isomorphism
if and only if it is non-equivariantly a $\pi_*$-isomorphism.
\end{Prop}

\begin{proof}
By definition, all cyclotomic $\pi_*$-isomorphisms are non-equivariant $\pi_*$-isomorphisms.

Assume that $f:T_1\rightarrow T_2$ is non-equivariantly a $\pi_*$-isomorphism.
Let $C$ be a finite normal subgroup of $S^1$ (or $O(2)$), and consider the diagram
\[
\begin{CD}
\rho_C^*\Phi^CT_1 @>{r_C}>> T_1\\
@V{\rho_C^*\Phi^C f}VV @VV{f}V\\
\rho_C^*\Phi^CT_2 @>{r_C}>> T_2
\end{CD}\quad.
\]
Since the $r_C$'s are cyclotomic $\pi_*$-isos, it follows that $\Phi^Cf$ is
non-equivariantly a $\pi_*$-isomorphism. Using proposition~\ref{Prop:geomfpdetectspiiso}
we recognize $f$ as a cyclotomic $\pi_*$-isomorphism.
\end{proof}

\begin{Rem}
The key ingredient in the proof of proposition~\ref{Prop:geomfpdetectspiiso}
was a homotopy cofiber sequence
\[
Q^{\cy}(L\wedge E\mathscr{F}_+)^C\rightarrow Q^{\cy}(L)^C \rightarrow \Phi^C L\quad,
\]
where $\mathscr{F}$ is a specific family of subgroups, and $L$ an orthogonal $S^1$- (or $O(2)$-) spectrum.
This sequence is a generalization of the ``fundamental cofibration sequence'', see formula~2.4.6
in~\cite{Madsen:94}, or theorem~2.2 in~\cite{HesselholtMadsen:97}.
\end{Rem}

In order to define $TC(T)$ we introduce the category $\mathbb{I}$.
It has the natural numbers, $\{1,2,3,\ldots\}$, as its objects, and
the set of all morphisms in $\mathbb{I}$ is generated by two classes of
morphisms $R_r:rm\rightarrow m$ and $F_r:rm\rightarrow m$, $m\geq1$, 
subject to the relations
\begin{align*}
R_1= F_1 &= \id_n\quad,\\
R_r R_s &= R_{rs}\quad,\\
F_r F_s &= F_{rs}\quad,\\
R_r F_s &= F_s R_r\quad.
\end{align*}

Given a cyclotomic spectrum $T$ we now construct a functor 
$\mathbb{I}\rightarrow\mathscr{IS}$ by sending $n$ to the categorical $C_n$-fixed points $T^{C_n}$.
The map $F_r:T^{C_n}\rightarrow T^{C_m}$ is given by inclusion of categorical fixed points.
To construct $R_r$ we recall that there is a natural map $T^C\rightarrow \Phi^C T$, 
see construction~\ref{Constr:cattogeomfp}. 
We define $R_r$ as the composition
\[
T^{C_n}=(T^{C_r})^{C_m}\rightarrow (\Phi^{C_r}T)^{C^m}\xrightarrow{r_{C_r}} T^{C_m}\quad.
\]

\begin{Def}
The topological cyclic homology of $T$, $TC(T)$, is the orthogonal spectrum defined as
\[
TC(T)=\holim_{n\in\mathbb{I}} T^{C_n}\quad.
\]
\end{Def}

\begin{Lem}
If $T$ is a cyclotomic spectrum with involution, then $TC(T)$ is an orthogonal $\Z/2$-spectrum. 
\end{Lem}

\begin{proof}
The dihedral group of order $2n$ is the subgroup of $O(2)$ spanned by $C_n$ and the matrix 
$\begin{pmatrix}0&1\\1&0\end{pmatrix}$.
Since $T$ is an orthogonal $O(2)$-spectrum, we can restrict the action getting an orthogonal $D_{2n}$-spectrum.
Taking categorical $C_n$-fixed points, we get an $\Z/2=D_{2n}/C_n$ action on each $T^{C_n}$. 
Clearly, both $F_r$ and $R_r$ become $\Z/2$-maps.
\end{proof}

\begin{Rem}
If $T$ is not an $\Omega$-$G$-spectrum, $G=S^1$ or $O(2)$, 
then the $C$-fixed points, $T^C$, might have the wrong homotopy groups. 
See warning~V.3.6 in~\cite{MandellMay:02}. Hence, one
should apply $Q^{\cy}$ to $T$ before calculating $TC$. 
\end{Rem}

\begin{Prop}
The fibrant replacement functor $Q^{\cy}$ preserves cyclotomic spectra.
\end{Prop}

\begin{proof}
By construction $Q^{\cy}$ comes with a natural 
acyclic q-cofibration $T\rightarrow Q^{\cy}T$.
By lemma~\ref{Lem:Qcyandgeomfp} there is a natural transformation 
$\rho_C^*\Phi^C Q^{\cy}T\rightarrow Q^{\cy}\rho_C^*\Phi^C T$. And the following
diagram commutes:
\[
\begin{CD}
\rho_C^*\Phi^C T @= \rho_C^*\Phi^C T @>{r_C}>{\simeq}> T\\
@V{\simeq}VV @VVV @VV{\simeq}V\\
\rho_C^*\Phi^C Q^{\cy}T @>>> Q^{\cy}\rho_C^*\Phi^C T @>{Q^{\cy}r_C}>> Q^{\cy} T
\end{CD}\quad.
\]
The left vertical map is a genuine $\pi_*$-isomorphism since $\rho^*_C\Phi^C$ preserves
acyclic q-cofibrations.
We take the composition of the two bottom maps as the definition of $r_C$ for $Q^{\cy}T$.
It is automatically a cyclotomic $\pi_*$-isomorphism since $r_C$ is.
\end{proof}

\begin{Prop}\label{Prop:TCiso}
If a cyclotomic map $f:T_1\rightarrow T_2$ is a cyclotomic $\pi_*$-isomorphism
between $\Omega$-$G$-spectra, $G=S^1$ or $O(2)$,
then the induced map $TC(T_1)\rightarrow TC(T_2)$ is a $\pi_*$-isomorphism.
\end{Prop}

\begin{proof}
Due to homotopy invariance of homotopy limits, see proposition~\ref{Prop:htpyinvofholim},
it is enough to show that for each $n$ the map
\[
T_1^{C_n}\rightarrow T^{C_n}_2
\]
is a $\pi_*$-iso between $\Omega$-spectra.

It follows directly from the definitions that the categorical $H$-fixed points of
an $\Omega$-$G$-spectrum is an $\Omega$-spectrum. Hence, $T_1^{C_n}$ and $T_2^{C_n}$
are $\Omega$-spectra.

Furthermore,
\[
\pi_*T_1^{C_n}=\pi_*^{C_n}T_1\xrightarrow{\cong}\pi_*^{C_n}T_2=\pi_*T_2^{C_n}\quad.
\]
Here the map in the middle is an isomorphism since $T_1\rightarrow T_2$ is a cyclotomic $\pi_*$-isomorphism.
\end{proof}

\begin{Rem}
Let $\mathscr{F}$ be the family of finite normal subgroups of $S^1$ or $O(2)$.
One can define the notion of an $\Omega$-$\mathscr{F}$-spectrum. All statements above
probably remain true if replacing $\Omega$-$G$-spectra, $G=S^1$ or $O(2)$, by
$\Omega$-$\mathscr{F}$-spectra. Furthermore, one can probably show that these spectra
are the fibrant objects of the stable cyclotomic model structure on $G\mathscr{IS}$.
\end{Rem}

\section{Topological Hochschild homology}

Since the time when B\"okstedt defined $THH$ based on an idea of Goodwillie, the technology of
spectra has evolved so much that we now can use Goodwillie's idea as definition, see~\cite{Shipley:00}.
What is needed is a symmetric smash product for spectra. We write out the definition for 
the category of orthogonal spectra. Furthermore, we show that $THH(L)$ is a
cyclotomic spectrum, when $L$ is cofibrant. 

We also consider the involutive case.

\begin{Def}
Let $L$ be an orthogonal ring spectrum. Define $THH_{\bullet}(L)$ to be the simplicial
orthogonal spectrum with $q$-simplices
\[
THH_q(L)= L^{\wedge(q+1)}=L\wedge L\wedge \cdots\wedge L\quad, 
\]
and face and degeneracy maps given by
\[
d_i=\begin{cases}
\id_L^{\wedge i}\wedge\mu\wedge\id_L^{\wedge (q-i-1)}&\text{ for $0\leq i<q$,}\\
\mu\wedge\id_L^{\wedge(q-1)}\circ\pi_{L^{\wedge q},L}&\text{ for $i=q$,}
\end{cases}
\]
and
\[
s_i=\id_L^{\wedge (i+1)}\wedge\eta\wedge\id_L^{\wedge (q-i)}\quad.
\]
We define $THH(L)$ to be the geometrical realization of $THH_{\bullet}(L)$.
\end{Def}

To clarify the definition of $d_q$ we write it as the composition
\[
L^{\wedge(q+1)}=(L^{\wedge q})\wedge L\xrightarrow{\text{twist}}L\wedge (L^{\wedge q})=L\wedge L \wedge (L^{\wedge (q-1)})
\xrightarrow{\mu\wedge\id}L \wedge (L^{\wedge (q-1)}) = L^{\wedge q}\quad.
\]

\begin{Rem}
If the unit of the orthogonal ring spectrum, $\eta:S\rightarrow L$, is not a q-cofibration, then there is no reason to
expect the homotopy of $L$ and the homotopy of $THH(L)$ to be related to each other.
Hence, we will often restrict attention to such orthogonal ring spectra, and we call them \textit{cofibrant}.

Given an arbitrary orthogonal ring spectrum $L$, it can often be checked directly that
the unit $\eta:S\rightarrow L$ is a closed inclusion. If this is the case, then we may apply
the cofibrant replacement functor $\Gamma$ from theorem~\ref{Thm:cofrepl}, to produce
a new orthogonal ring spectrum $\Gamma L$, which is cofibrant.

This replacement procedure also works when $L$ comes with an involution. We must then define the involution on
$\Gamma L$ as the composed map
\[
\Gamma L\xrightarrow{\iota_L} \Gamma L\xrightarrow{\Gamma \iota} \Gamma L\quad.
\]
Here the first $\iota$ comes from theorem~\ref{Thm:cofrepl}, while the second $\iota$ is the involution on $L$,
see definition~\ref{Def:ringLwithinvolution}. Because $\Gamma$ is a skew-symmetric functor, it follows that
$\Gamma L$ is an orthogonal ring spectrum with involution. Furthermore, $\Gamma L$ is $\Z/2$-equivariantly 
cofibrant by proposition~\ref{Prop:Gammaequivariantly}.
\end{Rem}

We now specify $S^1$- and $O(2)$-actions on $THH(L)$.

\begin{Prop}
$THH_{\bullet}(L)$ is a cyclic orthogonal spectrum. If $L$ has involution, then $THH_{\bullet}(L)$ is dihedral.
\end{Prop}

\begin{proof}
We define the cyclic operator $t_q:L^{\wedge(q+1)}\rightarrow L^{\wedge(q+1)}$ as
\[
L^{\wedge(q+1)}=(L^{\wedge q})\wedge L\xrightarrow{\text{twist}}L\wedge (L^{\wedge q})=L^{\wedge(q+1)}\quad.
\]
If $L$ has involution $\iota:L\rightarrow L$, we can define the 
involutive operator $r_q:L^{\wedge(q+1)}\rightarrow L^{\wedge(q+1)}$ as
\[
L^{\wedge(q+1)}\xrightarrow{\iota^{\wedge(q+1)}}
L^{\wedge(q+1)}\xrightarrow{\text{permute}}L^{\wedge(q+1)}\quad.
\]
The arrow labeled ``permute'' permutes the order of the factors in the smash product as follows:
We label the factors from $0$'th to $q$'th. The $0$'th factor maps to the $0$'th factor, while the $i$'th factor, $i>0$, 
maps to the $(q+1-i)$'th factor.
\end{proof}

\begin{Cor}
$THH(L)$ is an orthogonal $S^1$-spectrum. If $L$ has involution, then $THH(L)$ is an orthogonal $O(2)$-spectrum.
\end{Cor}

We now show:

\begin{Thm}\label{thm:THHstructure}
Let $L$ be an orthogonal ring spectrum (with involution).
\begin{itemize}
\item[i)] If $S\rightarrow L$ is a q-cofibration, then there is an $S^1$-isomorphism
\[
r_C:\rho_C^*\Phi^C THH(L)\cong THH(L)
\]
for every finite subgroup $C$ of $S^1$, and $THH(L)$ is a cyclotomic spectrum. 
When $L$ has involution the isomorphism is $O(2)$-equivariant, 
and in this case $THH(L)$ is a cyclotomic spectrum with involution.
\item[ii)] If $L\rightarrow K$ is a $\pi_*$-isomorphism between 
cofibrant orthogonal ring spectra (with involution), then
\[
THH(L)\rightarrow THH(K)
\]
is a cyclotomic $\pi_*$-isomorphism. 
\end{itemize}
\end{Thm}

Before giving a proof, let us define topological cyclic homology:

\begin{Def}
The \textit{topological cyclic homology} of a cofibrant orthogonal ring spectrum $L$ (with involution)
is defined as $TC(THH(L))$. We abbreviate this notation, and write $TC(L)$.
\end{Def}

\begin{proof}\nopagebreak \paragraph{Part i):} We first consider 
the case of orthogonal ring spectra $L$ without involution.
It is sufficient to prove the statement in the case where $S\rightarrow L$ is 
a relative $FI$-cellular map.
Let $C$ be the finite subgroup of $S^1$ of order $r$. We will now construct the
isomorphism $r_C$ of genuine orthogonal $S^1$-spectra
\[
\rho^{*}_C\Phi^C THH(L)\cong THH(L)\quad.
\]
By definition, $THH(L)$ is the geometric realization of a cyclic orthogonal spectrum $THH_{\bullet}(L)$.
Edgewise subdivision gives an $S^1$-isomorphism $|THH_{\bullet}(L)|\cong|\sd_C THH_{\bullet}(L)|$.
We will construct $r_C$ by computing $\rho_C^*\Phi^C|\sd_C THH_{\bullet}(L)|$.

The geometric realization is level-wise, so we can use the filtration of the geometric realization of 
$r$-cyclic spaces given in construction~\ref{Constr:filtercrossedsimplreal}. 
By induction we shall prove that 
\[
\rho^{*}_C\Phi^C F_q^{\catDeltaC_r}|\sd_C THH_{\bullet}(L)|\cong F_q^{\catDeltaC}|THH_{\bullet}(L)|
\]
for all $q\geq0$. Letting $q$ go to infinity, this statement yields part i) of the theorem.

To prove the induction step we begin with a few calculations.
The $q$-simplices of $\sd_C THH_{\bullet}(L)$ are $L^{\wedge rq}$. Here $r$ is the order of $C$.
We now have
\begin{align*}
\rho^{*}_C\Phi^C\left(L^{\wedge rq}\wedge_{C_{rq}}{\Delta C_r^q}_+\right)
&\cong \rho^{*}_C\Phi^C\left( (L^{\wedge rq}\wedge \Delta^q_+ )\wedge_{C_{rq}} S^1_+\right)\\
&\cong \rho^*_C\left( (\Phi^C(L^{\wedge rq})\wedge \Delta^q_+ )\wedge_{C_q} S^1/C_+\right)\\
&\cong (\Phi^C(L^{\wedge rq})\wedge \Delta^q_+ )\wedge_{C_q} S^1_+\\
&\cong \Phi^C(L^{\wedge rq})\wedge_{C_q}{\Delta C^q}_+\\
&\cong L^{\wedge q} \wedge_{C_q}{\Delta C^q}_+\quad.
\end{align*}
In this calculation we have used the 
following facts:
\begin{itemize}
\item[] The topological $r$-cyclic $q$-simplex, $\Delta C_r^q$ is defined as $\Delta^q\times S^1$.
\item[] $\Phi^C(K\wedge A)=(\Phi^C K)\wedge A^C$, when $K$ is an orthogonal $S^1$-spectrum and $A$ a based $S^1$-space,
see proposition~\ref{Prop:geomfpofsuspandprod}. 
\item[] $\Phi^C(K\wedge_{C_{rq}} S^1_+)\cong (\Phi^CK)\wedge_{C_q} S^1/C_+$, when $K$ is an orthogonal $C_{rq}$-spectrum,
see proposition~\ref{Prop:geomfpandindsp}.
\item[] The diagonal map $L^{\wedge q}\cong\Phi^C L^{\wedge rq}$ is an isomorphism for cofibrant orthogonal spectra $L$,
see proposition~\ref{Prop:diagonalisomorphism}.
\end{itemize}

Recall from remark~\ref{Rem:degeniteratedprod}
the notation $s^rL^{\wedge rq-r}$ for the orthogonal $C_{rq}$-spectrum
\[
s^rL^{\wedge rq-r}=\bigcup_i L^{\wedge i-1}\wedge S \wedge L^{\wedge q-1} \wedge S \wedge  L^{\wedge q-1} \wedge S \wedge \cdots
\wedge S \wedge L^{\wedge q-i}\quad.
\]
Observe that the degenerate $q$-simplices of $\sd_C THH_{\bullet}(L)$ are exactly $s^rL^{\wedge rq-r}$.
By a calculation similar to that above, we get
\[
\rho^{*}_C\Phi^C\left(s^rL^{\wedge rq-r}\wedge_{C_{rq}}{\Delta C_r^q}_+\right)\cong sL^{\wedge q-1} \wedge_{C_q}{\Delta C^q}_+\quad.
\]

Restricting to the boundary of the topological $r$-cyclic $q$-simplex, $\partial\Delta C_r^q$, we get
\begin{align*}
\rho^{*}_C\Phi^C\left(L^{\wedge rq}\wedge_{C_{rq}}{\partial\Delta C_r^q}_+\right)
&\cong L^{\wedge q} \wedge_{C_q}{\partial\Delta C^q}_+\quad\text{and}\\
\rho^{*}_C\Phi^C\left(s^rL^{\wedge rq-r}\wedge_{C_{rq}}{\partial\Delta C_r^q}_+\right)
&\cong sL^{\wedge q-1} \wedge_{C_q}{\partial\Delta C^q}_+\quad.
\end{align*}

Consider the diagram
\[
\begin{CD}
s^rL^{\wedge rq-r}\wedge_{C_{rq}}{\partial\Delta C_r^q}_+ 
@>>> L^{\wedge rq}\wedge_{C_{rq}}{\partial\Delta C_r^q}_+\\
@VVV @VVV\\
s^rL^{\wedge rq-r}\wedge_{C_{rq}}{\Delta C_r^q}_+
@>>>L^{\wedge rq}\wedge_{C_{rq}}{\Delta C_r^q}_+
\end{CD}\quad.
\]
As an orthogonal spectrum $L^{\wedge rq}\wedge_{C_{rq}}{\Delta C_r^q}_+$ has an $FI$-cellular
structure such that the three other orthogonal spectra are $FI$-cellular subspectra. It follows that
the map
\[
L^{\wedge rq}\wedge_{C_{rq}}{\partial\Delta C_r^q}_+\cup_{s^rL^{\wedge rq-r}\wedge_{C_{rq}}{\partial\Delta C_r^q}_+}
s^rL^{\wedge rq-r}\wedge_{C_{rq}}{\Delta C_r^q}_+
\rightarrow L^{\wedge rq}\wedge_{C_{rq}}{\Delta C_r^q}_+
\]
is a closed inclusion. Now consider the diagram
\[
L^{\wedge rq}\wedge_{C_{rq}}{\Delta C_r^q}_+ \leftarrow 
L^{\wedge rq}\wedge_{C_{rq}}{\partial\Delta C_r^q}_+\cup s^rL^{\wedge rq-r}\wedge_{C_{rq}}{\Delta C_r^q}_+
\rightarrow F_q^{\catDeltaC_r}|\sd_C THH_{\bullet}(L)|\quad.
\]
Since the left map is a closed inclusion, it follows by
proposition~\ref{Prop:fundresultofgeomfp} that $\rho_C^*\Phi^C$ of the pushout
is the pushout of $\rho_C^*\Phi^C$ applied to the diagram.
At last we look at the following diagram:
\footnotesize
\[
\begin{CD}
\rho_C^*\Phi^C\left(L^{\wedge rq}\wedge_{C_{rq}}{\Delta C_r^q}_+\right)
&{\cong}& L^{\wedge q} \wedge_{C_q}{\Delta C^q}_+\\
@AAA @AAA\\
\rho_C^*\Phi^C\left(L^{\wedge rq}\wedge_{C_{rq}}{\partial\Delta C_r^q}_+\cup s^rL^{\wedge rq-r}\wedge_{C_{rq}}{\Delta C_r^q}_+\right)
&{\cong}& L^{\wedge q} \wedge_{C_q}{\partial\Delta C^q}_+ \cup sL^{\wedge q-1} \wedge_{C_q}{\Delta C^q}_+\\
@VVV @VVV\\
\rho_C^*\Phi^CF_{q-1}^{\catDeltaC_r}|\sd_C THH_{\bullet}(L)| 
&{\cong}& F_{q-1}^{\catDeltaC}|THH_{\bullet}(L)|
\end{CD}\quad.
\]
\normalsize
By the calculations above, the top and the middle horizontal maps are isomorphisms.
By the induction hypothesis, the bottom horizontal map is an isomorphism. It follows that
the map of column-wise pushouts, 
\[
\rho^{*}_C\Phi^C F_q^{\catDeltaC_r}|\sd_C THH_{\bullet}(L)|\cong F_q^{\catDeltaC}|THH_{\bullet}(L)|\quad,
\]
is an isomorphism.

Now assume that $L$ is an orthogonal ring spectrum with involution.
We prove that $THH(L)$ is a cyclotomic spectrum with involution by an argument similar to that above.
Edgewise subdivision works also in the dihedral case, and 
we construct the $O(2)$-isomorphism
\[
r_C:\rho^{*}_C\Phi^C THH(L)\cong THH(L)
\]
by induction over the $O(2)$-equivariant filtration 
for the geometric realization of 
$r$-dihedral spaces provided by construction~\ref{Constr:filtercrossedsimplreal}.
Since the diagonal map is dihedral for $L$ with involution, see
proposition~\ref{Prop:diagonalisomorphisminvolutive},
the rest of the argument works exactly as before.

\paragraph{Part ii):} By proposition~\ref{Prop:piisoiscyclpiiso}, it is enough to show that the induced map
$THH(L)\rightarrow THH(K)$ is non-equivariantly a $\pi_*$-isomorphism.
Since both $S\rightarrow L$ and $S\rightarrow K$ are q-cofibrations, it 
follows that $THH_{\bullet}(L)$ and $THH_{\bullet}(K)$
are good simplicial orthogonal spectra. 
It remains to show that the induced map
$THH_{\bullet}(L)\rightarrow THH_{\bullet}(K)$ is a $\pi_*$-iso in each simplicial degree,
see proposition~\ref{Prop:piisoofrealiz}.
We can factor the map
$THH_q(L)=L^{\wedge (q+1)}\rightarrow K^{\wedge (q+1)}=THH_q(K)$ as
\[
L^{\wedge (q+1)}\rightarrow L^{\wedge q}\wedge K\rightarrow L^{\wedge (q-1)}\wedge K^{\wedge 2}\rightarrow\cdots\rightarrow L\wedge K^{\wedge q}
\rightarrow K^{\wedge (q+1)}\quad.
\]
The smash product of cofibrant orthogonal spectra is cofibrant by proposition~\ref{Prop:os:qcofsquare}, and
smashing with cofibrant orthogonal spectra preserves $\pi_*$-isomorphisms, see proposition~\ref{Prop:os:wedge}.
Hence, each map in the sequence above is a $\pi_*$-isomorphism. The result follows.
\end{proof}

\section{Matrices over an orthogonal ring spectrum}

An important ingredient when constructing a trace map from $K$-theory to $THH$ or $TC$,
is the definition of a matrix ring. Since we focus on orthogonal ring spectra with involution,
we need a construction which is involutive. In ordinary linear algebra we have such an involution,
namely the conjugate transposed matrix. However, if we consider the customary definition
of the matrix-FSP, $M_n(L)=F(\mathbf{n}_+,\mathbf{n}_+\wedge L)$, we see that transposition is
not well-defined. Hence a modification of the definition is required.

The purpose of this section is to provide a construction of $(n\times n)$-matrices 
for orthogonal ring spectra $L$, which also works when $L$ has involution.
Let me make a list of our hopes and needs regarding the construction:

\begin{itemize}
\item[] $M_n$ should be an endofunctor on orthogonal spectra.
\item[] Up to $\pi_*$-isomorphism $\bigvee_{n^2}L$, $M_n(L)$ and $L^{\times n^2}$ should be the same, 
at least when $L$ is cofibrant.
\item[] We want a matrix multiplication $M_n(L)\wedge M_n(L)\rightarrow M_n(L)$.
\item[] Direct sum of matrices should give a functor $M_{n_1}(L)\times M_{n_2}(L)\rightarrow M_{n_1+n_2}(L)$.
\item[] For $L$ with involution taking the transposed involuted matrix should be an involution on $M_n(L)$.
\item[] We want a trace map from $M_n(L)$ to some additive model for $L$.
\end{itemize}

Recall the concept of induced functors on orthogonal spectra, see subsection~\ref{subsect:indfunct}.
We first define a continuous endofunctor on $\Top_*$, which we also will denote by $M_n$. Then we define 
$M_n$ on orthogonal spectra as the induced functor.

\begin{Def}\label{Def:matrix}
For a based space $X$ let $M_n(X)$ be the 
subspace of $X^{\times n^2}$ consisting for those matrices
where each row contains at most one element 
different from $*$ and each column contains at most one element different form $*$.
\end{Def}

Using formulas we may write $M_n(X)$ as
\begin{multline*}
M_n(X)=\{(x_{i\,j})\in X^{\times n^2}\;|\;\\
 \text{if $x_{i_0\,j_0}\neq *$, 
then $x_{i_0\,j}=*$ and $x_{i\,j_0}=*$ for all $i\neq i_0$ and
$j\neq j_0$.}\}
\end{multline*}

It is clear that $M_n$ is a continuous functor. Since $M_n(*)=*$, there is a canonical right assembly
\[
\sigma_{X,Y}:M_n(X)\wedge Y\rightarrow M_n(X\wedge Y)\quad,
\]
see page~208 in~\cite{Madsen:94}. In our case we can easily write down a formula for $\sigma$:
Let $(x_{i\,j})$ be a matrix in $M_n(X)$ and $y$ a point in $Y$, then $\sigma((x_{i\,j}),y)$ is the
matrix $(z_{i\,j})$ in $M_n(X\wedge Y)$, where $z_{i\,j}=(x_{i\,j},y)\in X\wedge Y$.
Similarly, there is a left assembly $\bar{\sigma}$.

We now describe the structure of the functor $M_n$:

\begin{Lem}\label{Lem:Mnstr}
There are natural transformations
\begin{align*}
1_X&:X\rightarrow M_n(X)\quad,\\
\mu_{X,Y}&: M_n(X)\wedge M_n(Y)\rightarrow M_n(X\wedge Y)\quad,\text{ and}\\
\iota_X&:M_n(X)\rightarrow M_n(X)\quad,
\end{align*}
such that 
\begin{itemize}
\item[] the composition $M_n(X)\wedge Y\xrightarrow{\id\wedge 1_Y}M_n(X)\wedge M_n(Y)\xrightarrow{\mu_{X,Y}} M_n(X\wedge Y)$
is equal to the right assembly,
\item[] the composition $X \wedge M_n(Y)\xrightarrow{1_X\wedge \id}M_n(X)\wedge M_n(Y)\xrightarrow{\mu_{X,Y}} M_n(X\wedge Y)$
is equal to the left assembly,
\item[] $\mu$ is associative, 
\item[] $\iota^2=\id$, and
\item[] $\iota$ anti-commutes with $\mu$, this means that the following diagram commutes:
\[
\begin{CD}
M_n(X)\wedge M_n(Y) @>{\iota_X\wedge \iota_Y}>> M_n(X)\wedge M_n(Y) @>{\text{twist}}>> M_n(Y)\wedge M_n(X)\\
@V{\mu_{X,Y}}VV && @VV{\mu_{Y,X}}V\\
M_n(X\wedge Y) @>{\iota_{X\wedge Y}}>> M_n(X\wedge Y) @>{M_n(\text{twist})}>> M_n(Y\wedge X)
\end{CD}\quad.
\]
\end{itemize}
\end{Lem}

Another way to phrase this lemma is to say that $M_n$ is an FSP with involution.

\begin{proof}
The unit $1_X:X\rightarrow M_n(X)$ is defined by sending $x\in X$ to the $(n\times n)$-matrix
\[
\begin{bmatrix}
x & * & \cdots & *\\
* & x & \cdots & *\\
\vdots & \vdots & \ddots & \vdots\\
* & * & \cdots & x
\end{bmatrix}\quad,
\]
which has $x$ on the diagonal and $*$ elsewhere. Multiplication $\mu_{X,Y}:M_n(X)\wedge M_n(Y)\rightarrow M_n(X\wedge Y)$
is given as ordinary matrix multiplication. Explicitly this is given by sending $(x_{i\,j})$ and $(y_{i\,j})$
to $(z_{i\,j})$, where
\[
z_{i\,j}=\begin{cases}
(x_{i\,k},y_{k\,j})&\text{if $k$ is such that $x_{i\,k}\neq *$ and $y_{k\,j}\neq *$,}\\
*&\text{otherwise.}
\end{cases}
\]
By direct computation it is easily seen that $\mu\circ(1\wedge\id)$ and $\mu\circ(\id\wedge1)$ are the left and right
assemblies respectively. And easy calculations also show that $\mu$ is associative.

We define the involution $\iota$ by transposition. $\iota_X(x_{i\,j})$ is $(x'_{i\,j})$, where $x'_{i\,j}=x_{j\,i}$.
Clearly, $\iota^2=\id$. To check that $\iota$ anti-commutes with $\mu$, we take matrices $(x_{i\,j})$ and $(y_{i\,j})$
in $M_n(X)$ and $M_n(Y)$ respectively. Let $(z_{i\,j})$ be the matrix of $\mu(\iota(y_{i\,j}),\iota(x_{i\,j}))$ and
$(z'_{i\,j})=\iota(\mu((x_{i\,j}),(y_{i\,j})))$. Calculating we see that
\[
z_{i\,j}=\begin{cases}
(y_{k\,i},x_{j\,k})&\text{if $k$ is such that $x_{j\,k}\neq *$ and $y_{k\,i}\neq *$,}\\
*&\text{otherwise,}
\end{cases}
\]
and
\[
z'_{i\,j}=\begin{cases}
(x_{j\,k},y_{k\,i})&\text{if $k$ is such that $x_{j\,k}\neq *$ and $y_{k\,i}\neq *$,}\\
*&\text{otherwise.}
\end{cases}
\]
And we see that the diagram commutes.
\end{proof}

\begin{Cor}
If $L$ is an orthogonal ring spectrum $L$, then $M_n(L)$ is also an orthogonal ring spectrum, 
and if $L$ has involution, then
$M_n(L)$ has an induced involution.
\end{Cor}

\begin{proof}
We use the external description of the smash product, example~\ref{exa:extsmashdescr}. The unit for $M_n(L)$
is defined as the composition
\[
S^V\rightarrow L(V)\xrightarrow{1_{L(V)}} M_n(L(V))\quad.
\]
The multiplication for $M_n(L)$, defined externally, is given by
\footnotesize
\[
M_n(L(V_1))\wedge M_n(L(V_2))\xrightarrow{\mu_{L(V_1),L(V_2)}}M_n(L(V_1)\wedge L(V_2))\xrightarrow{M_n(\text{multiplication})} M_n(L(V_1\oplus V_2))\quad.
\]
\normalsize
And in the case $L$ has involution, the induced involution on $M_n(L)$ is given by
\[
M_n(L(V))\xrightarrow{\iota_{L(V)}}M_n(L(V))\xrightarrow{M_n(\text{involution})} M_n(L(V))\quad.
\]
The structure of $M_n$ described in lemma~\ref{Lem:Mnstr} ensures that $M_n(L)$ is an orthogonal ring spectrum (with involution).
\end{proof}

\begin{Exa}[Direct sum]
The direct sum $M_{n_1}(L)\times M_{n_2}(L)\rightarrow M_{n_1+n_2}(L)$ is easily defined. 
First observe that direct sum $M_{n_1}(X)\times M_{n_2}(X)\rightarrow M_{n_1+n_2}(X)$ is defined
for based spaces $X$ by the ordinary direct sum of matrices. Applying the concept of induced functors
we get the direct sum for matrices of orthogonal ring spectra.
\end{Exa}

Next we want to compare the weak homotopy type of $M_n(L)$ to $\bigvee_{n^2} L$ and $L^{\times n^2}$.
We clearly have maps
\[
\bigvee_{n^2} L\rightarrow M_n(L) \rightarrow L^{\times n^2}\quad,
\]
and when $L$ is cofibrant, the composition is a $\pi_*$-iso by proposition~\ref{Prop:veepiisotimes}.
Our strategy is to use corollary~\ref{Cor:indfunct} to show that the first map
also is a $\pi_*$-iso, given that $L$ is cofibrant. 

The third condition in corollary~\ref{Cor:indfunct} demands that the functors must commute with
colimit over sequences of cofibrations. We check this for the functors above.
Assume that $X_0\rightarrow X_1\rightarrow X_2\rightarrow \cdots$ is a sequence of 
cofibrations of spaces, and $X$ is the colimit. 
Clearly we have that $\colim_i \left(\bigvee_{n^2} X_i\right)=\bigvee_{n^2}X$. 
(We can describe the wedge as a colimit, and interchanging colimits does not affect the result.)
Theorem~10.3 in~\cite{Steenrod:67} also holds for the category of compactly generated spaces defined in~\cite{McCord:69}.
Therefore, we also have
\[
\colim_i \left( X_i^{\times n^2}\right)=X^{\times n^2}\quad.
\]
What we really are saying is that two a priori different topologies on the same set actually coincide.
It is easy to see that $\colim_i M_n(X_i)$ is equal to $M_n(X)$ as sets. But the topology of
$\colim_i M_n(X_i)$ is the subspace topology from $\colim_i \left( X_i^{\times n^2}\right)$, while
$M_n(X)$ has the subspace topology from $X^{\times n^2}$. However, the equality above implies that
\[
\colim_i M_n(X_i)=M_n(X)
\]
as topological spaces (=compactly generated spaces).

The following proposition checks the second condition for $M_n(-)$.

\begin{Prop}
If $A\rightarrow X$ is an unbased closed cofibration of spaces, then $M_n(A)\rightarrow M_n(X)$ is also
an unbased closed cofibration.
\end{Prop}

\begin{proof}
Represent the cofibration
$A\rightarrow X$ by a homotopy $H:X\times I\rightarrow X$ and $\phi:X\rightarrow I$, see remark~\ref{Rem:cofibdef}.
Define $\bar{H}:M_n(X)\times I\rightarrow M_n(X)$ by
\[
\bar{H}((x_{i\,j}),t)=(H(x_{i\,j},t))\quad,
\]
and $\bar{\phi}:M_n(X)\rightarrow I$ by
\[
\bar{\phi}(x_{i\,j})=\sup_{i,j}\phi(x_{i\,j})\quad.
\]
Clearly, $\bar{H}$ is a homotopy rel $M_n(A)$ with $\bar{H}(-,0)=\id_{M_n(X)}$, and $M_n(A)\subseteq \bar{\phi}^{-1}(0)$.
Assume that $t>\bar{\phi}(x_{i\,j})$, then for each $i$ and $j$ we have $t>\phi(x_{i\,j})$ and thus $H(x_{i\,j},t)\in A$.
It follows that $\bar{H}((x_{i\,j}),t)\in M_n(A)$.

This shows that $\bar{H}$ and $\bar{\phi}$ represent $M_n(A)\rightarrow M_n(X)$ as an unbased cofibration.
\end{proof}

We immediately get the following two corollaries:

\begin{Cor}
The induced functor $M_n$ on orthogonal spectra preserves l-cofibrations.
\end{Cor}

\begin{Cor}
If $L$ is a well-pointed orthogonal spectrum, then also
$M_n(L)$ is well-pointed.
\end{Cor}

\begin{Rem}
Observe that we do not claim that $M_n(L)$ is cofibrant. If we want a cofibrant version,
then we just apply the cofibrant replacement functor $\Gamma$ from theorem~\ref{Thm:cofrepl}.
Furthermore, whenever the unit $\eta:S\rightarrow L$ of an orthogonal ring spectrum is an 
l-cofibration, we have that $S\rightarrow M_n(L)$ is an l-cofibration and $S\rightarrow \Gamma M_n(L)$
is a q-cofibration.
\end{Rem}

To compare $M_n(X)$ and $\bigvee_{n^2}X$ up to homotopy, we now provide a filtration. 
Define $M_n^k(X)$ to be the subspace of
$M_n(X)$ consisting of those matrices with at most $k$ elements different form $*$. It is easily seen that
$M_n^1(X)$ is equal to $\bigvee_{n^2}X$, while $M_n^n(X)$ equals $M_n(X)$. 
The key lemma for analyzing this filtrations is:

\begin{Lem}\label{Lem:filtofMn}
For well-pointed $X$ there is a natural cofiber sequence
\[
M_n^{k-1}(X)\rightarrow M_n^k(X)\rightarrow \bigvee_{A}X^{\wedge k}\quad,
\]
where the wedge is indexed over a finite set $A$. 
\end{Lem}

\begin{proof}
Let $A$ be the set of maps, $f$, from $\{1,\ldots,k\}$ to $\{1,\ldots,n\}^2$ such that 
both
\[
\pr_1\circ f:\{1,\ldots,k\}\rightarrow \{1,\ldots,n\}
\]
is strictly increasing and
\[
\pr_2\circ f:\{1,\ldots,k\}\rightarrow \{1,\ldots,n\}
\]
is injective. For each such $f$ we construct a map $f_*:X^{\times k}\rightarrow M_n^k(X)$ by
\[
(x_1,\ldots,x_k)\mapsto (y_{i\,j})\quad\text{where}\quad y_{i\,j}=\begin{cases}
x_l&\text{if $f(l)=(i,j)$ and}\\
*&\text{otherwise.}
\end{cases}
\]
Let $sX^{\times k-1}$ be the subspace of $X^{\times k}$ consisting of the tuples $(x_1,\ldots,x_k)$
where at least one $x_l=*$. By Steenrod's product theorem for cofibrations we know that (by induction)
$sX^{\times k-1}\rightarrow X^{\times k}$ is a cofibration. Also observe that the image of $sX^{\times k-1}$ 
in $M_n^k(X)$ under $f_*$ actually lies in $M_n^{k-1}(X)$. Furthermore the diagram
\[
\begin{CD}
\coprod_A sX^{\times k-1} @>>> \coprod_A X^{\times k}\\
@V{f_*}VV @VV{f_*}V\\
M_n^{k-1}(X) @>>> M_n^{k}(X)
\end{CD}
\]
is pushout. The lemma follows by the observation that $\bigvee_{A}X^{\wedge k}$ is the cofiber of the top row.
\end{proof}

By counting one can check that $A$ contains $\binom{n}{k}^2 k!$ elements.

Using the filtration we prove the following result regarding the connectivity of the map
$\bigvee_{n^2}X\rightarrow M_n(X)$:

\begin{Prop}\label{Prop:WMhtpyequiv}
If $X$ is $r$-connected and well-pointed, then the map $\bigvee_{n^2}X\rightarrow M_n(X)$ is $2r$-connected.
\end{Prop}

\begin{proof}
We prove by induction on $r$ and $k$ that $M_n^k(X)$ is $r$-connected and the map
\[
M_n^{k-1}(X)\rightarrow M_n^k(X)\quad,k\geq2\quad,
\]
is $2r$-connected when $X$ is $r$-connected.

For $k=1$ observe that $M_n^1(X)=\bigvee_{n^2}X$. Therefore $M_n^1(X)$ is $r$-connected whenever $X$ is.

For $r=-1$ there is nothing to prove. 
For $r=0$ we will give a direct argument that shows that all $M_n^k(X)$ are $0$-connected.
Consequently, the maps $M_n^{k-1}(X)\rightarrow M_n^k(X)$ are all $0$-connected.
Let $(x_{i\,j})$ be a matrix in $M_n^k(X)$. Since $X$ is path-connected, we can for each $x_{i\,j}$
choose a path $\gamma_{i\,j}$ from $x_{i\,j}$ to $*$. When $x_{i\,j}=*$, we let the path be constant.
Then
\[
t\mapsto\big( \gamma_{i\,j}(t)\big)
\]
will be a path from $(x_{i\,j})$ to the base point in $M_n^k(X)$. Hence this space is also path-connected.

Let $r\geq1$ and $k\geq 2$ and assume that the induction hypothesis holds for smaller $r$ and $k$.
By lemma~\ref{Lem:filtofMn} above, the map
\[
M_n^{k-1}(X)\rightarrow M_n^k(X)
\]
is an unbased closed cofibration with cofiber $\bigvee_{A}X^{\wedge k}$. And it follows 
from proposition~\ref{Prop:connofsmash} that this cofiber is at least $(2r+1)$-connected.
The induction hypothesis says that $M_n^{k-1}(X)$ is $r$-connected and the pair $(M_n^k(X),M_n^{k-1}(X))$ is
$(2r-2)$-connected. We now apply proposition~\ref{Prop:homotopycofseq} and get that
\[
\pi_q(M_n^k(X),M_n^{k-1}(X))\rightarrow \pi_q(\bigvee_{A}X^{\wedge k})
\]
is an isomorphism for $q<3r-1$. If $r>1$, we immediately get that the map $M_n^{k-1}(X)\rightarrow M_n^k(X)$ is
$2r$-connected and $M_n^k(X)$ is $r$-connected.

If $r=1$, the statement above only says that the map $M_n^{k-1}(X)\rightarrow M_n^k(X)$ is $1$-connected. But in this
case we can apply the proposition again with the improved input that the pair $(M_n^k(X),M_n^{k-1}(X))$ is
$1$-connected. And we get that also $\pi_2(M_n^k(X),M_n^{k-1}(X))\rightarrow \pi_2(\bigvee_{A}X^{\wedge k})=0$
is an isomorphism. Thus $M_n^{k-1}(X)\rightarrow M_n^k(X)$ is actually $2$-connected and $M_n^k(X)$ is $1$-connected.
\end{proof}

\begin{Cor}
If $L$ is a cofibrant orthogonal spectrum, then $\bigvee_{n^2}L\rightarrow M_n(L)$ is a $\pi_*$-isomorphism.
\end{Cor}

\begin{proof}
The considerations and results above verify the conditions of corollary~\ref{Cor:indfunct}.
Hence, the map of induced functors is a $\pi_*$-isomorphism.
\end{proof}

\section{Cross product formula for $THH$}

In this section we will derive a cross product formula for $THH$. Unfortunately, 
our result, proposition~\ref{Prop:THHcrossprodformula} below, is not as strong as we would like.
The author suggests two ways to improve the conclusion, see remark~\ref{Rem:improvecrossprodformula}.
We begin the section by discussing the cross product of orthogonal ring spectra.

\begin{Lem}
If $L$ and $K$ are orthogonal ring spectra, then $L\times K$ is also an orthogonal ring spectra. If
$L$ and $K$ have involutions, then also $L\times K$ comes with an involution. 
The projections are maps of orthogonal ring spectra (with involution).
\end{Lem}

\begin{proof}
Observe that maps into a categorical product $L\times K$ correspond to pairs of maps, 
one into each factor. In other words, for any orthogonal spectrum $X$ there is a 
homeomorphism 
\[
\mathscr{IS}(X,L\times K)\cong \mathscr{IS}(X,L)\times \mathscr{IS}(X,K)\quad.
\]
When $L$ and $K$ are orthogonal ring spectra, we define 
the unit for $L\times K$ as the map determined by $\eta_L:S\rightarrow L$ and $\eta_K:S\rightarrow K$.
The multiplication is determined by the two maps
\begin{align*}
(L\times K)\wedge (L\times K)&\xrightarrow{\pr_L\wedge\pr_L} L\wedge L\xrightarrow{\mu_L} L\quad\text{and}\\
(L\times K)\wedge (L\times K)&\xrightarrow{\pr_K\wedge\pr_K} K\wedge K\xrightarrow{\mu_K} K\quad.
\end{align*}
In case $L$ and $K$ have involutions, we define an
involution on $L\times K$ by the map $\iota_L\times \iota_K$.

To check that $L\times K$ is an orthogonal ring spectrum (with involution), 
one needs to see that certain diagrams commute. This is an easy
computation, done by projecting the diagrams to $L$ and $K$, where they commute by assumption.
\end{proof}

However, it is not clear that the cross product of cofibrant orthogonal ring spectra is cofibrant.
So we provide a cofibrant replacement:

\begin{Lem}
If $L$ is an orthogonal ring spectrum and $\eta:S\rightarrow L$ is a closed inclusion, 
then $\Gamma L$ is a cofibrant orthogonal ring spectrum.
$\Gamma L$ is involutive whenever $L$ has involution. The natural map $\Gamma L\rightarrow L$ is a map of
orthogonal ring spectra (with involution).
\end{Lem}

\begin{proof}
The unit is defined as the composition
\[
S\rightarrow \Gamma S\xrightarrow{\Gamma\eta} \Gamma L\quad,
\]
the multiplication is given by
\[
\Gamma L\wedge \Gamma L\xrightarrow{\subseteq}\Gamma (L\wedge L)\xrightarrow{\Gamma\mu}\Gamma L\quad,
\]
and the involution is the composition
\[
\Gamma L \xrightarrow{\iota_L} \Gamma L \xrightarrow{\Gamma\iota} \Gamma L\quad,
\]
where the first map comes from the natural transformation $\iota:\Gamma\rightarrow\Gamma$,
and the second map is $\Gamma$ applied to the involution on $L$.

Commutativity of the required diagrams follows since $\Gamma$ is lax symmetric monoidal.
$S\rightarrow \Gamma L$ is a q-cofibration since $S\rightarrow \Gamma S$ is a q-cofibration and $\Gamma$ applied
to the closed inclusion $\eta:S\rightarrow L$ is a q-cofibration.
\end{proof}

We now state or cross product formula:

\begin{Prop}\label{Prop:THHcrossprodformula}
Assume that $L$ and $K$ are cofibrant orthogonal ring spectra (with involution).
The $S^1$-map ($O(2)$-map)
\[
THH(\Gamma(L\times K))\rightarrow THH(L)\times THH(K)
\]
induced by the projections $\Gamma(L\times K)\rightarrow L$ and 
$\Gamma(L\times K)\rightarrow K$ is non-equivariantly a $\pi_*$-isomorphism.
\end{Prop}

We will show this by comparing both 
$THH_{\bullet}(\Gamma(L\times K))$ and $THH_{\bullet}(L)\times THH_{\bullet}(K)$
to a third cyclic (dihedral) orthogonal spectrum, called $T_{\bullet}(L,K)$.
Together corollary~\ref{cor:THHprodfirst} and lemma~\ref{lem:THHprodsecond} 
below prove the proposition.

We begin by defining $T_{\bullet}(L,K)$. Let the $q$-simplices be the $\times$-product 
of all $(q+1)$-fold $\wedge$-products where each factor is either $L$ or $K$.
Let us write up explicitly what we get for small $q$:
\begin{align*}
T_0(L,K)&=L\times K\\
T_1(L,K)&=(L\wedge L)\times (L\wedge K)\times (K\wedge L)\times (K\wedge K)\\
T_2(L,K)&=(L\wedge L\wedge L)\times (L\wedge L\wedge K)\times (L\wedge K\wedge L)\times (L\wedge K\wedge K)\\
&\quad \times (K\wedge L\wedge L)\times (K\wedge L\wedge K)\times (K\wedge K\wedge L)\times (K\wedge K\wedge K)\quad.
\end{align*}
We define the cyclic (dihedral) structure by considering each $\times$-factor separately. Consider a factor
$X_0\wedge X_1\wedge\cdots \wedge X_q$, where each $X_i$ is $L$ or $K$. The face operator $d_i$ for $i<q$ tries to multiply
$X_i$ and $X_{i+1}$. We define
\small
\[
d_i:X_0\wedge\cdots \wedge X_q\rightarrow\begin{cases}
X_0\wedge\cdots \wedge X_{i-1}\wedge L\wedge X_{i+2}\wedge\cdots \wedge X_q&\text{if $X_i=X_{i+1}=L$,}\\
X_0\wedge\cdots \wedge X_{i-1}\wedge K\wedge X_{i+2}\wedge\cdots \wedge X_q&\text{if $X_i=X_{i+1}=K$, and}\\
*&\text{otherwise.}
\end{cases}
\]
\normalsize
Here the map is induced by $\mu_L$ in the first case and $\mu_K$ in the second.
To define $d_q$ we try to multiply $X_q$ with $X_0$, using $\mu_L$ if $X_q=X_0=L$, $\mu_K$ if $X_q=X_0=K$ and mapping to $*$ otherwise.
The degeneracy map is given by
\small
\[
s_i:X_0\wedge\cdots \wedge X_q\rightarrow(X_0\wedge\cdots \wedge X_{i}\wedge L\wedge X_{i+1}\wedge\cdots \wedge X_q)\times
(X_0\wedge\cdots \wedge X_{i}\wedge K\wedge X_{i+1}\wedge\cdots \wedge X_q)
\]
\normalsize
using $\eta_L$ into the first factor and $\eta_K$ into the second.
The cyclic operator permutes the factors:
\[
t_q:X_0\wedge\cdots \wedge X_q\rightarrow X_q\wedge X_0\wedge\cdots \wedge X_{q-1}\quad.
\]
If $L$ and $K$ come with involutions, we can define the involutive operator by applying $\iota$ to
each $X_i$ an then permute:
\[
r_q:X_0\wedge\cdots \wedge X_q\xrightarrow{\iota\wedge\cdots\wedge\iota}
X_0\wedge\cdots \wedge X_q\xrightarrow{\text{permute}} X_0\wedge X_q\wedge \cdots\wedge X_1\quad.
\]

There is a map
\[
\pr:T_{\bullet}(L,K)\rightarrow THH_{\bullet}(L)\times THH_{\bullet}(K)
\]
defined on factors by
\[
X_0\wedge\cdots \wedge X_q\rightarrow\begin{cases}
L^{\wedge (q+1)}=THH_q(L)&\text{if $X_0=\cdots=X_q=L$,}\\
K^{\wedge (q+1)}=THH_q(K)&\text{if $X_0=\cdots=X_q=K$, and}\\
*&\text{otherwise.}
\end{cases}
\]
In the opposite direction we have the inclusion
\[
\incl: THH_{\bullet}(L)\times THH_{\bullet}(K)\rightarrow T_{\bullet}(L,K)\quad.
\]

\begin{Lem}
The map $\pr$ is cyclic (dihedral), $\incl$ is presimplicial,
$\pr\incl=\id$, while $\incl\pr\simeq\id$ via a presimplicial homotopy.
\end{Lem}

\begin{proof}
The first three statements are obvious. 
To prove the last statement we define a presimplicial homotopy on factors as follows:
\small
\[
h_i:X_0\wedge\cdots\wedge X_q\rightarrow\begin{cases}
X_0\wedge\cdots\wedge X_i\wedge L \wedge X_{i+1}\wedge\cdots\wedge X_q&\text{if $X_0=\cdots=X_i=L$,}\\
X_0\wedge\cdots\wedge X_i\wedge K \wedge X_{i+1}\wedge\cdots\wedge X_q&\text{if $X_0=\cdots=X_i=K$, and}\\
*&\text{otherwise.}
\end{cases}
\]
\normalsize
Here $L$ is inserted using $\eta_L$ in the first case, in the second case $K$ is inserted using $\eta_K$.
We see that
\[
d_0 h_0=\id\quad\text{and}\quad d_{q+1}h_q=\incl\pr\quad.
\]
Assume $i<j$, then $d_ih_j$ and $h_{j-1}d_i$ are $*$ if not $X_0=\cdots=X_j$. When $X_0=\cdots=X_j$, it is easy to check that
$d_ih_j=h_{j-1}d_i$. Similarly, one shows that $d_ih_j=h_jd_{i-1}$ for $i>j+1$. At last we check that on the factor $X_0\wedge\cdots\wedge X_q$
we have
\[
d_ih_i=\begin{cases}\id&\text{if $X_0=\cdots=X_i$, and}\\
*&\text{otherwise,}\end{cases}
\]
and also
\[
d_ih_{i-1}=\begin{cases}\id&\text{if $X_0=\cdots=X_i$, and}\\
*&\text{otherwise.}\end{cases}
\]
Thus we have a presimplicial homotopy.
\end{proof}

\begin{Cor}\label{cor:THHprodfirst}
Assume that $L$ and $K$ are cofibrant orthogonal ring spectra (with involution). 
Then there is an $S^1$-map ($O(2)$-map) $\pr:|T_{\bullet}(L,K)|\rightarrow THH(L)\times THH(K)$,
which is non-equivariantly a $\pi_*$-isomorphism.
\end{Cor}

\begin{proof}
Since $L$ and $K$ are cofibrant, it follows that $T_{\bullet}(L,K)$ and $THH_{\bullet}(L)\times THH_{\bullet}(K)$
are good. The result now follows from the lemma above, 
proposition~\ref{Prop:sos:comparerealizarions}, and the fact that
a presimplicial homotopy induces a homotopy on presimplicial realization.
\end{proof}

Next, we compare $THH_{\bullet}(\Gamma(L\times K))$ and $T_{\bullet}(L,K)$.
There is a cyclic (dihedral) map 
$f:THH_{\bullet}(\Gamma(L\times K))\rightarrow T_{\bullet}(L,K)$ defined as follows:
Fix a simplicial degree $q$ and consider the factor 
$X_0\wedge\cdots\wedge X_q$ of $T_q(L,K)$. Let $p_i$ be the composition
\[
\Gamma(L\times K)\rightarrow L\times K\xrightarrow{\pr}X_i\quad.
\]
Here $\pr$ denotes the projection to the first 
factor, it is $L\times K\rightarrow L$ if $X_i=L$ and $L\times K\rightarrow K$ if $X_i=K$.
Now we map $THH_q(\Gamma(L\times K))$ 
into $X_0\wedge\cdots\wedge X_q$ by
\small
\[
THH_q(\Gamma(L\times K))=\Gamma(L\times K)^{\wedge (q+1)}=\Gamma(L\times K)\wedge\cdots\wedge \Gamma(L\times K)
\xrightarrow{p_0\wedge\cdots\wedge p_q} X_0\wedge\cdots\wedge X_q\quad.
\]
\normalsize
Our map $f_q$ is determined by the 
collection of all such maps when $X_0\wedge\cdots\wedge X_q$ runs through all
factors of $T_q(L,K)$.

\begin{Lem}\label{lem:THHprodsecond}
Assume that $L$ and $K$ are cofibrant orthogonal ring spectra (with involution).
The geometric realization of $f_{\bullet}$ 
is an $S^1$-map ($O(2)$-map) $f:THH(\Gamma(L\times K))\rightarrow |T_{\bullet}(L,K)|$,
which is non-equivariantly a $\pi_*$-isomorphism.
\end{Lem}

\begin{proof}
We have that $THH_{\bullet}(\Gamma(L\times K))$ and $T_{\bullet}(L,K)$ are both good
cyclic (dihedral) orthogonal spectra. 
Hence by proposition~\ref{Prop:piisoofrealiz}, 
it is sufficient to prove that in each simplicial degree 
\[
f_q:THH_q(\Gamma(L\times K))\rightarrow T_q(L,K)
\]
is non-equivariantly a $\pi_*$-isomorphism.
The clue to prove this is to replace $\times$ by $\vee$. Consider the diagram
\[
\begin{CD}
\Gamma(L\vee K)\wedge \cdots\wedge \Gamma(L\vee K) @>>> \bigvee X_0\wedge\cdots\wedge X_q\\
@VVV @VVV\\
\Gamma(L\times K)\wedge \cdots\wedge \Gamma(L\times K) @>{f_q}>> \prod X_0\wedge\cdots\wedge X_q
\end{CD}\quad.
\]
The lower left corner is $THH_q(\Gamma(L\times K))$, 
the lower right corner is $T_q(L,K)$, and the map at the bottom is $f_q$.

Look at the left vertical map.
Since $L$ and $K$ are cofibrant, we have that $L\vee K\rightarrow L\times K$ is a 
$\pi_*$-iso, see proposition~\ref{Prop:veepiisotimes}. 
$\Gamma$ preserves $\pi_*$-isomorphisms, and
$\Gamma(L\vee K)$ and $\Gamma(L\times K)$ are both cofibrant. 
Recall that the smash product of a $\pi_*$-iso with a cofibrant
orthogonal spectrum yields a new $\pi_*$-iso. From these 
considerations it follows that the left vertical map is a $\pi_*$-iso.

The right vertical map is also a $\pi_*$-iso. 
This follows from the fact that each $X_0\wedge\cdots\wedge X_q$ is cofibrant
and the fact that wedge and $\times$-products of 
cofibrant orthogonal spectra are $\pi_*$-isomorphic.

To see that the top map is a $\pi_*$-iso, we decompose it as
\[
\Gamma(L\vee K)\wedge \cdots\wedge \Gamma(L\vee K) \rightarrow
(L\vee K)^{\wedge (q+1)}
\rightarrow \bigvee X_0\wedge\cdots\wedge X_q\quad.
\]
The first map is a $\pi_*$-iso since $\Gamma(L\vee K)\rightarrow L\vee K$ is 
a $\pi_*$-iso between cofibrant orthogonal spectra.
Distributivity of $\vee$ over $\wedge$ shows that 
the second map is an isomorphism of orthogonal spectra.

Commutativity of the diagram implies that $f_q$ also is a $\pi_*$-iso.
\end{proof}

\begin{Rem}\label{Rem:improvecrossprodformula}
The conclusion of the cross product formula we just have derived, proposition~\ref{Prop:THHcrossprodformula} above,
is too weak. Due to proposition~\ref{Prop:TCiso}, we would like our map
\[
f:THH(\Gamma(L\times K))\rightarrow THH(L)\times THH(K)
\]
to be a cyclotomic $\pi_*$-isomorphism between cyclotomic spectra.
We know only that $f$ is non-equivariantly a $\pi_*$-isomorphism, and we do not know that
$THH(L)\times THH(K)$ is a cyclotomic spectrum. These two problems are closely connected:
\begin{itemize}
\item[] If we can show that $THH(L)\times THH(K)$ is a cyclotomic spectrum, then proposition~\ref{Prop:piisoiscyclpiiso}
would imply that $f$ is a cyclotomic $\pi_*$-isomorphism.
\item[] There are good candidates for the cyclotomic structure maps $r_C$ on $THH(L)\times THH(K)$, but we do not
know that these maps are cyclotomic $\pi_*$-isomorphisms. If
we knew that $f$ was a cyclotomic $\pi_*$-isomorphism, then we could show that the $r_C$'s also are
cyclotomic $\pi_*$-isomorphisms.
\end{itemize}
In both approaches we should allow ourselves to take cofibrant or fibrant replacements.
To study the first approach, one should give a more explicit description of a cyclotomic spectrum.
Lemma~2.2 in~\cite{HesselholtMadsen:97} can be the inspiration for such a description.
To study the second approach, one could try to transfer McCarthy's concept of a special homotopy
to the setting of cyclic orthogonal spectra. Analogies of propositions~1.5.12 and~1.6.15 in~\cite{DundasMcCarthy:96}
should then show that $f$ is a cyclotomic $\pi_*$-isomorphism.
\end{Rem}

\section{The Barratt-Eccles functor}

When trying to define a trace map $M_n(L)\rightarrow L$ one encounters the problem that 
there is a priori no notion of addition on $L$, i.e. one cannot add 
two points in $L(V)$ and get a new point. 
We will use the Barratt-Eccles
functor $\Gamma^+$, see~\cite{BarrattEccles:74}, to solve this problem. 
In the next section we will construct a trace map from $THH(\Gamma M_n(L))$.
As a target for this trace we now introduce a new model for $THH(L)$,
called $THH^+(L)$. The construction uses $\Gamma^+$.
This is an idea due to Christian Schlichtkrull, see~\cite{Schlichtkrull:98}.
The main result of this section is proposition~\ref{Prop:THH+}, 
but this result is not as strong as hoped for,
see remark~\ref{Rem:notoptimalTHH+}.

Let $G$ be a discrete group. Recall that  $E_{\bullet}G$ is the simplicial $G$-space given by
\[
E_qG= G^{q+1}
\]
where
\begin{align*}
d_i(g_0,\ldots,g_q)&=(g_0,\ldots,g_{i-1},g_{i+1},\ldots,g_q)\quad,\\
s_i(g_0,\ldots,g_q)&=(g_0,\ldots,g_{i},g_{i},\ldots,g_q)\quad\text{and}\\
(g_0,\ldots,g_q).g&=(g_0g,\ldots,g_qg)\quad.
\end{align*}
Let $EG$ be the geometric realization.

Write $\mathbf{n}$ for the set $\{1,2,\ldots,n\}$. Let $\mathscr{M}(\mathbf{m},\mathbf{n})$
be the set of all strictly increasing functions from $\mathbf{m}$ to $\mathbf{n}$. Given a permutation
$\sigma\in\Sigma_n$ and $\alpha\in \mathscr{M}(\mathbf{m},\mathbf{n})$, then there is a unique function 
in $\mathscr{M}(\mathbf{m},\mathbf{n})$ which has the same image as the composition $\sigma\alpha$.
Denote this by $\sigma_*(\alpha)$.

\begin{Def}
For $\alpha\in \mathscr{M}(\mathbf{m},\mathbf{n})$ define the restriction map $\alpha^*:\Sigma_n\rightarrow \Sigma_m$
by commutativity of the diagram
\[
\begin{CD}
\mathbf{m}@>{\alpha}>>\mathbf{n}\\
@V{\alpha^*(\sigma)}VV @VV{\sigma}V\\
\mathbf{m}@>{\sigma_*(\alpha)}>>\mathbf{n}
\end{CD}\quad.
\]
\end{Def}

On the Cartesian product $X^n$ we have a right action of $\Sigma_n$ given by
\[
(x_1,\ldots,x_n).\sigma= (x_{\sigma(1)},\ldots,x_{\sigma(n)})\quad,
\]
and given $\alpha\in\mathscr{M}(\mathbf{m},\mathbf{n})$ we have an induced map $\alpha^*:X^n\rightarrow X^m$
defined by the formula
\[
\alpha^*(x_1,\ldots,x_n)=(x_{\alpha(1)},\ldots,x_{\alpha(m)})\quad.
\]
We say that $\alpha$ is entire for $(x_1,\ldots,x_n)$ if $i\not\in\alpha(\mathbf{m})$ implies $x_i=*$.

Consider the equivalence relation on
\[
\coprod_{m\geq0} E_{\bullet}\Sigma_m\times X^m
\]
given by
\begin{itemize}
\item[(a)] $(\mathbf{c},\mathbf{x})\sim (\mathbf{c}.\sigma,\mathbf{x}.\sigma)$ for $\mathbf{c}\in E_{\bullet}\Sigma_m$, $\mathbf{x}\in X^m$ and $\sigma\in \Sigma_m$.
\item[(b)] $(\mathbf{c},\mathbf{x})\sim (\alpha^*\mathbf{c},\alpha^*\mathbf{x})$ if $\mathbf{c}\in E_{\bullet}\Sigma_m$, $\mathbf{x}\in X^m$ and 
$\alpha\in\mathscr{M}(\mathbf{m},\mathbf{n})$ is entire for $\mathbf{x}$.
\end{itemize}

\begin{Def}
Let $\Gamma_{\bullet}^+(X)$ be $(\coprod_{m\geq0} E_{\bullet}\Sigma_m\times X^m)/\sim$,
and $\Gamma^+(X)$ it's geometrical realization.
\end{Def}

\begin{Prop}
$\Gamma_{\bullet}^+(X)$ is a dihedral space. 
\end{Prop}

\begin{proof}
We must define the cyclic and involutive operators. Let
\[
t_q[(\sigma_0,\ldots,\sigma_q);(x_1,\ldots,x_m)]=[(\sigma_q,\sigma_0,\ldots,\sigma_{q-1});(x_1,\ldots,x_m)]
\]
and
\[
r_q[(\sigma_0,\ldots,\sigma_q);(x_1,\ldots,x_m)]=[(\sigma_q,\sigma_{q-1},\ldots,\sigma_1,\sigma_0);(x_1,\ldots,x_m)]\quad.
\]
\end{proof}

\begin{Prop}\label{Prop:gamma+addition}
$\Gamma^+(X)$ is a monoid with unit, the operation is commutative up to homotopy and the monoid is free.
\end{Prop}

\begin{proof}
A proof of this result can be found in~\cite{BarrattEccles:74}, corollary~3.10 and proposition~3.11.
However in their proofs $X$ is a based simplicial set. But the definition coincides with the one given here
when $X$ is a discrete based set. Therefore, the statements about the algebraic structure of $\Gamma^+(X)$
follow by applying Barratt and Eccles' proofs to $X$'s underlying discrete set $X^{\delta}$. Hence,
$\Gamma^+(X)$ is a free monoid with unit. Checking homotopy commutativity can be done as in the reference,
but in order to familiarize ourselves with the operation we write out the argument for homotopy commutativity here:

We write $+$ for the operation. To start we need a homomorphism $\Sigma_{m_1}\times \Sigma_{m_2}\rightarrow \Sigma_{m_1+m_2}$. 
We call this homomorphism $\amalg$ and it is defined by
\[
(\sigma\amalg\rho)(j)=\begin{cases}
\sigma(j)&\text{if $j\leq m_1$, and}\\
\rho(j-m_1)+m_1&\text{if $j> m_1$.}
\end{cases}
\]
Now given $[(\sigma_0,\ldots,\sigma_q);(x_1,\ldots,x_{m_1})]$ and $[(\rho_0,\ldots,\rho_q);(y_1,\ldots,y_{m_2})]$ in $\Gamma_q^+(X)$
we define their sum as:
\[
[(\sigma_0\amalg\rho_0,\ldots,\sigma_q\amalg\rho_q);(x_1,\ldots,x_{m_1},y_1,\ldots,y_{m_2})]\quad.
\]
This defines a simplicial map $\Gamma_{\bullet}^+(X)\times\Gamma_{\bullet}^+(X)\rightarrow \Gamma_{\bullet}^+(X)$, and
we define $+$ to be its geometric realization. 

To get homotopy commutativity we define a simplicial homotopy $h_i:\Gamma_{q}^+(X)\times\Gamma_{q}^+(X)\rightarrow \Gamma_{q+1}^+(X)$, $i=0,\ldots,q$.
Let $\tau\in\Sigma_{m_1+m_2}$ be the permutation defined by
\[
\tau(j)=
\begin{cases}
j+m_2&\text{if $j\leq m_1$, and}\\
j-m_1&\text{if $j> m_1$.}
\end{cases}
\]
Notice that $\tau(\sigma\amalg\rho)\tau^{-1}=\rho\amalg\sigma$. Now we define $h_i$ by:
\begin{multline*}
h_i([(\sigma_0,\ldots,\sigma_q);(x_1,\ldots,x_{m_1})],[(\rho_0,\ldots,\rho_q);(y_1,\ldots,y_{m_2})])=\\
[(\tau \sigma_0\amalg\rho_0,\ldots,\tau \sigma_i\amalg\rho_i, \sigma_i\amalg\rho_i ,\ldots,\sigma_q\amalg\rho_q);
(x_1,\ldots,x_{m_1},y_1,\ldots,y_{m_2})]\quad.
\end{multline*}
This is a simplicial homotopy between $x+y$ and $y+x$. 
\end{proof}

\begin{Prop}\label{Prop:Gamma+cof}
If $i:A\rightarrow X$ is an unbased cofibration of based topological spaces, 
then $\Gamma^+(A)\rightarrow \Gamma^+(X)$ is also an unbased cofibration.
\end{Prop}

\begin{proof}
Since $i$ is a cofibration, it is an inclusion and we view $A$ as a subspace of $X$. By Str\o{}m's criterion
there are maps $H:X\times I\rightarrow X$ and $\phi:X\rightarrow I$ with $A\subset\phi^{-1}(0)$, $H(x,0)=x$ for all $x\in X$,
$H(a,t)=a$ for all $a\in A$ and $t\in I$ and $H(x,t)\in A$ when $\phi(x)>t$. Now define $H^m:X^m\times I\rightarrow I$ and $\phi^m:X^m\rightarrow I$
by
\small
\[
H^m(x_1,\ldots,x_m,t)=(H(x_1,t),\ldots,H(x_m,t))\quad\text{and}\quad \phi^m(x_1,\ldots,x_m)=\max_i \phi(x_i)\quad.
\]
\normalsize
$H^m$ and $\phi^m$ satisfies Str\o{}m's criterion. In addition they are $\Sigma_m$-equivariant and if $\alpha:\mathbf{m}\rightarrow\mathbf{n}$
is entire for $(x_1,\ldots,x_m)$, then
\small
\[
\alpha^*H^n(x_1,\ldots,x_n,t)=H^m(\alpha^*(x_1,\ldots,x_n),t)\quad\text{and}\quad \phi^n(x_1,\ldots,x_n)=\phi^m\alpha^*(x_1,\ldots,x_n)\quad.
\]
\normalsize
Therefore, we get induced maps
\[
H'_q:(\Gamma^+_q(X))\times I\rightarrow \Gamma^+_q(X)\quad\text{and}\quad \phi'_q:\Gamma^+_q(X)\rightarrow I
\]
showing that $\Gamma^+_{\bullet}(A)\rightarrow \Gamma^+_{\bullet}(X)$ is a cofibration in each simplicial degree.
Moreover the $H'_q$'s and the $\phi'_q$'s respect the face and degeneracy maps. Thus by geometric realization we get 
maps $H':\Gamma^+(X)\times I\rightarrow \Gamma^+(X)$ and $\phi':\Gamma^+(X)\rightarrow I$,
showing that $\Gamma^+(A)\rightarrow\Gamma^+(X)$ is a cofibration.
\end{proof}

\begin{Prop}\label{Prop:Gamma+conn}
If $X$ is an $n$-connected well-pointed space, 
then the map $X\rightarrow\Gamma^+(X)$ is $(2n+1)$-connected.
\end{Prop}

\begin{proof}
First we observe from the definition of $s_i:\Gamma^+_q(X)\rightarrow \Gamma^+_{q+1}(X)$ that all 
degeneracy maps are cofibrations. Hence, $\Gamma^+_{\bullet}(X)$ is a good simplicial space for all $X$.
It follows that $\Gamma^+(-)$ preserves weak equivalences.

By the natural weak equivalence $|\Sing X|\rightarrow X$, we see that 
it is sufficient to prove the result when
$X$ is a simplicial set. Assume that $X$ is $n$-connected. If $n=-1$, 
there is nothing to prove. Therefore consider $n\geq0$.
Lemma~4.8 in~\cite{BarrattEccles:74} (see also corollary~5.4) says 
that $X\rightarrow \Gamma^+(X)$ is $(2n+1)$-connected. 
\end{proof}

\begin{Cor}\label{Cor:Gamma+piiso}
We may apply $\Gamma^+(-)$ level-wise to an orthogonal spectrum $L$.
If $L$ is cofibrant, then the natural map
\[
L\rightarrow\Gamma^+(L)
\]
is a $\pi_*$-iso.
\end{Cor}

\begin{proof}
The propositions~\ref{Prop:Gamma+conn} and~\ref{Prop:Gamma+cof} verify the conditions
required to apply corollary~\ref{Cor:indfunct}.
\end{proof}

We now define the model $THH^+(L)$:

\begin{Def}
Assume that $L$ is an orthogonal ring spectrum (with involution).
We define $THH^+(L)$ to be the geometric realization of the bicyclic (bidihedral)
orthogonal spectrum $\Gamma^+_{\bullet}(THH_{\bullet}(L))$.
\end{Def}

A bicyclic orthogonal spectrum is a functor $(\catDeltaC\times\catDeltaC)^{\op}\rightarrow\mathscr{IS}$.
Restricting via the diagonal $\catDeltaC\rightarrow\catDeltaC\times\catDeltaC$ we get a
cyclic orthogonal spectrum, whose geometric realization is an orthogonal $S^1$-spectrum.
Analogously, in the bidihedral case the geometric realization has $O(2)$-action.

\begin{Rem}
Observe that $THH^+(L)$ inherits an addition from $\Gamma^+$. Fixing a simplicial degree $p$
and a level $V$, we have an associative operation
\[
+:\Gamma^+(THH_p(L)(V))\times\Gamma^+(THH_p(L)(V))\rightarrow\Gamma^+(THH_p(L)(V))
\]
by proposition~\ref{Prop:gamma+addition}. Taking geometric realization in the $p$ direction,
we get addition
\[
+:THH^+(L)\times THH^+(L)\rightarrow THH^+(L)\quad.
\]
\end{Rem}

\begin{Prop}\label{Prop:THH+}
Assume that $L$ is a cofibrant orthogonal ring spectrum (with involution). 
Then the natural $S^1$-map ($O(2)$-map) $THH(L)\rightarrow THH^+(L)$ is non-equivariantly a $\pi_*$-isomorphism.
\end{Prop}

\begin{proof}
Forgetting about the action, we can consider $THH^+(L)$ as the geometric realization of the 
bisimplicial orthogonal spectrum $\Gamma^+_{\bullet}(THH_{\bullet}(L))$.
It is classical that the two geometric realizations
\[
\left|[q]\mapsto \Gamma^+_{q}(THH_{q}(L))\right|\quad\text{and}\quad 
\Big|[p]\mapsto \left|[q]\mapsto \Gamma^+_{q}(THH_{p}(L))\right|\;\Big|
\]
are homeomorphic. Notice that both $THH_{\bullet}(L)$ and $\Gamma^+(THH_{\bullet}(L))$ are good
simplicial orthogonal spectra.
Since $THH_{\bullet}(L)$ is cofibrant
in each simplicial degree, corollary~\ref{Cor:Gamma+piiso} says that
\[
THH_p(L)\rightarrow \left|[q]\mapsto \Gamma^+_{q}(THH_{p}(L))\right|=\Gamma^+(THH_p(L))
\]
is a $\pi_*$-isomorphism for all $p$. Furthermore, proposition~\ref{Prop:piisoofrealiz} yields that the map
\[
THH(L)=\left|[p]\mapsto THH_p(L)\right|\longrightarrow 
\Big|[p]\mapsto \left|[q]\mapsto \Gamma^+_{q}(THH_{p}(L))\right|\;\Big|\cong THH^+(L)
\]
is also a $\pi_*$-isomorphism.
\end{proof}

\begin{Rem}\label{Rem:notoptimalTHH+}
Again, our result is weaker than what we hoped for, namely that the natural map
$THH(L)\rightarrow THH^+(L)$ would be a cyclotomic $\pi_*$-isomorphism between
cyclotomic spectra. As before, there are two strategies for improving the result: Either
we could show that $THH^+(L)$ is a cyclotomic spectrum, or we could show that
the natural map is a cyclotomic $\pi_*$-isomorphism. In any case, the other part then 
should follow formally. Compare with remark~\ref{Rem:improvecrossprodformula}.
\end{Rem}

\section{Morita equivalence}\label{sect:morita}

In this section we show Morita equivalence for $THH$. We adopt
an approach by Christian Schlichtkrull, see theorem~3.6 in~\cite{Schlichtkrull:98}, 
to the setting of orthogonal ring spectra (with involution). 
In order to carry out the proof, we view orthogonal spectra as $\mathscr{J}$-spaces
and do our constructions externally. Unfortunately, our result is not as
strong as we would like, see remark~\ref{Rem:Moritanotoptimal}.

Let us start by stating the result, the proof spans the following subsections:

\begin{Prop}\label{Prop:Morita}
Assume that $L$ is a cofibrant orthogonal ring spectrum (with involution).
Then there is a natural $S^1$-map ($O(2)$-map)
\[
\Tr: THH(\Gamma M_n(L))\rightarrow THH^+(L)\quad,
\]
which is non-equivariantly a $\pi_*$-isomorphism. Furthermore, we have a commutative
diagram
\[
\begin{CD}
THH(\Gamma (M_{n_1}\times M_{n_2})(L)) @>{(\Tr\circ\pr_1,\Tr\circ\pr_2)}>> THH^+(L)\times THH^+(L)\\
@V{\oplus}VV @VV{+}V\\
THH(\Gamma M_{n_1+n_2}(L)) @>{\Tr}>> THH^+(L)
\end{CD}
\]
of orthogonal $S^1$-spectra (orthogonal $O(2)$-spectra).
\end{Prop}

\begin{Rem}\label{Rem:Moritanotoptimal}
The result above is not as strong as we could hope for. We have not shown that $\Tr$ is
a cyclotomic $\pi_*$-isomorphism, and we do not know that $THH^+(L)$ is a cyclotomic spectrum.
Again, it is probable that proving one of these wishes would yield the other as an easy
corollary. See also the remarks~\ref{Rem:improvecrossprodformula} and~\ref{Rem:notoptimalTHH+}.
\end{Rem}

\subsection{Internalizing}

Recall that orthogonal spectra can be described as diagram spaces over a topological category $\mathscr{J}$,
see theorem~\ref{thm:orthspasdiagr}. Furthermore, the category $\mathscr{J}$ has a symmetric operation, 
namely direct sum $\oplus$. Using left Kan extension, see theorem~X.4.1 in~\cite{MacLane:98},
we can therefore lift constructions on based topological spaces to constructions on orthogonal spectra.
To be more precise: 
Let $\mathscr{J}^{q+1}\Top_*$ denote the category of continuous functors $\mathscr{J}^{q+1}\rightarrow\Top_*$.
Assume that $F$ is a continuous 
functor $\Top_*^{q+1}\rightarrow\Top_*$. If we are given orthogonal spectra
$L^0$, $L^1$, $\ldots$, $L^q$, we can consider these as functors $\mathscr{J}\rightarrow\Top_*$ and take the 
composition
\[
\mathscr{J}^{q+1} \xrightarrow{L^0\times L^1\times\ldots\times L^q} \Top_*^{q+1} \xrightarrow{F} \Top_*\quad.
\]
This is an object in $\mathscr{J}^{q+1}\Top_*$.
The iterated direct sum is a functor $\mathscr{J}^{q+1}\rightarrow\mathscr{J}$, and left Kan extension
is a functor $\bbP:\mathscr{J}^{q+1}\Top_*\rightarrow\mathscr{J}\Top_*=\mathscr{IS}$.
$\bbP$ therefore turns the above composition into a functor
\[
\mathscr{J}\rightarrow\Top_*\quad,
\]
i.e. a new orthogonal spectrum.
This process of internalizing is natural 
with respect to natural transformations of functors
$\Top_*^{q+1}\rightarrow\Top_*$.

Let us look at some examples:

\begin{Exa}
If $F$ is given by
\[
F(X_0,\ldots,X_q)=X_0\wedge\cdots\wedge X_q\quad,
\]
then it follows from the definition of smash product of orthogonal spectra that the left Kan extension of
\[
F(L^0,\ldots,L^q)\quad\text{is}\quad L^0\wedge \cdots \wedge L^q\quad.
\]
\end{Exa}

\begin{Exa}
If $F$ is given by
\[
F(X_0,\ldots,X_q)=M_n(X_0)\wedge\cdots\wedge M_n(X_q)\quad,
\]
then the left Kan extension of
\[
F(L^0,\ldots,L^q)\quad\text{is}\quad M_n(L^0)\wedge \cdots \wedge M_n(L^q)\quad.
\]
\end{Exa}

\begin{Exa}\label{Exa:lKeG+}
Suppose that $F$ is given by
\[
F(X_0,\ldots,X_q)=\Gamma^+(X_0\wedge\cdots\wedge X_q)\quad.
\]
Let $G(L^0,\ldots,L^q)$ be the left Kan extension of 
$(F(L^0,\ldots,L^q))$. We observe that $G(L^0,\ldots,L^q)$
is an orthogonal spectrum, and there is a natural map
\[
G(L^0,\ldots,L^q)\rightarrow \Gamma^+(L^0\wedge \cdots\wedge L^q)\quad.
\]
To see how the map is defined, notice that the adjoint of the identity of 
$L^0\wedge \cdots\wedge L^q$ is a natural transformation
\[
L^0(V_0)\wedge \cdots\wedge L^q(V_q)\rightarrow (L^0\wedge \cdots\wedge L^q)(V_0\oplus\cdots\oplus V_q)\quad.
\]
Apply the functor $\Gamma^+(-)$ to get a natural transformation
\[
F(L^0(V_0),\ldots,L^q(V_q))\rightarrow \Gamma^+(L^0\wedge \cdots\wedge L^q)(V_0\oplus\cdots\oplus V_q)\quad.
\]
Its adjoint is the natural map we seek.
\end{Exa}

\subsection{Dihedral structure on functors $\Top_*^{q+1}\rightarrow\Top_*$}\label{subsect:dihedralcoll}

In this subsection we study two collections of functors $\Top_*^{q+1}\rightarrow\Top_*$, $q\geq0$,
these are
\[
(X_0,\ldots,X_q)\mapsto M_n(X_0)\wedge\cdots\wedge M_n(X_q)\quad\text{and}\quad (X_0,\ldots,X_q)\mapsto\Gamma_q^+(X_0\wedge\ldots\wedge X_q)\quad.
\]
A dihedral structure for such a collection consists of natural 
transformations $d_i$, $s_i$, $t_q$ and $r_q$ satisfying the dihedral identities.
For $M_n(X_0)\wedge\cdots\wedge M_n(X_q)$ we have:
\footnotesize
\begin{align*}
d_i&:M_n(X_0)\wedge\cdots\wedge M_n(X_q)\rightarrow 
M_n(X_0)\wedge\cdots \wedge M_n(X_i\wedge X_{i+1})\wedge\cdots\wedge M_n(X_q)\quad,
\quad\text{if $i< q$,}\\
d_q&:M_n(X_0)\wedge\cdots\wedge M_n(X_q)\rightarrow 
M_n(X_q\wedge X_0)\wedge M_n(X_1)\wedge\cdots \wedge M_n(X_{q-1})\quad,\\
s_i&:M_n(X_0)\wedge\cdots\wedge M_n(X_q)\rightarrow 
M_n(X_0)\wedge\cdots \wedge M_n(X_i)\wedge M_n(S^0)\wedge M_n(X_{i+1})\wedge\cdots\wedge M_n(X_q)
\quad,\\
t_q&:M_n(X_0)\wedge\cdots\wedge M_n(X_q)\rightarrow 
M_n(X_q)\wedge M_n(X_0)\wedge M_n(X_1)\wedge\cdots \wedge M_n(X_{q-1})\quad\text{, and}\\
r_q&:M_n(X_0)\wedge\cdots\wedge M_n(X_q)\rightarrow 
M_n(X_0)\wedge M_n(X_q)\wedge M_n(X_{q-1})\wedge\cdots\wedge M_n(X_1)\quad.
\end{align*}
\normalsize
Here the $d_i$'s are defined using the multiplication for $M_n$, the $s_i$'s are defined using the unit in $M_n(S^0)$,
cyclic operators $t_q$ are given by permuting the factors and the involutive operators $r_q$
are defined using the transposition on $M_n$ together with the order reversing permutation of the last $q$ factors.

For $\Gamma^+_q(X_0\wedge\cdots\wedge X_q)$ we have:
\small
\begin{align*}
d_i&:\Gamma^+_q(X_0\wedge\cdots\wedge X_q)\rightarrow 
\Gamma_{q-1}^+(X_0\wedge\cdots \wedge (X_i\wedge X_{i+1})\wedge\cdots\wedge X_q)\quad,
\quad\text{if $i< q$,}\\
d_q&:\Gamma^+_q(X_0\wedge\cdots\wedge X_q)\rightarrow 
\Gamma_{q-1}^+((X_q\wedge X_0)\wedge X_1\wedge \cdots \wedge X_{q-1})\quad,\\
s_i&:\Gamma^+_q(X_0\wedge\cdots\wedge X_q)\rightarrow 
\Gamma_{q+1}^+(X_0\wedge\cdots \wedge X_i\wedge S^0\wedge X_{i+1}\wedge\cdots\wedge X_q)\quad,\\
t_q&:\Gamma^+_q(X_0\wedge\cdots\wedge X_q)\rightarrow 
\Gamma_{q}^+(X_q \wedge X_0 \wedge X_1\wedge \cdots \wedge X_{q-1})\quad\text{, and}\\
r_q&:\Gamma^+_q(X_0\wedge\cdots\wedge X_q)\rightarrow 
\Gamma_{q}^+(X_0\wedge X_q\wedge X_{q-1}\wedge \cdots\wedge X_1)\quad.
\end{align*}
\normalsize
And the dihedral structure is directly inherited from $\Gamma^+_{\bullet}(-)$.

The purpose of such dihedral collections of functors is to construct cyclic (dihedral) orthogonal spectra
using left Kan extension. Our result is:

\begin{Prop}\label{Prop:dihedralcollection}
Assume that $L$ is an orthogonal ring spectrum (with involution). 
Given a collection of functors $F_q:\Top_*^{q+1}\rightarrow\Top_*$
with natural transformations $d_i$, $s_i$, $t_q$ and $r_q$ satisfying the dihedral identities, then the process 
of internalizing using $L$ in all factors yields a cyclic (dihedral) orthogonal spectrum. 
This construction is natural in $\{F_q\}$.
\end{Prop}

Applying this proposition to $M_n(X_0)\wedge\cdots\wedge M_n(X_q)$ we get precisely $THH_{\bullet}(M_n(L))$.
In the case $\Gamma^+_q(X_0\wedge\cdots\wedge X_q)$ we get a cyclic (dihedral) orthogonal spectrum $Y_{\bullet}$,
and by example~\ref{Exa:lKeG+} a map of cyclic (dihedral) orthogonal spectra
\[
Y_q\rightarrow \Gamma^+_q(THH_q(L))\quad.
\]
We end the subsection by proving the proposition.

\begin{proof}
Let $X_q$ be the left Kan extension of $F_q(L,\ldots,L)$. 
We will show that $X_{\bullet}$ is a cyclic orthogonal spectrum, and that $X_{\bullet}$
is dihedral whenever $L$ comes with an involution.

By definition of the left Kan extension there is a homeomorphism between the space of
\[
\text{continuous natural transformations}\quad F_q(L(V_0),\ldots,L(V_q))\rightarrow K(V_0\oplus\cdots\oplus V_q)
\]
and the space of
\[
\text{orthogonal spectrum maps}\quad X_q \rightarrow K
\]
for any orthogonal spectrum $K$. This is adjointness. In particular there is an
adjoint to the identity map of $X_q$, this means that we have a canonical map: 
\[
c_q:F_q(L(V_0),\ldots,L(V_q))\rightarrow X_q(V_0\oplus\cdots\oplus V_q)
\]
for all $q$. Now observe that all we have to do in order to define a map $X_q\rightarrow X_p$ is to specify a natural transformation
\[
F_q(L(V_0),\ldots,L(V_q))\rightarrow X_p(V_0\oplus\cdots\oplus V_q)\quad.
\]
This is how we are going to define the cyclic operators $d^X_i$, $s^X_i$ and $t^X_i$ of $X_{\bullet}$.

Recall that an orthogonal ring spectrum is the same as an $\mathscr{I}$-FSP,
see remark~\ref{Rem:extring}. Using this external description we have
unit $\eta:S^0\rightarrow L(0)$ and multiplication $\mu:L(V)\wedge L(W)\rightarrow L(V\oplus W)$.
If $L$ has involution, we also have a natural transformation $\iota:L(V)\rightarrow L(V)$.

Face maps of $X_{\bullet}$ are given as follows: For $i<q$ we consider the composition
\small
\begin{align*}
F_q(L(V_0),\ldots,L(V_q))
&\xrightarrow{d_i} F_{q-1}(L(V_0),\ldots,L(V_{i-1}),L(V_i)\wedge L(V_{i+1}),L(V_{i+2}),\ldots,L(V_q))\\
&\xrightarrow{\mu_*} F_{q-1}(L(V_0),\ldots,L(V_{i-1}),L(V_i\oplus V_{i+1}),L(V_{i+2}),\ldots,L(V_q))\\
&\xrightarrow{c_{q-1}}X_{q-1}(V_0\oplus\cdots\oplus V_q)\quad,
\end{align*}
\normalsize
and let the adjoint be $d^X_i:X_{q}\rightarrow X_{q-1}$. In the case $i=q$ we have:
\small
\begin{align*}
F_q(L(V_0),\ldots,L(V_q))
&\xrightarrow{d_q} F_{q-1}(L(V_q)\wedge L(V_0),L(V_1),\ldots,L(V_{q-1}))\\
&\xrightarrow{\mu_*} F_{q-1}(L(V_q\oplus V_0),L(V_{1}),\ldots,L(V_{q-1}))\\
&\xrightarrow{c_{q-1}}X_{q-1}(V_q\oplus V_0\oplus\cdots\oplus V_{q-1})\\
&\xrightarrow{\text{permute $V_i$'s}} X_{q-1}(V_0\oplus\cdots\oplus V_q)\quad,
\end{align*}
\normalsize
and let this define $d^X_q:X_{q}\rightarrow X_{q-1}$.
For degeneracy maps we consider the composition:
\small\begin{align*}
F_q(L(V_0),\ldots,L(V_q))
&\xrightarrow{s_i} F_{q+1}(L(V_0),\ldots,L(V_i),S^0,L(V_{i+1}),\ldots,L(V_q))\\
&\xrightarrow{\eta_*} F_{q+1}(L(V_0),\ldots,L(V_{i}),L(0),L(V_{i+1}),\ldots,L(V_q))\\
&\xrightarrow{c_{q+1}}X_{q+1}(V_0\oplus\cdots \oplus V_{i}\oplus 0\oplus V_{i+1}\oplus\cdots\oplus V_q)\\
&\xrightarrow{\text{$0$ is the unit for $\oplus$}} X_{q+1}(V_0\oplus\cdots\oplus V_q)\quad,
\end{align*}\normalsize
and we define $s^X_i:X_{q}\rightarrow X_{q+1}$ to be its adjoint.
The cyclic operator $t^X_q:X_{q}\rightarrow X_{q}$ is defined as the adjoint of the composition
\small\begin{align*}
F_q(L(V_0),\ldots,L(V_q))
&\xrightarrow{t_q} F_{q}(L(V_q),L(V_0),\ldots,L(V_{q-1}))\\
&\xrightarrow{c_{q}}X_{q}(V_q\oplus V_0\oplus\cdots\oplus V_{q-1})\\
&\xrightarrow{\text{permute $V_i$'s}} X_{q}(V_0\oplus\cdots\oplus V_q)\quad.
\end{align*}\normalsize

Now assume that $L$ has an involution $\iota:L\rightarrow L$, 
then we can give $X_{\bullet}$ dihedral structure by defining
the involutive operator $r^X_q:X_q\rightarrow X_q$ to be the adjoint of the composition
\small\begin{align*}
F_q(L(V_0),\ldots,L(V_q))
&\xrightarrow{r_q} F_{q}(L(V_0),L(V_q),\ldots,L(V_{1}))\\
&\xrightarrow{F(\iota,\ldots,\iota)} F_{q}(L(V_0),L(V_q),\ldots,L(V_{1}))\\
&\xrightarrow{c_{q}}X_{q}(V_0\oplus V_q\oplus \cdots \oplus V_{1})\\
&\xrightarrow{\text{permute $V_i$'s}} X_{q}(V_0\oplus\cdots\oplus V_q)\quad.
\end{align*}\normalsize

We are now supposed to verify a long list of identities 
involving the operators $d^X_i$, $s^X_i$ and $t^X_i$, and in the case
where $L$ has involution there are even more identities. 
We will skip this painful task with one exception: We will prove the dihedral identity
\[
d^X_i r^X_q = r_{q-1}^X d^X_{q-i}
\] 
in the case $0<i<q$. This case will illustrate the 
techniques used when proving the other cyclic and dihedral identities.
Moreover, we will also see how
the anti-commutativity of the involution plays a role.

By adjointness and the definitions we observe that $d^X_i r^X_q$ is the adjoint of
\footnotesize\begin{align*}
F_q(L(V_0)&,\ldots,L(V_q))
\xrightarrow{r_q} F_{q}(L(V_0),L(V_q),\ldots,L(V_{1}))\\
&\xrightarrow{F(\iota,\ldots,\iota)} F_{q}(L(V_0),L(V_q),\ldots,L(V_{1}))\\
&\xrightarrow{d_i} F_{q-1}(L(V_0),L(V_q),\ldots,L(V_{q-i+2}),L(V_{q-i+1})\wedge L(V_{q-i}),L(V_{q-i-1}),\ldots,L(V_1))\\
&\xrightarrow{\mu_*} F_{q-1}(L(V_0),L(V_q),\ldots,L(V_{q-i+2}),L(V_{q-i+1}\oplus V_{q-i}),L(V_{q-i-1}),\ldots,L(V_1))\\
&\xrightarrow{c_{q-1}}X_{q-1}(V_0\oplus V_q\oplus\cdots\oplus V_1)\\
&\xrightarrow{\text{permute}}X_{q-1}(V_0\oplus\cdots\oplus V_q)\quad,
\end{align*}\normalsize 
and $r_{q-1}^X d^X_{q-i}$ is the adjoint of 
\footnotesize\begin{align*}
F_q(L(V_0)&,\ldots,L(V_q))
\xrightarrow{d_{q-i}} F_{q-1}(L(V_0),\ldots,L(V_{q-i-1}),L(V_{q-i})\wedge L(V_{q-i+1}),L(V_{q-i+2}),\ldots,L(V_q))\\
&\xrightarrow{\mu_*} F_{q-1}(L(V_0),\ldots,L(V_{q-i-1}),L(V_{q-i}\oplus V_{q-i+1}),L(V_{q-i+2}),\ldots,L(V_q))\\
&\xrightarrow{r_{q-1}} F_{q-1}(L(V_0),L(V_q),\ldots,L(V_{q-i+2}),L(V_{q-i}\oplus V_{q-i+1}),L(V_{q-i-1}),\ldots,L(V_1))\\
&\xrightarrow{F(\iota,\ldots,\iota)} F_{q-1}(L(V_0),L(V_q),\ldots,L(V_{q-i+2}),L(V_{q-i}\oplus V_{q-i+1}),L(V_{q-i-1}),\ldots,L(V_1))\\
&\xrightarrow{c_{q-1}}X_{q-1}(V_0\oplus V_q\oplus\cdots\oplus V_{q-i+2}\oplus(V_{q-i}\oplus V_{q-i+1})\oplus V_{q-i-1}\oplus\cdots\oplus V_1)\\
&\xrightarrow{\text{permute}}X_{q-1}(V_0\oplus\cdots\oplus V_q)\quad.
\end{align*}\normalsize 
These two compositions can be compared and found to be equal using naturality of $d_i$ and $r_{q-1}$ and the fact that
\[
\begin{CD}
L(V)\wedge L(W) @>{\iota\wedge\iota}>> L(V)\wedge L(W) @>{\text{twist}}>> L(W)\wedge L(V)\\
@V{\mu}VV && @VV{\mu}V\\
L(V\oplus W) @>{\iota}>> L(V\oplus W) @>{L(\text{twist})}>> L(W\oplus V)
\end{CD}
\]
commutes.
\end{proof}

\subsection{The trace}

In this subsection we will define the trace map
\[
\Tr:THH(\Gamma M_nL)\rightarrow THH^+(L)\quad,
\]
and prove proposition~\ref{Prop:Morita}.
The definition of $\Tr$ will use the machinery developed above.
Below we will construct a collection of natural transformations 
\[
\Tr_q:M_n(X_0)\wedge\cdots\wedge M_n(X_q)\rightarrow \Gamma_q^+(X_0\wedge\ldots\wedge X_q)
\]
that commutes with the cyclic (dihedral) structure. Feeding this into proposition~\ref{Prop:dihedralcollection},
we get an $S^1$-map ($O(2)$-map) $THH(M_n(L))\rightarrow THH^+(L)$. 
We now define $\Tr$ as the composition
\[
THH(\Gamma M_n(L))\rightarrow THH(M_n(L))\rightarrow THH^+(L)\quad.
\]
Here $\Gamma$ denotes the cofibrant replacement functor of theorem~\ref{Thm:cofrepl}, and the first
map is induced by the natural transformation $\Gamma M_n(L)\rightarrow M_n(L)$.

To define $\Tr_q$, we would like to send a point $((x^0_{i\,j}),\ldots,(x^q_{i\,j}))$ in 
$M_n(X_0)\wedge\cdots\wedge M_n(X_q)$ to a sum of $(x^0_{j_q\,j_0},x^1_{j_0\,j_1},\ldots,x^q_{j_{q-1}\,j_q})$
taken in $\Gamma_q^+(X_0\wedge\ldots\wedge X_q)$. Here $j_0,\ldots,j_q$ run through $1,2,\ldots,n$.
Since $\Gamma^+$ is not a strictly commutative monoid,
this raises the question about how one should order the summands.

\paragraph{Orderings of $\mathbf{n}^{(q+1)}$:} An ordering is a 
bijection $\lambda:\mathbf{m}\rightarrow\mathbf{n}^{(q+1)}$. If we had to
choose a single ordering, the lexicographical ordering 
would be the natural choice. Denote this ordering by $\lambda_0$. It is defined as follows:
\[
\text{Assume that }(j_0,\ldots,j_q)=\lambda_0(k)\quad\text{and}\quad (j'_0,\ldots,j'_q)=\lambda_0(l)\quad.
\]
Then $k<l$ whenever there exists an $0\leq i\leq q$ such that
\[
j_0=j'_0,\quad j_1=j'_1,\quad\ldots,\quad j_{i-1}=j'_{i-1},\quad\text{and}\quad j_i<j'_i\quad. 
\]
However, in our situation we must also take other 
orderings into account. The reason is that we want the trace to commute with
the cyclic actions. Hence, we consider what happens when we permute the factors of $\mathbf{n}^{(q+1)}$ cyclically.
Let $\lambda_s$ be the ordering defined as:
\[
\text{Assume that }(j_0,\ldots,j_q)=\lambda_s(k)\quad\text{and}\quad (j'_0,\ldots,j'_q)=\lambda_s(l)\quad.
\]
Then $k<l$ whenever there exists an $s\leq i\leq q$ such that
\[
j_s=j'_s,\quad j_{s+1}=j'_{s+1},\quad\ldots,\quad j_{i-1}=j'_{i-1},\quad\text{and}\quad j_i<j'_i\quad, 
\]
or an $0\leq i\leq s-1$ such that
\[
j_s=j'_s,\quad\ldots,\quad j_{q}=j'_{q},\quad j_{0}=j'_{0},\quad\ldots,\quad j_{i-1}=j'_{i-1},\quad\quad\text{and}\quad j_i<j'_i\quad. 
\]

\paragraph{Definition of $\Tr$:} Given $((x^0_{i\,j}),\ldots,(x^q_{i\,j}))$ in $M_n(X_0)\wedge\cdots\wedge M_n(X_q)$, we
are now ready to write down the formula for the trace: First set
\[
\mathbf{x}(j_0,\ldots,j_q)=(x^0_{j_q\,j_0},x^1_{j_0\,j_1},\ldots,x^q_{j_{q-1}\,j_q})\in X_0\wedge\cdots\wedge X_q\quad.
\]

\begin{Def}
Let
\[
\Tr_q((x^0_{i\,j}),\ldots,(x^q_{i\,j}))=
[(\lambda_0^{-1}\lambda_0,\lambda_1^{-1}\lambda_0,\ldots,\lambda_q^{-1}\lambda_0);
\mathbf{x}(\lambda_0(1)),\ldots,\mathbf{x}(\lambda_0(m))]\quad.
\]
\end{Def}

To prove that $\Tr_q$ is involutive we need some facts about the lexicographical orderings. We let $\rho$ denote the
bijection on $\mathbf{n}^{(q+1)}$ defined by $\rho(j_0,\ldots,j_q)=(j_q,\ldots,j_0)$. We also need a technical
definition:

\begin{Def}
Let $\beta:\mathbf{k}\rightarrow\mathbf{m}$ be ordering preserving. We say that $\beta$ is \textit{special} if
the composition
\[
\mathbf{k}\xrightarrow{\beta}\mathbf{m}\xrightarrow{\lambda_0}\mathbf{n}^{(q+1)}\xrightarrow{\pr_i}\mathbf{n}
\]
is injective for all $i$.
\end{Def}

Our technical lemma is:

\begin{Lem}\label{Lem:special}
If $\beta$ is special, then $\beta^*(\lambda_{q-i}^{-1}\rho\lambda_0)=\beta^*(\lambda_{i}^{-1}\lambda_0)$ for all $i$.
\end{Lem}

\begin{proof}
Consider the two diagrams:
\[
\begin{CD}
\mathbf{k}@>{\beta}>> \mathbf{m}@>{\lambda_0}>> \mathbf{n}^{(q+1)} @>{\pr_i}>>\mathbf{n}\\
@V{\beta^*(\lambda_{i}^{-1}\lambda_0)}VV @V{\lambda_{i}^{-1}\lambda_0}VV @VV{=}V @VV{=}V\\
\mathbf{k}@>\text{order preserving}>> \mathbf{m}@>{\lambda_i}>> \mathbf{n}^{(q+1)} @>{\pr_i}>>\mathbf{n}
\end{CD}
\]
and
\[
\begin{CD}
\mathbf{k}@>{\beta}>> \mathbf{m}@>{\lambda_0}>> \mathbf{n}^{(q+1)} @>{\pr_i}>>\mathbf{n}\\
@V{\beta^*(\lambda_{q-i}^{-1}\rho\lambda_0)}VV @V{\lambda_{q-i}^{-1}\rho\lambda_0}VV @VV{\rho}V @VV{=}V\\
\mathbf{k}@>\text{order preserving}>> \mathbf{m}@>{\lambda_{q-i}}>> \mathbf{n}^{(q+1)} @>{\pr_{q-i}}>>\mathbf{n}
\end{CD}\quad.
\]
The left squares of the first and second diagram define $\beta^*(\lambda_{i}^{-1}\lambda_0)$ and $\beta^*(\lambda_{q-i}^{-1}\rho\lambda_0)$
respectively. Now notice that the compositions at the bottom of both diagrams are order preserving, while the maps at the top are
the same for both diagrams. Since there is a unique factorization of the injective map
\[
\pr_i\circ\lambda_0\circ\beta:\mathbf{k}\rightarrow \mathbf{n}
\]
as the composition of a permutation $\mathbf{k}\rightarrow\mathbf{k}$ and an order preserving map $\mathbf{k}\rightarrow \mathbf{n}$,
we get that
\[
\beta^*(\lambda_{q-i}^{-1}\rho\lambda_0)=\beta^*(\lambda_{i}^{-1}\lambda_0)\quad.
\]
\end{proof}

\begin{Lem}
$\Tr_q$ is a dihedral collection of natural transformations.
\end{Lem}

\begin{proof}
We have to verify that $\Tr_{\bullet}$ commutes with the operators $d_i$, $s_i$, $t_q$ and $r_q$.
These operators are specified in subsection~\ref{subsect:dihedralcoll}.
Let $((x^0_{a\,b}),\ldots,(x^q_{a\,b}))$ be a point in $M_n(X_0)\wedge\cdots\wedge M_n(X_q)$.

\paragraph{Face operators:}
Set $((y^0_{a\,b}),\ldots,(y^{q-1}_{a\,b}))$ to be equal to $d_i((x^0_{a\,b}),\ldots,(x^q_{a\,b}))$. Then
we have
\[
y^s_{a\,b}=\begin{cases}
x^s_{a\,b}&\text{if $s<i$,}\\
(x^i_{a\,c(a)},x^{i+1}_{c(a),\,b})&\text{if $s=i$, and}\\
x^s_{a\,b}&\text{if $s>i$.}
\end{cases}
\]
Here $c(a)$ is a choice of index such that $x^i_{a\,d}=*$ whenever $d\neq c(a)$.
Define $\alpha:\mathbf{n}^{q}\rightarrow \mathbf{n}^{q+1}$ by
\[
\alpha(j_0,\ldots,j_{q-1})=(j_0,\ldots,j_{i-1},c(j_{i-1}),j_i,\ldots,j_{q-1})\quad.
\]
Let $\lambda'_s:\mathbf{m'}\rightarrow \mathbf{n}^q$ be the cycled lexicographical orderings on $q-1$ factors. Define
$\beta:\mathbf{m'}\rightarrow \mathbf{m}$ to be $\lambda_0\alpha\lambda_0'$. Then one can easily prove that
\[
\beta^*(\lambda_s^{-1}\lambda_0)=\begin{cases}
{\lambda'}_s^{-1}\lambda'_0&\text{if $s<i$, and}\\
{\lambda'}_{s-1}^{-1}\lambda'_0&\text{if $s>i$.}
\end{cases}
\]
Now we see that
\begin{align*}
\Tr_{q-1}&d_i((x^0_{a\,b}),\ldots,(x^q_{a\,b}))\\
&=\Tr_{q-1}((y^0_{a\,b}),\ldots,(y^{q-1}_{a\,b}))\\
&=[({\lambda'}_0^{-1}{\lambda'}_0,{\lambda'}_1^{-1}{\lambda'}_0,\ldots,{\lambda'}_{q-1}^{-1}{\lambda'}_0);
\mathbf{y}(\lambda'_0(1)),\ldots,\mathbf{y}(\lambda'_0(m'))]\\
&=[(\beta^*(\lambda_0^{-1}\lambda_0),\ldots,\beta^*(\lambda_{i-1}^{-1}\lambda_0),
\beta^*(\lambda_{i+1}^{-1}\lambda_0),\ldots,\beta^*(\lambda_q^{-1}\lambda_0));\\
&\quad\quad\quad\quad\quad\quad\quad\quad\quad\quad\quad\quad\quad\quad\quad\quad\quad\quad\quad
\beta^*(\mathbf{x}(\lambda_0(1)),\ldots,\mathbf{x}(\lambda_0(m)))]\\
&=[(\lambda_0^{-1}\lambda_0,\ldots,\lambda_{i-1}^{-1}\lambda_0,\lambda_{i+1}^{-1}\lambda_0,\ldots,\lambda_q^{-1}\lambda_0);
\mathbf{x}(\lambda_0(1)),\ldots,\mathbf{x}(\lambda_0(m))]\\
&=d_i[(\lambda_0^{-1}\lambda_0,\ldots,\lambda_q^{-1}\lambda_0);\mathbf{x}(\lambda_0(1)),\ldots,\mathbf{x}(\lambda_0(m))]\\
&=d_i\Tr_q((x^0_{a\,b}),\ldots,(x^q_{a\,b}))\quad.
\end{align*}

\paragraph{Degeneracy operators:}
Set $((y^0_{a\,b}),\ldots,(y^{q+1}_{a\,b}))$ to be equal to $s_i((x^0_{a\,b}),\ldots,(x^q_{a\,b}))$. Then
we have
\[
y^s_{a\,b}=\begin{cases}
x^s_{a\,b}&\text{if $s<i$,}\\
1\in S^0&\text{if $s=i$ and $a=b$,}\\
*\in S^0&\text{if $s=i$ and $a\neq b$, and}\\
x^{s-1}_{a\,b}&\text{if $s>i$.}
\end{cases}
\]
Define $\alpha:\mathbf{n}^{q+2}\rightarrow \mathbf{n}^{q+1}$ by
\[
\alpha(j_0,\ldots,j_{q-1})=(j_0,\ldots,j_{i-1},j_i,j_i,j_{i+1},\ldots,j_{q})\quad.
\]
Let $\lambda'_s:\mathbf{m'}\rightarrow \mathbf{n}^{q+2}$ be the cycled lexicographical orderings on $q+1$ factors. Define
$\beta:\mathbf{m}\rightarrow \mathbf{m'}$ to be $\lambda'_0\alpha\lambda_0$. Then one can easily prove that
\[
\beta^*({\lambda'}_s^{-1}\lambda'_0)=\begin{cases}
{\lambda}_s^{-1}\lambda_0&\text{if $s\leq i$, and}\\
{\lambda}_{s-1}^{-1}\lambda_0&\text{if $s>i$.}
\end{cases}
\]
Now we see that
\begin{align*}
\Tr_{q+1}&s_i((x^0_{a\,b}),\ldots,(x^q_{a\,b}))\\
&= \Tr_{q+1}((y^0_{a\,b}),\ldots,(y^{q+1}_{a\,b}))\\
&=[({\lambda'}_0^{-1}\lambda'_0,{\lambda'}_1^{-1}\lambda'_0,\ldots,{\lambda'}_{q+1}^{-1}\lambda'_0);
\mathbf{y}(\lambda'_0(1)),\ldots,\mathbf{y}(\lambda'_0(m'))]\\
&=[(\beta^*({\lambda}_0^{-1}\lambda_0),\ldots,\beta^*({\lambda}_i^{-1}\lambda_0),
\beta^*({\lambda}_i^{-1}\lambda_0),\ldots,\beta^*({\lambda}_{q}^{-1}\lambda_0));\\
&\quad\quad\quad\quad\quad\quad\quad\quad\quad\quad\quad\quad\quad\quad\quad\quad\quad\quad\quad
 \beta^*(\mathbf{x}(\lambda_0(1)),\ldots,\mathbf{x}(\lambda_0(m)))]\\
&=[({\lambda}_0^{-1}\lambda_0,\ldots,{\lambda}_i^{-1}\lambda_0,{\lambda}_i^{-1}\lambda_0,\ldots,{\lambda}_{q}^{-1}\lambda_0);
\mathbf{x}(\lambda_0(1)),\ldots,\mathbf{x}(\lambda_0(m))]\\
&=s_i[(\lambda_0^{-1}\lambda_0,\ldots,\lambda_q^{-1}\lambda_0);\mathbf{x}(\lambda_0(1)),\ldots,\mathbf{x}(\lambda_0(m))]\\
&=s_i\Tr_q((x^0_{a\,b}),\ldots,(x^q_{a\,b}))\quad.
\end{align*}

\paragraph{The cyclic operator:}
Set $((y^0_{a\,b}),\ldots,(y^{q}_{a\,b}))$ to be equal to $t_q((x^0_{a\,b}),\ldots,(x^q_{a\,b}))$. Then
we have
\[
y^s_{a\,b}=\begin{cases}
x^{q}_{a\,b}&\text{if $s=0$, and}\\
x^{s-1}_{a\,b}&\text{if $s>0$.}
\end{cases}
\]
Define $\alpha:\mathbf{n}^{q+2}\rightarrow \mathbf{n}^{q+1}$ by
\[
\alpha(j_0,\ldots,j_{q-1})=(j_{q},j_0,\ldots,j_{q-1})\quad.
\]
Notice that $\alpha\lambda_s=\lambda_{s+1}$ for $s<q$ and $\alpha\lambda_q=\lambda_0$. Thus 
\[
\lambda_s^{-1}\lambda_0\lambda_q^{-1}\lambda_0=\begin{cases}
\lambda_q^{-1}\lambda_0&\text{if $s=0$, and}\\
\lambda_{s-1}^{-1}\lambda_0&\text{if $s>0$.}
\end{cases}
\]
Let $\tau:X_0\wedge\cdots\wedge X_q\rightarrow X_q\wedge X_0\wedge\cdots\wedge X_{q-1}$ be the homeomorphism which permutes the factors.
We have
\[
\mathbf{y}(\alpha(j_0,\ldots,j_q))=\tau(\mathbf{x}(j_0,\ldots,j_q))\quad.
\]
Now we see that
\begin{align*}
\Tr_{q}&t_q((x^0_{a\,b}),\ldots,(x^q_{a\,b}))\\
&= \Tr_{q}((y^0_{a\,b}),\ldots,(y^{q}_{a\,b}))\\
&=[({\lambda}_0^{-1}\lambda_0,{\lambda}_1^{-1}\lambda_0,\ldots,{\lambda}_{q}^{-1}\lambda_0);
\mathbf{y}(\lambda_0(1)),\ldots,\mathbf{y}(\lambda_0(m))]\\
&=[({\lambda}_0^{-1}\lambda_0,{\lambda}_1^{-1}\lambda_0,\ldots,{\lambda}_{q}^{-1}\lambda_0);
\tau\mathbf{x}(\lambda_q(1)),\ldots,\tau\mathbf{x}(\lambda_q(m))]\\
&=[({\lambda}_0^{-1}\lambda_0,{\lambda}_1^{-1}\lambda_0,\ldots,{\lambda}_{q}^{-1}\lambda_0);
(\tau\mathbf{x}(\lambda_0(1)),\ldots,\tau\mathbf{x}(\lambda_0(m))).(\lambda_0^{-1}\lambda_q)]\\
&=[({\lambda}_0^{-1}\lambda_0\lambda_q^{-1}\lambda_0,{\lambda}_1^{-1}\lambda_0\lambda_q^{-1}\lambda_0,\ldots,
{\lambda}_{q}^{-1}\lambda_0\lambda_q^{-1}\lambda_0);
\tau\mathbf{x}(\lambda_0(1)),\ldots,\tau\mathbf{x}(\lambda_0(m))]\\
&=[({\lambda}_q^{-1}\lambda_0,{\lambda}_0^{-1}\lambda_0,\ldots,{\lambda}_{q-1}^{-1}\lambda_0);
\tau\mathbf{x}(\lambda_0(1)),\ldots,\tau\mathbf{x}(\lambda_0(m))]\\
&=t_q[(\lambda_0^{-1}\lambda_0,\ldots,\lambda_q^{-1}\lambda_0);\mathbf{x}(\lambda_0(1)),\ldots,\mathbf{x}(\lambda_0(m))]\\
&=t_q\Tr_q((x^0_{a\,b}),\ldots,(x^q_{a\,b}))\quad.
\end{align*}

\paragraph{The involutive operator:}
Set $((y^0_{a\,b}),\ldots,(y^{q}_{a\,b}))$ to be equal to $r_q((x^0_{a\,b}),\ldots,(x^q_{a\,b}))$. Then
we have
\[
y^s_{a\,b}=\begin{cases}
x^{0}_{b\,a}&\text{if $s=0$, and}\\
x^{q+1-s}_{b\,a}&\text{if $s>0$.}
\end{cases}
\]
Recall that  $\rho:\mathbf{n}^{q+2}\rightarrow \mathbf{n}^{q+1}$ was defined by
\[
\rho(j_0,\ldots,j_{q-1})=(j_{q},j_{q-1},\ldots,j_{1},j_0)\quad.
\]
Recall also the fact that 
$\beta^*(\lambda_{q-s}^{-1}\rho\lambda_0)=\beta^*(\lambda_{s}^{-1}\lambda_0)$ if $\beta$ 
is special, see lemma~\ref{Lem:special}.
Let $\tau:X_0\wedge\cdots\wedge X_q\rightarrow X_q\wedge X_{q-1}\wedge\cdots\wedge X_{1}\wedge X_0$ 
be the homeomorphism which permutes the factors.
We have
\[
\mathbf{y}(\rho(j_0,\ldots,j_q))=\tau(\mathbf{x}(j_0,\ldots,j_q))\quad.
\]
Observe that there exists a $\beta:\mathbf{k}\rightarrow\mathbf{m}$ such that $\beta$ is entire for 
$(\mathbf{y}(\lambda_0(1)),\ldots,\mathbf{y}(\lambda_0(m)))$ and $\beta$ is special. This is a consequence of that
each column of a matrix in $M_n(X)$ contains at most one element different from $*$.
Now we see that
\begin{align*}
\Tr_{q}&r_q((x^0_{a\,b}),\ldots,(x^q_{a\,b}))\\
&= \Tr_{q}((y^0_{a\,b}),\ldots,(y^{q}_{a\,b}))\\
&=[({\lambda}_0^{-1}\lambda_0,{\lambda}_1^{-1}\lambda_0,\ldots,{\lambda}_{q}^{-1}\lambda_0);
\mathbf{y}(\lambda_0(1)),\ldots,\mathbf{y}(\lambda_0(m))]\\
&=[(\beta^*({\lambda}_0^{-1}\lambda_0),\beta^*({\lambda}_1^{-1}\lambda_0),\ldots,\beta^*({\lambda}_{q}^{-1}\lambda_0));
\beta^*(\mathbf{y}(\lambda_0(1)),\ldots,\mathbf{y}(\lambda_0(m)))]\\
&=[(\beta^*({\lambda}_q^{-1}\rho\lambda_0),\beta^*({\lambda}_{q-1}^{-1}\rho\lambda_0),\ldots,\beta^*({\lambda}_{0}^{-1}\rho\lambda_0));
\beta^*(\mathbf{y}(\lambda_0(1)),\ldots,\mathbf{y}(\lambda_0(m)))]\\
&=[({\lambda}_q^{-1}\rho\lambda_0,{\lambda}_{q-1}^{-1}\rho\lambda_0,\ldots,{\lambda}_{0}^{-1}\rho\lambda_0);
\mathbf{y}(\lambda_0(1)),\ldots,\mathbf{y}(\lambda_0(m))]\\
&=[({\lambda}_q^{-1}\rho\lambda_0,{\lambda}_{q-1}^{-1}\rho\lambda_0,\ldots,{\lambda}_{0}^{-1}\rho\lambda_0);
\tau\mathbf{x}(\rho\lambda_0(1)),\ldots,\mathbf{x}(\rho\lambda_0(m))]\\
&=[({\lambda}_q^{-1}\lambda_0,{\lambda}_{q-1}^{-1}\lambda_0,\ldots,{\lambda}_{0}^{-1}\lambda_0);
\tau\mathbf{x}(\lambda_0(1)),\ldots,\mathbf{x}(\lambda_0(m))]\\
&=r_q[(\lambda_0^{-1}\lambda_0,\ldots,\lambda_q^{-1}\lambda_0);\mathbf{x}(\lambda_0(1)),\ldots,\mathbf{x}(\lambda_0(m))]\\
&=r_q\Tr_q((x^0_{a\,b}),\ldots,(x^q_{a\,b}))\quad.
\end{align*}
\end{proof}

This lemma shows that $\Tr:THH(\Gamma M_n(L))\rightarrow THH^+(L)$ is an $S^1$-map,
or $O(2)$-map when $L$ comes with an involution. 
Now the next statement of proposition~\ref{Prop:Morita} is that
$\Tr:THH(\Gamma M_n(L))\rightarrow THH^+(L)$
is non-equivariantly a $\pi_*$-isomorphism. Let us now prove this:

\begin{proof}
We are assuming that $L$ is a cofibrant orthogonal ring spectrum (with involution).
Let $W_n$ be the endofunctor on $\Top_*$ defined by $X\mapsto \bigvee_{n^2}X$. 
We write elements of $W_n(X)$ on the form
\[
(a,x,b)\quad\text{where $a,b\in\{1,\ldots,n\}$ and $x\in X$.}
\]
We think about $W_n$ as an ``FSP without unit''. 
Multiplication $\mu:W_n(X)\wedge W_n(Y)\rightarrow W_n(X\wedge Y)$ is given by
\[
\mu((a,x,b);(c,y,d))=\begin{cases}
(a,(x,y),d)&\text{if $b=c$, and}\\
*&\text{otherwise.}
\end{cases}
\]
We have involution given by $(a,x,b)\mapsto(b,x,a)$. 
Furthermore, there is a natural 
transformation from $W_n$ to $M_n$ which respects both multiplication and
involution. This natural transformation is given by sending $(a,x,b)$ to the matrix $(x_{i\,j})$,
where $x_{a\,b}=x$ and $x_{i\,j}=*$ otherwise.

$THH_{\bullet}(W_n(L))$ is a presimplicial orthogonal spectrum, 
the $q$-simplices are $(W_n(L))^{\wedge(q+1)}$, and face maps are given by the usual formulas.
Moreover, we have a natural presimplicial map $THH_{\bullet}(W_n(L))\rightarrow THH_{\bullet}(M_n(L))$. 
However, we do not know that  
$THH_{\bullet}(M_n(L))$ is good as a simplicial orthogonal spectrum. 
The problem is that $M_n(L)$ might not be cofibrant.
Therefore we replace $M_n(L)$ by $\Gamma M_n(L)$, and consider
\[
THH_{\bullet}(\Gamma W_n(L))\rightarrow THH_{\bullet}(\Gamma M_n(L))\quad.
\]
By proposition~\ref{Prop:WMhtpyequiv}
and results about induced functors, we see that the presimplicial realization of this map is a $\pi_*$-iso.
Notice that $W_n(L)$ is a wedge of copies of $L$, thus cofibrant. Hence, the natural map
\[
THH_{\bullet}(\Gamma W_n(L))\rightarrow THH_{\bullet}(W_n(L))
\]
is a $\pi_*$-iso in each simplicial degree.

There is a trace map $\Tr:W_n(X_0)\wedge\cdots\wedge W_n(X_q)\rightarrow X_0\wedge\cdots\wedge X_q$ given by
\[
\Tr((a_0,x_0,b_0),\ldots,(a_q,x_q,b_q))=\begin{cases}
(x_0,\ldots,x_q)&\text{if $b_0=a_1$, $b_1=a_2$, $\ldots$, $b_q=a_0$, and}\\
*&\text{otherwise.}
\end{cases}
\]
We get an induced presimplicial map $THH_{\bullet}(W_n(L))\rightarrow THH_{\bullet}(L)$, which fits into the diagram
\[
\begin{CD}
THH_{\bullet}(\Gamma W_n(L)) @>{\simeq}>> THH_{\bullet}(W_n(L)) @>{\Tr}>> THH_{\bullet}(L)\\
@V{\simeq}VV @VVV @VV{\simeq}V\\
THH_{\bullet}(\Gamma M_n(L)) @>>> THH_{\bullet}(M_n(L)) @>{\Tr}>> THH^+_{\bullet}(L)
\end{CD}\quad.
\]
We already know that after presimplicial realization the left and right vertical maps become $\pi_*$-isomorphisms.
Hence, it remains only to show that the geometric realization of the top $\Tr$ is also a $\pi_*$-isomorphism.

To show this we construct a presimplicial homotopy inverse:
There is a presimplicial map $\incl:THH_{\bullet}(L)\rightarrow THH_{\bullet}(W_n(L))$ 
defined by the natural transformation
\[
X\rightarrow W_n(X)\quad,\text{ which sends $x$ to $(1,x,1)$.}
\]
It is easily seen that $\Tr\circ\incl$ is the identity. 
We complete the proof by constructing a presimplicial homotopy from
$\incl\circ\Tr$ to the identity on $THH_{\bullet}(W_n(L))$. 

Let the natural transformations
\small
\[
h_i:W_n(X_0)\wedge\cdots\wedge W_n(X_q)\rightarrow W_n(X_0)\wedge \cdots\wedge W_n(X_{i})\wedge W_n(S^0)\wedge W_n(X_{i+1})\wedge\cdots\wedge W_n(X_q)
\]
\normalsize
be given by
\scriptsize
\[
h_i((a_0,x_0,b_0),\ldots,(a_q,x_q,b_q))=\begin{cases}
((a_0,x_0,1),\ldots&,(1,x_{i},1),(1,1,b_{i}),(a_{i+1},x_{i+1},b_{i+1}),\ldots,(a_q,x_q,b_q))\\
&\text{if $b_0=a_1$, $b_1=a_2$, $\ldots$, $b_{i-1}=a_i$, and}\\
*&\text{otherwise.}
\end{cases}
\]
\normalsize
We see that
\begin{align*}
d_0h_0 &= \id\quad,\\
d_ih_j &= h_{j-1} d_i \quad\text{for $i<j$,}\\
d_ih_i &= d_ih_{i-1} \quad,\\
d_ih_j &= h_j d_{i-1} \quad\text{for $i>j+1$, and}\\
d_{q+1}h_q &= \incl\Tr_q \quad.
\end{align*}
Hence, $h$ is a presimplicial homotopy between $\id$ and $\incl\circ\Tr$.
\end{proof}

To finish the proof of proposition~\ref{Prop:Morita}, it remains only to check that direct
sum of matrices corresponds to the addition induced by the Barratt-Eccles functor, $\Gamma^+$.
To check this, we first consider the $q$-simplices via the external viewpoint. 
The following lemma is the key:

\begin{Lem}
The diagram
\scriptsize
\[
\begin{CD}
(M_{n_1}\times M_{n_2})(X_0)\wedge\cdots\wedge (M_{n_1}\times M_{n_2})(X_q) @>{(\Tr\circ\pr_1,\Tr\circ\pr_2)}>> 
\Gamma^+(X_0\wedge \cdots \wedge X_q)\times \Gamma^+(X_0\wedge \cdots \wedge X_q)\\
@V{\oplus}VV @VV{+}V\\
(M_{n_1+n_2})(X_0)\wedge\cdots\wedge (M_{n_1+n_2})(X_q) @>{\Tr}>> \Gamma^+(X_0\wedge \cdots \wedge X_q)
\end{CD}
\]
\normalsize
commutes.
\end{Lem}

\begin{proof}
Set $((x^0_{a\,b}),\ldots,(x^{q}_{a\,b}))$ to be equal to 
$((y^0_{a\,b})\oplus (z^0_{a\,b}),\ldots,(y^q_{a\,b})\oplus (z^q_{a\,b}))$. Then
we have
\[
x^s_{a\,b}=\begin{cases}
y^s_{a\,b}&\text{if $a\leq n_1$ and $b\leq n_1$,}\\
z^s_{(a-n_1)\,(b-n_1)} &\text{if $a> n_1$ and $b> n_1$, and}\\
*&\text{otherwise.}
\end{cases}
\]
Define $\alpha:\mathbf{n_1}^{q+1}\amalg\mathbf{n_2}^{q+1}\rightarrow \mathbf{(n_1+n_2)}^{q+1}$ by
\[
\alpha(j_0,\ldots,j_q)=\begin{cases}
(j_0,\ldots,j_q)&\text{for $(j_0,\ldots,j_q)\in \mathbf{n_1}^{q+1}$, and}\\
(j_0+n_1,\ldots,j_q+n_1)&\text{for $(j_0,\ldots,j_q)\in \mathbf{n_2}^{q+1}$.}
\end{cases}
\]
Denote by $\lambda_s^1$, $\lambda_s^2$ and $\lambda_s$ the cycled lexicographical 
orderings of $\mathbf{n_1}^{q+1}$, $\mathbf{n_2}^{q+1}$
and $\mathbf{(n_1+n_2)}^{q+1}$ respectively. 
Define $\beta_s$ to be the unique map such that the diagram below commutes:
\[
\begin{CD}
\mathbf{m_1+m_2} @>{\lambda_s^1\amalg\lambda^2_s}>> \mathbf{n_1}^{q+1}\amalg\mathbf{n_2}^{q+1}\\
@V{\beta_s}VV @VV{\alpha}V\\
\mathbf{m} @>{\lambda_s}>> \mathbf{(n_1+n_2)}^{q+1}
\end{CD}\quad.
\]
By the definition of $\alpha$ and the $\lambda$'s we see that 
$\beta_s$ is order preserving. And from the commutative diagram
\[
\begin{CD}
\mathbf{m_1+m_2} @>{\beta_0}>> \mathbf{m}\\
@V{(\lambda_s^1\amalg\lambda_s^2)^{-1}(\lambda_0^1\amalg\lambda_0^2)}VV @V{\lambda^{-1}_s\lambda_0}VV\\
\mathbf{m_1+m_2} @>{\beta_s}>> \mathbf{m}
\end{CD}
\]
it follows that
\[
\beta_0^*(\lambda^{-1}_s\lambda_0)=((\lambda^{1}_s)^{-1}\lambda_0^1)\amalg((\lambda^{2}_s)^{-1}\lambda_0^2)\quad.
\]
Now we see that
\small
\begin{align*}
\Tr\circ\pr_1&((y^0_{a\,b})\oplus (z^0_{a\,b}),\ldots,(y^q_{a\,b})\oplus (z^q_{a\,b}))+
\Tr\circ\pr_2((y^0_{a\,b})\oplus (z^0_{a\,b}),\ldots,(y^q_{a\,b})\oplus (z^q_{a\,b}))\\
&=\Tr((y^0_{a\,b}),\ldots,(y^q_{a\,b}))+\Tr((z^0_{a\,b}),\ldots,(z^q_{a\,b}))\\
&=[((\lambda^1)_0^{-1}\lambda^1_0,\ldots,(\lambda^1)_{q}^{-1}\lambda^1_0);
\mathbf{y}(\lambda^1_0(1)),\ldots,\mathbf{y}(\lambda^1_0(m_1))]\quad\quad\quad\\
&\quad\quad\quad\quad\quad\quad+[((\lambda^2)_0^{-1}\lambda^2_0,\ldots,(\lambda^2)_{q}^{-1}\lambda^2_0);
\mathbf{z}(\lambda^2_0(1)),\ldots,\mathbf{z}(\lambda^2_0(m_2))]\\
&=[(((\lambda^1)_0^{-1}\lambda^1_0)\amalg ((\lambda^2)_0^{-1}\lambda^2_0),\ldots,
((\lambda^1)_{q}^{-1}\lambda^1_0))\amalg ((\lambda^2)_{q}^{-1}\lambda^2_0));\quad\quad\quad\\
&\quad\quad\quad\quad\quad\quad\quad\mathbf{y}(\lambda^1_0(1)),\ldots,\mathbf{y}(\lambda^1_0(m_1)),\mathbf{z}(\lambda^2_0(1)),\ldots,\mathbf{z}(\lambda^2_0(m_2))]\\
&=[(\beta_0^*({\lambda}_0^{-1}\lambda_0),\ldots,\beta_0^*({\lambda}_{q}^{-1}\lambda_0));
\mathbf{x}(\lambda_0\beta_0(1)),\ldots,\mathbf{x}(\lambda_0\beta_0(m))]\\
&=[({\lambda}_0^{-1}\lambda_0,\ldots,{\lambda}_{q}^{-1}\lambda_0);
\mathbf{x}(\lambda_0(1)),\ldots,\mathbf{x}(\lambda_0(m))]\\
&=\Tr((x^0_{a\,b}),\ldots,(x^{q}_{a\,b}))\quad.
\end{align*}
\normalsize
This concludes the proof of the lemma.
\end{proof}

Let us now finish the proof of proposition~\ref{Prop:Morita}:

\begin{proof}
By internalizing the diagram of the lemma, we get that
\[
\begin{CD}
(M_{n_1}\times M_{n_2})(L)^{\wedge(q+1)} @>{(\Tr\circ\pr_1,\Tr\circ\pr_2)}>> 
\Gamma^+(L^{\wedge(q+1)})\times \Gamma^+(L^{\wedge(q+1)})\\
@V{\oplus}VV @VV{+}V\\
(M_{n_1+n_2})(L)^{\wedge(q+1)} @>{\Tr}>> \Gamma^+(L^{\wedge(q+1)})
\end{CD}
\]
commutes for any $L$ orthogonal ring spectrum (with involution). 
We identify the corners with $THH_q((M_{n_1}\times M_{n_2})(L))$, $THH_q^+(L)\times THH_q^+(L)$,
$THH_q((M_{n_1+n_2})(L))$ and $THH_q^+(L)$.
Take geometric realization and use the natural transformation $\Gamma X\rightarrow X$ to get
the commutative diagram
\footnotesize
\[
\begin{CD}
THH(\Gamma (M_{n_1}\times M_{n_2})(L)) &\rightarrow& THH((M_{n_1}\times M_{n_2})(L)) @>{(\Tr\circ\pr_1,\Tr\circ\pr_2)}>> THH^+(L)\times THH^+(L)\\
@V{\oplus}VV @V{\oplus}VV @VV{+}V\\
THH(\Gamma M_{n_1+n_2}(L)) &\rightarrow& THH(M_{n_1+n_2}(L)) @>{\Tr}>> THH^+(L)
\end{CD}\quad.
\]
\normalsize
The outer square is the diagram we are interested in.
\end{proof}

\appendix
\chapter{Useful results}

\section{Homotopy theory}

Let us first recall some results about CW-approximations:

\begin{Prop}
For every space $X$ there exists a CW-complex $Z$ and a weak homotopy equivalence
$f:Z\rightarrow X$. Furthermore, $Z$ is unique up to homotopy equivalence.
\end{Prop}

See theorem~7.8.1 in~\cite{Spanier:91} for a proof.
In the relative case we will use the following:

\begin{Prop}\label{Prop:CWapprox}
For every cofibration $i:A\rightarrow X$ there exists a CW-pair $(Z,C)$ and a commutative diagram
\[
\begin{CD}
C@>>> Z\\
@V{f_0}VV @VV{f}V\\
A@>{i}>> X
\end{CD}\quad,
\]
where $f$ and $f_0$ are weak equivalences. Furthermore, $(Z,C)$ is unique up to homotopy equivalence.
\end{Prop}

\begin{proof}
First choose a CW-approximation $f_0:C\rightarrow A$. Let $M$ be the mapping cylinder of $f_0$.
Since $i$ is a cofibration, the map $M\cup_A X\rightarrow X$ induced by the projection $M\rightarrow A$
is a homotopy equivalence. By proposition~4.13 in~\cite{Hatcher:02} there exists a CW-space $Z$ containing $C$ as
a subcomplex and a weak equivalence of pairs $f':(Z,C)\rightarrow (M\cup_A X,C)$. Composing with the
projection $(M\cup_A X,C)\rightarrow (X,A)$ we get our CW-approximation.

Uniqueness follows by applying corollary~4.19 in~\cite{Hatcher:02} twice,
first to show that the choice of $C$ is unique up to homotopy,
then to show that for fixed $C$ the choice of $Z$ is unique up to homotopy rel $C$.
\end{proof}

In order to efficiently apply CW-approximation we need a gluing theorem for weak equivalences. The proof
is formal once the following lemma is established:

\begin{Lem}
Let $Y$ be the pushout of $X\overset{f}{\leftarrow} A\overset{i}{\rightarrow} B$, where $f$ is a weak equivalence
and $i$ a cofibration. Then $g:B\rightarrow Y$ is also a weak equivalence.
\end{Lem}

\begin{proof}
Since $i$ is a cofibration, we have that $Y$ is homotopic to the homotopy pushout. 
Hence, both the van Kampen theorem and
Mayer-Vietoris sequences can be applied. 
By elementary considerations it is seen that $g$ induces a bijection of path components.
For $\pi_1$ we consider each path-component of $B$ separately, so we may just as well assume that $B$ is path-connected.
If $A$ also is path-connected, then the van Kampen theorem applies 
to $Y=X\cup B$ and shows that $\pi_1 B\rightarrow \pi_1 Y$
is an isomorphism. When $A$ has more than one path-component, we write
\[
A=\bigcup A_{\alpha}\quad\text{and}\quad X=\bigcup X_{\alpha}
\]
where all $A_{\alpha}$ and $X_{\alpha}$ are path connected. Then we apply van Kampen to the union
\[
Y=\bigcup(X_{\alpha}\cup B)\quad.
\]
Since each $X_{\alpha}\cup B$ is $\pi_1$-isomorphic to $B$, it follows that $g$ is a $\pi_1$-iso.

For the higher homotopy groups we use Mayer-Vietoris sequences and the Hurewicz theorem.
\end{proof}

Form this lemma it is a formal argument due to Thomas Gunnarsson, see the proof of lemma~8.8 
in~\cite{GoerssJardine:99}, to show the gluing theorem:

\begin{Prop}\label{Prop:gluewe}
If we have a commutative diagram
\[
\begin{CD}
X @<<< A @>{i}>> B\\
@V{\simeq}VV @VV{\simeq}V @VV{\simeq}V\\
X' @<<< A' @>{i'}>> B'
\end{CD}
\]
where $i$ and $i'$ are cofibrations and the vertical maps are weak equivalences, then the map of pushouts
is also a weak equivalence.
\end{Prop}

We have several times used the Blakers-Massey homotopy excision theorem.
The following form of the theorem is suitable for our purposes:

\begin{Thm}\label{Thm:BlakersMassey}
Suppose that $X$ is a pointed space and that $A$ and $B$ are pointed subspaces of $X$ such that
\begin{itemize}
\item[] $X=A\cup B$ and
\item[] the inclusions $A\cap B\rightarrow A$ and $A\cap B\rightarrow B$ are cofibrations.
\end{itemize}
If the pair $(A,A\cap B)$ is $m$-connected and the pair $(B,A\cap B)$ is $n$-connected,
$m\geq 0$, $n\geq 0$, then the homomorphism induced by the inclusion, namely $i_*:\pi_q(A,A\cap B)\rightarrow \pi_q(X,B)$,
is an isomorphism for $q< m+n$ and is surjective for $q=m+n$.
\end{Thm}

By proposition~\ref{Prop:gluewe} we can apply CW-approximation. Then
the result follows from theorem~4.23 in~\cite{Hatcher:02}.
Two useful consequences of homotopy excision are:

\begin{Prop}\label{Prop:homotopycofseq}
Suppose that $A\rightarrow X$ is a cofibration, that the pair $(X,A)$ is $(r-1)$-connected, and that the subspace
$A$ is $(s-1)$-connected, $r\geq 1$, $s\geq1$. Then the homomorphism induced by the quotient map, namely
\[
\pi_q(X,A)\rightarrow \pi_q(X/A)\quad,
\]
is an isomorphism for $q<r+s-1$ and surjective for $q=r+s-1$.
\end{Prop}

\begin{Thm}
Let $X$ be an $(n-1)$-connected well-pointed space. Then the suspension map $\pi_q(X)\rightarrow \pi_{q+1}(S^1 \wedge X)$
is an isomorphism for $q<2n-1$ and surjective for $q=2n-1$.
\end{Thm}

The last result is known as the Freudenthal suspension theorem. Proofs of both results can be found in~\cite{Hatcher:02}.
Just observe that his conditions concern CW-pairs instead of cofibrations, but the only places where he uses these
conditions in his proofs, are when applying the Blakers-Massey homotopy extension theorem.

As a corollary of the proposition we have:

\begin{Cor}\label{Cor:collweakcontr}
If $A\rightarrow X$ is a cofibration and $A$ is weakly contractible, then
$X\rightarrow X/A$ is a weak equivalence.
\end{Cor}

\begin{proof}
We can assume that $X$ is path-connected without loss of generality. By the proposition, the maps
\[
\pi_q(X,A)\rightarrow \pi_q(X/A)
\]
are isomorphisms for all $q$. Furthermore, the maps $\pi_q(X)\rightarrow \pi_q(X,A)$
are also isomorphisms for all $q$ since $A$ is weakly contractible.
\end{proof}

\begin{Prop}\label{Prop:connofsmash}
If $X$ and $Y$ are well-pointed spaces $(r-1)$- and $(s-1)$-connected respectively, then
$X\wedge Y$ is $(r+s-1)$-connected.
\end{Prop}

\begin{proof}
We can approximate $X$ and $Y$ by CW-complexes $X'$ and $Y'$ such that all cells except $*$ has
dimension greater than $(r-1)$ and $(s-1)$ respectively. By proposition~\ref{Prop:gluewe}
we have a weak equivalence $X'\wedge Y'\rightarrow X\wedge Y$. Since all cells of $X'\wedge Y'$
have dimension greater than $(r+s-1)$, the smash product is $(r+s-1)$-connected.
\end{proof}

\section{Monoidal categories}\label{sect:MonCat}

Here we specify the language used for monoidal categories. The standard reference is~\cite{MacLane:98}. 
See also~\S{}20 in~\cite{MandellMaySchwedeShipley:01}.

\begin{Def}
A \textit{monoidal category} $M$ is a category with a bifunctor, $\square:M\times M\rightarrow M$,
a unit $e\in M$ and natural isomorphisms 
\begin{align*}
\alpha:a\square(b\square c) &\cong (a\square b)\square c\quad,\\
\lambda:e\square a &\cong a\quad,\text{ and}\\
\rho:a\square e&\cong a\quad,
\end{align*}
such that the diagrams (i), (ii) and (iii) commute.
\end{Def}

\begin{equation}
\begin{CD}
a\square(b\square(c\square d)) @>{\alpha}>> (a \square b )\square (c\square d) @>{\alpha}>> ((a\square b)\square c)\square d\\
@V{\id\square\alpha}VV && @AA{\alpha\square\id}A\\
a\square((b\square c)\square d) & @>{\alpha}>> & (a\square(b\square c))\square d
\end{CD}
\tag{i}\end{equation}

\begin{equation}
\begin{CD}
a\square( e\square c) @>{\alpha}>> (a \square e) \square c\\
@V{\id\square\lambda}VV @VV{\rho\square\id}V\\
a\square c @= a\square c
\end{CD}
\tag{ii}\end{equation}

\begin{equation}
\begin{CD}
e\square e @= e\square e\\
@V{\lambda}VV @VV{\rho}V\\ 
e @= e
\end{CD}
\tag{iii}\end{equation}

\begin{Def}
A \textit{symmetric monoidal category} $M$ is a monoidal category $M$ with a natural isomorphism
\[
\gamma: a\square b\cong b\square a
\]
such that $\gamma^2=\id$ and the diagrams (iv) and (v) commute.
\end{Def}

\begin{equation}
\begin{CD}
a\square e @>{\gamma}>> e \square a\\
@V{\rho}VV @VV{\lambda}V\\
a @= a
\end{CD}
\tag{iv}\end{equation}

\begin{equation}
\begin{CD}
(a\square b)\square c @>{\gamma}>> c\square (a\square b) @>{\alpha}>> (c\square a)\square b\\
@A{\alpha}AA && @V{\gamma\square\id}VV\\
a\square (b\square c) @>{\id\square\gamma}>> a\square (c\square b) @>{\alpha}>> (a\square c)\square b
\end{CD}
\tag{v}\end{equation}

\begin{Def}
A functor $F:M\rightarrow B$ between monoidal categories is \textit{lax monoidal}
if there is a map $\eta:e_B\rightarrow F(e_M)$ and a natural transformation
\[
\phi:F(a)\square F(b)\rightarrow F(a\square b)
\]
such that the diagrams (vi), (vii) and (viii) commute.
\end{Def}

\begin{equation}
\begin{CD}
F(a)\square(F(b)\square F(c)) @>{\id\square\phi}>> F(a)\square F(b\square c) @>{\phi}>> F(a\square(b\square c))\\
@V{\alpha}VV && @VV{F(\alpha)}V\\
(F(a)\square F(b))\square F(c) @>{\phi\square\id}>> F(a \square b)\square F(c) @>{\phi}>> F((a\square b)\square c)
\end{CD}
\tag{vi}\end{equation}

\begin{equation}
\begin{CD}
F(a)\square e_B @>{\rho}>> F(a)\\
@V{\id\square\eta}VV @AA{F(\rho)}A\\
F(a)\square F(e_M) @>{\phi}>> F(a\square e_M)
\end{CD}
\tag{vii}\end{equation}

\begin{equation}
\begin{CD}
e_B\square F(b) @>{\lambda}>> F(b)\\
@V{\eta\square\id}VV @AA{F(\lambda)}A\\
F(e_M)\square F(b) @>{\phi}>> F(e_M\square b)
\end{CD}
\tag{viii}\end{equation}

\begin{Def}
A functor $F:M\rightarrow B$ between monoidal categories is \textit{lax comonoidal}
if there is a map $\eta:F(e_M)\rightarrow e_B$ and a natural transformation
\[
\phi:F(a\square b)\rightarrow F(a)\square F(b)
\]
such that the diagrams (ix), (x) and (xi) commute.
\end{Def}

\begin{equation}
\begin{CD}
F(a)\square(F(b)\square F(c)) @<{\id\square\phi}<< F(a)\square F(b\square c) @<{\phi}<< F(a\square(b\square c))\\
@V{\alpha}VV && @VV{F(\alpha)}V\\
(F(a)\square F(b))\square F(c) @<{\phi\square\id}<< F(a \square b)\square F(c) @<{\phi}<< F((a\square b)\square c)
\end{CD}
\tag{ix}\end{equation}

\begin{equation}
\begin{CD}
F(a)\square e_B @>{\rho}>> F(a)\\
@A{\id\square\eta}AA @AA{F(\rho)}A\\
F(a)\square F(e_M) @<{\phi}<< F(a\square e_M)
\end{CD}
\tag{x}\end{equation}

\begin{equation}
\begin{CD}
e_B\square F(b) @>{\lambda}>> F(b)\\
@A{\eta\square\id}AA @AA{F(\lambda)}A\\
F(e_M)\square F(b) @<{\phi}<< F(e_M\square b)
\end{CD}
\tag{xi}\end{equation}

\begin{Def}
A lax monoidal functor $F:M\rightarrow B$ between symmetric monoidal categories
is \textit{lax symmetric monoidal} if diagram (xii) commutes.
\end{Def}

\begin{equation}
\begin{CD}
F(a)\square F(b) @>{\gamma}>> F(b)\square F(a)\\
@V{\phi}VV @VV{\phi}V\\
F(a\square b) @>{F(\gamma)}>> F(b\square a)
\end{CD}
\tag{xii}\end{equation}

\begin{Def}
A lax comonoidal functor $F:M\rightarrow B$ between symmetric monoidal categories
is \textit{lax symmetric comonoidal} if diagram (xiii) commutes.
\end{Def}

\begin{equation}
\begin{CD}
F(a\square b) @>{F(\gamma)}>> F(b\square a)\\
@V{\phi}VV @VV{\phi}V\\
F(a)\square F(b) @>{\gamma}>> F(b)\square F(a)
\end{CD}
\tag{xiii}\end{equation}

\begin{Def}
A lax monoidal functor $F:M\rightarrow B$ between monoidal categories
is \textit{strong monoidal} if $\eta$ and $\phi$ are isomorphisms.
$F$ is \textit{strong symmetric monoidal} if $F$ is both strong monoidal and lax symmetric monoidal.
\end{Def}

\begin{Rem}
We could also define \textit{strong comonoidal}, but this definition is redundant since demanding that
$\eta$ and $\phi$ are isomorphisms, for a lax comonoidal functor $F$, would imply that
$F$ together with $\eta^{-1}$ and $\phi^{-1}$ is strong monoidal.
\end{Rem}

\begin{Lem}\label{Lem:catprodimpllaxsymfunct}
If $F:M\rightarrow B$ is a functor between symmetric monoidal categories where $\square_M$ and $\square_B$ are categorical products
for $M$ and $B$, then $F$ is lax symmetric comonoidal.
\end{Lem}

\begin{proof}
Recall that $\square$ on $B$ is a categorical product if there are natural transformations 
\[
a\overset{p_1}{\leftarrow} a\square b\overset{p_2}{\rightarrow} b\quad,
\]
such that the induced map
\[
B(c,a\square b)\rightarrow B(c,a)\times B(c,b)
\]
is a bijection. 
Using that $\lambda$ is an isomorphism, it immediately follows that $e_B$ is a
terminal object. We define $\eta:F(e_M)\rightarrow e_B$ to be the unique map.

It is implicitly understood when saying that ``$\square$ is the categorical product''
that $\alpha$, $\rho$, $\lambda$ and $\gamma$ are related to $p_1$ and $p_2$.
We require that $p_1=\rho$ when $b=e_B$ and $p_2=\lambda$ when $a=e_B$.
Furthermore, the following diagrams must commute:
\[
\begin{CD}
a\square (b\square c) @>{\alpha}>> (a\square b)\square c\\
@V{p_1}VV @VV{p_1}V\\
a @<{p_1}<< a\square b
\end{CD}\quad,\quad
\begin{CD}
a\square (b\square c) @>{\alpha}>> (a\square b)\square c\\
@V{p_2}VV @VV{p_2}V\\
b\square c @<{p_2}<< c
\end{CD}\quad\text{and}
\]
\[
\begin{CD}
a\square (b\square c) & @>{\alpha}>> & (a\square b)\square c\\
@V{p_2}VV && @VV{p_1}V\\
b\square c @>{p_1}>> b @<{p_2}<< a\square b
\end{CD}\quad,
\]
and $\gamma:a\square b\rightarrow b\square a$ is the unique map such that
\[
\begin{CD}
a @<{p_1}<< a\square b @>{p_2}>> b\\
@V{=}VV @V{\gamma}VV @VV{=}V\\
a @<{p_2}<< b\square a @>{p_1}>> b
\end{CD}
\]
commute.

To define $\phi$ we apply $F$ to $p_1$ and $p_2$. Under the bijection 
\[
B(F(a\square b),F(a))\times B(F(a\square b),F(b))\cong B(F(a\square b),F(a)\square F(b))
\]
the pair $(F(p_1),F(p_2))$ corresponds to $\phi$.
It is now an easy exercise to check that diagrams (ix), (x), (xi) and (xiii) commutes.
\end{proof}

\begin{Rem}
There is a dual lemma: If $\square_M$ and $\square_B$ are categorical coproducts, then
any functor $F:M\rightarrow B$ is symmetric monoidal.
\end{Rem}

\section{Arithmetics for May's operad $\mathcal{M}$}

May's operad $\mathcal{M}$ encodes the structure
of a monoid with unit. In this section we will 
derive formulas for the composition operation for
this operad. Note the following fact: The definition of the action of
$\mathcal{M}$ on arbitrary monoids forces the definition of the $\circ_i$'s.

We begin by defining the spaces of $\mathcal{M}$. 
We let
\[
\mathcal{M}(j)=\Sigma_j\quad,
\]
where $\Sigma_j$ denotes the permutation group on the integers $1,2,\ldots, j$. 
When needed
we will write permutations $\rho$ in $\Sigma_j$ as $2\times j$-matrices:
\[
\rho=\begin{pmatrix}
1 & 2 & \ldots & j\\
\rho(1) & \rho(2) &\ldots &\rho(j)
\end{pmatrix}\quad.
\]

The way that $\mathcal{M}$ encodes the structure of a monoid $G$ with unit
is that for every $j$ there is an action
\[
\theta_j:\mathcal{M}(j)\times G^j\rightarrow G\quad.
\]
This is defined by sending $(\rho;g_1,\ldots,g_j)$
to the product $g_{\rho^{-1}(1)}\cdots g_{\rho^{-1}(j)}$.

A main part of an operad is the composition operations $\circ_i$.
The idea behind the $\circ_i$'s is that they describe how to act iteratively.
Assume given elements $\rho\in\mathcal{M}(k)$ and $\upsilon\in \mathcal{M}(j)$.
First use
\[
\theta_j(\upsilon;-)
\]
to multiply $(g_1,\ldots,g_j)$. Let $g'_i$ be the result, and insert it as the $i$'th
factor in $(g'_1,\ldots,g'_k)$.
Next multiply using
\[
\theta_k(\rho;-)\quad.
\]
Now we can hope that there exists some element $\mu\in\mathcal{M}(k+j-1)$
such that
\[
\theta_{k+j-1}(\mu;g'_1,\ldots,g'_{i-1},g_1,\ldots,g_j,g'_{i+1},\ldots,g'_k)
\]
is equal to the result of the two step process above. The composition operation
$\circ_i$ is defined so that $\rho\circ_i\upsilon$ is such a $\mu$.

\begin{Def}
Let $\rho\in\Sigma_k$ and $\upsilon\in\Sigma_j$ be 
permutations and $1\leq i\leq k$.
We define the composition operation
\[
\circ_i:\Sigma_k\times\Sigma_j\rightarrow \Sigma_{k+j-1}
\]
by the formula
\[
(\rho\circ_i\upsilon)(t)=
\begin{cases}
\rho(t)&\text{if $t<i$ and $\rho(t)<\rho(i)$,}\\
\rho(t)+j-1&\text{if $t<i$ and $\rho(t)>\rho(i)$,}\\
\upsilon(t-i+1)+\rho(i)-1&\text{if $i\leq t< j+i$,}\\
\rho(t-j+1)&\text{if $j+i\leq t$ and $\rho(t-j+1)<\rho(i)$ and}\\
\rho(t-j+1)+j-1&\text{if $j+i\leq t$ and $\rho(t-j+1)>\rho(i)$.}
\end{cases}
\]
\end{Def}

\begin{Exa}
We now look at some explicit examples:
For instance, if 
\[
\rho=
\begin{pmatrix}1&2&3&4\\ 2&4&1&3\end{pmatrix}
\quad\text{and}\quad\upsilon=
\begin{pmatrix}1&2&3\\3&2&1\end{pmatrix}\quad,
\]
then
\[
\rho\circ_1\upsilon =\begin{pmatrix}
1&2&3&4&5&6\\
4&3&2&6&1&5
\end{pmatrix}\quad,
\]
while
\[
\rho\circ_3\upsilon =\begin{pmatrix}
1&2&3&4&5&6\\
4&6&3&2&1&3
\end{pmatrix}.
\]
\end{Exa}

There is a ``box''-model that can be helpful when trying to visualize this operation.
Given $i$, $\rho$ and $\upsilon$, we put boxes around the integers from $1$ to $k+j-1$
as follows:
\[
\boxed{1},\ldots,\boxed{i-1},\boxed{i, i+1,\ldots, i+j-1},\boxed{i+j},\ldots,\boxed{k+j-1}\quad.
\]
We now use $\rho$ to permute the boxes, while we use $\upsilon$ to permute the
elements in the $i$'th box. Removing the boxes one get the permutation $\rho\circ_i\upsilon$.

To see that our definition of the composition operation is correct, we prove the
following lemma:

\begin{Lem}
If $g'_i=\theta_j(\upsilon;g_1,\ldots,g_j)$ then
\[
\theta_k(\rho;g'_1,\ldots,g'_k)=
\theta_{k+j-1}(\rho\circ_i\upsilon;g'_1,\ldots,g'_{i-1},g_1,\ldots,g_j,g'_{i+1},\ldots,g'_k)\quad.
\]
\end{Lem}

\begin{proof}
By the definition we have that
\[
g'_i=\theta_j(\upsilon;g_1,\ldots,g_j)=g_{\upsilon^{-1}(1)}\cdots g_{\upsilon^{-1}(j)}
\]
and
\[
\theta_k(\rho;g'_1,\ldots,g'_k)=g'_{\rho^{-1}(1)}\cdots g'_{\rho^{-1}(k)}\quad.
\]
If we let $s=\rho(i)$, then the $s$'th factor in the last product is $g'_{\rho^{-1}(s)}=g'_i$
and we have that
\[
\theta_k(\rho;g'_1,\ldots,g'_k)=g'_{\rho^{-1}(1)}\cdots g'_{\rho^{-1}(s-1)}
g_{\upsilon^{-1}(1)}\cdots g_{\upsilon^{-1}(j)}
g'_{\rho^{-1}(s+1)}\cdots g'_{\rho^{-1}(k)}\quad.
\]
To evaluate $\theta_{k+j-1}(\rho\circ_i\upsilon;g'_1,\ldots,g'_{i-1},g_1,\ldots,g_j,g'_{i+1},\ldots,g'_k)$
using the definition of $\theta$ we need an explicit expression for $(\rho\circ_i\upsilon)^{-1}$.
We use the definition of $\circ_i$ to deduce the formula:
\footnotesize
\[
(\rho\circ_i\upsilon)^{-1}(r)=
\begin{cases}
\rho^{-1}(r)&\text{if $\rho^{-1}(r)<i$ and $r<\rho(i)$,}\\
\rho^{-1}(r)+j-1&\text{if $\rho^{-1}(r)>i$ and $r<\rho(i)$,}\\
\upsilon^{-1}(r-\rho(i)+1)+i-1&\text{if $\rho(i)\leq r< \rho(i)+j$,}\\
\rho^{-1}(t-j+1)&\text{if $\rho^{-1}(r-j+1)< i$ and $r-j+1>\rho(i)$ and}\\
\rho^{-1}(t-j+1)+j-1&\text{if $\rho^{-1}(r-j+1)> i$ and $r-j+1>\rho(i)$.}
\end{cases}
\]
\normalsize
Now we see that
\begin{multline*}
\theta_{k+j-1}(\rho\circ_i\upsilon;g'_1,\ldots,g'_{i-1},g_1,\ldots,g_j,g'_{i+1},\ldots,g'_k)\\
=g'_{\rho^{-1}(1)}\cdots g'_{\rho^{-1}(s-1)}
g_{\upsilon^{-1}(1)}\cdots g_{\upsilon^{-1}(j)}
g'_{\rho^{-1}(s+1)}\cdots g'_{\rho^{-1}(k)}\quad.
\end{multline*}
This concludes the proof.
\end{proof}

Let us now deduce a couple of formulas telling us how to calculate using the 
composition operations:

\begin{Lem}\label{Lem:Mformulas}
If $\rho,\rho'\in\Sigma_k$ and $\upsilon,\upsilon'\in \Sigma_j$, then
the following formulas hold:
\begin{itemize}
\item[i)] $\rho\circ_i\upsilon=(\rho\circ_i\id_j)(\id_k\circ_i\upsilon)$.
\item[ii)] $\rho\circ_i\upsilon=(\id_k\circ_{\rho(i)}\upsilon)(\rho\circ_i\id_j)$.
\item[iii)] $\id_k\circ_i(\upsilon\upsilon')=(\id_k\circ_i\upsilon)(\id_k\circ_i\upsilon')$.
\item[iv)] $(\rho\rho')\circ_i\id_j=(\rho\circ_{\rho'(i)}\id_j)(\rho\circ_i\id_j)$.
\item[v)] $(\rho\rho')\circ_i(\upsilon\upsilon')=(\rho\circ_{\rho'(i)}\upsilon)(\rho'\circ_i\upsilon')$.
\end{itemize}
\end{Lem}

\begin{proof}
To check i) we pick $t$ and calculate $(\rho\circ_i\id_j)(\id_k\circ_i\upsilon)(t)$.
First we have
\[
(\id_k\circ_i\upsilon)(t)=\begin{cases}
t &\text{if $t<i$,}\\
\upsilon(t-i+1)+i-1 &\text{if $i\leq t <j+i$, and}\\
t &\text{if $t\geq j+i$.}
\end{cases}
\]
We also have that
\[
(\rho\circ_i\id_j)(t)=
\begin{cases}
\rho(t)&\text{if $t<i$ and $\rho(t)<\rho(i)$,}\\
\rho(t)+j-1&\text{if $t<i$ and $\rho(t)>\rho(i)$,}\\
t-i+\rho(i)&\text{if $i\leq t< j+i$,}\\
\rho(t-j+1)&\text{if $j+i\leq t$ and $\rho(t-j+1)<\rho(i)$ and}\\
\rho(t-j+1)+j-1&\text{if $j+i\leq t$ and $\rho(t-j+1)>\rho(i)$.}
\end{cases}
\]
Putting these together we get that $(\rho\circ_i\id_j)(\id_k\circ_i\upsilon)(t)=(\rho\circ_i\upsilon)(t)$.

For ii) we use the formulas above to compute that
$(\id_k\circ_{\rho(i)}\upsilon)(\rho\circ_i\id_j)(t)=(\rho\circ_i\upsilon)(t)$.
Let us verify this in the case $i\leq t<j+i$. Then
\footnotesize
\[
(\id_k\circ_{\rho(i)}\upsilon)(\rho\circ_i\id_j)(t)=
(\id_k\circ_{\rho(i)}\upsilon)(t-i+\rho(i))=
\upsilon(t-i+\rho(i)-\rho(i)+1)+\rho(i)-1=(\rho\circ_i\upsilon)(t)\quad.
\]
\normalsize

The case iii) is obvious from the formula for $(\id_k\circ_i\upsilon)(t)$.
Also the case iv) is easy. Now formula v) follows from the other cases:
\begin{multline*}
(\rho\rho')\circ_i(\upsilon\upsilon')\\
=((\rho\rho')\circ_i\id_j)(\id_k\circ_i(\upsilon\upsilon'))\\
=(\rho\circ_{\rho'(i)}\id_j)(\rho'\circ_i\id_j)(\id_k\circ_i\upsilon)(\id_k\circ_i\upsilon')\\
=(\rho\circ_{\rho'(i)}\id_j)(\id_k\circ_{\rho'(i)}\upsilon)(\rho'\circ_i\id_j)(\id_k\circ_i\upsilon')\\
=(\rho\circ_{\rho'(i)}\upsilon)(\rho'\circ_i\upsilon')\quad.
\end{multline*}
\end{proof}

We interpret case v) as a formula for the $\Sigma$-equivariance for the operad
$\mathcal{M}$. There are also formulas for iterated compositions. These are:

\begin{Lem}\label{Lem:circonM}
If $\rho\in\Sigma_k$, $\upsilon\in\Sigma_j$ and $\mu\in\Sigma_l$,
then
\begin{itemize}
\item[i)] $(\rho\circ_a\upsilon)\circ_b\mu=(\rho\circ_b\mu)\circ_{a+l-1}\upsilon$ for $b<a$,
\item[ii)] $(\rho\circ_a\upsilon)\circ_b\mu=\rho\circ_a(\upsilon\circ_{b-a+1}\mu)$ for $a\leq b<a+j$, and
\item[iii)] $(\rho\circ_a\upsilon)\circ_b\mu=(\rho\circ_{b-j+1}\mu)\circ_{a}\upsilon$ for $a+j\leq b$.
\end{itemize}
\end{Lem}

\begin{proof}
This is most easily verified using the box model.
\end{proof}

\newpage
\bibliographystyle{amsalpha}
\bibliography{main2,Ekstrabib}

\end{document}